%% file: paper_2d_fs_qtt_solver.tex
\documentclass[review]{siamart1116}

\usepackage[cal=boondoxo]{mathalfa}
\usepackage{tabularx,ragged2e,booktabs,caption}
\usepackage{amssymb,amsmath,amsfonts,amsopn,mathtools}
\usepackage[algo2e,linesnumbered,boxed,ruled,vlined]{algorithm2e}
\usepackage{subcaption}
\usepackage{hyperref}
\usepackage{pgf}
\usepackage{bm}
\usepackage{ulem}
\usepackage{verbatim}
\usepackage{fix-cm}

\usepackage{tikz}
\usetikzlibrary{calc}
\usetikzlibrary{plotmarks}

\newcolumntype{C}[1]{>{\Centering}m{#1}}

\DeclareMathOperator{\diag}{diag}
\DeclareMathOperator{\vdiag}{\mathbf{diag}}
\DeclareMathOperator{\myvec}{\mathbf{vec}}

\newcommand{\set}[1]{\mathbb{#1}}		
\newcommand{\myfunc}[1]{\textsf{#1}}	
\newcommand{\tensor}[1]{\mathcal{#1}}	
\newcommand{\qttvec}[1]{\mathbcal{#1}}	
\newcommand{\qttmat}[1]{\mathcal{#1}}	

\newcommand{\TheTitle}{ROBUST DISCRETIZATION IN QUANTIZED TENSOR TRAIN FORMAT FOR ELLIPTIC PROBLEMS IN TWO DIMENSIONS} 
\newcommand{\TheAuthors}{A. V. Chertkov, I. V. Oseledets, M. V. Rakhuba}
\headers{\TheTitle}{\TheAuthors}

\title{{\TheTitle}\thanks{
	Submitted to the editors December 6, 2016.
	\funding{
		The work was supported by the RFBR grant 16-31-60095 mol\_a\_ dk.
	}
}}

\author{
	Andrei V. Chertkov
	\thanks{
		Skolkovo Institute of Science and Technology, Skolkovo Innovation Center, Building 3, Moscow 143026, Russia
		(\email{andrei.chertkov@skoltech.ru}, \email{i.oseledets@skoltech.ru}, \email{m.rakhuba@skoltech.ru}).
	}
	\and
	Ivan V. Oseledets
	\footnotemark[2]
	\thanks{
		Institute of Numerical Mathematics, Russian Academy of Sciences. Gubkina St. 8, Moscow 119333, Russia.
	}
	\and
	Maxim V. Rakhuba
	\footnotemark[2]
}

\ifpdf
\hypersetup{
  pdftitle={\TheTitle},
  pdfauthor={\TheAuthors}
}
\fi

\begin{document}

\maketitle

\begin{abstract}
	In this work we propose an efficient black-box solver for two-dimensional stationary diffusion equations, which is based on a new robust discretization scheme.
	The idea is to formulate an equation in a certain form without derivatives with a non-local stencil, which leads us to a linear system of equations with dense matrix.
	This matrix and a right-hand side are represented in a low-rank parametric representation -- the quantized tensor train (QTT-) format, and then all operations are performed with logarithmic complexity and memory consumption.
	Hence very fine grids can be used, and very accurate solutions with extremely high spatial resolution can be obtained.
	Numerical experiments show that this formulation gives accurate results and can be used up to $2^{60}$ grid points with no problems with conditioning, while total computational time is around several seconds.
\end{abstract}
	


\section{Introduction} \label{section:intro} 

It has been recently proved in~\cite{ks-2dqtt-2015} that for a class of partial differential equations (PDEs) with piecewise-analytic coefficients the exact solution $u$ can be approximated with accuracy $\epsilon$ in the energy norm with $N=\mathcal{O}\left(\log^{\alpha} \epsilon^{-1}\right)$ degrees of freedom for $\alpha \leq 5$ in the QTT-FEM (quantized tensor train - finite element method) approach.
In the QTT-FEM approach the coefficient vector of the FEM solution is represented in the QTT-format~\cite{osel-2d2d-2010, khor-qtt-2011}.
The key idea of the QTT-format is to reshape the discretized solution into a tensor, which is then compressed in the low-rank tensor train (TT-) format~\cite{osel-tt-2011}.

The computation of the approximate solution in the QTT-format looks straightforward.
We discretize the PDE on a very fine virtual mesh with $2^d$ points and use a finite element/difference method with constraint on the number of parameters in the QTT-representation.
Efficient solvers in the QTT- and TT-format are available for this task~\cite{DoOs-dmrg-solve-2011, ds-amen-2014}.
However, if the equation is discretized using standard low-order FEM methods, it is not possible to get to very fine meshes, and this becomes a key issue.
Let us illustrate it for an example of the one-dimensional Poisson equation $-\nabla^2 u = f(x)$, $x \in [0, 1]$ with homogeneous Dirichlet boundary conditions.
The simplest discretization scheme reads
\begin{equation}
	\frac{u_{i+1} - 2 u_i + u_{i-1}}{h^2} = - f(x_i), \quad
	i = 1, 2, \ldots, n-1, \quad 
	u_0 = u_n = 0, \quad 
	h = \frac{1}{n},
\end{equation}
and it is well known that $|u_i - u(x_i)| = \mathcal{O}\left(h^2\right)$ for smooth enough solution, i.e. the smaller the grid size $h$, the better the approximation is.
However, in numerical computations we can not take $h$ too small. 
Indeed, let $\tau$ be the rounding error introduced by arithmetic operations. 
Then the approximation error of the action of discrete second-order derivative operator can be estimated as $\mathcal{O}\left( \tau / h^2 \right)$, and the total error is $\mathcal{O}\left( \tau / h^2 + h^2 \right)$, which means that the error reaches its minimal value at $h_{min} \sim \tau^{1/4}$, and for the double precision $\tau \approx 10^{-16}$ we have $h_{min} \sim 10^{-4}$.
Smaller values of $h$ will lead to larger errors.
Therefore, one can not treat straightforwardly problems where smaller $h$ is required, e.g. multiscale and high-frequency problems.
		
We propose to solve this problem by using a discretization scheme with a non-local stencil, that is robust for any $h$, but still can be efficiently implemented  in a low-rank tensor format.
The full list of our contributions is the following.
\begin{itemize}
	\item We derived a new robust discretization scheme for the two-dimensional diffusion-type equation (\cref{section:scheme}). The main idea is to rewrite the initial problem in a derivative-free formulation, which leads us to a linear system of equations with a dense matrix.
	\item We formulated this discretization scheme in the QTT-format, so all operations are performed with logarithmic complexity and memory (\cref{section:qtt_solver}), and we developed a robust Finite Sum (FS-QTT-) solver, based on this scheme.
	\item We proved that proposed scheme has a second order convergence with respect to the grid size (\cref{section:connection}), and that TT-ranks of the matrix of the corresponding linear system are bounded by a certain number, which depends on the TT-ranks of the right-hand side and of inverse of the PDE coefficients (\cref{section:qtt_solver}).
	\item We performed numerical experiments (\cref{section:numex}), that demonstrate robustness of the proposed scheme. FS-QTT-solver can handle up to $2^{60}$ virtual grid points, while total computational time is around several seconds.	
\end{itemize} 

\section{Robust discretization scheme} \label{section:scheme}

\subsection{Model problem}
		
Our problem is a diffusion-type PDE with homogeneous Dirichlet boundary conditions of the form
\begin{equation}\label{scheme:model_pde}
	- \nabla \left( \mathcal{K} \nabla u \right) = f, \quad u |_{\partial \Omega} = 0,
\end{equation}
in a bounded two-dimensional domain $\Omega = [0, 1]^2$, where $f$ is a function in $\Omega$, $\mathcal{K}$ is a diagonal diffusion tensor with components $k_x$ and $k_y$, and $u$ is an unknown scalar field.
	
In this section we present an informal derivation of a new discretization scheme for PDEs of the form~\eqref{scheme:model_pde} from the variational principle, and in \cref{section:connection} it will be proved that the obtained scheme is equivalent in the exact arithmetic to the second order finite-difference scheme.
It should be noted, that the proposed method can be generalized to other forms of coefficient $\mathcal{K}$ and types of boundary conditions. 
				
\subsection{Variational formulation}
		
The problem~\eqref{scheme:model_pde} with $k_x, k_y, f \in \set{L}_2(\Omega)$ is equivalent to the minimization of the functional\footnote{
	Functions are denoted hereinafter by lower-case letters ($v$ or $v(x,y)$, $\ldots$), functionals and operators are denoted by upper case letters with function in square brackets ($F[v]$, $\ldots$).
}
\begin{equation}\label{scheme:functional_def}
	u = \arg \min_{v \in {H}_0^1(\Omega)} F[v], \quad 
	F[v] = \int_{\Omega} k_x \left( \frac{\partial v}{\partial x} \right) ^2 + 
	       \int_{\Omega} k_y \left( \frac{\partial v}{\partial y} \right) ^2 - 
         2 \int_{\Omega} v f.
\end{equation}	
To transform~\eqref{scheme:functional_def} into a derivative-free form, we introduce two new variables 
\begin{equation}\label{scheme:vx_vy_def}
	v_x (x,y) = \frac{\partial v}{\partial x} (x,y), \quad
	v_y (x,y) = \frac{\partial v}{\partial y} (x,y).
\end{equation}
Using~\eqref{scheme:vx_vy_def} and taking into account the homogeneous boundary conditions for $x = 0$ and $y = 0$, we can write
\begin{equation}\label{scheme:v_from_vx}
	B_x [v_x] \equiv  \int_0^x v_x(t, y) \, dt = v(x,y),
\end{equation}
\begin{equation}\label{scheme:v_from_vy}
	B_y [v_y] \equiv  \int_0^y v_y(x, t) \, dt = v(x,y).
\end{equation}
To enforce the homogeneous boundary conditions for $x = 1$ and $y = 1$ we have to set the following constraints
\begin{equation}\label{scheme:vx_constraint}
	S_x [v_x] \equiv \int_0^1 v_x(t, y) \, dt = 0, \quad y \in [0, 1],
\end{equation}
\begin{equation}\label{scheme:vy_constraint}
	S_y [v_y] \equiv \int_0^1 v_y(x, t) \, dt = 0, \quad x \in [0, 1],
\end{equation}
and from equations~\eqref{scheme:v_from_vx},~\eqref{scheme:v_from_vy} we have one more constraint on the new variables $v_x$, $v_y$
\begin{equation}\label{scheme:vx_vs_vy_constraint}
	B_x [v_x] = B_y [v_y].
\end{equation}
If we substitute~\eqref{scheme:v_from_vy} to~\eqref{scheme:functional_def}, we come to the following optimization problem for the derivatives of the solution of equation~\eqref{scheme:model_pde}
\begin{equation}\label{scheme:functional}
	u_x, u_y = \arg \min  F [v_x, v_y], \quad
    F [v_x, v_y] = \int_{\Omega} k_x v_x ^2 + \int_{\Omega} k_y v_y ^2 - 2 \int_{\Omega} B_y [v_y] f,
\end{equation}
with constraints~\eqref{scheme:vx_constraint},~\eqref{scheme:vy_constraint} and~\eqref{scheme:vx_vs_vy_constraint}.

\subsection{Discretization on the spatial grid} \label{subsection:scheme_discr}
	
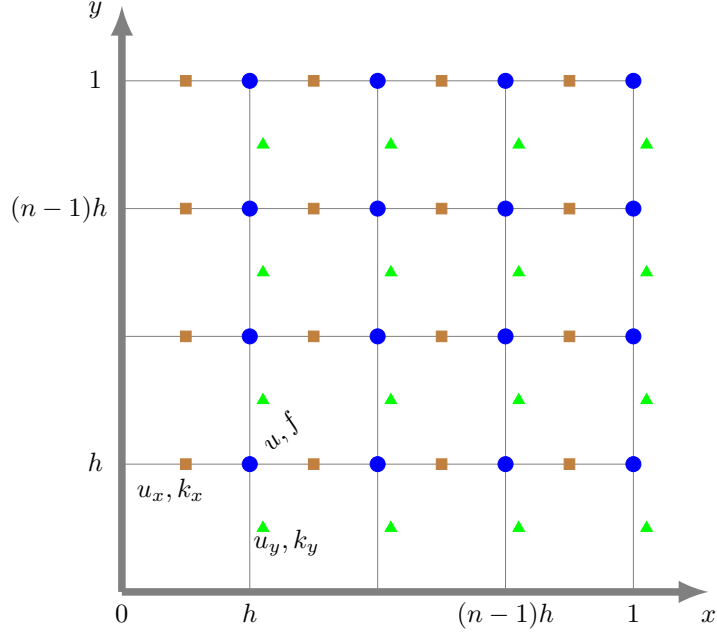
\begin{figure}[ht]\centering
	\begin{tikzpicture}[
		axtxt/.style={  
			draw=none,
			text height=1.5ex,
			text depth=.25ex,
			minimum height=3em,
		}]
		\def\h{1.7}; \def\n{4}; \def\L{\h*\n}; \def\dL{1.}
	    \clip (-2,-1) rectangle (\L+2*\dL,\L+2*\dL);
	    \draw [line width=.1cm, gray,-latex] (0,0) -- (\L+\dL,0);
	    \draw [line width=.1cm, gray,-latex] (0,0) -- (0,\L+\dL);
	    \draw[style=help lines] (0,0) grid[line width=.4,step=\h] (\L,\L);
	       
	    \foreach \i in {1,...,\n}{
	    	\foreach \j in {1,...,\n}{
	      		\node[draw,circle,inner sep=2pt,fill,blue] at (\h*\i,\h*\j) {};
	        	\node[mark size=2.5pt,color=green] at (\h*\i,\h*\j-\h*0.5) {\pgfuseplotmark{triangle*}};
	        	\node[mark size=2pt,color=brown] at (\h*\i-\h*0.5,\h*\j) {\pgfuseplotmark{square*}};
	      		}
		}
	      			
		\def\dy{-0.3}
		\node[axtxt] at (0,\dy) {$0$};
		\node[axtxt] at (\h,\dy) {$h$};
		\node[axtxt] at (\n*\h-\h,\dy) {$(n-1)h$};
		\node[axtxt] at (\L,\dy) {$1$};
		\node[axtxt] at (\L+\dL,\dy) {$x$};
	    			
		\def\dx{-0.35}
		\node[axtxt] at (\dx,\h) {$h$};
		\node[axtxt] at (\dx*2.4,\n*\h-\h) {$(n-1)h$};
		\node[axtxt] at (\dx,\L) {$1$};
		\node[axtxt] at (\dx,\L+\dL) {$y$};
	
		\node[draw=none,text width=1.9cm,rotate=45] at (\h+0.9,\h+0.9) {$u,f$};
		\node[draw=none,text width=1.9cm] at (\h+1.0,\h*0.5-0.2) {$u_y,k_y$};
		\node[draw=none,text width=1.9cm] at (\h*0.5+0.3,\h-0.35) {$u_x,k_x$}; 
	\end{tikzpicture}
  	\caption{
  		Example of the spatial grid with the step $h=2^{-d}$, where $d=2$.
  		Upper right cell corners (blue circles) are used for the discretization of the solution and the right-hand side.
  		The $x$-derivative and $k_x$ coefficient are discretized on midpoints of top edges of the cells (brown squares), and the $y$-derivative and $k_y$ coefficient are discretized on midpoints of right edges of the cells (green triangles).
  	}
  	\label{figure:grid}
\end{figure}
			
Now we discretize the functional~\eqref{scheme:functional} and constraints~\eqref{scheme:vx_constraint},~\eqref{scheme:vy_constraint},~\eqref{scheme:vx_vs_vy_constraint} on a tensor-product square uniform grid with $n^2=2^{2d}$ nodes and with grid step $h=1 / n = 2^{-d}$, where $d \in \set{N}$.
The grid for the case $d=2$ is presented in \cref{figure:grid}.
The solution $u$ and the right-hand side $f$ are discretized in upper right cell corners (are marked by blue circles in \cref{figure:grid}) and its values are collected in vectors\footnote{
	Vectors and matrices are denoted hereinafter by lower case bold letters ($\bm{a}, \bm{b}, \bm{c}, \ldots$) and upper case letters ($A, B, C, \ldots$) respectively.
	We denote the (i,j)th element of an $n \times n$ matrix $A$ as $A[i, j]$ and assume that numeration starts from the zero, so $0 \leq i, j < n$.
	For vectors we use the same notation: $\bm{a}[i]$ (for $i=0,1,\ldots,n-1$) is the $i$th element of the vector $\bm{a}$.
}
$\bm{u}$ and $\bm{f}$ respectively
$$
	\bm{u}[\alpha] = u\left(x_{i(\alpha)}, y_{j(\alpha)} \right),
$$
\begin{equation}\label{scheme:descr_f}
	\bm{f}[\alpha] = f\left(x_{i(\alpha)}, y_{j(\alpha)}\right),		
\end{equation}
where $\alpha = 0, 1, \ldots, n^2-1$, and
$$
    x_{i(\alpha)} = \left(i(\alpha)+1\right)h, \quad y_{j(\alpha)} = \left(j(\alpha)+1\right)h,
$$
$$
    i(\alpha) = \alpha \text{ mod } n, \quad j(\alpha) = \frac{ \alpha-i(\alpha)}{n}.
$$
Derivative $u_x$ and coefficient $k_x$ are discretized on midpoints of top edges of the cells (are marked by brown squares in \cref{figure:grid}) and are represented as vector $\bm{u_x}$ and an $n^2 \times n^2$ diagonal\footnote{
	We denote by $\diag (\bm{\cdot})$ an operation that constructs diagonal matrix for a given vector, and an operation $\vdiag (\cdot)$ constructs a vector, that is a diagonal of a given matrix.
}
matrix $K_x$ respectively
$$
	\bm{u}_x[\alpha] = \frac{\partial u}{\partial x} \left(x_{i(\alpha)-1/2}, y_{j(\alpha)}\right),
$$
\begin{equation}\label{scheme:descr_kx} 
	K_x = \diag \left(\bm{k}_x\right), \quad \bm{k}_x[\alpha] = k_x\left(x_{i(\alpha)-1/2}, y_{j(\alpha)}\right),
\end{equation}
where
$$
    x_{i(\alpha)-1/2} = \left(i(\alpha)+ \frac 12 \right)h.
$$
Midpoints of right edges are used for discretization  of $u_y$ and $k_y$ (are marked by green triangles in \cref{figure:grid}), and in the discrete setting these quantities are represented as $\bm{u_y}$ and $K_y$ respectively
$$
	\bm{u}_y[\alpha] = \frac{\partial u}{\partial y} \left(x_{i(\alpha)}, y_{j(\alpha)-1/2}\right),
$$
\begin{equation}\label{scheme:descr_ky} 
	K_y = \diag (\bm{k}_y), \quad \bm{k}_y[\alpha] =  k_y \left(x_{i(\alpha)}, y_{j(\alpha)-1/2}\right),
\end{equation}
where
$$
    y_{j(\alpha)-1/2} = \left(j(\alpha)+\frac 12 \right)h.
$$

Integrals in~\eqref{scheme:v_from_vx} and~\eqref{scheme:v_from_vy} are approximated by a simple rectangular quadrature formula with nodes in upper right cell corners (blue circles in \cref{figure:grid}), and then the discretized operators $B_x [\cdot]$ and $B_y [\cdot]$ take the form
\begin{equation}\label{scheme:Bx_By}
	B_x = I \otimes B, \quad B_y = B \otimes I,
\end{equation}
where $I$ is an $n \times n$ identity matrix and $B$ is an $n \times n$ matrix given as
\begin{equation}\label{scheme:B}
   	B [i, j] = 
   	\begin{cases}
		h, & i \geq j, \\
		0, & \mbox{otherwise},
    \end{cases}
\end{equation} 
for $i,j=0, 1, \ldots, n-1$. Integrals in~\eqref{scheme:vx_constraint} and~\eqref{scheme:vy_constraint} are approximated in the same manner, and for the discretized operators $S_x [\cdot]$ and $S_y [\cdot]$ we have
\begin{equation}\label{scheme:Sx_Sy}
	S_x = I \otimes \bm{e}^{\top}, \quad S_y = \bm{e}^{\top} \otimes I,
\end{equation}
where $\bm{e}$ is a vector of ones of length $n$.

Then, in the discrete setting, we can approximate the functional~\eqref{scheme:functional} by a discrete functional
\begin{equation}\label{scheme:functional_discr}
	F [\bm{v}_x, \bm{v}_y] = (K_x \bm{v}_x, \bm{v}_x) + (K_y \bm{v_y}, \bm{v}_y) - 2 (B_y \bm{v}_y, \bm{f}),
\end{equation}
and constraints~\eqref{scheme:vx_constraint},~\eqref{scheme:vy_constraint},~\eqref{scheme:vx_vs_vy_constraint} take the form
\begin{equation}\label{scheme:vx_vs_vy_constraint_discr}
	B_x \bm{v}_x = B_y \bm{v}_y,
\end{equation}
\begin{equation}\label{scheme:vx_and_vy_constraint_discr}
	S_x \bm{v}_x = 0, \quad S_y \bm{v}_y = 0.
\end{equation}

\subsection{Minimization of the functional}

\begin{theorem} \label{scheme:th:result}
	Minimization of the functional~\eqref{scheme:functional_discr}
	$$
		\bm{u}_x, \bm{u}_y = \arg \min F [\bm{v}_x, \bm{v}_y] 
	$$
	with constraints~\eqref{scheme:vx_vs_vy_constraint_discr} and~\eqref{scheme:vx_and_vy_constraint_discr} give the following formulas for the approximated derivatives $\bm{u}_x$ and $\bm{u}_y$ of the solution of PDE~\eqref{scheme:model_pde}
	\begin{equation}\label{scheme:result_ux_uy}
	   	\bm{u}_x = R_x \bm{\mu}, \quad
	   	\bm{u}_y = R_y \left( \bm{f} - \bm{\mu} \right),
	\end{equation}
	and the approximated solution $\bm{u}$ can be recovered as 
	\begin{equation}\label{scheme:result_u}
		\bm{u}   = H_x \bm{\mu},
	\end{equation}
	where $\bm{\mu}$ is a solution of a linear system
	\begin{equation}\label{scheme:eq_for_mu}
	  	\left( H_x + H_y \right) \bm{\mu} =  H_y \bm{f},
	\end{equation}
	and the following definitions for the matrices are used
	\begin{equation}\label{scheme:HRW_x_def}
	  	H_x = B_x R_x, \quad
	  	R_x = K_x^{-1} \left( I \otimes I - W_x K_x^{-1} \right) B_x^T, \quad
	  	W_x = Q_x \otimes E, 
	\end{equation}
	\begin{equation}\label{scheme:HRW_y_def}
	  	H_y = B_y R_y, \quad
	  	R_y = K_y^{-1} \left( I \otimes I - W_y K_y^{-1} \right) B_y^T, \quad
	  	W_y = E \otimes Q_y,
	\end{equation}
	where $I$ is an $n \times n$ identity matrix, $E$ is an $n \times n$ matrix of ones, and\footnote{
		Vector $\bm{a}^{-1}$ denotes elementwise inversion of the vector $\bm{a}$.
	}
	\begin{equation}\label{scheme:Qx_explicit}
		\begin{split}
	  		& Q_x = \diag (\bm{q}_x), \quad \bm{q}_x^{-1}[i] = \sum_{\alpha=0}^{n-1} \bm{k}_x^{-1} [in+\alpha], \\
	  		& \quad \bm{k}_x^{-1} = \vdiag (K_x^{-1}), \quad i=0,\ldots,n-1,
		\end{split}
	\end{equation}
	\begin{equation}\label{scheme:Qy_explicit}
		\begin{split}
			& Q_y = \diag (\bm{q}_y), \quad \bm{q}_y^{-1}[i] = \sum_{\alpha=0}^{n-1} \bm{k}_y^{-1} [i+n\alpha], \\
	  		& \quad \bm{k}_y^{-1} = \vdiag (K_y^{-1}), \quad i=0,\ldots,n-1.
		\end{split}
	\end{equation}
\end{theorem}
			
\begin{proof}	
	Let us introduce Lagrange multipliers for constraints~\eqref{scheme:vx_vs_vy_constraint_discr} and~\eqref{scheme:vx_and_vy_constraint_discr} as a vector $2\bm{\mu}$ of length $n^2$ and vectors $2\bm{\phi}_x$  and $2\bm{\phi}_y$ of length $n$, then minimization problem for the functional from~\eqref{scheme:functional_discr} is equivalent to the minimization of the unconstrained functional 
	\begin{equation}
		\begin{split}
			\widetilde{F} & \left[ \bm{v}_x, \bm{v}_y, \bm{\mu}, \bm{\phi}_x, \bm{\phi}_y \right] = 
			  (K_x \bm{v}_x, \bm{v}_x) + (K_y \bm{v}_y, \bm{v}_y) - 2 (B_y \bm{v}_y, \bm{f}) + \\
			& (2\bm{\mu}, B_y \bm{v}_y - B_x \bm{v}_x) + (2\bm{\phi}_x, S_x \bm{v}_x) + (2\bm{\phi}_y, S_y \bm{v}_y)
			\rightarrow \min.  
		\end{split}
	\end{equation}
	Optimality conditions give us
	\begin{equation}\label{scheme:opt_cond}
		\begin{cases}
	   		K_x \bm{u}_x - B_x^T \bm{\mu} + S_x^T \bm{\phi}_x  = 0, \\
	   		K_y \bm{u}_y + B_y^T \bm{\mu} + S_y^T \bm{\phi}_y  = B_y^T \bm{f},
		\end{cases}
	\end{equation}
	with constraints~\eqref{scheme:vx_vs_vy_constraint_discr} and~\eqref{scheme:vx_and_vy_constraint_discr}.
	From~\eqref{scheme:opt_cond} we can express $\bm{u}_x$ and $\bm{u}_y$
	\begin{equation}\label{scheme:system_elim_ux_uy}
		\begin{cases}
	   		\bm{u}_x = K_x^{-1} B_x^T \bm{\mu} - K_x^{-1} S_x^{\top} \bm{\phi}_x, \\
	   		\bm{u}_y =-K_y^{-1} B_y^T \bm{\mu} - K_y^{-1} S_y^{\top} \bm{\phi}_y + K_y^{-1} B_y^T \bm{f},
		\end{cases}
	\end{equation}
	and substituting it into~\eqref{scheme:vx_and_vy_constraint_discr} we have
	\begin{equation}\label{scheme:system_phix_phiy}
		\begin{cases}
	   		S_x K_x^{-1} B_x^T \bm{\mu} - S_x K_x^{-1} S_x^{\top} \bm{\phi}_x = 0, \\
	   	   -S_y K_y^{-1} B_y^T \bm{\mu} - S_y K_y^{-1} S_y^{\top} \bm{\phi}_y + S_y K_y^{-1} B_y^T \bm{f} = 0.
		\end{cases}
	\end{equation}
	Fortunately, the matrices 
	\begin{equation}\label{scheme:Q_def}
	  	Q_x^{-1} = S_x K_x^{-1} S_x^{\top}, \quad Q_y^{-1} = S_y K_y^{-1} S_y^{\top},
	\end{equation}
	are $n \times n$ diagonal matrices, and it can be easily shown that they have the form~\eqref{scheme:Qx_explicit},~\eqref{scheme:Qy_explicit}.
	We can express $\phi_x$ and $\phi_y$ from~\eqref{scheme:system_phix_phiy}, using $Q_x$ and $Q_y$ matrices from~\eqref{scheme:Q_def}
	$$
		\begin{cases}
	   		\bm{\phi}_x = Q_x S_x K_x^{-1} B_x^T \bm{\mu}, \\
	   		\bm{\phi}_y =-Q_y S_y K_y^{-1} B_y^T \bm{\mu} +  Q_y S_y K_y^{-1} B_y^T \bm{f}.
		\end{cases}
	$$
	Introduce intermediate $n^2 \times n^2$ matrices 
	\begin{equation}\label{scheme:W_def}
	  	W_x = S_x^{\top} Q_x S_x, \quad W_y = S_y^{\top} Q_y S_y,
	\end{equation}
	which can be represented also in a more compact form~\eqref{scheme:HRW_x_def},~\eqref{scheme:HRW_y_def}.
	Using $W_x$ and $W_y$ matrices, we can substitute the expressions for $\bm{\phi}_x$ and $\bm{\phi}_y$ into~\eqref{scheme:system_elim_ux_uy}
	$$
		\begin{cases}
	   		\bm{u}_x = K_x^{-1} B_x^T \bm{\mu} - K_x^{-1} W_x K_x^{-1} B_x^T \bm{\mu}, \\
	   		\bm{u}_y =-K_y^{-1} B_y^T \bm{\mu} + K_y^{-1} W_y K_y^{-1} B_y^T \bm{\mu} 
	   		          -K_y^{-1} W_y K_y^{-1} B_y^T \bm{f} + K_y^{-1} B_y^T \bm{f},
		\end{cases}
	$$
	or in a more compact form
	\begin{equation}\label{scheme:eq_for_u_der_pre}
		\begin{cases}
	   		\bm{u}_x = K_x^{-1} \left( I \otimes I - W_x K_x^{-1} \right) B_x^T \bm{\mu} 
	   		= R_x \bm{\mu}, \\
	   		\bm{u}_y = K_y^{-1} \left( I \otimes I - W_y K_y^{-1} \right) B_y^T \left( \bm{f} - \bm{\mu} \right) 
	   		= R_y \left( \bm{f} - \bm{\mu} \right),
		\end{cases}
	\end{equation}
	where $R_x$ and $R_y$ are came from~\eqref{scheme:HRW_x_def} and~\eqref{scheme:HRW_y_def}.
	Putting~\eqref{scheme:eq_for_u_der_pre} into~\eqref{scheme:vx_vs_vy_constraint_discr}, we get equation for $\bm{\mu}$
	$$
		B_x R_x \bm{\mu} = B_y R_y (\bm{f}-\bm{\mu}),
	$$
	which can be rewritten in the form~\eqref{scheme:eq_for_mu} if we introduce $H_x$ and $H_y$ matrices from~\eqref{scheme:HRW_x_def} and~\eqref{scheme:HRW_y_def}. Since $\bm{u} = B_x \bm{u}_x = B_y \bm{u}_y$, and using the first equation in~\eqref{scheme:eq_for_u_der_pre}, we immediately obtain $ \bm{u} = H_x \bm{\mu}$.
\end{proof}

\section{Connection with finite difference scheme} \label{section:connection}
	
Let consider the same spatial grid structure as in \cref{subsection:scheme_discr}, and write a second order finite-difference scheme for the model PDE~\eqref{scheme:model_pde} in the following form
\begin{equation}\label{connection:fd_fin_dif_eq_full}
\begin{split}
	& k_x (x_{i+1/2}, y_{j}) u_{i+1, j} - \left( k_x (x_{i+1/2}, y_{j}) + k_x (x_{i-1/2}, y_{j}) \right) u_{i, j} + \\
	& k_x (x_{i-1/2}, y_{j}) u_{i-1, j} + k_y (x_{i}, y_{j+1/2}) u_{i, j+1} + k_y (x_{i}, y_{j-1/2}) u_{i, j-1} - \\
	& \left( k_y (x_{i}, y_{j+1/2}) + k_y (x_{i}, y_{j-1/2}) \right) u_{i, j}
    = - h^2 f (x_{i}, y_{j}),
\end{split}
\end{equation}
for
$$
	i = 1, \ldots, n-2, \quad j = 1, \ldots, n-2,
$$
with boundary conditions
\begin{equation}\label{connection:fd_fin_dif_eq_full_bc}
\begin{split}
	& u_{i, j} = 0, \quad i = 0, n-1, \quad j = 1, \ldots, n-2, \\  
    & u_{i, j} = 0, \quad j = 0, n-1, \quad i = 1, \ldots, n-2, \\
    & u_{i, j} = 0, \quad i = 0, n-1, \quad j = 0, n-1,
\end{split}
\end{equation}
where we use notation $u_{i, j} = u ( x_{i}, y_{j} )$ for the nodal values of the unknown solution $u$.

At the same time, if we note that the inverse of the matrix $B$ from~\eqref{scheme:B} is a finite-difference matrix
\begin{equation}\label{connection:iB}
	B = h
	\begin{pmatrix}
		1      &       &        &   \\
 		1      & 1     &        &   \\
  		\vdots &       & \ddots &   \\
  	  	1      & \dots & \dots	& 1 
	\end{pmatrix},
	\quad
	B^{-1} = \frac{1}{h}
	\begin{pmatrix}
		1  &        &        &    \\
   		-1 & \ddots &        &    \\
   		   & \ddots & \ddots &    \\
  	  	   &        & -1     & 1
	\end{pmatrix},
\end{equation}
then, taking into account homogeneous Dirichlet boundary conditions on the bottom and left boundaries of the domain, we can write for the derivatives
$$
	\bm{u}_x = (I \otimes B^{-1}) (I \otimes J) \bm{u}, \quad 
	\bm{u}_y = (B^{-1} \otimes I) (J \otimes I) \bm{u},
$$
where $I$ is an $n \times n$ identity matrix and $J$ is an $n \times n$ matrix of the form $J=\diag (\bm{[1, 1, \ldots, 1, 0]}^T)$.
Using homogeneous Dirichlet boundary conditions on the top and right boundaries, we can write the matrix formulation of equation~\eqref{connection:fd_fin_dif_eq_full} as follows
\begin{equation}\label{connection:fin_dif_eq}
	\left( A_x + A_y \right) \bm{u} = \left( I \otimes Z_{n-1} \right) \left( Z_{n-1} \otimes I \right) \bm{f}.
\end{equation}
Matrices $A_x$ and $A_y$ in~\eqref{connection:fin_dif_eq} have the form
\begin{equation}\label{connection:Ax}
	A_x = (I \otimes J) (I \otimes B^{-T}) K_x (I \otimes B^{-1}) (I \otimes J) + (I \otimes Z_{n-1}),
\end{equation}
\begin{equation}\label{connection:Ay}
	A_y = (J \otimes I) (B^{-T} \otimes I) K_y (B^{-1} \otimes I) (J \otimes I) + (Z_{n-1} \otimes I),
\end{equation}
where
\begin{equation}\label{connection:Z_def}
   	Z_{\alpha}[i, j] = 
   	\begin{cases}
    	1, & i = j = \alpha, \\
       	0, & \mbox{otherwise},
    \end{cases}
\end{equation} 
is an $n \times n$ matrix for each $\alpha=0,1,\ldots,(n-1)$, and $K_x$, $K_y$ are discretized PDE coefficients $k_x$ and $k_y$ respectively.
	
\begin{lemma}\label{connection:lemma:A_compact}
	The matrices $A_x$ and $A_y$ from~\eqref{connection:Ax} and~\eqref{connection:Ay} can be rewritten in the form
	\begin{equation}\label{connection:Ax_compact}
		A_x = \sum_{\alpha=0}^{n-1} Z_{\alpha} \otimes \left( J B^{-T} K_{x,\alpha} B^{-1} J + Z_{n-1} \right),
	\end{equation}
	\begin{equation}\label{connection:Ay_compact}
		A_y = \sum_{\alpha=0}^{n-1} \left( J B^{-T} K_{y,\alpha} B^{-1} J + Z_{n-1} \right) \otimes Z_{\alpha},
	\end{equation}
	where $K_{x,\alpha}$ and $K_{y,\alpha}$ for each $\alpha=0,1,\ldots,n-1$ are partitions of matrices $K_x$ and $K_y$ respectively
	\begin{equation}\label{connection:Kx_part_def}
		K_{x,\alpha}[i, j] = K_x[\alpha n + i, \alpha n + j], \quad i,j=0,1,\ldots,n-1,
	\end{equation}
	\begin{equation}\label{connection:Ky_part_def}
		K_{y,\alpha}[i, j] = K_y[\alpha + i n, \alpha + j n], \quad i,j=0,1,\ldots,n-1.
	\end{equation}	
\end{lemma}

\begin{proof}
	Using partition~\eqref{connection:Kx_part_def}, we can rewrite diagonal matrix $K_x$ as a sum of Kronecker products
	\begin{equation}\label{connection:Kx_part}
		K_x = \sum_{\alpha=0}^{n-1} Z_{\alpha} \otimes K_{x,\alpha},
	\end{equation}
	and if we substitute it to~\eqref{connection:Ax}, we obtain
	$$
		A_x = (I \otimes J) \left(I \otimes B^{-T}\right)
		\left( \sum_{\alpha=0}^{n-1} Z_{\alpha} \otimes K_{x,\alpha} \right)
		\left(I \otimes B^{-1}\right) (I \otimes J) + (I \otimes Z_{n-1}).
	$$
	Since
	$$
		\left(I \otimes J\right) \left(I \otimes B^{-T}\right) = \left(I \otimes J B^{-T} \right),
	$$
	$$
		\left(I \otimes B^{-1}\right) \left(I \otimes J\right) = \left(I \otimes B^{-1} J \right),
	$$
	we can rewrite $A_x$ as follows
	$$
		A_x = \sum_{\alpha=0}^{n-1} Z_{\alpha} \otimes \left( J B^{-T} K_{x,\alpha} B^{-1} J \right) + (I \otimes Z_{n-1}),
	$$
	and, finally, if we notice that $I = \sum_{\alpha=0}^{n-1} Z_{\alpha}$, then we obtain the form~\eqref{connection:Ax_compact}.
	The proof for the matrix $A_y$ can be done by analogy.
\end{proof}
	
\begin{lemma}\label{connection:lemma:H_compact}
	The matrices $H_x$ and $H_y$ from~\eqref{scheme:HRW_x_def} and~\eqref{scheme:HRW_y_def} can be rewritten in the form
	\begin{equation}\label{connection:Hx_compact}
		H_x = \sum_{\alpha=0}^{n-1}
		      Z_{\alpha} \otimes \left( B K^{-1}_{x,\alpha} \left( I - E K^{-1}_{x,\alpha} \bm{q}_x[\alpha] \right) B^T \right),
	\end{equation}
	\begin{equation}\label{connection:Hy_compact}
		H_y = \sum_{\alpha=0}^{n-1} 
	  		  \left( B K^{-1}_{y,\alpha} \left( I - E K^{-1}_{y,\alpha} \bm{q}_y[\alpha] \right) B^T \right) \otimes Z_{\alpha},
	\end{equation}
	where $K^{-1}_{x,\alpha}$ and $K^{-1}_{y,\alpha}$ for each $\alpha=0,1,\ldots,n-1$ are partitions of $K^{-1}_x$ and $K^{-1}_y$ matrices respectively
	\begin{equation}\label{connection:iKx_part_def}
		K^{-1}_{x,\alpha}[i, j] = K^{-1}_x[\alpha n + i, \alpha n + j], \quad i,j=0,1,\ldots,n-1,
	\end{equation}
	\begin{equation}\label{connection:iKy_part_def}
		K^{-1}_{y,\alpha}[i, j] = K^{-1}_y[\alpha + i n, \alpha + j n], \quad i,j=0,1,\ldots,n-1.
	\end{equation}
\end{lemma}	
		
\begin{proof}
	Consider $H_x$ matrix that is defined by~\eqref{scheme:HRW_x_def} and substitute the formula~\eqref{scheme:HRW_x_def} for $W_x$  and~\eqref{scheme:Bx_By} for $B_x$
	$$
		H_x = \left(I \otimes B\right) K_x^{-1} \left( I \otimes I - (Q_x \otimes E) K_x^{-1} \right) \left(I \otimes B^T\right).
	$$
	With partition~\eqref{connection:iKx_part_def} we have
	$$
		H_x = \left(I \otimes B\right) \sum_{\beta=0}^{n-1} (Z_{\beta} \otimes K^{-1}_{x,\beta})
	  		  \left( I \otimes I - (Q_x \otimes E) \sum_{\alpha=0}^{n-1} Z_{\alpha} \otimes K^{-1}_{x,\alpha} \right) 
	  		  \left(I \otimes B^T\right),
	$$
	$$
		H_x = \left( \sum_{\beta=0}^{n-1} Z_{\beta} \otimes \left(B K^{-1}_{x,\beta} \right) \right)
	  		  \left( I \otimes I - \sum_{\alpha=0}^{n-1} \left( Q_x Z_{\alpha} \right) \otimes \left( E K^{-1}_{x,\alpha} \right) \right) 
	  		  \left(I \otimes B^T\right).
	$$
	Using definitions~\eqref{scheme:Qx_explicit}, and~\eqref{connection:Z_def} of $Q_x$ and $Z_{\alpha}$ matrices, and equality $I = \sum_{\alpha=0}^{n-1} Z_{\alpha}$, we obtain	
	$$
		H_x = \left( \sum_{\beta=0}^{n-1} Z_{\beta} \otimes \left(B K^{-1}_{x,\beta} \right) \right)
	  		  \left( \sum_{\alpha=0}^{n-1} Z_{\alpha} \otimes \left( I - E K^{-1}_{x,\alpha} \bm{q}_x[\alpha] \right) \right)
	  		  \left(I \otimes B^T\right),
	$$
	$$
		H_x = \sum_{\beta=0}^{n-1} \sum_{\alpha=0}^{n-1} 
	  		  \left( Z_{\beta} Z_{\alpha} \right) \otimes 
	  		  \left( B K^{-1}_{x,\beta} \left( I - E K^{-1}_{x,\alpha} \bm{q}_x[\alpha] \right) B^T \right).
	$$
	According to the forms of $Z_{\alpha}$ and $Z_{\beta}$ matrices, we can remove one summation, and finally arrive at the formula~\eqref{connection:Hx_compact}.
	For the matrix $H_y$ the proof is similar.
\end{proof}
		
\begin{theorem} \label{connection:th:H_vs_A}
	For the matrix $A_x$ from~\eqref{connection:Ax} and matrix $H_x$ from~\eqref{scheme:HRW_x_def} the equality $A_x H_x = I \otimes J$ holds.
	For the matrix $A_y$ from~\eqref{connection:Ay} and matrix $H_y$ from~\eqref{scheme:HRW_y_def} the equality $A_y H_y = J \otimes I$ holds.
	Here $I$ is an $n \times n$ identity matrix, and $J$ is an $n \times n$ matrix of the form $J=\diag (\bm{[1, 1, \ldots, 1, 0]}^T)$.
\end{theorem}
		
\begin{proof}
	We will use compact representations~\eqref{connection:Ax_compact} and~\eqref{connection:Hx_compact} for matrices $A_x$ and $H_x$ respectively.
	If we multiply these matrices, we have
	\begin{equation}\label{connection:AxHx}
		A_x H_x = \sum_{\alpha=0}^{n-1} Z_{\alpha} \otimes M_{\alpha},
	\end{equation}
	where we introduced a matrix
	$$
		M_{\alpha} = \left( J B^{-T} K_{x,\alpha} B^{-1} J + Z_{n-1} \right)
	                        B K^{-1}_{x,\alpha} \left( I - E K^{-1}_{x,\alpha} \bm{q}_x[\alpha] \right) B^T.
	$$
	Since $J=I-Z_{n-1}$ we can rewrite
	\begin{equation}
		\begin{split}
 			M_{\alpha} = & \left( J B^{-T} K_{x,\alpha} B^{-1} + 
			               (I - J B^{-T} K_{x,\alpha} B^{-1}) Z_{n-1} \right) \\
	    	             &  B K^{-1}_{x,\alpha} \left( I - E K^{-1}_{x,\alpha} \bm{q}_x[\alpha] \right) B^T.
		\end{split}               
	\end{equation}
	Using $B^{-1} B = I$ and $K_{x,\alpha} K^{-1}_{x,\alpha} = I$ we obtain
	\begin{equation}
		\begin{split}
			M_{\alpha} = & J B^{-T} \left( I - E K^{-1}_{x,\alpha} \bm{q}_x[\alpha] \right) B^T + \\
		             	 & (I - J B^{-T} K_{x,\alpha} B^{-1}) Z_{n-1} B K^{-1}_{x,\alpha} 
		               	   \left( I - E K^{-1}_{x,\alpha} \bm{q}_x[\alpha] \right) B^T,
		\end{split}               
	\end{equation}
	and therefore
	$$
		M_{\alpha} = J - J B^{-T} E K^{-1}_{x,\alpha} \bm{q}_x[\alpha] B^T +
			        (I - J B^{-T} K_{x,\alpha} B^{-1}) Z_{n-1} N_{\alpha},
	$$
	where
	\begin{equation}\label{connection:N_matrix}
		N_{\alpha} = B K^{-1}_{x,\alpha}\left( I - E K^{-1}_{x,\alpha} \bm{q}_x[\alpha] \right) B^T.
	\end{equation}
	It can be shown that $J B^{-T} E$ is a zero matrix, and then we have
	\begin{equation}\label{connection:M_by_N_matrix}
		M_{\alpha} = J + (I - J B^{-T} K_{x,\alpha} B^{-1}) Z_{n-1} N_{\alpha}.
	\end{equation}

	Let us write out matrix multiplications for $N_{\alpha}$ matrix \eqref{connection:N_matrix} explicitly and consider the $j$th element of its last row
	\begin{equation}
		\begin{split}
			N_{\alpha}[n-1, j] = & \sum_{\gamma_1=0}^{n-1} \sum_{\gamma_2=0}^{n-1} \sum_{\gamma_3=0}^{n-1}
		                           B[n-1, \gamma_1] K^{-1}_{x,\alpha} [\gamma_1, \gamma_2] \cdot \\
		                         & \left( I - E K^{-1}_{x,\alpha} \bm{q}_x[\alpha] \right) [\gamma_2, \gamma_3] 
		                           B [j, \gamma_3].
		\end{split}	
	\end{equation}
	Using that $K^{-1}_{x,\alpha}=\diag({\bm{k}^{-1}_{x,\alpha}})$, $\left( E K^{-1}_{x,\alpha} \right) [i,j]=\bm{k}^{-1}_{x,\alpha}[j]$, that the last row of the matrix $B$ contains all values equal to $h$, and that $B [j, \gamma_3]$ is non-zero (and is equal to $h$) only for $\gamma_3 \leq j$, we can rewrite
	$$
		N_{\alpha}[n-1, j] = h^2 \sum_{\gamma_2=0}^{n-1} \sum_{\gamma_3=0}^{j} 
			\bm{k}^{-1}_{x,\alpha} [\gamma_2] 
		    \left( \delta_{\gamma_2, \gamma_3} - \bm{k}^{-1}_{x,\alpha}[\gamma_3] \bm{q}_x[\alpha] \right),
	$$
	or
	\begin{equation}\label{connection:N_matrix_last_row}
		N_{\alpha}[n-1, j] =
			h^2                  \sum_{\gamma_3=0}^{j} \bm{k}^{-1}_{x,\alpha} [\gamma_3] -
			h^2 \bm{q}_x[\alpha] \sum_{\gamma_2=0}^{n-1} \bm{k}^{-1}_{x,\alpha} [\gamma_2] \sum_{\gamma_3=0}^{j} \bm{k}^{-1}_{x,\alpha}[\gamma_3].
	\end{equation}
	According to~\eqref{scheme:Qx_explicit} we have $\bm{q}_x^{-1}[\alpha]=\sum_{\beta=0}^{n-1} \bm{k}^{-1}_{x,\alpha} [\beta]$, then from~\eqref{connection:N_matrix_last_row} we have $N_{\alpha}[n-1, j] = 0$ for $j=0,1,\ldots,n-1$, and then $Z_{n-1} N_{\alpha}$ is a zero matrix.
	Therefore we conclude from~\eqref{connection:M_by_N_matrix} that $M_{\alpha} = J$, and finally, using~\eqref{connection:AxHx} and equality $\sum_{\alpha=0}^{n-1} Z_{\alpha} = I$, we obtain the required formula
	$$
		A_x H_x = \sum_{\alpha=0}^{n-1} Z_{\alpha} \otimes M_{\alpha} = \sum_{\alpha=0}^{n-1} Z_{\alpha} \otimes J = I \otimes J.
	$$
	
	If we use compact representations~\eqref{connection:Ay_compact} and~\eqref{connection:Hy_compact} for matrices $A_y$ and $H_y$ respectively, then an idea of the proof for the product $A_y H_y$ will be almost the same, and hence it is not presented here.
\end{proof}
	
The corollaries of the Theorem~\ref{connection:th:H_vs_A} are the following.
	
\begin{corollary} \label{connection:concl:equiv}
	In exact arithmetic solution, obtained according to the proposed discretization scheme~\eqref{scheme:result_u}, is equivalent to the solution, obtained by the second order finite-difference discretization scheme on a uniform grid.
\end{corollary}

\begin{corollary} \label{connection:concl:second_order}
	The solution $\bm{u}$, obtained from the formula~\eqref{scheme:result_u}, converges to the exact solution, when $h$ goes to zero, with the second order in the exact arithmetic under standard conditions for the PDE coefficients.	
\end{corollary}

\section{Solver in the QTT-format} \label{section:qtt_solver}

\subsection{Vectors in a low-rank format}

The proposed in \cref{section:scheme} discretization scheme is robust for small grid steps $h$.
Small $h$ are required in the case when standard finite element/difference scheme converges slowly (corners, point singularities in the spatial domain, etc.).
Storage of the solution when $h$ is small can be prohibitive, we hence go to the main algorithmic contribution of the paper, and propose a "tensor-based" version of the scheme. 

The right-hand side $f$ and coefficients $k_x$, $k_y$ of the model PDE~\eqref{scheme:model_pde} in the discrete setting are considered as the vectors $\bm{f}$, $\bm{k}_x$ and $\bm{k}_y$ from~\eqref{scheme:descr_f},~\eqref{scheme:descr_kx} and~\eqref{scheme:descr_ky} respectively.
These vectors will be represented in the memory-efficient QTT-format~\cite{osel-2d2d-2010, khor-qtt-2011}. 
	
The concept of the QTT-decomposition looks as follows.
Consider a vector $\bm{x} \in \set{R}^{I}$ of the size $I=I_1 I_2 \ldots I_d$, where $I_\alpha \in \set{N}$ for $\alpha=1, 2, \ldots, d$.
We can treat it as a $d$-dimensional tensor\footnote{
	By \emph{tensors} we mean multidimensional arrays with a number of dimensions $d$.
	A two-dimensional tensor ($d=2$) is a matrix, and when $d=1$ it is a vector.
	For tensors with $d>2$ we use upper case calligraphic letters ($\tensor{A}, \tensor{B}, \tensor{C}, \ldots$).
	The $(i_1, i_2, \ldots, i_d)$th entry of a $d$-dimensional tensor $\tensor{X} \in \set{R}^{I_1 \times I_2 \times \ldots \times I_d}$ is denoted by $\tensor{X}[i_1, i_2, \ldots, i_d]$, where $i_\alpha = 0, 1, \ldots, I_\alpha-1$ ($\alpha=1, 2, \ldots, d$) and $I_\alpha$ is a size of the $\alpha$th mode.
	Mode-$\alpha$ slice of such tensor is denoted by $\tensor{X}[i_1, \ldots, i_{\alpha-1}, :, i_{\alpha+1}, \ldots, i_d]$, and an operation $\myvec(\cdot)$ constructs a vector $\myvec(\tensor{X})\in \set{R}^{I_1 I_2 \ldots I_d}$ from the given tensor $\tensor{X}$ by a reshaping procedure from the formula~\eqref{qtt_solver:multi_index}, where Fortran style is used for index ordering.
}
$\tensor{X} \in \set{R}^{I_1 \times I_2 \times \ldots \times I_d}$, such that $\myvec(\tensor{X}) = \bm{x}$, or equivalently
$$
	\tensor{X} [i_1, \ldots, i_d] = \bm{x} [i], \quad i_\alpha = 0, 1, \ldots, I_\alpha-1 \quad (\alpha=1, 2, \ldots, d),
$$
where
\begin{equation}\label{qtt_solver:multi_index}	
	i = i_1 + i_2 I_1 + \ldots + i_d I_1 I_2 \ldots I_{d-1}.
\end{equation}

After this transformation, that is also called \emph{quantization}, we represent the tensor $\tensor{X}$ in the low-rank TT-format~\cite{osel-tt-2011}.
A tensor $\tensor{X} \in \set{R}^{I_1 \times I_2 \times \ldots \times I_d}$ is said to be in the TT-format, if its elements are represented by the formula
\begin{equation}\label{qtt_solver:tt_tensor_sum_form}	
	\tensor{X} [i_1, i_2, \ldots, i_d] = \sum_{r_1=0}^{R_1-1} \ldots \sum_{r_{d-1}=0}^{R_{d-1}-1}
	\tensor{G}_1 [0, i_1, r_1] \tensor{G}_2 [r_1, i_2, r_2] \ldots \tensor{G}_d [r_{d-1}, i_d, 0],
\end{equation}
where $i_\alpha = 0, 1, \ldots, I_\alpha-1$ for $\alpha = 1, 2, \ldots, d$, $\tensor{G}_\alpha \in \set{R}^{R_{\alpha-1} \times I_\alpha \times R_\alpha}$ are three-dimensional tensors, that are named TT-cores, and integers $R_{0}, R_{1}, \ldots, R_{d}$ (with convention $R_{0}=R_{d}=1$) are named TT-ranks.
The last formula can be also rewritten in a more compact form
\begin{equation}\label{qtt_solver:tt_tensor_matrix_form}
	\tensor{X} [i_1, i_2, \ldots, i_d] = G_1(i_1) G_2(i_2) \ldots G_d(i_d),
\end{equation}
where $G_\alpha(i_\alpha) = \tensor{G}_\alpha [:, i_\alpha, :]$ is an $R_{\alpha-1} \times R_\alpha$ matrix for each fixed $i_\alpha$ (since $R_{0}=R_{d}=1$, the result of matrix multiplications in~\eqref{qtt_solver:tt_tensor_matrix_form} is a scalar).
A vector form of the TT-decomposition looks like
\begin{equation}\label{qtt_solver:tt_tensor_kron_form}
	\bm{x} = \myvec(\tensor{X}) = \sum_{r_1=0}^{R_1-1} \ldots \sum_{r_{d-1}=0}^{R_d-1}
	\tensor{G}_1 [0, :, r_1] \otimes \tensor{G}_2 [r_1, :, r_2] \otimes \ldots \otimes \tensor{G}_d [r_{d-1}, :, 0],
\end{equation}
where the slices of the TT-cores $\tensor{G}_\alpha [r_{\alpha-1}, :, r_{\alpha}]$ are vectors of length $I_\alpha$ for $\alpha = 1, 2, \ldots, d$.

The benefit of the TT-decomposition is the following.
Storage of the TT-cores $\tensor{G}_1, \tensor{G}_2, \ldots, \tensor{G}_d$ requires less or equal than
$$
	 d \times \max{I_\alpha} \times (\max{R_\alpha})^2
$$
memory cells (instead of $I_1 I_2 \ldots I_d$ cells for the uncompressed tensor), and hence the TT-decomposition is free from the curse of dimensionality\footnote{
	By \emph{the full format tensor representation} or \emph{uncompressed tensor} we mean the case, when one calculates and saves in the memory all tensor elements.
	The number of elements of an uncompressed tensor (hence, the memory required to store it) and the amount of operations required to perform basic operations with such tensor grows exponentially in the dimensionality, and this problem is called \emph{the curse of dimensionality}.
}
if the TT-ranks are bounded.
For representation of a vector in the QTT-format, the most convenient way is to set its size as $I=2^d$ ($d \in \set{N}$) and, correspondingly, $I_1=I_2=\ldots=I_d=2$, hence only less or equal than
$$
	 2 \log_2 I \times (\max{R_\alpha})^2
$$
memory cells are required for the storage of the vector in the QTT-format (instead of $I$ cells for the uncompressed representation).

For a given tensor $\tensor{\tilde{X}}$ in the full format, the TT-decomposition (compression) can be performed by a stable TT-SVD algorithm~\cite{osel-tt-2011}.
This algorithm constructs an approximation $\tensor{X}$ in the TT-format to the given tensor $\tensor{\tilde{X}}$ with a prescribed accuracy $\tau$ in the Frobenius norm\footnote{
	An exact TT-representation  exists for the given full tensor $\tensor{\tilde{X}}$, and TT-ranks of such representation are bounded by ranks of the corresponding unfolding matrices~\cite{osel-tt-2011}.
	Nevertheless, in practical applications it is more useful to construct TT-approximation with a prescribed accuracy $\tau$, and then carry out all operations (summations, products, etc) in the TT-format, maintaining the same accuracy $\tau$ of the result.
}
\begin{equation}\label{qtt_solver:tt_tensor_appr_acc}
	|| \tensor{X} - \tensor{\tilde{X}} ||_F \leq \tau || \tensor{\tilde{X}} ||_F,
\end{equation}
but a procedure of the tensor approximation in the full format is too costly, and is even impossible for large dimensions due to the curse of dimensionality.
A TT-cross method~\cite{ot-ttcross-2010} can be used instead.
This method is a generalization of the cross approximation method for matrices~\cite{tee-cross-2000} to higher dimensions.
TT-cross method constructs a TT-approximation $\tensor{X}$ to the tensor $\tensor{\tilde{X}}$, given as a function $f(i_1, i_2, \ldots, i_d)$, that returns the $(i_1, i_2, \ldots, i_d)$th entry of $\tensor{\tilde{X}}$ for a given set of indices.
In a more general form this method can be formulated as follows.
For a given function $f(\tensor{\tilde{A}}, \tensor{\tilde{B}}, \ldots)$, where arguments $\tensor{\tilde{A}}, \tensor{\tilde{B}}, \ldots \in \set{R}^{I_1 \times I_2 \times \ldots \times I_d}$ are given tensors of the equal shape, the TT-cross method constructs a TT-approximation $\tensor{X}$ to the tensor
$$
	\tensor{\tilde{X}} = f(\tensor{\tilde{A}}, \tensor{\tilde{B}}, \ldots) \in \set{R}^{I_1 \times I_2 \times \ldots \times I_d},
$$
using TT-approximations $\tensor{A}, \tensor{B}, \ldots$ to the tensors $\tensor{\tilde{A}}, \tensor{\tilde{B}}, \ldots$.
This method requires only
$$
	\mathcal{O} \left( d \times \max{I_\alpha} \times (\max{R_\alpha^{\tensor{X}}})^3 \right),
$$
operations for the construction of the approximation with a prescribed accuracy $\tau$.
	
\begin{figure}[t] \begin{center} \begin{algorithm}[H]
	\SetAlgoLined 	
	\KwData  {functions $f^{func}$, $k_x^{func}$ and $k_y^{func}$, grid factor $d$, accuracy $\tau$.}
	\KwResult{vectors $\qttvec{f}$, $\qttvec{k}_x^{-1}$ and $\qttvec{k}_y^{-1}$ of length $2^{2d}$  in the QTT-format.}

	$h = 2^{-d}$
	
	$\qttvec{e} = \myfunc{tt\_ones}(2, d)$
	
	$\qttvec{p} = \myfunc{tt\_xfun}(2, d)$

	$\qttvec{x}^{(c)}  = h \cdot \myfunc{tt\_round}( \qttvec{e} \otimes ( \qttvec{p} + 0.5 \qttvec{e} ), \tau)$

	$\qttvec{y}^{(c)}  = h \cdot \myfunc{tt\_round}( ( \qttvec{p} + 0.5 \qttvec{e})  \otimes \qttvec{e}, \tau)$

	$\qttvec{x}^{(r)}  = h \cdot \myfunc{tt\_round}( \qttvec{e} \otimes ( \qttvec{p} + \qttvec{e} ), \tau)$

	$\qttvec{y}^{(r)}  = h \cdot \myfunc{tt\_round}( ( \qttvec{p} + \qttvec{e} ) \otimes \qttvec{e}, \tau)$
  			
	$\qttvec{f}        = \myfunc{tt\_cross}(f^{func},     [\qttvec{x}^{(r)}, \qttvec{y}^{(r)}], \tau)$

	$\qttvec{k}_x^{-1} = \myfunc{tt\_cross}(1/k_x^{func}, [\qttvec{x}^{(c)}, \qttvec{y}^{(r)}], \tau)$
	
	$\qttvec{k}_y^{-1} = \myfunc{tt\_cross}(1/k_y^{func}, [\qttvec{x}^{(r)}, \qttvec{y}^{(c)}], \tau)$

	\caption{discretization of the PDE coefficients and the right-hand side on a spatial grid in the QTT-format.}
	\label{qtt_solver:alg:pde_construction}
\end{algorithm} \end{center} \end{figure}
			
Construction of the vectors $\bm{f}$, $\bm{k}_x^{-1}$ and $\bm{k}_y^{-1}$ (from formulas~\eqref{scheme:descr_f},~\eqref{scheme:descr_kx} and~\eqref{scheme:descr_ky} respectively) in the QTT-format\footnote{
	For vectors in the QTT-format hereinafter we use bold lower case calligraphic letters ($\qttvec{x}, \qttvec{y}, \qttvec{f}, \ldots$) to emphasize that in operations of linear algebra they play the same role as ordinary vectors.
	At the same time, for matrices in the QTT-format (see \cref{subsection:qtt_solver_matrices}) we will use upper case calligraphic letters.
}
is described in \cref{qtt_solver:alg:pde_construction}.
First we prepare the spatial grid from \cref{subsection:scheme_discr} in the QTT-format, then the TT-cross method (function $\myfunc{tt\_cross}$) is used for construction of the QTT-vectors on the spatial grid from functions $f^{func}$, $k_x^{func}$, $k_y^{func}$ that calculate the right-hand side $f$ and the coefficients $k_x$, $k_y$ of the model PDE~\eqref{scheme:model_pde} for a given point $(x, y)$.
We use\footnote{
	All basic operations in the QTT- and TT-format are implemented in the \emph{ttpy} package (see \cref{qtt_solver:tbl:tt_operations} with a list of the main used operations and functions with complexity estimations and TT-ranks of the result).
	This package is publicly available from \url{https://github.com/oseledets/ttpy} (python programming language) and from \url{https://github.com/oseledets/TT-Toolbox} (MATLAB version).
}
function $\myfunc{tt\_ones}(s, d)$ for the construction of rank-1 QTT-vector of ones of the length $s^d$, function $\myfunc{tt\_xfun}(s, d)$ for construction of rank-1 QTT-vector of the form $\bm{[0, 1, \ldots, s^d-1]}^T$,  and function $\myfunc{tt\_round}$ for the rounding of the given QTT-vector to the prescribed accuracy $\tau$. 

\subsection{Matrices in a low-rank format} \label{subsection:qtt_solver_matrices}
		
\begin{figure}[t] \begin{center} \begin{algorithm}[H]
	\SetAlgoLined 	
	\KwData  {vector $\qttvec{k}_x^{-1}$ of length $2^{2d}$ in the QTT-format, accuracy $\tau$.}
	\KwResult{an $2^{2d} \times 2^{2d}$  matrix $\qttmat{W_x}$ in the QTT-format from formula~\eqref{scheme:HRW_x_def}.}

	Let $G$ be a list of TT-cores of $\qttvec{k}_x^{-1}$ (it has $2d$ TT-cores)

	Introduce matrix $M = G[0][:, 0, :] + G[0][:, 1, :]$

	\For{$\alpha = 1, 2, \ldots, d-1$}{

		Introduce intermediate matrix $P = G[\alpha][:, 0, :] + G[\alpha][:, 1, :]$
 				
 		Calculate matrix-by-matrix product $M = M P$
	}

	Set $\tensor{G} = G[d]$ and calculate tensor $\tensor{\hat G} = M  \times_1 \tensor{G}$
		
	Construct QTT-vector $\qttvec{q}_x^{-1}$ from cores $\left( \tensor{\hat G}, G[d+1], \ldots, G[2d-1] \right)$

	$\qttvec{q}_x = \myfunc{tt\_cross}(f, \qttvec{q}_x^{-1}, \tau)$  // Function $f(t)$ is $1/t$ here
		
	$\qttmat{Q_x} = \myfunc{tt\_diag}(\qttvec{q}_x)$
        	
	$\qttmat{E} = \myfunc{tt\_ones\_mat}(2, d)$

	$\qttmat{W_x} = \myfunc{tt\_round}(\qttmat{Q_x} \otimes \qttmat{E}, \tau)$
		
	\caption{construction of the $W_x$ matrix in the QTT-format.}
	\label{qtt_solver:alg:constructWx}
\end{algorithm} \end{center} \end{figure}

\begin{figure}[t] \begin{center} \begin{algorithm}[H]
	\SetAlgoLined 	
	\KwData  {vector $\qttvec{k}_y^{-1}$ of length $2^{2d}$ in the QTT-format, accuracy $\tau$.}
	\KwResult{an $2^{2d} \times 2^{2d}$  matrix $\qttmat{W_y}$ in the QTT-format from formula~\eqref{scheme:HRW_y_def}.}

	Let $G$ be a list of TT-cores of $\qttvec{k}_y^{-1}$ (it has $2d$ TT-cores)

	Introduce matrix $M = G[2d-1][:, 0, :] + G[2d-1][:, 1, :]$

	\For{$\alpha = 2d-2, 2d-3, \ldots, d+1$}{

		Introduce intermediate matrix $P = G[\alpha][:, 0, :] + G[\alpha][:, 1, :]$
 				
 		Calculate matrix-by-matrix product $M = P M$
	}
		
	Set $\tensor{G} = G[d]$ and calculate tensor $\tensor{\hat G} = \tensor{G} \times_1 M$
		
	Construct QTT-vector $\qttvec{q}_y^{-1}$ from cores $\left( G[0], \ldots, G[d-1], \tensor{\hat G} \right)$

	$\qttvec{q}_y = \myfunc{tt\_cross}(f, \qttvec{q}_y^{-1}, \tau)$ // Function $f(t)$ is $1/t$ here

	$\qttmat{Q_y} = \myfunc{tt\_diag}(\qttvec{q}_y)$
        	
	$\qttmat{E} = \myfunc{tt\_ones\_mat}(2, d)$

	$\qttmat{W_y} = \myfunc{tt\_round}(\qttmat{E} \otimes \qttmat{Q_y}, \tau)$
		
	\caption{construction of the $W_y$ matrix in the QTT-format.}
	\label{qtt_solver:alg:constructWy}
\end{algorithm} \end{center} \end{figure}

The next step is to represent matrices.
Even though these matrices are in general of full matrix rank, they also admit low-rank QTT-approximation, and we will use it for the effective construction of the matrices~\eqref{scheme:HRW_x_def},~\eqref{scheme:HRW_y_def},~\eqref{scheme:Qx_explicit} and~\eqref{scheme:Qy_explicit}, that are involved in the formulation of the proposed discretization scheme.

The QTT-representation of matrices is defined as follows.
Consider a matrix $A \in \set{R}^{I \times J}$, where 
$$
	I=I_1 I_2 \ldots I_d, \quad J=J_1 J_2 \ldots J_d, \quad I_\alpha, J_\alpha \in \set{N} \quad (\alpha=1, 2, \ldots, d).
$$
We can treat this matrix as a $d$-dimensional tensor
$$
	\qttmat{A} \in \set{R}^{(I_1 J_1) \times (I_2 J_2) \times \ldots \times (I_d J_d)}, \quad
	\qttmat{A} [\overline{i_1 j_1}, \overline{i_2 j_2}, \ldots, \overline{i_d j_d}] = \myvec (A) [\overline{i j}],
$$
where
$$
	\overline{i_\alpha j_\alpha} = i_\alpha + j_\alpha I_\alpha, \quad
	i_\alpha = 0, 1, \ldots, I_\alpha-1, \quad 
	j_\alpha = 0, 1, \ldots, J_\alpha-1,
$$
for $\alpha=1,2,\ldots,d$, and multi-indices are constructed in the same manner, as in~\eqref{qtt_solver:multi_index} 
$$
	\overline{i j} = \overline{i_1 j_1} + \overline{i_2 j_2} (I_1 J_1) + \ldots + 
	                 \overline{i_d j_d} (I_1 J_1) (I_2 J_2) \ldots (I_{d-1} J_{d-1}).
$$
Then we can apply the TT-decomposition in the form~\eqref{qtt_solver:tt_tensor_matrix_form} to the tensor $\qttmat{A}$
$$
	\qttmat{A} [\overline{i_1 j_1}, \overline{i_2 j_2}, \ldots, \overline{i_d j_d}] = 
	\tilde{G_1}(\overline{i_1 j_1}) \tilde{G_2}(\overline{i_2 j_2}) \ldots \tilde{G_d}(\overline{i_d j_d}),
$$
where three-dimensional TT-cores $\tilde{\tensor{G}_\alpha} \in \set{R}^{R_{\alpha-1} \times I_\alpha J_\alpha \times R_{\alpha}}$ are represented as matrices $\tilde{G_\alpha}$ depending on a multi-index $\overline{i_\alpha j_\alpha}$ for $\alpha=1, 2, \ldots, d$.
If we reshape $\tilde{\tensor{G}_\alpha}$ to four-dimensional tensors $\tensor{G}_\alpha \in \set{R}^{R_{\alpha-1} \times I_\alpha \times J_\alpha \times R_\alpha}$, then the QTT-decomposition for the matrix $A$ can be written as
\begin{equation}\label{qtt_solver:tt_matrix_matrix_form}
	\qttmat{A} [i_1, i_2, \ldots, i_d; j_1, j_2, \ldots, j_d] = G_1(i_1, j_1) G_2(i_2, j_2) \ldots G_d(i_d, j_d),
\end{equation}
or in the form like~\eqref{qtt_solver:tt_tensor_sum_form}
\begin{equation}\label{qtt_solver:tt_matrix_sum_form}
	\begin{split}
		\qttmat{A} & [i_1, i_2, \ldots, i_d; j_1, j_2, \ldots, j_d] = \\
		           & \sum_{r_1=0}^{R_1-1} \ldots \sum_{r_{d-1}=0}^{R_{d-1}-1}
		             \tensor{G} [0, i_1, j_1, r_1] 
		             \tensor{G} [r_1, i_2, j_2, r_2] \ldots 
		             \tensor{G} [r_{d-1}, i_d, j_d, 0].
	\end{split}
\end{equation}

The QTT-representation~\eqref{qtt_solver:tt_matrix_matrix_form} or~\eqref{qtt_solver:tt_matrix_sum_form} for the matrices\footnote{
	We use special notation $[i_1, i_2, \ldots, i_d; j_1, j_2, \ldots, j_d]$ for the indices in~\eqref{qtt_solver:tt_matrix_matrix_form} and~\eqref{qtt_solver:tt_matrix_sum_form} to emphasize that $\qttmat{A}$ acts as a matrix with rows given by the multi-index $(i_1, i_2, \ldots, i_d)$ and columns given by $(j_1, j_2, \ldots, j_d)$.
} 
makes it possible (see~\cite{osel-tt-2011, osel-2d2d-2010, ds-amen-2014}) to formulate standard linear algebra operations, like matrix-by-vector and matrix-by-matrix product, and to solve linear systems in the QTT/TT-format.
For example, matrix-by-vector multiplication $\qttvec{x} = \qttmat{A} \qttvec{y}$ can be efficiently implemented in the QTT/TT-format.
It has complexity 
$$
	\mathcal{O} \left( d \times (\max{(I_\alpha, J_\alpha)})^2 \times 
	(\max{R^{\qttmat{A}}_\alpha})^2 \times (\max{R^{\qttvec{y}}_\alpha})^2 \right),
$$
and TT-ranks of the resulting QTT-vector $\qttvec{x}$ are equal to the product of the ranks of the QTT-matrix $\qttmat{A}$ and the QTT-vector $\qttvec{y}$: $R^{\qttvec{x}}_\alpha = R^{\qttmat{A}}_\alpha R^{\qttvec{y}}_\alpha$ for $\alpha=1,2,\ldots,d$.
It should be noted, that in many cases ranks of the QTT-vector $\qttvec{x}$ are overestimated, and to avoid rank growth one has to reduce the ranks of the product, while maintaining accuracy.
Robust TT-round algorithm~\cite{osel-tt-2011} is available for this purpose.
It has complexity
$$
	\mathcal{O} \left( d \times \max{I_\alpha} \times (\max{R^{\qttvec{x}}_\alpha})^3 \right)
$$
and should be used after such operations in the QTT/TT-format that leads to the growing of the ranks.

In \cref{qtt_solver:alg:constructWx} and \cref{qtt_solver:alg:constructWy} we present a pseudo code for the construction of the matrices $\qttmat{W_x}$ and $\qttmat{W_y}$ in the QTT-format from the formulas~\eqref{scheme:HRW_x_def} and~\eqref{scheme:HRW_y_def} respectively.
Function $\myfunc{tt\_diag}(\cdot)$ constructs a diagonal matrix in the QTT-format from the given QTT-vector, function $\myfunc{tt\_ones\_mat}(s, d)$ is used for construction of an $s^d \times s^d$ QTT-matrix of all ones, and $\myfunc{tt\_round}$ function is used for rounding of the given QTT-matrix to the prescribed accuracy.
Inversions of vectors $\qttvec{q}_x$ and $\qttvec{q}_y$ from~\eqref{scheme:Qx_explicit} and~\eqref{scheme:Qy_explicit} are performed by the TT-cross method (function $\myfunc{tt\_cross}$).
Notation "$\times_1$" in \cref{qtt_solver:alg:constructWx} and \cref{qtt_solver:alg:constructWy}  means tensor-by-matrix multiplication.
Given a tensor $\tensor{X} \in \set{R}^{I_1 \times I_2 \times \ldots I_d}$ and matrices $A \in \set{R}^{I_d \times I_A}$ and $B \in \set{R}^{I_B \times I_1}$, we define right-tensor-by-matrix multiplication operation as follows
$$
	\tensor{Y} = \tensor{X} \times_1 A, \quad
	\tensor{Y} \in \set{R}^{I_1 \times \ldots I_{d-1} \times I_A},
$$
$$
	\tensor{Y}[i_1, \ldots, i_{d-1}, \alpha] = 
	\sum_{i_d=0}^{I_d-1} \tensor{X}[i_1, \ldots, i_{d-1}, i_d] A[i_d, \alpha],
	\quad \alpha = 0, 1, \ldots, I_A-1,
$$
and left-tensor-by-matrix multiplication
$$
	\tensor{Z} = B \times_1 \tensor{X}, \quad
	\tensor{Z} \in \set{R}^{I_B \times I_2 \times \ldots \times I_d},
$$
$$
	\tensor{Z}[\alpha, i_2 \ldots, i_d] =
	\sum_{i_1=0}^{I_1-1} B[\alpha, i_1] \tensor{X}[i_1, i_2, \ldots, i_d], 
	\quad \alpha = 0, 1, \ldots, I_B-1.
$$

\subsection{Linear system construction and solution}

\begin{figure}[t] \begin{center} \begin{algorithm}[H]
	\SetAlgoLined 	
	\KwData  {functions $k_x^{func}$, $k_y^{func}$, $f^{func}$, grid factor $d$, accuracy $\tau$.}
	\KwResult{approximate solution $\qttvec{u}$ of the PDE~\eqref{scheme:model_pde} and its derivatives $\qttvec{u}_x$, $\qttvec{u}_y$ in the QTT-format.}

	Construct QTT-vectors $\qttvec{f}$, $\qttvec{k}_x^{-1}$ and $\qttvec{k}_y^{-1}$ according to \cref{qtt_solver:alg:pde_construction}
  	
  	Construct QTT-matrix $\qttmat{W_x}$ according to \cref{qtt_solver:alg:constructWx}

	Construct QTT-matrix $\qttmat{W_y}$ according to \cref{qtt_solver:alg:constructWy}

	$\qttmat{I} = \myfunc{tt\_eye}(2, d), \quad \qttmat{I}_2   = \myfunc{tt\_eye}(2, 2d)$
	
	$\qttmat{B} = \myfunc{tt\_volterra}(2, d)$
	
	$\qttmat{B_x} = \qttmat{I} \otimes \qttmat{B}$

	$\qttmat{B_y} = \qttmat{B} \otimes \qttmat{I}$
	
	$\qttmat{K_x^{-1}} = \myfunc{tt\_diag}(\qttvec{k}_x^{-1})$

	$\qttmat{K_y^{-1}} = \myfunc{tt\_diag}(\qttvec{k}_y^{-1})$

	$\qttmat{R_x} = \qttmat{K_x^{-1}} \left( \qttmat{I}_2 - \qttmat{W_x} \qttmat{K_x^{-1}} \right) \qttmat{B_x}^T$

	$\qttmat{R_y} = \qttmat{K_y^{-1}} \left( \qttmat{I}_2 - \qttmat{W_y} \qttmat{K_y^{-1}} \right) \qttmat{B_y}^T$
	
	$\qttmat{H_x} = \qttmat{B_x} \qttmat{R_x}$

	$\qttmat{H_y} = \qttmat{B_y} \qttmat{R_y}$

	$\qttvec{\mu} = \myfunc{tt\_amen}(\qttmat{H_x} + \qttmat{H_y}, \qttmat{H_y} \qttvec{f}, \tau)$

	$\qttvec{u}_x = \qttmat{R_x} \qttvec{\mu}$,
	$\quad$
	$\qttvec{u}_y = \qttmat{R_y} ( \qttvec{f}-\qttvec{\mu} )$,
	$\quad$
	$\qttvec{u}   = \qttmat{B_x} \qttvec{u}_x$
	
	\caption{FS-QTT-solver.}
	\label{qtt_solver:alg:fs_qtt}
\end{algorithm} \end{center} \end{figure}

A practical implementation\footnote{
	The proposed discretization scheme is implemented as a \emph{qttpdesolver} python package.
	The code is publicly available from \url{https://github.com/AndreChertkov/qttpdesolver}.
}
of the new descritization scheme in the QTT-format is described in \cref{qtt_solver:alg:fs_qtt}. 
The main steps of the computational process were formulated in Theorem~\ref{scheme:th:result} from \cref{section:scheme}.
To obtain the solution $u$ of the model PDE~\eqref{scheme:model_pde} and its derivatives $u_x$, $u_y$ in the QTT-format, we have to select a grid factor $d$, construct QTT-vectors $\qttvec{f}$, $\qttvec{k}_x^{-1}$ and $\qttvec{k}_y^{-1}$, that are the discrete versions of the right-hand side and inverses to the PDE coefficients $k_x$ and $k_y$ respectively.
Then we calculate QTT-matrices $\qttmat{W_x}$ and $\qttmat{W_y}$ from~\eqref{scheme:HRW_x_def} and~\eqref{scheme:HRW_y_def}, using $\qttmat{Q_x}$ and $\qttmat{Q_y}$ QTT-matrices from~\eqref{scheme:Qx_explicit} and~\eqref{scheme:Qy_explicit}.
After that, we construct QTT-matrix $\qttmat{B}$ from~\eqref{scheme:B} by a function $\myfunc{tt\_volterra}$ and QTT-matrices $\qttmat{B_x}$, $\qttmat{B_y}$ by Kronecker products in the QTT-format (function $\myfunc{tt\_eye}(s, d)$ is used for construction of an $s^d \times s^d$ diagonal QTT-matrix of ones), and then QTT-matrices $\qttmat{R_x}$, $\qttmat{H_x}$, $\qttmat{R_y}$ and $\qttmat{H_y}$ from~\eqref{scheme:HRW_x_def} and~\eqref{scheme:HRW_y_def} are calculated.

With the matrix-by-vector operation the problem of solving linear systems in the QTT-format can be formulated.
For a given square matrix $A \in \set{R}^{I \times I}$ and a vector $\bm{f} \in \set{R}^{I}$, that are given in the QTT-format as a QTT-matrix $\qttmat{A}$ and a QTT-vector $\qttvec{f}$ respectively, one has to find a QTT-vector $\qttvec{x}$, that is an approximation of solution $\bm{x} \in \set{R}^{I}$ of a linear system $A \bm{x} = \bm{f}$.
Efficient iterative solver with step-complexity
$$	
	\mathcal{O} \left(
	\max{I_\alpha}     \times (\max{R^\qttmat{A}_\alpha})   (\max{R^\qttvec{x}_\alpha})^3 +
	(\max{I_\alpha})^2 \times (\max{R^\qttmat{A}_\alpha})^2 (\max{R^\qttvec{x}_\alpha})^2 \right),
$$ 
named AMEn-solver, exists~\cite{ds-amen-2014} for such class of problems.
We use this solver (function $\myfunc{tt\_amen}$ in \cref{qtt_solver:alg:fs_qtt}) for approximation of the solution of the linear system~\eqref{scheme:eq_for_mu} in the QTT-format.
And, finally, we construct PDE solution and derivatives from explicit formulas~\eqref{scheme:result_u} and~\eqref{scheme:result_ux_uy} respectively.
It should be noted, that after each operation in the QTT/TT-format, TT-round procedure with accuracy $\tau$ should be performed (it is omitted in \cref{qtt_solver:alg:fs_qtt} for the compactness of the presentation).

\subsection{Ranks and complexity estimation}
		
\begin{figure}[t] \begin{center}
	\captionof{table}{basic operations in the QTT/TT-format.}
	\label{qtt_solver:tbl:tt_operations} 
	\begin{tabular}{ l l l } \toprule[1.5pt]
		\bf Operation & \bf Maximum TT-rank & \bf Complexity 
		\\ $\tensor{C} = a + \tensor{B}$
				&  $R_\tensor{C} \leq R_\tensor{B} + 1$
				&  -
		\\ $\tensor{C} = \tensor{A}+\tensor{B}$
				&  $R_\tensor{C} \leq R_\tensor{A}+R_\tensor{B}$
				&  -
		\\ $\tensor{C} = \tensor{A} \otimes \tensor{B}$
				&  $R_\tensor{C} \leq \max(R_\tensor{A}, R_\tensor{B})$   
				&  -
		\\ $\tensor{C} = a \tensor{B}$
				&  $R_\tensor{C} = R_\tensor{B}$
				&  $\overline{I} R_\tensor{C}$
		\\ $\qttvec{x} = \qttmat{A} \qttvec{y}$
				&  $R_\qttvec{x} \leq R_\qttmat{A} R_\qttvec{y}$  
				&  $d \overline{I}^2 R_\qttmat{A}^2 R_\qttvec{y}^2$
		\\ $\tensor{C} = \tensor{A} \tensor{B}$
				&  $R_\tensor{C} \leq R_\tensor{A} R_\tensor{B}$  
				&  $d \overline{I}^3 R_\qttmat{A}^2 R_\qttmat{B}^2$
		\\ $\tensor{C} = \myfunc{tt\_round}  (\tensor{A}, \tau)$
				&  $R_\tensor{C} \leq R_\tensor{A}$
				&  $ d \overline{I} R_\tensor{A}^3$
		\\ $\qttvec{x} = \myfunc{tt\_cross}  (\qttvec{y}, \tau)$
				&  $R_\qttvec{x}$
				&  $ d \overline{I} R_\qttvec{x}^3$
		\\ $\qttvec{x} = \myfunc{tt\_amen}  (\qttmat{A}, \qttvec{y}, \tau)$
				&  $R_\qttvec{x}$
				&  $ d \overline{I}   R_\qttmat{A}   R_\qttvec{x}^3 + 
				     d \overline{I}^2 R_\qttmat{A}^2 R_\qttvec{x}^2$  
		\\ \bottomrule[1.25pt]
	\end {tabular}
\end{center} \end{figure}

Linear algebra operations in the QTT/TT-format are implemented in linear in $d$, polynomial in $R_\tensor{\tensor{X}}=\max_{\alpha=1,2,\ldots,d-1}{R^{\tensor{X}}_\alpha}$ and $\overline{I}=\max_{\alpha=1,2,\ldots,d}{I_\alpha}$  complexity with the result $\tensor{X} \in \set{R}^{I_1 \times I_2 \times \ldots \times I_d}$ also in the QTT/TT-format (see~\cite{osel-tt-2011} for more details).
In \cref{qtt_solver:tbl:tt_operations} we present main operations with estimates for maximum TT-rank and complexity\footnote{
	To construct a sum of two TT-tensors $\tensor{C} = \tensor{A}+\tensor{B}$ (or a sum of a constant and a TT-tensor $\tensor{C} = a+\tensor{B}$), we only need to put each TT-core of $\tensor{A}$ and $\tensor{B}$ into diagonal of the corresponding TT-core of $\tensor{C}$, and for the case of Kronecker product $\tensor{C} = \tensor{A} \otimes \tensor{B}$ only concatenation of TT-cores is performed, hence these operations have, formally, the zero complexity.
	Nevertheless, as was mentioned above, the TT-rounding procedure with complexity $ d \overline{I} R_\tensor{C}^3$ should be done after such operations to avoid rank growth.
}
. Using these estimates, we can derive the rank bounds for the matrix and the right-hand side of the linear system~\eqref{scheme:eq_for_mu}.

\begin{theorem} \label{qtt_solver:th:hxy_ranks}
	The QTT-matrix $\qttmat{H_x}+\qttmat{H_y}$ in the linear system~\eqref{scheme:eq_for_mu} has TT-ranks, that are bounded by the value
	\begin{equation}\label{qtt_solver:hx_plus_hy_rank}
		4 \left( 1 + r_{1/k,x} r_{q,x} \right) r_{1/k,x} +
		4 \left( 1 + r_{1/k,y} r_{q,y} \right) r_{1/k,y},
	\end{equation}
	where
	$$
		r_{1/k,x}=\max_{\alpha=1,2,\ldots,2d-1}{R^{\qttmat{K_x^{-1}}}_\alpha}, \quad 
		r_{1/k,y}=\max_{\alpha=1,2,\ldots,2d-1}{R^{\qttmat{K_y^{-1}}}_\alpha},
	$$
	are the maximum TT-ranks of the inverse of discretized coefficients $k_x$ and $k_y$ of the model PDE~\eqref{scheme:model_pde},
	$$
		r_{q,x}=\max_{\alpha=1,2,\ldots,d-1}{R^{\qttvec{q}_x}_\alpha}, \quad
		r_{q,y}=\max_{\alpha=1,2,\ldots,d-1}{R^{\qttvec{q}_y}_\alpha},
	$$
	are the maximum TT-ranks of the vectors $\bm{q}_x$ and $\bm{q}_y$ from~\eqref{scheme:Qx_explicit} and~\eqref{scheme:Qy_explicit} respectively\footnote{
		TT-ranks of the QTT-vectors $\qttvec{q}_x$ and $\qttvec{q}_y$ can be expressed in terms of the TT-ranks of discretized PDE coefficients $k_x$ and $k_y$ under some additional restrictions on their smoothness, using an approach similar to the one described in~\cite{ks-2dqtt-2015}, but this work will be reported elsewhere.
	}
	. The right-hand side $\qttmat{H_y} \qttvec{f}$ in the linear system~\eqref{scheme:eq_for_mu} has TT-ranks, that are bounded by the value
	\begin{equation}\label{qtt_solver:rhs_rank}
		4 \left( 1 + r_{1/k,y} r_{q,y} \right) r_{1/k,y} r_{f},
	\end{equation}
	where
	$$
		r_{f}=\max_{\alpha=1,2,\ldots,2d-1}{R^{\qttvec{f}}_\alpha},
	$$
	is the maximum TT-rank of the discretized right-hand side of the model PDE.
\end{theorem}
			
\begin{proof}
	As presented in \cref{qtt_solver:alg:fs_qtt}, to construct matrix $R_x$ in the QTT-format, we use formula~\eqref{scheme:HRW_x_def}, where all operations are performed in the QTT-format.
	Given that $\textit{r}(\qttmat{B})=2$, $\textit{r}(\qttmat{I})=1$, $\textit{r}(\qttmat{I}_2)=1$, and with rank estimates from \cref{qtt_solver:tbl:tt_operations}, we have
	$$
		\textit{r}(\qttmat{R_x}) \leq 
		2 \left( 1 + r_{1/k,x} \textit{r}(\qttmat{W_x}) \right) r_{1/k,x}, 
	$$
	where $\textit{r}(\qttmat{W_x})$ is the maximum TT-rank of the QTT-matrix $\qttmat{W_x}$.
	Due to \cref{qtt_solver:alg:constructWx} we can conclude that $\textit{r}(\qttmat{W_x}) \leq r_{q,x}$, since Kronecker product does not increase TT-ranks, then
	$$
		\textit{r}(\qttmat{R_x}) \leq 
		2 \left( 1 + r_{1/k,x} r_{q,x} \right) r_{1/k,x}, 
	$$
	and since $\qttmat{H_x} = \qttmat{B_x} \qttmat{R_x}$, we obtain
	\begin{equation}\label{qtt_solver:hx_rank}
		\textit{r}(\qttmat{H_x}) \leq 
		4 \left( 1 + r_{1/k,x} r_{q,x} \right) r_{1/k,x}.
	\end{equation}
	Using the same idea, we can obtain a similar estimate for $\qttmat{H_y}$
	\begin{equation}\label{qtt_solver:hy_rank}
		\textit{r}(\qttmat{H_y}) \leq 
		4 \left( 1 + r_{1/k,y} r_{q,y} \right) r_{1/k,y}.
	\end{equation}
	Summing~\eqref{qtt_solver:hx_rank} and~\eqref{qtt_solver:hy_rank}, we immediately obtain~\eqref{qtt_solver:hx_plus_hy_rank}.
	Using~\eqref{qtt_solver:hy_rank} and the estimate of the maximum TT-rank for the matrix-by-vector product  from \cref{qtt_solver:tbl:tt_operations}, we obtain~\eqref{qtt_solver:rhs_rank}.
\end{proof}
			
\section{Numerical exaples} \label{section:numex}
			
In this section we illustrate the theoretical results presented above with numerical experiments.
We compare three different solvers:
\begin{itemize}
	\item FS-QTT-solver, that is based on the new scheme from \cref{qtt_solver:alg:fs_qtt},
	\item FD-solver, that is based on the finite difference discretization scheme in the standard sparse format, which is described in \cref{section:connection},
	\item FD-QTT-solver, that is the QTT-version of FD-Solver.
\end{itemize}
		
\subsection{PDE with known analytic solution}

\begin{figure}[t]
	\begin{center}
		\captionof{table}{Effective TT-ranks of the main vectors and matrices from FS-QTT-solver, applied to the model PDE with known analytic solution for different grid factors $d$.}
		\label{numex:tbl:res_analyt_eranks}
		\begin{tabular}{ | l | l | l | l | l | l | l | l | l | l | l | l | l | l |}
			\hline
			$d$ & $\qttmat{K_x^{-1}}$ & $\qttmat{K_y^{-1}}$ & $\qttvec{q_x}$ & $\qttvec{q_y}$ & $\qttmat{H_x}$ & $\qttmat{H_y}$ & $\qttmat{A}$ & $\qttvec{u}$ 
		\\ \hline\bf 5  & 7.2   & 7.2   & 3.1   & 3.1   & 16.6  & 20.8  & 30.1  & 7.1   
			\\ \bf 10 & 8.9   & 8.9   & 4.3   & 3.9   & 18.5  & 37.8  & 43.5  & 6.7   
			\\ \bf 15 & 9.2   & 9.2   & 4.0   & 3.6   & 16.6  & 36.6  & 41.5  & 4.5   
			\\ \bf 20 & 9.1   & 9.1   & 3.7   & 3.3   & 16.1  & 34.1  & 39.3  & 6.4   
			\\ \bf 25 & 9.0   & 9.0   & 3.4   & 3.0   & 15.5  & 32.3  & 37.7  & 9.7   
			\\ \bf 30 & 8.8   & 8.8   & 3.2   & 2.9   & 15.1  & 30.5  & 36.2  & 12.2  
			\\ \hline
		\end{tabular}
	\end{center}
\end{figure}

\begin{figure}[t] \begin{center}
	\resizebox{0.3\textwidth}{!}{
		\input{res_analyt_all_time.pgf}
	}
	\resizebox{0.29\textwidth}{!}{
		\input{res_analyt_all_u_calc_erank.pgf}
	}
	\resizebox{0.3\textwidth}{!}{
		\input{res_analyt_all_uf_conv.pgf}
	}
	\caption{
		Total computational time (on the left plot), effective TT-rank of the calculated solution $u$ (on the middle plot) and $(u, f)-(u, f)_{RE}$ value (on the right plot), where $(u, f)_{RE}$ value is obtained from Richardson extrapolation, w.r.t. the mesh size factor $d$ for the model PDE with known analytic solution.
		Results are presented for FS-QTT-solver (blue line with circle marker), FD-QTT-solver (green line with square marker) and FD-solver (red line with triangle marker).
	}
	\label{figure:res_analyt_all_time_u_calc_erank}
\end{center} \end{figure}

\begin{figure}[t] \begin{center}
	\resizebox{0.45\textwidth}{!}{
		\input{res_analyt_all_u_err.pgf}
	}
	\resizebox{0.45\textwidth}{!}{
		\input{res_analyt_all_u_der_err.pgf}
	}
	\caption{
		Error of the calculated solution $u$ (on the left plot) and errors of the calculated $x$- and $y$-derivatives (on the right plot) w.r.t. the mesh size factor $d$ for the model PDE with known analytic solution.
		Results are presented for FS-QTT-solver (blue line with circle marker), FD-QTT-solver (green line with square marker) and FD-solver (red line with triangle marker). For the $y$-derivatives we use dotted lines with the same color and marker shape.
	}
	\label{figure:res_analyt_all_u_and_der_err}
\end{center} \end{figure}

Firstly, we consider a PDE with known analytic solution and homogeneous Dirichlet boundary conditions
\begin{equation}\label{numex:pde_analyt}
	- \nabla \cdot \left( k (x, y)  \nabla u (x, y) \right) = f(x, y), \quad 
	(x, y) \in \Omega = [0, 1]^2, \quad 
	u |_{\partial \Omega} = 0,
\end{equation}
with a scalar coefficient $k$
\begin{equation}\label{numex:pde_analyt_k}
	k(x, y) = 1 + x y^2,
\end{equation}
and the right-hand side
\begin{equation}\label{numex:pde_analyt_f}
\begin{split}
	&f(x, y) = (4 w_1^2 x^2 + w_2^2)(1+xy^2) \sin(w_1 x^2) \sin(w_2 y) - \\
	          & 2 w_1 (1 + 2xy^2) \cos(w_1 x^2) \sin(w_2 y) - 2 w_2 xy \sin(w_1 x^2) \cos(w_2 y).
\end{split}
\end{equation}
It can be shown, that the problem~\eqref{numex:pde_analyt},~\eqref{numex:pde_analyt_k},~\eqref{numex:pde_analyt_f} has exact analytic solution
\begin{equation}\label{numex:pde_analyt_u_real}
	u(x, y) = \sin(w_1 x^2) \sin(w_2 y).
\end{equation}

We select $w_1=\pi$, $w_2=2\pi$ in~\eqref{numex:pde_analyt_f} and~\eqref{numex:pde_analyt_u_real}, and perform calculations for different grid factor values: 
\begin{itemize}
	\item $d=4, 5, \ldots 30$ for FS-QTT-Solver,
	\item $d=4, 5, \ldots 15$ for FD-QTT-Solver,
	\item $d=4, 5, \ldots 10$ for FD-Solver.
\end{itemize}

Accuracy of AMEn-solver (function $\myfunc{tt\_amen}$ in \cref{qtt_solver:alg:fs_qtt}) was selected as $10^{-10}$, at the same time, for TT-round and TT-cross operations we use a higher accuracy, that is equal to $10^{-12}$ (since it provides better stability).
For error estimation we construct the analytic solution and its derivatives in the QTT-format with accuracy $10^{-14}$.
As a measure of solution (and derivatives) approximation error we use the the following quantity
$$
\epsilon = \frac{|| \qttvec{u}_\text{real} - \qttvec{u}_\text{calc} ||_F }{|| \qttvec{u}_\text{real} ||_F },
$$
for FS-QTT- and FD-QTT-solver, where $\qttvec{u}_\text{real}$ and $\qttvec{u}_\text{calc}$ are discretized analytic solution and calculated approximation in the QTT-format respectively, and $||\cdot||_F$ is Frobenius norm of the QTT-vector (for FD-solver we operate with vectors in the full format).
		
Total computational times for all solvers and effective TT-ranks\footnote{
	\emph{Effective TT-rank} $\hat{R}$ of a TT-tensor $\tensor{X} \in \set{R}^{I_1 \times I_2 \times \ldots \times I_d}$ with TT-ranks $R_0, R_1, \ldots, R_d$ ($R_0=R_d=1$) is a solution of quadratic equation
$$
	I_1 \hat{R} + \sum_{\alpha=2}^{d-1} I_\alpha \hat{R}^2 + I_d \hat{R} = \sum_{\alpha=1}^{d} I_\alpha R_{\alpha-1} R_{\alpha}.
$$
The representation with a constant TT-rank $\hat{R}$ ($\hat{R}_0=1$, $\hat{R}_1=\hat{R}_2 \ldots =\hat{R}_{d-1}=\hat{R}$,  $\hat{R}_d=1$) yields the same total number of parameters as in the original decomposition of the tensor $\tensor{X}$.
}
of obtained solutions for QTT-based solvers are shown in \cref{figure:res_analyt_all_time_u_calc_erank}.
Complexity (and hence, computational time) for FD-solver depends  exponentially on the grid factor $d$, hence this solver can efficiently operate only on moderate $d$ values.
For FD-QTT-solver we have a fast rank growth due to rounding errors, hence this solver can be also applied only for moderate grids.
At the same time, new FS-QTT-solver has almost linear dependence of computation time w.r.t. $d$ and bounded TT-ranks of solution even for huge grids (for example, for grid factor $d=30$, that means $2^{60}$ for the total number of grid cells).
Effective TT-ranks of discretized inverse to PDE coefficient $k$ (QTT-matrices $\qttmat{K_x^{-1}}$ and $\qttmat{K_y^{-1}}$), of intermediate QTT-vectors $\qttvec{q}_x$ and $\qttvec{q}_y$, of QTT-matrices $\qttmat{H_x}$, $\qttmat{H_y}$, of their sum $\qttmat{A}$ and of obtained solution $\qttvec{u}$ are presented in \cref{numex:tbl:res_analyt_eranks}.
As we can see from the table, all ranks are bounded and only slightly depend on the grid factor $d$. 
			
All solvers have almost the same accuracy for moderate grids according to \cref{figure:res_analyt_all_u_and_der_err}, where relative errors of solution and its derivatives are presented.
But for the grid factor $d > 12$ FD-QTT-solver becomes unstable, and the error grows.
FS-QTT-solver keeps second order accuracy for PDE solution, until $d=18$, where accuracy reaches the value of the selected accuracy of AMEn-solver ($10^{-10}$).
The same conclusion can be made from the analysis of the energy functional $(u, f)$, that is a scalar product of the calculated solution and the right-hand side of the PDE.
On the right plot in \cref{figure:res_analyt_all_time_u_calc_erank} we present dependence of the value $(u, f) - (u, f)_{RE}$ on the grid factor $d$, where $(u, f)_{RE}$ is Richardson extrapolation under assumption
$$
	(u, f)_{exact} = (u, f) + C h^{2} + \mathcal{O}(h^3),
$$
where $h=2^{-d}$ is a grid step, and, as can be seen from the figure, the value of $(u, f) - (u, f)_{RE}$ for FS-QTT-solver tends to zero with the second order convergence until $d=18$.

\subsection{PDE with constant right-hand side}
			
\begin{figure}[t]
	\begin{center}
		\captionof{table}{Effective TT-ranks of the main vectors and matrices from FS-QTT-solver, applied to the model PDE with right hand side, that is equal to one for different grid factors $d$.}
		\label{numex:tbl:res_rhs1_eranks}
		\begin{tabular}{ | l | l | l | l | l | l | l | l | l | l | l | l | l | l |}
			\hline
			$d$ & $\qttmat{K_x^{-1}}$ & $\qttmat{K_y^{-1}}$ & $\qttvec{q_x}$ & $\qttvec{q_y}$ & $\qttmat{H_x}$ & $\qttmat{H_y}$ & $\qttmat{A}$ & $\qttvec{u}$ 
		\\ \hline\bf 5  & 5.8   & 5.8   & 3.1   & 2.8   & 11.6  & 15.2  & 21.1  & 7.9   
			\\ \bf 10 & 6.4   & 6.4   & 3.4   & 3.0   & 11.8  & 21.4  & 25.6  & 11.0  
			\\ \bf 15 & 6.3   & 6.3   & 3.0   & 2.7   & 10.8  & 19.7  & 24.1  & 12.9  
			\\ \bf 20 & 6.2   & 6.2   & 2.8   & 2.5   & 10.6  & 18.5  & 23.0  & 12.8  
			\\ \bf 25 & 6.1   & 6.1   & 2.6   & 2.3   & 9.6   & 17.2  & 21.5  & 11.1  
			\\ \bf 30 & 5.9   & 5.9   & 2.4   & 2.2   & 9.5   & 16.1  & 20.7  & 10.9  
			\\ \hline
		\end{tabular}
	\end{center}
\end{figure}

\begin{figure}[t] \begin{center}
	\resizebox{0.3\textwidth}{!}{
		\input{res_rhs1_all_time.pgf}
	}
	\resizebox{0.29\textwidth}{!}{
		\input{res_rhs1_all_u_calc_erank.pgf}
	}
	\resizebox{0.3\textwidth}{!}{
		\input{res_rhs1_all_uf_conv.pgf}
	}
	\caption{
		Total computational time (on the left plot), effective TT-rank of the calculated solution $u$ (on the middle) and $(u, f)-(u, f)_{RE}$ value, where $(u, f)_{RE}$ value (on the right plot) is obtained from Richardson extrapolation, w.r.t. the mesh size factor $d$ for the model PDE with right-hand side, that is equal to one.
		Results are presented for FS-QTT-solver (blue line with circle marker), FD-QTT-solver (green line with square marker) and FD-solver (red line with triangle marker).
	}
	\label{figure:res_rhs1_all_time_u_calc_erank}
\end{center} \end{figure}
			
Next, we consider one more example, that is similar to the problem~\eqref{numex:pde_analyt},~\eqref{numex:pde_analyt_k}, but instead of~\eqref{numex:pde_analyt_f}, we select the case of the uniform constant source
\begin{equation}\label{numex:pde_rhs1_f}
    f(x, y) = 1, \quad \mbox{for all} \quad (x,y) \in \Omega.
\end{equation}
			
We select $10^{-6}$ as accuracy of AMEn-solver, and $10^{-8}$ as accuracy for TT-round and TT-cross operations and perform calculations for the following grid factor values: 
\begin{itemize}
	\item $d=4, 5, \ldots 30$ for FS-QTT-Solver,
	\item $d=4, 5, \ldots 15$ for FD-QTT-Solver,
	\item $d=4, 5, \ldots 10$ for FD-Solver.
\end{itemize}		
			
\cref{figure:res_rhs1_all_time_u_calc_erank} and \cref{numex:tbl:res_rhs1_eranks} represent computational results for three solvers, that were described above.
FD-QTT-Solver, as in the previous example, becomes unstable for $d>12$, at the same time, FS-QTT-Solver keeps second order accuracy until $d = 15$, and the approximate solution has bounded TT-ranks even for the case $d = 30$.

\subsection{PDE with point source}

\begin{figure}[t]
	\begin{center}
		\captionof{table}{Effective TT-ranks of the main vectors and matrices from FS-QTT-solver, applied to the model PDE with four point sources for different grid factors $d$.}
		\label{numex:tbl:res_ps4_eranks}
		\begin{tabular}{ | l | l | l | l | l | l | l | l | l | l | l | l | l | l |}
			\hline
			$d$ & $\qttmat{K_x^{-1}}$ & $\qttmat{K_y^{-1}}$ & $\qttvec{q_x}$ & $\qttvec{q_y}$ & $\qttmat{H_x}$ & $\qttmat{H_y}$ & $\qttmat{A}$ & $\qttvec{u}$ 
		\\ \hline\bf 5  & 5.8   & 5.8   & 3.1   & 2.8   & 11.6  & 15.2  & 21.1  & 11.1  
			\\ \bf 10 & 6.4   & 6.4   & 3.4   & 3.0   & 11.8  & 21.4  & 25.6  & 24.3  
			\\ \bf 15 & 6.3   & 6.3   & 3.0   & 2.7   & 10.8  & 19.7  & 24.1  & 38.2  
			\\ \bf 20 & 6.2   & 6.2   & 2.8   & 2.5   & 10.6  & 18.5  & 23.0  & 53.6  
			\\ \hline
		\end{tabular}
	\end{center}
\end{figure}

\begin{figure}[t] \begin{center}
	\resizebox{0.3\textwidth}{!}{
		\input{res_ps4_all_time.pgf}
	}
	\resizebox{0.29\textwidth}{!}{
		\input{res_ps4_all_u_calc_erank.pgf}
	}
	\resizebox{0.3\textwidth}{!}{
		\input{res_ps4_all_maxres.pgf}
	}
	\caption{
		Total computational time (on the left plot), effective TT-rank of the calculated solution $u$ (on the middle) and relative error for residual of linear system, that is solved by AMEn solver (on the right plot) w.r.t. the mesh size factor $d$ for the model PDE with four point sources.
		Results are presented for FS-QTT-solver (blue line with circle marker), FD-QTT-solver (green line with square marker) and FD-solver (red line with triangle marker).
	}
	\label{figure:res_ps4_all_time_u_calc_erank}
\end{center} \end{figure}

\begin{figure}[t!] 

	\begin{center}
		\resizebox{0.95\textwidth}{!}{
			\input{res_ps4_u_calc_from_d_8.pgf}
		}
		\caption{
			Calculated solution of the model PDE with four point sources for the grid factor value $d=8$.
		}
		\label{figure:res_ps4_u_calc_d8}
	\end{center}

	\begin{center}
		\resizebox{0.95\textwidth}{!}{
			\input{res_ps4_u_calc_from_d_10.pgf}
		}
		\caption{
			Calculated solution of the model PDE with four point sources for the grid factor value $d=10$, that is restricted to the grid with factor $d=8$.
		}
		\label{figure:res_ps4_u_calc_d10}
	\end{center}

	\begin{center}
		\resizebox{0.95\textwidth}{!}{
			\input{res_ps4_u_calc_from_d_15.pgf}
		}
		\caption{
			Calculated solution of the model PDE with four point sources for the grid factor value $d=15$ (on the left plot and on the middle), and for $d=20$ (on the right plot), that is restricted to the grid with factor $d=8$.
		}
		\label{figure:res_ps4_u_calc_d15}
	\end{center}
	
\end{figure}
		
We consider one more model problem of the form~\eqref{numex:pde_analyt} with coefficient $k$ from~\eqref{numex:pde_analyt_k} and with a right-hand side that is a model of four point sources
\begin{equation}\label{numex:pde_ps4_f}
	\begin{split}
    	f(x, y) = & \delta(x-0.2)\delta(y-0.2) + \delta(x-0.8)\delta(y-0.2) + \\
                  & \delta(x-0.2)\delta(y-0.8) + \delta(x-0.8)\delta(y-0.8).
    \end{split}
\end{equation}
			
We select $10^{-6}$ as accuracy of AMEn-solver and $10^{-8}$ as accuracy for TT-round and TT-cross operations, and perform calculations for the following grid factor values: 
\begin{itemize}
	\item $d=4, 5, \ldots 20$ for FS-QTT-Solver,
	\item $d=4, 5, \ldots 15$ for FD-QTT-Solver,
	\item $d=4, 5, \ldots 10$ for FD-Solver.
\end{itemize}		
			
\cref{figure:res_ps4_all_time_u_calc_erank} and \cref{numex:tbl:res_ps4_eranks} represent computational results for three solvers, that were described above, and on ~\cref{figure:res_ps4_u_calc_d8,figure:res_ps4_u_calc_d10,figure:res_ps4_u_calc_d15} we present calculated solution for the grid factor values $d=8, 10, 15, 20$.
For QTT-based solvers we transform QTT-vectors of the solution to the full format, and for $d$ values higher than $d=8$, we use simple restriction procedure and plot solution on the coarser grid with $d=8$.
As we can see from the figures, FS-QTT-Solver gives accurate solution for all considered values of the grid factor $d$.
					
\section{Related work} \label{section:related}

Tensor numerical methods are becoming increasingly popular in various fields of science (see, for example, review~\cite{larskres-survey-2013} and book~\cite{hackbusch-2012}).
Applications of low-rank tensor based techniques to PDEs are typically limited to multidimensional equations. In~\cite{khos-dmrg-2010} the QTT-decomposition was applied for the problem of quantum molecular computations.
Molecular Schrodinger equation was represented in the QTT-format, and the corresponding high-dimensional eigenvalue problem was efficiently solved by DMRG method~\cite{white-dmrg-1993}.
Low-rank tensor techniques were also applied to parametric and stochastic PDEs, that arise in uncertainty quantification and optimization problems.
Tensor decompositions can be used for solving these problems, for instance combined with a finite element discretisation~\cite{Schneider-HT-SFEM-2015, dklm-tt-pce-2015,khos-pde-2010}.
The total solution error is then influenced by the finite element discretisation, the truncation of coefficient expansions and polynomial degrees, and by the tensor approximation ranks.

Tensor based approach can be applied to a low-dimensional PDEs via QTT-decomposition technique~\cite{osel-2d2d-2010, khor-qtt-2011}.
In~\cite{ks-2dqtt-2015} it was proved, that approximate solution of elliptic problem with piecewise-analytic coefficients admit a compact QTT-representation, where a number of parameters depends polylogarithmic on the accuracy of approximation. 
		
Huge grids are required in a list of practical applications.
For example, in multiscale modeling~\cite{hoang-mscfem-2005, hoang-mscsfem-2013} one has to construct such a grid, that makes it possible to resolve the smallest spatial scale and the grids with $h \sim 2^{-20}$ or even finer grids should be used for a number of problems.
Formalism from the work~\cite{ks-2dqtt-2015} was successfully applied in~\cite{ksro-multiscale-2016} to two classes of the one-dimensional multiscale problems: two-scaled diffusion equation and Helmholtz equation with high wave numbers.
It was proved that solutions of such problems can be represented in the QTT-format with the polylogarithmic dependency of parameters number under prescribed tolerance of the solution, and in~\cite{khoromskij-quasiperiodic-2015} a preconditioner in the QTT-format was suggested for such class of problems.
		
However, standard discretization schemes (for example, finite  element or difference scheme) become inefficient for the fine grids due to conditioning problem, as was mentioned in \cref{section:intro}.
Hence, in the work~\cite{cor-eccomas-2016} we proposed a new derivative-free discretization scheme for solution of one-dimensional diffusion type PDEs, which is based on explicit formula for the PDE solution in the QTT-format.
It was also shown in~\cite{cor-eccomas-2016} that such scheme is effective for the multiscale modeling, and can handle up to $2^{50}$ virtual grid points, without problems with conditioning.
			
\section{Conclusions} \label{section:concl}

In this paper we proposed the efficient robust FS-QTT-solver for equations of diffusion type in two dimensions that is implemented in the low-rank QTT-format and resolves solution of the equation with high accuracy on very fine grids.
Presented numerical examples illustrate its efficiency.
Derivative-free discretization scheme, that is used in FS-QTT-solver, can be naturally generalized to the three-dimensional case and to the other forms of equations.
The multiscale problems and oscillatory problems are promising areas of application of the proposed solver, since in such problems one needs a very fine grid to resolve all the scales.

\section*{Acknowledgments}

Authors would like to thank Prof. Dr. Christoph Schwab and Dr. Vladimir Kazeev for providing helpful suggestions
on the topic of the manuscript.

\bibliographystyle{siamplain} 
\bibliography{./../../../bibtex_my/general,./../../../bibtex/our,./../../../bibtex/tensor,./../../../bibtex/stochastic,./../../../bibtex/dmrg}

\end{document}

%% file: res_analyt_all_time.pgf
\begingroup%
\makeatletter%
\begin{pgfpicture}%
\pgfpathrectangle{\pgfpointorigin}{\pgfqpoint{7.668116in}{7.355067in}}%
\pgfusepath{use as bounding box, clip}%
\begin{pgfscope}%
\pgfsetbuttcap%
\pgfsetmiterjoin%
\definecolor{currentfill}{rgb}{1.000000,1.000000,1.000000}%
\pgfsetfillcolor{currentfill}%
\pgfsetlinewidth{0.000000pt}%
\definecolor{currentstroke}{rgb}{1.000000,1.000000,1.000000}%
\pgfsetstrokecolor{currentstroke}%
\pgfsetdash{}{0pt}%
\pgfpathmoveto{\pgfqpoint{0.000000in}{0.000000in}}%
\pgfpathlineto{\pgfqpoint{7.668116in}{0.000000in}}%
\pgfpathlineto{\pgfqpoint{7.668116in}{7.355067in}}%
\pgfpathlineto{\pgfqpoint{0.000000in}{7.355067in}}%
\pgfpathclose%
\pgfusepath{fill}%
\end{pgfscope}%
\begin{pgfscope}%
\pgfsetbuttcap%
\pgfsetmiterjoin%
\definecolor{currentfill}{rgb}{1.000000,1.000000,1.000000}%
\pgfsetfillcolor{currentfill}%
\pgfsetlinewidth{0.000000pt}%
\definecolor{currentstroke}{rgb}{0.000000,0.000000,0.000000}%
\pgfsetstrokecolor{currentstroke}%
\pgfsetstrokeopacity{0.000000}%
\pgfsetdash{}{0pt}%
\pgfpathmoveto{\pgfqpoint{1.209638in}{0.899712in}}%
\pgfpathlineto{\pgfqpoint{7.409638in}{0.899712in}}%
\pgfpathlineto{\pgfqpoint{7.409638in}{7.099712in}}%
\pgfpathlineto{\pgfqpoint{1.209638in}{7.099712in}}%
\pgfpathclose%
\pgfusepath{fill}%
\end{pgfscope}%
\begin{pgfscope}%
\pgfpathrectangle{\pgfqpoint{1.209638in}{0.899712in}}{\pgfqpoint{6.200000in}{6.200000in}} %
\pgfusepath{clip}%
\pgfsetbuttcap%
\pgfsetroundjoin%
\pgfsetlinewidth{0.803000pt}%
\definecolor{currentstroke}{rgb}{0.800000,0.800000,0.800000}%
\pgfsetstrokecolor{currentstroke}%
\pgfsetdash{{1.000000pt}{3.000000pt}}{0.000000pt}%
\pgfpathmoveto{\pgfqpoint{1.209638in}{0.899712in}}%
\pgfpathlineto{\pgfqpoint{1.209638in}{7.099712in}}%
\pgfusepath{stroke}%
\end{pgfscope}%
\begin{pgfscope}%
\definecolor{textcolor}{rgb}{0.150000,0.150000,0.150000}%
\pgfsetstrokecolor{textcolor}%
\pgfsetfillcolor{textcolor}%
\pgftext[x=1.209638in,y=0.821934in,,top]{\color{textcolor}\sffamily\fontsize{24.000000}{28.800000}\selectfont \(\displaystyle 0\)}%
\end{pgfscope}%
\begin{pgfscope}%
\pgfpathrectangle{\pgfqpoint{1.209638in}{0.899712in}}{\pgfqpoint{6.200000in}{6.200000in}} %
\pgfusepath{clip}%
\pgfsetbuttcap%
\pgfsetroundjoin%
\pgfsetlinewidth{0.803000pt}%
\definecolor{currentstroke}{rgb}{0.800000,0.800000,0.800000}%
\pgfsetstrokecolor{currentstroke}%
\pgfsetdash{{1.000000pt}{3.000000pt}}{0.000000pt}%
\pgfpathmoveto{\pgfqpoint{2.242972in}{0.899712in}}%
\pgfpathlineto{\pgfqpoint{2.242972in}{7.099712in}}%
\pgfusepath{stroke}%
\end{pgfscope}%
\begin{pgfscope}%
\definecolor{textcolor}{rgb}{0.150000,0.150000,0.150000}%
\pgfsetstrokecolor{textcolor}%
\pgfsetfillcolor{textcolor}%
\pgftext[x=2.242972in,y=0.821934in,,top]{\color{textcolor}\sffamily\fontsize{24.000000}{28.800000}\selectfont \(\displaystyle 5\)}%
\end{pgfscope}%
\begin{pgfscope}%
\pgfpathrectangle{\pgfqpoint{1.209638in}{0.899712in}}{\pgfqpoint{6.200000in}{6.200000in}} %
\pgfusepath{clip}%
\pgfsetbuttcap%
\pgfsetroundjoin%
\pgfsetlinewidth{0.803000pt}%
\definecolor{currentstroke}{rgb}{0.800000,0.800000,0.800000}%
\pgfsetstrokecolor{currentstroke}%
\pgfsetdash{{1.000000pt}{3.000000pt}}{0.000000pt}%
\pgfpathmoveto{\pgfqpoint{3.276305in}{0.899712in}}%
\pgfpathlineto{\pgfqpoint{3.276305in}{7.099712in}}%
\pgfusepath{stroke}%
\end{pgfscope}%
\begin{pgfscope}%
\definecolor{textcolor}{rgb}{0.150000,0.150000,0.150000}%
\pgfsetstrokecolor{textcolor}%
\pgfsetfillcolor{textcolor}%
\pgftext[x=3.276305in,y=0.821934in,,top]{\color{textcolor}\sffamily\fontsize{24.000000}{28.800000}\selectfont \(\displaystyle 10\)}%
\end{pgfscope}%
\begin{pgfscope}%
\pgfpathrectangle{\pgfqpoint{1.209638in}{0.899712in}}{\pgfqpoint{6.200000in}{6.200000in}} %
\pgfusepath{clip}%
\pgfsetbuttcap%
\pgfsetroundjoin%
\pgfsetlinewidth{0.803000pt}%
\definecolor{currentstroke}{rgb}{0.800000,0.800000,0.800000}%
\pgfsetstrokecolor{currentstroke}%
\pgfsetdash{{1.000000pt}{3.000000pt}}{0.000000pt}%
\pgfpathmoveto{\pgfqpoint{4.309638in}{0.899712in}}%
\pgfpathlineto{\pgfqpoint{4.309638in}{7.099712in}}%
\pgfusepath{stroke}%
\end{pgfscope}%
\begin{pgfscope}%
\definecolor{textcolor}{rgb}{0.150000,0.150000,0.150000}%
\pgfsetstrokecolor{textcolor}%
\pgfsetfillcolor{textcolor}%
\pgftext[x=4.309638in,y=0.821934in,,top]{\color{textcolor}\sffamily\fontsize{24.000000}{28.800000}\selectfont \(\displaystyle 15\)}%
\end{pgfscope}%
\begin{pgfscope}%
\pgfpathrectangle{\pgfqpoint{1.209638in}{0.899712in}}{\pgfqpoint{6.200000in}{6.200000in}} %
\pgfusepath{clip}%
\pgfsetbuttcap%
\pgfsetroundjoin%
\pgfsetlinewidth{0.803000pt}%
\definecolor{currentstroke}{rgb}{0.800000,0.800000,0.800000}%
\pgfsetstrokecolor{currentstroke}%
\pgfsetdash{{1.000000pt}{3.000000pt}}{0.000000pt}%
\pgfpathmoveto{\pgfqpoint{5.342972in}{0.899712in}}%
\pgfpathlineto{\pgfqpoint{5.342972in}{7.099712in}}%
\pgfusepath{stroke}%
\end{pgfscope}%
\begin{pgfscope}%
\definecolor{textcolor}{rgb}{0.150000,0.150000,0.150000}%
\pgfsetstrokecolor{textcolor}%
\pgfsetfillcolor{textcolor}%
\pgftext[x=5.342972in,y=0.821934in,,top]{\color{textcolor}\sffamily\fontsize{24.000000}{28.800000}\selectfont \(\displaystyle 20\)}%
\end{pgfscope}%
\begin{pgfscope}%
\pgfpathrectangle{\pgfqpoint{1.209638in}{0.899712in}}{\pgfqpoint{6.200000in}{6.200000in}} %
\pgfusepath{clip}%
\pgfsetbuttcap%
\pgfsetroundjoin%
\pgfsetlinewidth{0.803000pt}%
\definecolor{currentstroke}{rgb}{0.800000,0.800000,0.800000}%
\pgfsetstrokecolor{currentstroke}%
\pgfsetdash{{1.000000pt}{3.000000pt}}{0.000000pt}%
\pgfpathmoveto{\pgfqpoint{6.376305in}{0.899712in}}%
\pgfpathlineto{\pgfqpoint{6.376305in}{7.099712in}}%
\pgfusepath{stroke}%
\end{pgfscope}%
\begin{pgfscope}%
\definecolor{textcolor}{rgb}{0.150000,0.150000,0.150000}%
\pgfsetstrokecolor{textcolor}%
\pgfsetfillcolor{textcolor}%
\pgftext[x=6.376305in,y=0.821934in,,top]{\color{textcolor}\sffamily\fontsize{24.000000}{28.800000}\selectfont \(\displaystyle 25\)}%
\end{pgfscope}%
\begin{pgfscope}%
\pgfpathrectangle{\pgfqpoint{1.209638in}{0.899712in}}{\pgfqpoint{6.200000in}{6.200000in}} %
\pgfusepath{clip}%
\pgfsetbuttcap%
\pgfsetroundjoin%
\pgfsetlinewidth{0.803000pt}%
\definecolor{currentstroke}{rgb}{0.800000,0.800000,0.800000}%
\pgfsetstrokecolor{currentstroke}%
\pgfsetdash{{1.000000pt}{3.000000pt}}{0.000000pt}%
\pgfpathmoveto{\pgfqpoint{7.409638in}{0.899712in}}%
\pgfpathlineto{\pgfqpoint{7.409638in}{7.099712in}}%
\pgfusepath{stroke}%
\end{pgfscope}%
\begin{pgfscope}%
\definecolor{textcolor}{rgb}{0.150000,0.150000,0.150000}%
\pgfsetstrokecolor{textcolor}%
\pgfsetfillcolor{textcolor}%
\pgftext[x=7.409638in,y=0.821934in,,top]{\color{textcolor}\sffamily\fontsize{24.000000}{28.800000}\selectfont \(\displaystyle 30\)}%
\end{pgfscope}%
\begin{pgfscope}%
\definecolor{textcolor}{rgb}{0.150000,0.150000,0.150000}%
\pgfsetstrokecolor{textcolor}%
\pgfsetfillcolor{textcolor}%
\pgftext[x=4.309638in,y=0.441780in,,top]{\color{textcolor}\sffamily\fontsize{26.400000}{31.680000}\selectfont d}%
\end{pgfscope}%
\begin{pgfscope}%
\pgfpathrectangle{\pgfqpoint{1.209638in}{0.899712in}}{\pgfqpoint{6.200000in}{6.200000in}} %
\pgfusepath{clip}%
\pgfsetbuttcap%
\pgfsetroundjoin%
\pgfsetlinewidth{0.803000pt}%
\definecolor{currentstroke}{rgb}{0.800000,0.800000,0.800000}%
\pgfsetstrokecolor{currentstroke}%
\pgfsetdash{{1.000000pt}{3.000000pt}}{0.000000pt}%
\pgfpathmoveto{\pgfqpoint{1.209638in}{0.899712in}}%
\pgfpathlineto{\pgfqpoint{7.409638in}{0.899712in}}%
\pgfusepath{stroke}%
\end{pgfscope}%
\begin{pgfscope}%
\definecolor{textcolor}{rgb}{0.150000,0.150000,0.150000}%
\pgfsetstrokecolor{textcolor}%
\pgfsetfillcolor{textcolor}%
\pgftext[x=1.131861in,y=0.899712in,right,]{\color{textcolor}\sffamily\fontsize{24.000000}{28.800000}\selectfont \(\displaystyle 10^{-2}\)}%
\end{pgfscope}%
\begin{pgfscope}%
\pgfpathrectangle{\pgfqpoint{1.209638in}{0.899712in}}{\pgfqpoint{6.200000in}{6.200000in}} %
\pgfusepath{clip}%
\pgfsetbuttcap%
\pgfsetroundjoin%
\pgfsetlinewidth{0.803000pt}%
\definecolor{currentstroke}{rgb}{0.800000,0.800000,0.800000}%
\pgfsetstrokecolor{currentstroke}%
\pgfsetdash{{1.000000pt}{3.000000pt}}{0.000000pt}%
\pgfpathmoveto{\pgfqpoint{1.209638in}{2.449712in}}%
\pgfpathlineto{\pgfqpoint{7.409638in}{2.449712in}}%
\pgfusepath{stroke}%
\end{pgfscope}%
\begin{pgfscope}%
\definecolor{textcolor}{rgb}{0.150000,0.150000,0.150000}%
\pgfsetstrokecolor{textcolor}%
\pgfsetfillcolor{textcolor}%
\pgftext[x=1.131861in,y=2.449712in,right,]{\color{textcolor}\sffamily\fontsize{24.000000}{28.800000}\selectfont \(\displaystyle 10^{-1}\)}%
\end{pgfscope}%
\begin{pgfscope}%
\pgfpathrectangle{\pgfqpoint{1.209638in}{0.899712in}}{\pgfqpoint{6.200000in}{6.200000in}} %
\pgfusepath{clip}%
\pgfsetbuttcap%
\pgfsetroundjoin%
\pgfsetlinewidth{0.803000pt}%
\definecolor{currentstroke}{rgb}{0.800000,0.800000,0.800000}%
\pgfsetstrokecolor{currentstroke}%
\pgfsetdash{{1.000000pt}{3.000000pt}}{0.000000pt}%
\pgfpathmoveto{\pgfqpoint{1.209638in}{3.999712in}}%
\pgfpathlineto{\pgfqpoint{7.409638in}{3.999712in}}%
\pgfusepath{stroke}%
\end{pgfscope}%
\begin{pgfscope}%
\definecolor{textcolor}{rgb}{0.150000,0.150000,0.150000}%
\pgfsetstrokecolor{textcolor}%
\pgfsetfillcolor{textcolor}%
\pgftext[x=1.131861in,y=3.999712in,right,]{\color{textcolor}\sffamily\fontsize{24.000000}{28.800000}\selectfont \(\displaystyle 10^{0}\)}%
\end{pgfscope}%
\begin{pgfscope}%
\pgfpathrectangle{\pgfqpoint{1.209638in}{0.899712in}}{\pgfqpoint{6.200000in}{6.200000in}} %
\pgfusepath{clip}%
\pgfsetbuttcap%
\pgfsetroundjoin%
\pgfsetlinewidth{0.803000pt}%
\definecolor{currentstroke}{rgb}{0.800000,0.800000,0.800000}%
\pgfsetstrokecolor{currentstroke}%
\pgfsetdash{{1.000000pt}{3.000000pt}}{0.000000pt}%
\pgfpathmoveto{\pgfqpoint{1.209638in}{5.549712in}}%
\pgfpathlineto{\pgfqpoint{7.409638in}{5.549712in}}%
\pgfusepath{stroke}%
\end{pgfscope}%
\begin{pgfscope}%
\definecolor{textcolor}{rgb}{0.150000,0.150000,0.150000}%
\pgfsetstrokecolor{textcolor}%
\pgfsetfillcolor{textcolor}%
\pgftext[x=1.131861in,y=5.549712in,right,]{\color{textcolor}\sffamily\fontsize{24.000000}{28.800000}\selectfont \(\displaystyle 10^{1}\)}%
\end{pgfscope}%
\begin{pgfscope}%
\pgfpathrectangle{\pgfqpoint{1.209638in}{0.899712in}}{\pgfqpoint{6.200000in}{6.200000in}} %
\pgfusepath{clip}%
\pgfsetbuttcap%
\pgfsetroundjoin%
\pgfsetlinewidth{0.803000pt}%
\definecolor{currentstroke}{rgb}{0.800000,0.800000,0.800000}%
\pgfsetstrokecolor{currentstroke}%
\pgfsetdash{{1.000000pt}{3.000000pt}}{0.000000pt}%
\pgfpathmoveto{\pgfqpoint{1.209638in}{7.099712in}}%
\pgfpathlineto{\pgfqpoint{7.409638in}{7.099712in}}%
\pgfusepath{stroke}%
\end{pgfscope}%
\begin{pgfscope}%
\definecolor{textcolor}{rgb}{0.150000,0.150000,0.150000}%
\pgfsetstrokecolor{textcolor}%
\pgfsetfillcolor{textcolor}%
\pgftext[x=1.131861in,y=7.099712in,right,]{\color{textcolor}\sffamily\fontsize{24.000000}{28.800000}\selectfont \(\displaystyle 10^{2}\)}%
\end{pgfscope}%
\begin{pgfscope}%
\definecolor{textcolor}{rgb}{0.150000,0.150000,0.150000}%
\pgfsetstrokecolor{textcolor}%
\pgfsetfillcolor{textcolor}%
\pgftext[x=0.441780in,y=3.999712in,,bottom,rotate=90.000000]{\color{textcolor}\sffamily\fontsize{26.400000}{31.680000}\selectfont Total time (seconds)}%
\end{pgfscope}%
\begin{pgfscope}%
\pgfpathrectangle{\pgfqpoint{1.209638in}{0.899712in}}{\pgfqpoint{6.200000in}{6.200000in}} %
\pgfusepath{clip}%
\pgfsetroundcap%
\pgfsetroundjoin%
\pgfsetlinewidth{2.007500pt}%
\definecolor{currentstroke}{rgb}{0.298039,0.447059,0.690196}%
\pgfsetstrokecolor{currentstroke}%
\pgfsetdash{}{0pt}%
\pgfpathmoveto{\pgfqpoint{2.036305in}{3.211103in}}%
\pgfpathlineto{\pgfqpoint{2.242972in}{3.380280in}}%
\pgfpathlineto{\pgfqpoint{2.449638in}{3.653283in}}%
\pgfpathlineto{\pgfqpoint{2.656305in}{3.878885in}}%
\pgfpathlineto{\pgfqpoint{2.862972in}{3.948931in}}%
\pgfpathlineto{\pgfqpoint{3.069638in}{4.080890in}}%
\pgfpathlineto{\pgfqpoint{3.276305in}{4.182627in}}%
\pgfpathlineto{\pgfqpoint{3.482972in}{4.198001in}}%
\pgfpathlineto{\pgfqpoint{3.689638in}{4.226733in}}%
\pgfpathlineto{\pgfqpoint{3.896305in}{4.495096in}}%
\pgfpathlineto{\pgfqpoint{4.102972in}{4.359018in}}%
\pgfpathlineto{\pgfqpoint{4.309638in}{4.392572in}}%
\pgfpathlineto{\pgfqpoint{4.516305in}{4.456513in}}%
\pgfpathlineto{\pgfqpoint{4.722972in}{4.529438in}}%
\pgfpathlineto{\pgfqpoint{4.929638in}{4.555747in}}%
\pgfpathlineto{\pgfqpoint{5.136305in}{4.779251in}}%
\pgfpathlineto{\pgfqpoint{5.342972in}{4.643508in}}%
\pgfpathlineto{\pgfqpoint{5.549638in}{4.750337in}}%
\pgfpathlineto{\pgfqpoint{5.756305in}{4.814531in}}%
\pgfpathlineto{\pgfqpoint{5.962972in}{4.886377in}}%
\pgfpathlineto{\pgfqpoint{6.169638in}{4.937949in}}%
\pgfpathlineto{\pgfqpoint{6.376305in}{5.058718in}}%
\pgfpathlineto{\pgfqpoint{6.582972in}{5.055392in}}%
\pgfpathlineto{\pgfqpoint{6.789638in}{5.052600in}}%
\pgfpathlineto{\pgfqpoint{6.996305in}{5.233966in}}%
\pgfpathlineto{\pgfqpoint{7.202972in}{5.217919in}}%
\pgfpathlineto{\pgfqpoint{7.409638in}{5.319332in}}%
\pgfusepath{stroke}%
\end{pgfscope}%
\begin{pgfscope}%
\pgfpathrectangle{\pgfqpoint{1.209638in}{0.899712in}}{\pgfqpoint{6.200000in}{6.200000in}} %
\pgfusepath{clip}%
\pgfsetbuttcap%
\pgfsetroundjoin%
\definecolor{currentfill}{rgb}{0.298039,0.447059,0.690196}%
\pgfsetfillcolor{currentfill}%
\pgfsetlinewidth{0.000000pt}%
\definecolor{currentstroke}{rgb}{0.000000,0.000000,0.000000}%
\pgfsetstrokecolor{currentstroke}%
\pgfsetdash{}{0pt}%
\pgfsys@defobject{currentmarker}{\pgfqpoint{-0.038889in}{-0.038889in}}{\pgfqpoint{0.038889in}{0.038889in}}{%
\pgfpathmoveto{\pgfqpoint{0.000000in}{-0.038889in}}%
\pgfpathcurveto{\pgfqpoint{0.010313in}{-0.038889in}}{\pgfqpoint{0.020206in}{-0.034791in}}{\pgfqpoint{0.027499in}{-0.027499in}}%
\pgfpathcurveto{\pgfqpoint{0.034791in}{-0.020206in}}{\pgfqpoint{0.038889in}{-0.010313in}}{\pgfqpoint{0.038889in}{0.000000in}}%
\pgfpathcurveto{\pgfqpoint{0.038889in}{0.010313in}}{\pgfqpoint{0.034791in}{0.020206in}}{\pgfqpoint{0.027499in}{0.027499in}}%
\pgfpathcurveto{\pgfqpoint{0.020206in}{0.034791in}}{\pgfqpoint{0.010313in}{0.038889in}}{\pgfqpoint{0.000000in}{0.038889in}}%
\pgfpathcurveto{\pgfqpoint{-0.010313in}{0.038889in}}{\pgfqpoint{-0.020206in}{0.034791in}}{\pgfqpoint{-0.027499in}{0.027499in}}%
\pgfpathcurveto{\pgfqpoint{-0.034791in}{0.020206in}}{\pgfqpoint{-0.038889in}{0.010313in}}{\pgfqpoint{-0.038889in}{0.000000in}}%
\pgfpathcurveto{\pgfqpoint{-0.038889in}{-0.010313in}}{\pgfqpoint{-0.034791in}{-0.020206in}}{\pgfqpoint{-0.027499in}{-0.027499in}}%
\pgfpathcurveto{\pgfqpoint{-0.020206in}{-0.034791in}}{\pgfqpoint{-0.010313in}{-0.038889in}}{\pgfqpoint{0.000000in}{-0.038889in}}%
\pgfpathclose%
\pgfusepath{fill}%
}%
\begin{pgfscope}%
\pgfsys@transformshift{2.036305in}{3.211103in}%
\pgfsys@useobject{currentmarker}{}%
\end{pgfscope}%
\begin{pgfscope}%
\pgfsys@transformshift{2.242972in}{3.380280in}%
\pgfsys@useobject{currentmarker}{}%
\end{pgfscope}%
\begin{pgfscope}%
\pgfsys@transformshift{2.449638in}{3.653283in}%
\pgfsys@useobject{currentmarker}{}%
\end{pgfscope}%
\begin{pgfscope}%
\pgfsys@transformshift{2.656305in}{3.878885in}%
\pgfsys@useobject{currentmarker}{}%
\end{pgfscope}%
\begin{pgfscope}%
\pgfsys@transformshift{2.862972in}{3.948931in}%
\pgfsys@useobject{currentmarker}{}%
\end{pgfscope}%
\begin{pgfscope}%
\pgfsys@transformshift{3.069638in}{4.080890in}%
\pgfsys@useobject{currentmarker}{}%
\end{pgfscope}%
\begin{pgfscope}%
\pgfsys@transformshift{3.276305in}{4.182627in}%
\pgfsys@useobject{currentmarker}{}%
\end{pgfscope}%
\begin{pgfscope}%
\pgfsys@transformshift{3.482972in}{4.198001in}%
\pgfsys@useobject{currentmarker}{}%
\end{pgfscope}%
\begin{pgfscope}%
\pgfsys@transformshift{3.689638in}{4.226733in}%
\pgfsys@useobject{currentmarker}{}%
\end{pgfscope}%
\begin{pgfscope}%
\pgfsys@transformshift{3.896305in}{4.495096in}%
\pgfsys@useobject{currentmarker}{}%
\end{pgfscope}%
\begin{pgfscope}%
\pgfsys@transformshift{4.102972in}{4.359018in}%
\pgfsys@useobject{currentmarker}{}%
\end{pgfscope}%
\begin{pgfscope}%
\pgfsys@transformshift{4.309638in}{4.392572in}%
\pgfsys@useobject{currentmarker}{}%
\end{pgfscope}%
\begin{pgfscope}%
\pgfsys@transformshift{4.516305in}{4.456513in}%
\pgfsys@useobject{currentmarker}{}%
\end{pgfscope}%
\begin{pgfscope}%
\pgfsys@transformshift{4.722972in}{4.529438in}%
\pgfsys@useobject{currentmarker}{}%
\end{pgfscope}%
\begin{pgfscope}%
\pgfsys@transformshift{4.929638in}{4.555747in}%
\pgfsys@useobject{currentmarker}{}%
\end{pgfscope}%
\begin{pgfscope}%
\pgfsys@transformshift{5.136305in}{4.779251in}%
\pgfsys@useobject{currentmarker}{}%
\end{pgfscope}%
\begin{pgfscope}%
\pgfsys@transformshift{5.342972in}{4.643508in}%
\pgfsys@useobject{currentmarker}{}%
\end{pgfscope}%
\begin{pgfscope}%
\pgfsys@transformshift{5.549638in}{4.750337in}%
\pgfsys@useobject{currentmarker}{}%
\end{pgfscope}%
\begin{pgfscope}%
\pgfsys@transformshift{5.756305in}{4.814531in}%
\pgfsys@useobject{currentmarker}{}%
\end{pgfscope}%
\begin{pgfscope}%
\pgfsys@transformshift{5.962972in}{4.886377in}%
\pgfsys@useobject{currentmarker}{}%
\end{pgfscope}%
\begin{pgfscope}%
\pgfsys@transformshift{6.169638in}{4.937949in}%
\pgfsys@useobject{currentmarker}{}%
\end{pgfscope}%
\begin{pgfscope}%
\pgfsys@transformshift{6.376305in}{5.058718in}%
\pgfsys@useobject{currentmarker}{}%
\end{pgfscope}%
\begin{pgfscope}%
\pgfsys@transformshift{6.582972in}{5.055392in}%
\pgfsys@useobject{currentmarker}{}%
\end{pgfscope}%
\begin{pgfscope}%
\pgfsys@transformshift{6.789638in}{5.052600in}%
\pgfsys@useobject{currentmarker}{}%
\end{pgfscope}%
\begin{pgfscope}%
\pgfsys@transformshift{6.996305in}{5.233966in}%
\pgfsys@useobject{currentmarker}{}%
\end{pgfscope}%
\begin{pgfscope}%
\pgfsys@transformshift{7.202972in}{5.217919in}%
\pgfsys@useobject{currentmarker}{}%
\end{pgfscope}%
\begin{pgfscope}%
\pgfsys@transformshift{7.409638in}{5.319332in}%
\pgfsys@useobject{currentmarker}{}%
\end{pgfscope}%
\end{pgfscope}%
\begin{pgfscope}%
\pgfpathrectangle{\pgfqpoint{1.209638in}{0.899712in}}{\pgfqpoint{6.200000in}{6.200000in}} %
\pgfusepath{clip}%
\pgfsetroundcap%
\pgfsetroundjoin%
\pgfsetlinewidth{2.007500pt}%
\definecolor{currentstroke}{rgb}{0.333333,0.658824,0.407843}%
\pgfsetstrokecolor{currentstroke}%
\pgfsetdash{}{0pt}%
\pgfpathmoveto{\pgfqpoint{2.036305in}{2.655854in}}%
\pgfpathlineto{\pgfqpoint{2.242972in}{2.741539in}}%
\pgfpathlineto{\pgfqpoint{2.449638in}{2.945407in}}%
\pgfpathlineto{\pgfqpoint{2.656305in}{3.121639in}}%
\pgfpathlineto{\pgfqpoint{2.862972in}{3.233915in}}%
\pgfpathlineto{\pgfqpoint{3.069638in}{3.419692in}}%
\pgfpathlineto{\pgfqpoint{3.276305in}{3.638324in}}%
\pgfpathlineto{\pgfqpoint{3.482972in}{4.446039in}}%
\pgfpathlineto{\pgfqpoint{3.689638in}{4.722603in}}%
\pgfpathlineto{\pgfqpoint{3.896305in}{4.712190in}}%
\pgfpathlineto{\pgfqpoint{4.102972in}{4.767548in}}%
\pgfpathlineto{\pgfqpoint{4.309638in}{4.828556in}}%
\pgfusepath{stroke}%
\end{pgfscope}%
\begin{pgfscope}%
\pgfpathrectangle{\pgfqpoint{1.209638in}{0.899712in}}{\pgfqpoint{6.200000in}{6.200000in}} %
\pgfusepath{clip}%
\pgfsetbuttcap%
\pgfsetmiterjoin%
\definecolor{currentfill}{rgb}{0.333333,0.658824,0.407843}%
\pgfsetfillcolor{currentfill}%
\pgfsetlinewidth{0.000000pt}%
\definecolor{currentstroke}{rgb}{0.000000,0.000000,0.000000}%
\pgfsetstrokecolor{currentstroke}%
\pgfsetdash{}{0pt}%
\pgfsys@defobject{currentmarker}{\pgfqpoint{-0.038889in}{-0.038889in}}{\pgfqpoint{0.038889in}{0.038889in}}{%
\pgfpathmoveto{\pgfqpoint{-0.038889in}{-0.038889in}}%
\pgfpathlineto{\pgfqpoint{0.038889in}{-0.038889in}}%
\pgfpathlineto{\pgfqpoint{0.038889in}{0.038889in}}%
\pgfpathlineto{\pgfqpoint{-0.038889in}{0.038889in}}%
\pgfpathclose%
\pgfusepath{fill}%
}%
\begin{pgfscope}%
\pgfsys@transformshift{2.036305in}{2.655854in}%
\pgfsys@useobject{currentmarker}{}%
\end{pgfscope}%
\begin{pgfscope}%
\pgfsys@transformshift{2.242972in}{2.741539in}%
\pgfsys@useobject{currentmarker}{}%
\end{pgfscope}%
\begin{pgfscope}%
\pgfsys@transformshift{2.449638in}{2.945407in}%
\pgfsys@useobject{currentmarker}{}%
\end{pgfscope}%
\begin{pgfscope}%
\pgfsys@transformshift{2.656305in}{3.121639in}%
\pgfsys@useobject{currentmarker}{}%
\end{pgfscope}%
\begin{pgfscope}%
\pgfsys@transformshift{2.862972in}{3.233915in}%
\pgfsys@useobject{currentmarker}{}%
\end{pgfscope}%
\begin{pgfscope}%
\pgfsys@transformshift{3.069638in}{3.419692in}%
\pgfsys@useobject{currentmarker}{}%
\end{pgfscope}%
\begin{pgfscope}%
\pgfsys@transformshift{3.276305in}{3.638324in}%
\pgfsys@useobject{currentmarker}{}%
\end{pgfscope}%
\begin{pgfscope}%
\pgfsys@transformshift{3.482972in}{4.446039in}%
\pgfsys@useobject{currentmarker}{}%
\end{pgfscope}%
\begin{pgfscope}%
\pgfsys@transformshift{3.689638in}{4.722603in}%
\pgfsys@useobject{currentmarker}{}%
\end{pgfscope}%
\begin{pgfscope}%
\pgfsys@transformshift{3.896305in}{4.712190in}%
\pgfsys@useobject{currentmarker}{}%
\end{pgfscope}%
\begin{pgfscope}%
\pgfsys@transformshift{4.102972in}{4.767548in}%
\pgfsys@useobject{currentmarker}{}%
\end{pgfscope}%
\begin{pgfscope}%
\pgfsys@transformshift{4.309638in}{4.828556in}%
\pgfsys@useobject{currentmarker}{}%
\end{pgfscope}%
\end{pgfscope}%
\begin{pgfscope}%
\pgfpathrectangle{\pgfqpoint{1.209638in}{0.899712in}}{\pgfqpoint{6.200000in}{6.200000in}} %
\pgfusepath{clip}%
\pgfsetroundcap%
\pgfsetroundjoin%
\pgfsetlinewidth{2.007500pt}%
\definecolor{currentstroke}{rgb}{0.768627,0.305882,0.321569}%
\pgfsetstrokecolor{currentstroke}%
\pgfsetdash{}{0pt}%
\pgfpathmoveto{\pgfqpoint{2.036305in}{0.993619in}}%
\pgfpathlineto{\pgfqpoint{2.242972in}{1.159378in}}%
\pgfpathlineto{\pgfqpoint{2.449638in}{1.726230in}}%
\pgfpathlineto{\pgfqpoint{2.656305in}{2.701749in}}%
\pgfpathlineto{\pgfqpoint{2.862972in}{3.613228in}}%
\pgfpathlineto{\pgfqpoint{3.069638in}{4.964659in}}%
\pgfpathlineto{\pgfqpoint{3.276305in}{6.350394in}}%
\pgfusepath{stroke}%
\end{pgfscope}%
\begin{pgfscope}%
\pgfpathrectangle{\pgfqpoint{1.209638in}{0.899712in}}{\pgfqpoint{6.200000in}{6.200000in}} %
\pgfusepath{clip}%
\pgfsetbuttcap%
\pgfsetmiterjoin%
\definecolor{currentfill}{rgb}{0.768627,0.305882,0.321569}%
\pgfsetfillcolor{currentfill}%
\pgfsetlinewidth{0.000000pt}%
\definecolor{currentstroke}{rgb}{0.000000,0.000000,0.000000}%
\pgfsetstrokecolor{currentstroke}%
\pgfsetdash{}{0pt}%
\pgfsys@defobject{currentmarker}{\pgfqpoint{-0.038889in}{-0.038889in}}{\pgfqpoint{0.038889in}{0.038889in}}{%
\pgfpathmoveto{\pgfqpoint{0.000000in}{0.038889in}}%
\pgfpathlineto{\pgfqpoint{-0.038889in}{-0.038889in}}%
\pgfpathlineto{\pgfqpoint{0.038889in}{-0.038889in}}%
\pgfpathclose%
\pgfusepath{fill}%
}%
\begin{pgfscope}%
\pgfsys@transformshift{2.036305in}{0.993619in}%
\pgfsys@useobject{currentmarker}{}%
\end{pgfscope}%
\begin{pgfscope}%
\pgfsys@transformshift{2.242972in}{1.159378in}%
\pgfsys@useobject{currentmarker}{}%
\end{pgfscope}%
\begin{pgfscope}%
\pgfsys@transformshift{2.449638in}{1.726230in}%
\pgfsys@useobject{currentmarker}{}%
\end{pgfscope}%
\begin{pgfscope}%
\pgfsys@transformshift{2.656305in}{2.701749in}%
\pgfsys@useobject{currentmarker}{}%
\end{pgfscope}%
\begin{pgfscope}%
\pgfsys@transformshift{2.862972in}{3.613228in}%
\pgfsys@useobject{currentmarker}{}%
\end{pgfscope}%
\begin{pgfscope}%
\pgfsys@transformshift{3.069638in}{4.964659in}%
\pgfsys@useobject{currentmarker}{}%
\end{pgfscope}%
\begin{pgfscope}%
\pgfsys@transformshift{3.276305in}{6.350394in}%
\pgfsys@useobject{currentmarker}{}%
\end{pgfscope}%
\end{pgfscope}%
\begin{pgfscope}%
\pgfsetrectcap%
\pgfsetmiterjoin%
\pgfsetlinewidth{1.254687pt}%
\definecolor{currentstroke}{rgb}{0.150000,0.150000,0.150000}%
\pgfsetstrokecolor{currentstroke}%
\pgfsetdash{}{0pt}%
\pgfpathmoveto{\pgfqpoint{1.209638in}{0.899712in}}%
\pgfpathlineto{\pgfqpoint{7.409638in}{0.899712in}}%
\pgfusepath{stroke}%
\end{pgfscope}%
\begin{pgfscope}%
\pgfsetrectcap%
\pgfsetmiterjoin%
\pgfsetlinewidth{1.254687pt}%
\definecolor{currentstroke}{rgb}{0.150000,0.150000,0.150000}%
\pgfsetstrokecolor{currentstroke}%
\pgfsetdash{}{0pt}%
\pgfpathmoveto{\pgfqpoint{1.209638in}{0.899712in}}%
\pgfpathlineto{\pgfqpoint{1.209638in}{7.099712in}}%
\pgfusepath{stroke}%
\end{pgfscope}%
\begin{pgfscope}%
\pgfsetbuttcap%
\pgfsetmiterjoin%
\definecolor{currentfill}{rgb}{1.000000,1.000000,1.000000}%
\pgfsetfillcolor{currentfill}%
\pgfsetlinewidth{0.240900pt}%
\definecolor{currentstroke}{rgb}{0.150000,0.150000,0.150000}%
\pgfsetstrokecolor{currentstroke}%
\pgfsetdash{}{0pt}%
\pgfpathmoveto{\pgfqpoint{3.839358in}{5.400917in}}%
\pgfpathlineto{\pgfqpoint{7.242972in}{5.400917in}}%
\pgfpathlineto{\pgfqpoint{7.242972in}{6.933045in}}%
\pgfpathlineto{\pgfqpoint{3.839358in}{6.933045in}}%
\pgfpathclose%
\pgfusepath{stroke,fill}%
\end{pgfscope}%
\begin{pgfscope}%
\pgfsetroundcap%
\pgfsetroundjoin%
\pgfsetlinewidth{2.007500pt}%
\definecolor{currentstroke}{rgb}{0.298039,0.447059,0.690196}%
\pgfsetstrokecolor{currentstroke}%
\pgfsetdash{}{0pt}%
\pgfpathmoveto{\pgfqpoint{3.972692in}{6.674842in}}%
\pgfpathlineto{\pgfqpoint{4.639358in}{6.674842in}}%
\pgfusepath{stroke}%
\end{pgfscope}%
\begin{pgfscope}%
\pgfsetbuttcap%
\pgfsetroundjoin%
\definecolor{currentfill}{rgb}{0.298039,0.447059,0.690196}%
\pgfsetfillcolor{currentfill}%
\pgfsetlinewidth{0.000000pt}%
\definecolor{currentstroke}{rgb}{0.000000,0.000000,0.000000}%
\pgfsetstrokecolor{currentstroke}%
\pgfsetdash{}{0pt}%
\pgfsys@defobject{currentmarker}{\pgfqpoint{-0.038889in}{-0.038889in}}{\pgfqpoint{0.038889in}{0.038889in}}{%
\pgfpathmoveto{\pgfqpoint{0.000000in}{-0.038889in}}%
\pgfpathcurveto{\pgfqpoint{0.010313in}{-0.038889in}}{\pgfqpoint{0.020206in}{-0.034791in}}{\pgfqpoint{0.027499in}{-0.027499in}}%
\pgfpathcurveto{\pgfqpoint{0.034791in}{-0.020206in}}{\pgfqpoint{0.038889in}{-0.010313in}}{\pgfqpoint{0.038889in}{0.000000in}}%
\pgfpathcurveto{\pgfqpoint{0.038889in}{0.010313in}}{\pgfqpoint{0.034791in}{0.020206in}}{\pgfqpoint{0.027499in}{0.027499in}}%
\pgfpathcurveto{\pgfqpoint{0.020206in}{0.034791in}}{\pgfqpoint{0.010313in}{0.038889in}}{\pgfqpoint{0.000000in}{0.038889in}}%
\pgfpathcurveto{\pgfqpoint{-0.010313in}{0.038889in}}{\pgfqpoint{-0.020206in}{0.034791in}}{\pgfqpoint{-0.027499in}{0.027499in}}%
\pgfpathcurveto{\pgfqpoint{-0.034791in}{0.020206in}}{\pgfqpoint{-0.038889in}{0.010313in}}{\pgfqpoint{-0.038889in}{0.000000in}}%
\pgfpathcurveto{\pgfqpoint{-0.038889in}{-0.010313in}}{\pgfqpoint{-0.034791in}{-0.020206in}}{\pgfqpoint{-0.027499in}{-0.027499in}}%
\pgfpathcurveto{\pgfqpoint{-0.020206in}{-0.034791in}}{\pgfqpoint{-0.010313in}{-0.038889in}}{\pgfqpoint{0.000000in}{-0.038889in}}%
\pgfpathclose%
\pgfusepath{fill}%
}%
\begin{pgfscope}%
\pgfsys@transformshift{4.306025in}{6.674842in}%
\pgfsys@useobject{currentmarker}{}%
\end{pgfscope}%
\end{pgfscope}%
\begin{pgfscope}%
\definecolor{textcolor}{rgb}{0.150000,0.150000,0.150000}%
\pgfsetstrokecolor{textcolor}%
\pgfsetfillcolor{textcolor}%
\pgftext[x=4.906025in,y=6.558175in,left,base]{\color{textcolor}\sffamily\fontsize{24.000000}{28.800000}\selectfont FS-QTT-solver}%
\end{pgfscope}%
\begin{pgfscope}%
\pgfsetroundcap%
\pgfsetroundjoin%
\pgfsetlinewidth{2.007500pt}%
\definecolor{currentstroke}{rgb}{0.333333,0.658824,0.407843}%
\pgfsetstrokecolor{currentstroke}%
\pgfsetdash{}{0pt}%
\pgfpathmoveto{\pgfqpoint{3.972692in}{6.197466in}}%
\pgfpathlineto{\pgfqpoint{4.639358in}{6.197466in}}%
\pgfusepath{stroke}%
\end{pgfscope}%
\begin{pgfscope}%
\pgfsetbuttcap%
\pgfsetmiterjoin%
\definecolor{currentfill}{rgb}{0.333333,0.658824,0.407843}%
\pgfsetfillcolor{currentfill}%
\pgfsetlinewidth{0.000000pt}%
\definecolor{currentstroke}{rgb}{0.000000,0.000000,0.000000}%
\pgfsetstrokecolor{currentstroke}%
\pgfsetdash{}{0pt}%
\pgfsys@defobject{currentmarker}{\pgfqpoint{-0.038889in}{-0.038889in}}{\pgfqpoint{0.038889in}{0.038889in}}{%
\pgfpathmoveto{\pgfqpoint{-0.038889in}{-0.038889in}}%
\pgfpathlineto{\pgfqpoint{0.038889in}{-0.038889in}}%
\pgfpathlineto{\pgfqpoint{0.038889in}{0.038889in}}%
\pgfpathlineto{\pgfqpoint{-0.038889in}{0.038889in}}%
\pgfpathclose%
\pgfusepath{fill}%
}%
\begin{pgfscope}%
\pgfsys@transformshift{4.306025in}{6.197466in}%
\pgfsys@useobject{currentmarker}{}%
\end{pgfscope}%
\end{pgfscope}%
\begin{pgfscope}%
\definecolor{textcolor}{rgb}{0.150000,0.150000,0.150000}%
\pgfsetstrokecolor{textcolor}%
\pgfsetfillcolor{textcolor}%
\pgftext[x=4.906025in,y=6.080799in,left,base]{\color{textcolor}\sffamily\fontsize{24.000000}{28.800000}\selectfont FD-QTT-solver}%
\end{pgfscope}%
\begin{pgfscope}%
\pgfsetroundcap%
\pgfsetroundjoin%
\pgfsetlinewidth{2.007500pt}%
\definecolor{currentstroke}{rgb}{0.768627,0.305882,0.321569}%
\pgfsetstrokecolor{currentstroke}%
\pgfsetdash{}{0pt}%
\pgfpathmoveto{\pgfqpoint{3.972692in}{5.720090in}}%
\pgfpathlineto{\pgfqpoint{4.639358in}{5.720090in}}%
\pgfusepath{stroke}%
\end{pgfscope}%
\begin{pgfscope}%
\pgfsetbuttcap%
\pgfsetmiterjoin%
\definecolor{currentfill}{rgb}{0.768627,0.305882,0.321569}%
\pgfsetfillcolor{currentfill}%
\pgfsetlinewidth{0.000000pt}%
\definecolor{currentstroke}{rgb}{0.000000,0.000000,0.000000}%
\pgfsetstrokecolor{currentstroke}%
\pgfsetdash{}{0pt}%
\pgfsys@defobject{currentmarker}{\pgfqpoint{-0.038889in}{-0.038889in}}{\pgfqpoint{0.038889in}{0.038889in}}{%
\pgfpathmoveto{\pgfqpoint{0.000000in}{0.038889in}}%
\pgfpathlineto{\pgfqpoint{-0.038889in}{-0.038889in}}%
\pgfpathlineto{\pgfqpoint{0.038889in}{-0.038889in}}%
\pgfpathclose%
\pgfusepath{fill}%
}%
\begin{pgfscope}%
\pgfsys@transformshift{4.306025in}{5.720090in}%
\pgfsys@useobject{currentmarker}{}%
\end{pgfscope}%
\end{pgfscope}%
\begin{pgfscope}%
\definecolor{textcolor}{rgb}{0.150000,0.150000,0.150000}%
\pgfsetstrokecolor{textcolor}%
\pgfsetfillcolor{textcolor}%
\pgftext[x=4.906025in,y=5.603423in,left,base]{\color{textcolor}\sffamily\fontsize{24.000000}{28.800000}\selectfont FD-solver}%
\end{pgfscope}%
\end{pgfpicture}%
\makeatother%
\endgroup%

%% file: res_analyt_all_u_calc_erank.pgf
\begingroup%
\makeatletter%
\begin{pgfpicture}%
\pgfpathrectangle{\pgfpointorigin}{\pgfqpoint{7.364436in}{7.355067in}}%
\pgfusepath{use as bounding box, clip}%
\begin{pgfscope}%
\pgfsetbuttcap%
\pgfsetmiterjoin%
\definecolor{currentfill}{rgb}{1.000000,1.000000,1.000000}%
\pgfsetfillcolor{currentfill}%
\pgfsetlinewidth{0.000000pt}%
\definecolor{currentstroke}{rgb}{1.000000,1.000000,1.000000}%
\pgfsetstrokecolor{currentstroke}%
\pgfsetdash{}{0pt}%
\pgfpathmoveto{\pgfqpoint{0.000000in}{0.000000in}}%
\pgfpathlineto{\pgfqpoint{7.364436in}{0.000000in}}%
\pgfpathlineto{\pgfqpoint{7.364436in}{7.355067in}}%
\pgfpathlineto{\pgfqpoint{0.000000in}{7.355067in}}%
\pgfpathclose%
\pgfusepath{fill}%
\end{pgfscope}%
\begin{pgfscope}%
\pgfsetbuttcap%
\pgfsetmiterjoin%
\definecolor{currentfill}{rgb}{1.000000,1.000000,1.000000}%
\pgfsetfillcolor{currentfill}%
\pgfsetlinewidth{0.000000pt}%
\definecolor{currentstroke}{rgb}{0.000000,0.000000,0.000000}%
\pgfsetstrokecolor{currentstroke}%
\pgfsetstrokeopacity{0.000000}%
\pgfsetdash{}{0pt}%
\pgfpathmoveto{\pgfqpoint{0.905958in}{0.899712in}}%
\pgfpathlineto{\pgfqpoint{7.105958in}{0.899712in}}%
\pgfpathlineto{\pgfqpoint{7.105958in}{7.099712in}}%
\pgfpathlineto{\pgfqpoint{0.905958in}{7.099712in}}%
\pgfpathclose%
\pgfusepath{fill}%
\end{pgfscope}%
\begin{pgfscope}%
\pgfpathrectangle{\pgfqpoint{0.905958in}{0.899712in}}{\pgfqpoint{6.200000in}{6.200000in}} %
\pgfusepath{clip}%
\pgfsetbuttcap%
\pgfsetroundjoin%
\pgfsetlinewidth{0.803000pt}%
\definecolor{currentstroke}{rgb}{0.800000,0.800000,0.800000}%
\pgfsetstrokecolor{currentstroke}%
\pgfsetdash{{1.000000pt}{3.000000pt}}{0.000000pt}%
\pgfpathmoveto{\pgfqpoint{0.905958in}{0.899712in}}%
\pgfpathlineto{\pgfqpoint{0.905958in}{7.099712in}}%
\pgfusepath{stroke}%
\end{pgfscope}%
\begin{pgfscope}%
\definecolor{textcolor}{rgb}{0.150000,0.150000,0.150000}%
\pgfsetstrokecolor{textcolor}%
\pgfsetfillcolor{textcolor}%
\pgftext[x=0.905958in,y=0.821934in,,top]{\color{textcolor}\sffamily\fontsize{24.000000}{28.800000}\selectfont \(\displaystyle 0\)}%
\end{pgfscope}%
\begin{pgfscope}%
\pgfpathrectangle{\pgfqpoint{0.905958in}{0.899712in}}{\pgfqpoint{6.200000in}{6.200000in}} %
\pgfusepath{clip}%
\pgfsetbuttcap%
\pgfsetroundjoin%
\pgfsetlinewidth{0.803000pt}%
\definecolor{currentstroke}{rgb}{0.800000,0.800000,0.800000}%
\pgfsetstrokecolor{currentstroke}%
\pgfsetdash{{1.000000pt}{3.000000pt}}{0.000000pt}%
\pgfpathmoveto{\pgfqpoint{1.939292in}{0.899712in}}%
\pgfpathlineto{\pgfqpoint{1.939292in}{7.099712in}}%
\pgfusepath{stroke}%
\end{pgfscope}%
\begin{pgfscope}%
\definecolor{textcolor}{rgb}{0.150000,0.150000,0.150000}%
\pgfsetstrokecolor{textcolor}%
\pgfsetfillcolor{textcolor}%
\pgftext[x=1.939292in,y=0.821934in,,top]{\color{textcolor}\sffamily\fontsize{24.000000}{28.800000}\selectfont \(\displaystyle 5\)}%
\end{pgfscope}%
\begin{pgfscope}%
\pgfpathrectangle{\pgfqpoint{0.905958in}{0.899712in}}{\pgfqpoint{6.200000in}{6.200000in}} %
\pgfusepath{clip}%
\pgfsetbuttcap%
\pgfsetroundjoin%
\pgfsetlinewidth{0.803000pt}%
\definecolor{currentstroke}{rgb}{0.800000,0.800000,0.800000}%
\pgfsetstrokecolor{currentstroke}%
\pgfsetdash{{1.000000pt}{3.000000pt}}{0.000000pt}%
\pgfpathmoveto{\pgfqpoint{2.972625in}{0.899712in}}%
\pgfpathlineto{\pgfqpoint{2.972625in}{7.099712in}}%
\pgfusepath{stroke}%
\end{pgfscope}%
\begin{pgfscope}%
\definecolor{textcolor}{rgb}{0.150000,0.150000,0.150000}%
\pgfsetstrokecolor{textcolor}%
\pgfsetfillcolor{textcolor}%
\pgftext[x=2.972625in,y=0.821934in,,top]{\color{textcolor}\sffamily\fontsize{24.000000}{28.800000}\selectfont \(\displaystyle 10\)}%
\end{pgfscope}%
\begin{pgfscope}%
\pgfpathrectangle{\pgfqpoint{0.905958in}{0.899712in}}{\pgfqpoint{6.200000in}{6.200000in}} %
\pgfusepath{clip}%
\pgfsetbuttcap%
\pgfsetroundjoin%
\pgfsetlinewidth{0.803000pt}%
\definecolor{currentstroke}{rgb}{0.800000,0.800000,0.800000}%
\pgfsetstrokecolor{currentstroke}%
\pgfsetdash{{1.000000pt}{3.000000pt}}{0.000000pt}%
\pgfpathmoveto{\pgfqpoint{4.005958in}{0.899712in}}%
\pgfpathlineto{\pgfqpoint{4.005958in}{7.099712in}}%
\pgfusepath{stroke}%
\end{pgfscope}%
\begin{pgfscope}%
\definecolor{textcolor}{rgb}{0.150000,0.150000,0.150000}%
\pgfsetstrokecolor{textcolor}%
\pgfsetfillcolor{textcolor}%
\pgftext[x=4.005958in,y=0.821934in,,top]{\color{textcolor}\sffamily\fontsize{24.000000}{28.800000}\selectfont \(\displaystyle 15\)}%
\end{pgfscope}%
\begin{pgfscope}%
\pgfpathrectangle{\pgfqpoint{0.905958in}{0.899712in}}{\pgfqpoint{6.200000in}{6.200000in}} %
\pgfusepath{clip}%
\pgfsetbuttcap%
\pgfsetroundjoin%
\pgfsetlinewidth{0.803000pt}%
\definecolor{currentstroke}{rgb}{0.800000,0.800000,0.800000}%
\pgfsetstrokecolor{currentstroke}%
\pgfsetdash{{1.000000pt}{3.000000pt}}{0.000000pt}%
\pgfpathmoveto{\pgfqpoint{5.039292in}{0.899712in}}%
\pgfpathlineto{\pgfqpoint{5.039292in}{7.099712in}}%
\pgfusepath{stroke}%
\end{pgfscope}%
\begin{pgfscope}%
\definecolor{textcolor}{rgb}{0.150000,0.150000,0.150000}%
\pgfsetstrokecolor{textcolor}%
\pgfsetfillcolor{textcolor}%
\pgftext[x=5.039292in,y=0.821934in,,top]{\color{textcolor}\sffamily\fontsize{24.000000}{28.800000}\selectfont \(\displaystyle 20\)}%
\end{pgfscope}%
\begin{pgfscope}%
\pgfpathrectangle{\pgfqpoint{0.905958in}{0.899712in}}{\pgfqpoint{6.200000in}{6.200000in}} %
\pgfusepath{clip}%
\pgfsetbuttcap%
\pgfsetroundjoin%
\pgfsetlinewidth{0.803000pt}%
\definecolor{currentstroke}{rgb}{0.800000,0.800000,0.800000}%
\pgfsetstrokecolor{currentstroke}%
\pgfsetdash{{1.000000pt}{3.000000pt}}{0.000000pt}%
\pgfpathmoveto{\pgfqpoint{6.072625in}{0.899712in}}%
\pgfpathlineto{\pgfqpoint{6.072625in}{7.099712in}}%
\pgfusepath{stroke}%
\end{pgfscope}%
\begin{pgfscope}%
\definecolor{textcolor}{rgb}{0.150000,0.150000,0.150000}%
\pgfsetstrokecolor{textcolor}%
\pgfsetfillcolor{textcolor}%
\pgftext[x=6.072625in,y=0.821934in,,top]{\color{textcolor}\sffamily\fontsize{24.000000}{28.800000}\selectfont \(\displaystyle 25\)}%
\end{pgfscope}%
\begin{pgfscope}%
\pgfpathrectangle{\pgfqpoint{0.905958in}{0.899712in}}{\pgfqpoint{6.200000in}{6.200000in}} %
\pgfusepath{clip}%
\pgfsetbuttcap%
\pgfsetroundjoin%
\pgfsetlinewidth{0.803000pt}%
\definecolor{currentstroke}{rgb}{0.800000,0.800000,0.800000}%
\pgfsetstrokecolor{currentstroke}%
\pgfsetdash{{1.000000pt}{3.000000pt}}{0.000000pt}%
\pgfpathmoveto{\pgfqpoint{7.105958in}{0.899712in}}%
\pgfpathlineto{\pgfqpoint{7.105958in}{7.099712in}}%
\pgfusepath{stroke}%
\end{pgfscope}%
\begin{pgfscope}%
\definecolor{textcolor}{rgb}{0.150000,0.150000,0.150000}%
\pgfsetstrokecolor{textcolor}%
\pgfsetfillcolor{textcolor}%
\pgftext[x=7.105958in,y=0.821934in,,top]{\color{textcolor}\sffamily\fontsize{24.000000}{28.800000}\selectfont \(\displaystyle 30\)}%
\end{pgfscope}%
\begin{pgfscope}%
\definecolor{textcolor}{rgb}{0.150000,0.150000,0.150000}%
\pgfsetstrokecolor{textcolor}%
\pgfsetfillcolor{textcolor}%
\pgftext[x=4.005958in,y=0.441780in,,top]{\color{textcolor}\sffamily\fontsize{26.400000}{31.680000}\selectfont d}%
\end{pgfscope}%
\begin{pgfscope}%
\pgfpathrectangle{\pgfqpoint{0.905958in}{0.899712in}}{\pgfqpoint{6.200000in}{6.200000in}} %
\pgfusepath{clip}%
\pgfsetbuttcap%
\pgfsetroundjoin%
\pgfsetlinewidth{0.803000pt}%
\definecolor{currentstroke}{rgb}{0.800000,0.800000,0.800000}%
\pgfsetstrokecolor{currentstroke}%
\pgfsetdash{{1.000000pt}{3.000000pt}}{0.000000pt}%
\pgfpathmoveto{\pgfqpoint{0.905958in}{0.899712in}}%
\pgfpathlineto{\pgfqpoint{7.105958in}{0.899712in}}%
\pgfusepath{stroke}%
\end{pgfscope}%
\begin{pgfscope}%
\definecolor{textcolor}{rgb}{0.150000,0.150000,0.150000}%
\pgfsetstrokecolor{textcolor}%
\pgfsetfillcolor{textcolor}%
\pgftext[x=0.828181in,y=0.899712in,right,]{\color{textcolor}\sffamily\fontsize{24.000000}{28.800000}\selectfont \(\displaystyle 2\)}%
\end{pgfscope}%
\begin{pgfscope}%
\pgfpathrectangle{\pgfqpoint{0.905958in}{0.899712in}}{\pgfqpoint{6.200000in}{6.200000in}} %
\pgfusepath{clip}%
\pgfsetbuttcap%
\pgfsetroundjoin%
\pgfsetlinewidth{0.803000pt}%
\definecolor{currentstroke}{rgb}{0.800000,0.800000,0.800000}%
\pgfsetstrokecolor{currentstroke}%
\pgfsetdash{{1.000000pt}{3.000000pt}}{0.000000pt}%
\pgfpathmoveto{\pgfqpoint{0.905958in}{1.933045in}}%
\pgfpathlineto{\pgfqpoint{7.105958in}{1.933045in}}%
\pgfusepath{stroke}%
\end{pgfscope}%
\begin{pgfscope}%
\definecolor{textcolor}{rgb}{0.150000,0.150000,0.150000}%
\pgfsetstrokecolor{textcolor}%
\pgfsetfillcolor{textcolor}%
\pgftext[x=0.828181in,y=1.933045in,right,]{\color{textcolor}\sffamily\fontsize{24.000000}{28.800000}\selectfont \(\displaystyle 4\)}%
\end{pgfscope}%
\begin{pgfscope}%
\pgfpathrectangle{\pgfqpoint{0.905958in}{0.899712in}}{\pgfqpoint{6.200000in}{6.200000in}} %
\pgfusepath{clip}%
\pgfsetbuttcap%
\pgfsetroundjoin%
\pgfsetlinewidth{0.803000pt}%
\definecolor{currentstroke}{rgb}{0.800000,0.800000,0.800000}%
\pgfsetstrokecolor{currentstroke}%
\pgfsetdash{{1.000000pt}{3.000000pt}}{0.000000pt}%
\pgfpathmoveto{\pgfqpoint{0.905958in}{2.966379in}}%
\pgfpathlineto{\pgfqpoint{7.105958in}{2.966379in}}%
\pgfusepath{stroke}%
\end{pgfscope}%
\begin{pgfscope}%
\definecolor{textcolor}{rgb}{0.150000,0.150000,0.150000}%
\pgfsetstrokecolor{textcolor}%
\pgfsetfillcolor{textcolor}%
\pgftext[x=0.828181in,y=2.966379in,right,]{\color{textcolor}\sffamily\fontsize{24.000000}{28.800000}\selectfont \(\displaystyle 6\)}%
\end{pgfscope}%
\begin{pgfscope}%
\pgfpathrectangle{\pgfqpoint{0.905958in}{0.899712in}}{\pgfqpoint{6.200000in}{6.200000in}} %
\pgfusepath{clip}%
\pgfsetbuttcap%
\pgfsetroundjoin%
\pgfsetlinewidth{0.803000pt}%
\definecolor{currentstroke}{rgb}{0.800000,0.800000,0.800000}%
\pgfsetstrokecolor{currentstroke}%
\pgfsetdash{{1.000000pt}{3.000000pt}}{0.000000pt}%
\pgfpathmoveto{\pgfqpoint{0.905958in}{3.999712in}}%
\pgfpathlineto{\pgfqpoint{7.105958in}{3.999712in}}%
\pgfusepath{stroke}%
\end{pgfscope}%
\begin{pgfscope}%
\definecolor{textcolor}{rgb}{0.150000,0.150000,0.150000}%
\pgfsetstrokecolor{textcolor}%
\pgfsetfillcolor{textcolor}%
\pgftext[x=0.828181in,y=3.999712in,right,]{\color{textcolor}\sffamily\fontsize{24.000000}{28.800000}\selectfont \(\displaystyle 8\)}%
\end{pgfscope}%
\begin{pgfscope}%
\pgfpathrectangle{\pgfqpoint{0.905958in}{0.899712in}}{\pgfqpoint{6.200000in}{6.200000in}} %
\pgfusepath{clip}%
\pgfsetbuttcap%
\pgfsetroundjoin%
\pgfsetlinewidth{0.803000pt}%
\definecolor{currentstroke}{rgb}{0.800000,0.800000,0.800000}%
\pgfsetstrokecolor{currentstroke}%
\pgfsetdash{{1.000000pt}{3.000000pt}}{0.000000pt}%
\pgfpathmoveto{\pgfqpoint{0.905958in}{5.033045in}}%
\pgfpathlineto{\pgfqpoint{7.105958in}{5.033045in}}%
\pgfusepath{stroke}%
\end{pgfscope}%
\begin{pgfscope}%
\definecolor{textcolor}{rgb}{0.150000,0.150000,0.150000}%
\pgfsetstrokecolor{textcolor}%
\pgfsetfillcolor{textcolor}%
\pgftext[x=0.828181in,y=5.033045in,right,]{\color{textcolor}\sffamily\fontsize{24.000000}{28.800000}\selectfont \(\displaystyle 10\)}%
\end{pgfscope}%
\begin{pgfscope}%
\pgfpathrectangle{\pgfqpoint{0.905958in}{0.899712in}}{\pgfqpoint{6.200000in}{6.200000in}} %
\pgfusepath{clip}%
\pgfsetbuttcap%
\pgfsetroundjoin%
\pgfsetlinewidth{0.803000pt}%
\definecolor{currentstroke}{rgb}{0.800000,0.800000,0.800000}%
\pgfsetstrokecolor{currentstroke}%
\pgfsetdash{{1.000000pt}{3.000000pt}}{0.000000pt}%
\pgfpathmoveto{\pgfqpoint{0.905958in}{6.066379in}}%
\pgfpathlineto{\pgfqpoint{7.105958in}{6.066379in}}%
\pgfusepath{stroke}%
\end{pgfscope}%
\begin{pgfscope}%
\definecolor{textcolor}{rgb}{0.150000,0.150000,0.150000}%
\pgfsetstrokecolor{textcolor}%
\pgfsetfillcolor{textcolor}%
\pgftext[x=0.828181in,y=6.066379in,right,]{\color{textcolor}\sffamily\fontsize{24.000000}{28.800000}\selectfont \(\displaystyle 12\)}%
\end{pgfscope}%
\begin{pgfscope}%
\pgfpathrectangle{\pgfqpoint{0.905958in}{0.899712in}}{\pgfqpoint{6.200000in}{6.200000in}} %
\pgfusepath{clip}%
\pgfsetbuttcap%
\pgfsetroundjoin%
\pgfsetlinewidth{0.803000pt}%
\definecolor{currentstroke}{rgb}{0.800000,0.800000,0.800000}%
\pgfsetstrokecolor{currentstroke}%
\pgfsetdash{{1.000000pt}{3.000000pt}}{0.000000pt}%
\pgfpathmoveto{\pgfqpoint{0.905958in}{7.099712in}}%
\pgfpathlineto{\pgfqpoint{7.105958in}{7.099712in}}%
\pgfusepath{stroke}%
\end{pgfscope}%
\begin{pgfscope}%
\definecolor{textcolor}{rgb}{0.150000,0.150000,0.150000}%
\pgfsetstrokecolor{textcolor}%
\pgfsetfillcolor{textcolor}%
\pgftext[x=0.828181in,y=7.099712in,right,]{\color{textcolor}\sffamily\fontsize{24.000000}{28.800000}\selectfont \(\displaystyle 14\)}%
\end{pgfscope}%
\begin{pgfscope}%
\definecolor{textcolor}{rgb}{0.150000,0.150000,0.150000}%
\pgfsetstrokecolor{textcolor}%
\pgfsetfillcolor{textcolor}%
\pgftext[x=0.441780in,y=3.999712in,,bottom,rotate=90.000000]{\color{textcolor}\sffamily\fontsize{26.400000}{31.680000}\selectfont Solution effective TT-rank}%
\end{pgfscope}%
\begin{pgfscope}%
\pgfpathrectangle{\pgfqpoint{0.905958in}{0.899712in}}{\pgfqpoint{6.200000in}{6.200000in}} %
\pgfusepath{clip}%
\pgfsetroundcap%
\pgfsetroundjoin%
\pgfsetlinewidth{2.007500pt}%
\definecolor{currentstroke}{rgb}{0.298039,0.447059,0.690196}%
\pgfsetstrokecolor{currentstroke}%
\pgfsetdash{}{0pt}%
\pgfpathmoveto{\pgfqpoint{1.732625in}{2.966379in}}%
\pgfpathlineto{\pgfqpoint{1.939292in}{3.545958in}}%
\pgfpathlineto{\pgfqpoint{2.145958in}{3.739679in}}%
\pgfpathlineto{\pgfqpoint{2.352625in}{3.543333in}}%
\pgfpathlineto{\pgfqpoint{2.559292in}{3.646397in}}%
\pgfpathlineto{\pgfqpoint{2.765958in}{3.300143in}}%
\pgfpathlineto{\pgfqpoint{2.972625in}{3.314380in}}%
\pgfpathlineto{\pgfqpoint{3.179292in}{2.774762in}}%
\pgfpathlineto{\pgfqpoint{3.385958in}{2.821213in}}%
\pgfpathlineto{\pgfqpoint{3.592625in}{2.466736in}}%
\pgfpathlineto{\pgfqpoint{3.799292in}{2.163851in}}%
\pgfpathlineto{\pgfqpoint{4.005958in}{2.191379in}}%
\pgfpathlineto{\pgfqpoint{4.212625in}{2.392628in}}%
\pgfpathlineto{\pgfqpoint{4.419292in}{2.457722in}}%
\pgfpathlineto{\pgfqpoint{4.625958in}{2.788667in}}%
\pgfpathlineto{\pgfqpoint{4.832625in}{2.952060in}}%
\pgfpathlineto{\pgfqpoint{5.039292in}{3.184366in}}%
\pgfpathlineto{\pgfqpoint{5.245958in}{3.565762in}}%
\pgfpathlineto{\pgfqpoint{5.452625in}{4.234530in}}%
\pgfpathlineto{\pgfqpoint{5.659292in}{4.213037in}}%
\pgfpathlineto{\pgfqpoint{5.865958in}{4.644712in}}%
\pgfpathlineto{\pgfqpoint{6.072625in}{4.885945in}}%
\pgfpathlineto{\pgfqpoint{6.279292in}{5.420245in}}%
\pgfpathlineto{\pgfqpoint{6.485958in}{5.398937in}}%
\pgfpathlineto{\pgfqpoint{6.692625in}{5.736806in}}%
\pgfpathlineto{\pgfqpoint{6.899292in}{6.227054in}}%
\pgfpathlineto{\pgfqpoint{7.105958in}{6.173677in}}%
\pgfusepath{stroke}%
\end{pgfscope}%
\begin{pgfscope}%
\pgfpathrectangle{\pgfqpoint{0.905958in}{0.899712in}}{\pgfqpoint{6.200000in}{6.200000in}} %
\pgfusepath{clip}%
\pgfsetbuttcap%
\pgfsetroundjoin%
\definecolor{currentfill}{rgb}{0.298039,0.447059,0.690196}%
\pgfsetfillcolor{currentfill}%
\pgfsetlinewidth{0.000000pt}%
\definecolor{currentstroke}{rgb}{0.000000,0.000000,0.000000}%
\pgfsetstrokecolor{currentstroke}%
\pgfsetdash{}{0pt}%
\pgfsys@defobject{currentmarker}{\pgfqpoint{-0.038889in}{-0.038889in}}{\pgfqpoint{0.038889in}{0.038889in}}{%
\pgfpathmoveto{\pgfqpoint{0.000000in}{-0.038889in}}%
\pgfpathcurveto{\pgfqpoint{0.010313in}{-0.038889in}}{\pgfqpoint{0.020206in}{-0.034791in}}{\pgfqpoint{0.027499in}{-0.027499in}}%
\pgfpathcurveto{\pgfqpoint{0.034791in}{-0.020206in}}{\pgfqpoint{0.038889in}{-0.010313in}}{\pgfqpoint{0.038889in}{0.000000in}}%
\pgfpathcurveto{\pgfqpoint{0.038889in}{0.010313in}}{\pgfqpoint{0.034791in}{0.020206in}}{\pgfqpoint{0.027499in}{0.027499in}}%
\pgfpathcurveto{\pgfqpoint{0.020206in}{0.034791in}}{\pgfqpoint{0.010313in}{0.038889in}}{\pgfqpoint{0.000000in}{0.038889in}}%
\pgfpathcurveto{\pgfqpoint{-0.010313in}{0.038889in}}{\pgfqpoint{-0.020206in}{0.034791in}}{\pgfqpoint{-0.027499in}{0.027499in}}%
\pgfpathcurveto{\pgfqpoint{-0.034791in}{0.020206in}}{\pgfqpoint{-0.038889in}{0.010313in}}{\pgfqpoint{-0.038889in}{0.000000in}}%
\pgfpathcurveto{\pgfqpoint{-0.038889in}{-0.010313in}}{\pgfqpoint{-0.034791in}{-0.020206in}}{\pgfqpoint{-0.027499in}{-0.027499in}}%
\pgfpathcurveto{\pgfqpoint{-0.020206in}{-0.034791in}}{\pgfqpoint{-0.010313in}{-0.038889in}}{\pgfqpoint{0.000000in}{-0.038889in}}%
\pgfpathclose%
\pgfusepath{fill}%
}%
\begin{pgfscope}%
\pgfsys@transformshift{1.732625in}{2.966379in}%
\pgfsys@useobject{currentmarker}{}%
\end{pgfscope}%
\begin{pgfscope}%
\pgfsys@transformshift{1.939292in}{3.545958in}%
\pgfsys@useobject{currentmarker}{}%
\end{pgfscope}%
\begin{pgfscope}%
\pgfsys@transformshift{2.145958in}{3.739679in}%
\pgfsys@useobject{currentmarker}{}%
\end{pgfscope}%
\begin{pgfscope}%
\pgfsys@transformshift{2.352625in}{3.543333in}%
\pgfsys@useobject{currentmarker}{}%
\end{pgfscope}%
\begin{pgfscope}%
\pgfsys@transformshift{2.559292in}{3.646397in}%
\pgfsys@useobject{currentmarker}{}%
\end{pgfscope}%
\begin{pgfscope}%
\pgfsys@transformshift{2.765958in}{3.300143in}%
\pgfsys@useobject{currentmarker}{}%
\end{pgfscope}%
\begin{pgfscope}%
\pgfsys@transformshift{2.972625in}{3.314380in}%
\pgfsys@useobject{currentmarker}{}%
\end{pgfscope}%
\begin{pgfscope}%
\pgfsys@transformshift{3.179292in}{2.774762in}%
\pgfsys@useobject{currentmarker}{}%
\end{pgfscope}%
\begin{pgfscope}%
\pgfsys@transformshift{3.385958in}{2.821213in}%
\pgfsys@useobject{currentmarker}{}%
\end{pgfscope}%
\begin{pgfscope}%
\pgfsys@transformshift{3.592625in}{2.466736in}%
\pgfsys@useobject{currentmarker}{}%
\end{pgfscope}%
\begin{pgfscope}%
\pgfsys@transformshift{3.799292in}{2.163851in}%
\pgfsys@useobject{currentmarker}{}%
\end{pgfscope}%
\begin{pgfscope}%
\pgfsys@transformshift{4.005958in}{2.191379in}%
\pgfsys@useobject{currentmarker}{}%
\end{pgfscope}%
\begin{pgfscope}%
\pgfsys@transformshift{4.212625in}{2.392628in}%
\pgfsys@useobject{currentmarker}{}%
\end{pgfscope}%
\begin{pgfscope}%
\pgfsys@transformshift{4.419292in}{2.457722in}%
\pgfsys@useobject{currentmarker}{}%
\end{pgfscope}%
\begin{pgfscope}%
\pgfsys@transformshift{4.625958in}{2.788667in}%
\pgfsys@useobject{currentmarker}{}%
\end{pgfscope}%
\begin{pgfscope}%
\pgfsys@transformshift{4.832625in}{2.952060in}%
\pgfsys@useobject{currentmarker}{}%
\end{pgfscope}%
\begin{pgfscope}%
\pgfsys@transformshift{5.039292in}{3.184366in}%
\pgfsys@useobject{currentmarker}{}%
\end{pgfscope}%
\begin{pgfscope}%
\pgfsys@transformshift{5.245958in}{3.565762in}%
\pgfsys@useobject{currentmarker}{}%
\end{pgfscope}%
\begin{pgfscope}%
\pgfsys@transformshift{5.452625in}{4.234530in}%
\pgfsys@useobject{currentmarker}{}%
\end{pgfscope}%
\begin{pgfscope}%
\pgfsys@transformshift{5.659292in}{4.213037in}%
\pgfsys@useobject{currentmarker}{}%
\end{pgfscope}%
\begin{pgfscope}%
\pgfsys@transformshift{5.865958in}{4.644712in}%
\pgfsys@useobject{currentmarker}{}%
\end{pgfscope}%
\begin{pgfscope}%
\pgfsys@transformshift{6.072625in}{4.885945in}%
\pgfsys@useobject{currentmarker}{}%
\end{pgfscope}%
\begin{pgfscope}%
\pgfsys@transformshift{6.279292in}{5.420245in}%
\pgfsys@useobject{currentmarker}{}%
\end{pgfscope}%
\begin{pgfscope}%
\pgfsys@transformshift{6.485958in}{5.398937in}%
\pgfsys@useobject{currentmarker}{}%
\end{pgfscope}%
\begin{pgfscope}%
\pgfsys@transformshift{6.692625in}{5.736806in}%
\pgfsys@useobject{currentmarker}{}%
\end{pgfscope}%
\begin{pgfscope}%
\pgfsys@transformshift{6.899292in}{6.227054in}%
\pgfsys@useobject{currentmarker}{}%
\end{pgfscope}%
\begin{pgfscope}%
\pgfsys@transformshift{7.105958in}{6.173677in}%
\pgfsys@useobject{currentmarker}{}%
\end{pgfscope}%
\end{pgfscope}%
\begin{pgfscope}%
\pgfpathrectangle{\pgfqpoint{0.905958in}{0.899712in}}{\pgfqpoint{6.200000in}{6.200000in}} %
\pgfusepath{clip}%
\pgfsetroundcap%
\pgfsetroundjoin%
\pgfsetlinewidth{2.007500pt}%
\definecolor{currentstroke}{rgb}{0.333333,0.658824,0.407843}%
\pgfsetstrokecolor{currentstroke}%
\pgfsetdash{}{0pt}%
\pgfpathmoveto{\pgfqpoint{1.732625in}{2.734519in}}%
\pgfpathlineto{\pgfqpoint{1.939292in}{3.616588in}}%
\pgfpathlineto{\pgfqpoint{2.145958in}{3.435436in}}%
\pgfpathlineto{\pgfqpoint{2.352625in}{3.582002in}}%
\pgfpathlineto{\pgfqpoint{2.559292in}{3.454228in}}%
\pgfpathlineto{\pgfqpoint{2.765958in}{3.340819in}}%
\pgfpathlineto{\pgfqpoint{2.972625in}{2.823274in}}%
\pgfpathlineto{\pgfqpoint{3.179292in}{1.874828in}}%
\pgfpathlineto{\pgfqpoint{3.385958in}{1.642973in}}%
\pgfpathlineto{\pgfqpoint{3.592625in}{1.306363in}}%
\pgfpathlineto{\pgfqpoint{3.799292in}{1.183245in}}%
\pgfpathlineto{\pgfqpoint{4.005958in}{1.102281in}}%
\pgfusepath{stroke}%
\end{pgfscope}%
\begin{pgfscope}%
\pgfpathrectangle{\pgfqpoint{0.905958in}{0.899712in}}{\pgfqpoint{6.200000in}{6.200000in}} %
\pgfusepath{clip}%
\pgfsetbuttcap%
\pgfsetmiterjoin%
\definecolor{currentfill}{rgb}{0.333333,0.658824,0.407843}%
\pgfsetfillcolor{currentfill}%
\pgfsetlinewidth{0.000000pt}%
\definecolor{currentstroke}{rgb}{0.000000,0.000000,0.000000}%
\pgfsetstrokecolor{currentstroke}%
\pgfsetdash{}{0pt}%
\pgfsys@defobject{currentmarker}{\pgfqpoint{-0.038889in}{-0.038889in}}{\pgfqpoint{0.038889in}{0.038889in}}{%
\pgfpathmoveto{\pgfqpoint{-0.038889in}{-0.038889in}}%
\pgfpathlineto{\pgfqpoint{0.038889in}{-0.038889in}}%
\pgfpathlineto{\pgfqpoint{0.038889in}{0.038889in}}%
\pgfpathlineto{\pgfqpoint{-0.038889in}{0.038889in}}%
\pgfpathclose%
\pgfusepath{fill}%
}%
\begin{pgfscope}%
\pgfsys@transformshift{1.732625in}{2.734519in}%
\pgfsys@useobject{currentmarker}{}%
\end{pgfscope}%
\begin{pgfscope}%
\pgfsys@transformshift{1.939292in}{3.616588in}%
\pgfsys@useobject{currentmarker}{}%
\end{pgfscope}%
\begin{pgfscope}%
\pgfsys@transformshift{2.145958in}{3.435436in}%
\pgfsys@useobject{currentmarker}{}%
\end{pgfscope}%
\begin{pgfscope}%
\pgfsys@transformshift{2.352625in}{3.582002in}%
\pgfsys@useobject{currentmarker}{}%
\end{pgfscope}%
\begin{pgfscope}%
\pgfsys@transformshift{2.559292in}{3.454228in}%
\pgfsys@useobject{currentmarker}{}%
\end{pgfscope}%
\begin{pgfscope}%
\pgfsys@transformshift{2.765958in}{3.340819in}%
\pgfsys@useobject{currentmarker}{}%
\end{pgfscope}%
\begin{pgfscope}%
\pgfsys@transformshift{2.972625in}{2.823274in}%
\pgfsys@useobject{currentmarker}{}%
\end{pgfscope}%
\begin{pgfscope}%
\pgfsys@transformshift{3.179292in}{1.874828in}%
\pgfsys@useobject{currentmarker}{}%
\end{pgfscope}%
\begin{pgfscope}%
\pgfsys@transformshift{3.385958in}{1.642973in}%
\pgfsys@useobject{currentmarker}{}%
\end{pgfscope}%
\begin{pgfscope}%
\pgfsys@transformshift{3.592625in}{1.306363in}%
\pgfsys@useobject{currentmarker}{}%
\end{pgfscope}%
\begin{pgfscope}%
\pgfsys@transformshift{3.799292in}{1.183245in}%
\pgfsys@useobject{currentmarker}{}%
\end{pgfscope}%
\begin{pgfscope}%
\pgfsys@transformshift{4.005958in}{1.102281in}%
\pgfsys@useobject{currentmarker}{}%
\end{pgfscope}%
\end{pgfscope}%
\begin{pgfscope}%
\pgfsetrectcap%
\pgfsetmiterjoin%
\pgfsetlinewidth{1.254687pt}%
\definecolor{currentstroke}{rgb}{0.150000,0.150000,0.150000}%
\pgfsetstrokecolor{currentstroke}%
\pgfsetdash{}{0pt}%
\pgfpathmoveto{\pgfqpoint{0.905958in}{0.899712in}}%
\pgfpathlineto{\pgfqpoint{7.105958in}{0.899712in}}%
\pgfusepath{stroke}%
\end{pgfscope}%
\begin{pgfscope}%
\pgfsetrectcap%
\pgfsetmiterjoin%
\pgfsetlinewidth{1.254687pt}%
\definecolor{currentstroke}{rgb}{0.150000,0.150000,0.150000}%
\pgfsetstrokecolor{currentstroke}%
\pgfsetdash{}{0pt}%
\pgfpathmoveto{\pgfqpoint{0.905958in}{0.899712in}}%
\pgfpathlineto{\pgfqpoint{0.905958in}{7.099712in}}%
\pgfusepath{stroke}%
\end{pgfscope}%
\begin{pgfscope}%
\pgfsetbuttcap%
\pgfsetmiterjoin%
\definecolor{currentfill}{rgb}{1.000000,1.000000,1.000000}%
\pgfsetfillcolor{currentfill}%
\pgfsetlinewidth{0.240900pt}%
\definecolor{currentstroke}{rgb}{0.150000,0.150000,0.150000}%
\pgfsetstrokecolor{currentstroke}%
\pgfsetdash{}{0pt}%
\pgfpathmoveto{\pgfqpoint{1.072625in}{5.878293in}}%
\pgfpathlineto{\pgfqpoint{4.476238in}{5.878293in}}%
\pgfpathlineto{\pgfqpoint{4.476238in}{6.933045in}}%
\pgfpathlineto{\pgfqpoint{1.072625in}{6.933045in}}%
\pgfpathclose%
\pgfusepath{stroke,fill}%
\end{pgfscope}%
\begin{pgfscope}%
\pgfsetroundcap%
\pgfsetroundjoin%
\pgfsetlinewidth{2.007500pt}%
\definecolor{currentstroke}{rgb}{0.298039,0.447059,0.690196}%
\pgfsetstrokecolor{currentstroke}%
\pgfsetdash{}{0pt}%
\pgfpathmoveto{\pgfqpoint{1.205958in}{6.674842in}}%
\pgfpathlineto{\pgfqpoint{1.872625in}{6.674842in}}%
\pgfusepath{stroke}%
\end{pgfscope}%
\begin{pgfscope}%
\pgfsetbuttcap%
\pgfsetroundjoin%
\definecolor{currentfill}{rgb}{0.298039,0.447059,0.690196}%
\pgfsetfillcolor{currentfill}%
\pgfsetlinewidth{0.000000pt}%
\definecolor{currentstroke}{rgb}{0.000000,0.000000,0.000000}%
\pgfsetstrokecolor{currentstroke}%
\pgfsetdash{}{0pt}%
\pgfsys@defobject{currentmarker}{\pgfqpoint{-0.038889in}{-0.038889in}}{\pgfqpoint{0.038889in}{0.038889in}}{%
\pgfpathmoveto{\pgfqpoint{0.000000in}{-0.038889in}}%
\pgfpathcurveto{\pgfqpoint{0.010313in}{-0.038889in}}{\pgfqpoint{0.020206in}{-0.034791in}}{\pgfqpoint{0.027499in}{-0.027499in}}%
\pgfpathcurveto{\pgfqpoint{0.034791in}{-0.020206in}}{\pgfqpoint{0.038889in}{-0.010313in}}{\pgfqpoint{0.038889in}{0.000000in}}%
\pgfpathcurveto{\pgfqpoint{0.038889in}{0.010313in}}{\pgfqpoint{0.034791in}{0.020206in}}{\pgfqpoint{0.027499in}{0.027499in}}%
\pgfpathcurveto{\pgfqpoint{0.020206in}{0.034791in}}{\pgfqpoint{0.010313in}{0.038889in}}{\pgfqpoint{0.000000in}{0.038889in}}%
\pgfpathcurveto{\pgfqpoint{-0.010313in}{0.038889in}}{\pgfqpoint{-0.020206in}{0.034791in}}{\pgfqpoint{-0.027499in}{0.027499in}}%
\pgfpathcurveto{\pgfqpoint{-0.034791in}{0.020206in}}{\pgfqpoint{-0.038889in}{0.010313in}}{\pgfqpoint{-0.038889in}{0.000000in}}%
\pgfpathcurveto{\pgfqpoint{-0.038889in}{-0.010313in}}{\pgfqpoint{-0.034791in}{-0.020206in}}{\pgfqpoint{-0.027499in}{-0.027499in}}%
\pgfpathcurveto{\pgfqpoint{-0.020206in}{-0.034791in}}{\pgfqpoint{-0.010313in}{-0.038889in}}{\pgfqpoint{0.000000in}{-0.038889in}}%
\pgfpathclose%
\pgfusepath{fill}%
}%
\begin{pgfscope}%
\pgfsys@transformshift{1.539292in}{6.674842in}%
\pgfsys@useobject{currentmarker}{}%
\end{pgfscope}%
\end{pgfscope}%
\begin{pgfscope}%
\definecolor{textcolor}{rgb}{0.150000,0.150000,0.150000}%
\pgfsetstrokecolor{textcolor}%
\pgfsetfillcolor{textcolor}%
\pgftext[x=2.139292in,y=6.558175in,left,base]{\color{textcolor}\sffamily\fontsize{24.000000}{28.800000}\selectfont FS-QTT-solver}%
\end{pgfscope}%
\begin{pgfscope}%
\pgfsetroundcap%
\pgfsetroundjoin%
\pgfsetlinewidth{2.007500pt}%
\definecolor{currentstroke}{rgb}{0.333333,0.658824,0.407843}%
\pgfsetstrokecolor{currentstroke}%
\pgfsetdash{}{0pt}%
\pgfpathmoveto{\pgfqpoint{1.205958in}{6.197466in}}%
\pgfpathlineto{\pgfqpoint{1.872625in}{6.197466in}}%
\pgfusepath{stroke}%
\end{pgfscope}%
\begin{pgfscope}%
\pgfsetbuttcap%
\pgfsetmiterjoin%
\definecolor{currentfill}{rgb}{0.333333,0.658824,0.407843}%
\pgfsetfillcolor{currentfill}%
\pgfsetlinewidth{0.000000pt}%
\definecolor{currentstroke}{rgb}{0.000000,0.000000,0.000000}%
\pgfsetstrokecolor{currentstroke}%
\pgfsetdash{}{0pt}%
\pgfsys@defobject{currentmarker}{\pgfqpoint{-0.038889in}{-0.038889in}}{\pgfqpoint{0.038889in}{0.038889in}}{%
\pgfpathmoveto{\pgfqpoint{-0.038889in}{-0.038889in}}%
\pgfpathlineto{\pgfqpoint{0.038889in}{-0.038889in}}%
\pgfpathlineto{\pgfqpoint{0.038889in}{0.038889in}}%
\pgfpathlineto{\pgfqpoint{-0.038889in}{0.038889in}}%
\pgfpathclose%
\pgfusepath{fill}%
}%
\begin{pgfscope}%
\pgfsys@transformshift{1.539292in}{6.197466in}%
\pgfsys@useobject{currentmarker}{}%
\end{pgfscope}%
\end{pgfscope}%
\begin{pgfscope}%
\definecolor{textcolor}{rgb}{0.150000,0.150000,0.150000}%
\pgfsetstrokecolor{textcolor}%
\pgfsetfillcolor{textcolor}%
\pgftext[x=2.139292in,y=6.080799in,left,base]{\color{textcolor}\sffamily\fontsize{24.000000}{28.800000}\selectfont FD-QTT-solver}%
\end{pgfscope}%
\end{pgfpicture}%
\makeatother%
\endgroup%

%% file: res_analyt_all_uf_conv.pgf
\begingroup%
\makeatletter%
\begin{pgfpicture}%
\pgfpathrectangle{\pgfpointorigin}{\pgfqpoint{7.781360in}{7.355067in}}%
\pgfusepath{use as bounding box, clip}%
\begin{pgfscope}%
\pgfsetbuttcap%
\pgfsetmiterjoin%
\definecolor{currentfill}{rgb}{1.000000,1.000000,1.000000}%
\pgfsetfillcolor{currentfill}%
\pgfsetlinewidth{0.000000pt}%
\definecolor{currentstroke}{rgb}{1.000000,1.000000,1.000000}%
\pgfsetstrokecolor{currentstroke}%
\pgfsetdash{}{0pt}%
\pgfpathmoveto{\pgfqpoint{0.000000in}{0.000000in}}%
\pgfpathlineto{\pgfqpoint{7.781360in}{0.000000in}}%
\pgfpathlineto{\pgfqpoint{7.781360in}{7.355067in}}%
\pgfpathlineto{\pgfqpoint{0.000000in}{7.355067in}}%
\pgfpathclose%
\pgfusepath{fill}%
\end{pgfscope}%
\begin{pgfscope}%
\pgfsetbuttcap%
\pgfsetmiterjoin%
\definecolor{currentfill}{rgb}{1.000000,1.000000,1.000000}%
\pgfsetfillcolor{currentfill}%
\pgfsetlinewidth{0.000000pt}%
\definecolor{currentstroke}{rgb}{0.000000,0.000000,0.000000}%
\pgfsetstrokecolor{currentstroke}%
\pgfsetstrokeopacity{0.000000}%
\pgfsetdash{}{0pt}%
\pgfpathmoveto{\pgfqpoint{1.322882in}{0.899712in}}%
\pgfpathlineto{\pgfqpoint{7.522882in}{0.899712in}}%
\pgfpathlineto{\pgfqpoint{7.522882in}{7.099712in}}%
\pgfpathlineto{\pgfqpoint{1.322882in}{7.099712in}}%
\pgfpathclose%
\pgfusepath{fill}%
\end{pgfscope}%
\begin{pgfscope}%
\pgfpathrectangle{\pgfqpoint{1.322882in}{0.899712in}}{\pgfqpoint{6.200000in}{6.200000in}} %
\pgfusepath{clip}%
\pgfsetbuttcap%
\pgfsetroundjoin%
\pgfsetlinewidth{0.803000pt}%
\definecolor{currentstroke}{rgb}{0.800000,0.800000,0.800000}%
\pgfsetstrokecolor{currentstroke}%
\pgfsetdash{{1.000000pt}{3.000000pt}}{0.000000pt}%
\pgfpathmoveto{\pgfqpoint{1.322882in}{0.899712in}}%
\pgfpathlineto{\pgfqpoint{1.322882in}{7.099712in}}%
\pgfusepath{stroke}%
\end{pgfscope}%
\begin{pgfscope}%
\definecolor{textcolor}{rgb}{0.150000,0.150000,0.150000}%
\pgfsetstrokecolor{textcolor}%
\pgfsetfillcolor{textcolor}%
\pgftext[x=1.322882in,y=0.821934in,,top]{\color{textcolor}\sffamily\fontsize{24.000000}{28.800000}\selectfont \(\displaystyle 0\)}%
\end{pgfscope}%
\begin{pgfscope}%
\pgfpathrectangle{\pgfqpoint{1.322882in}{0.899712in}}{\pgfqpoint{6.200000in}{6.200000in}} %
\pgfusepath{clip}%
\pgfsetbuttcap%
\pgfsetroundjoin%
\pgfsetlinewidth{0.803000pt}%
\definecolor{currentstroke}{rgb}{0.800000,0.800000,0.800000}%
\pgfsetstrokecolor{currentstroke}%
\pgfsetdash{{1.000000pt}{3.000000pt}}{0.000000pt}%
\pgfpathmoveto{\pgfqpoint{2.356215in}{0.899712in}}%
\pgfpathlineto{\pgfqpoint{2.356215in}{7.099712in}}%
\pgfusepath{stroke}%
\end{pgfscope}%
\begin{pgfscope}%
\definecolor{textcolor}{rgb}{0.150000,0.150000,0.150000}%
\pgfsetstrokecolor{textcolor}%
\pgfsetfillcolor{textcolor}%
\pgftext[x=2.356215in,y=0.821934in,,top]{\color{textcolor}\sffamily\fontsize{24.000000}{28.800000}\selectfont \(\displaystyle 5\)}%
\end{pgfscope}%
\begin{pgfscope}%
\pgfpathrectangle{\pgfqpoint{1.322882in}{0.899712in}}{\pgfqpoint{6.200000in}{6.200000in}} %
\pgfusepath{clip}%
\pgfsetbuttcap%
\pgfsetroundjoin%
\pgfsetlinewidth{0.803000pt}%
\definecolor{currentstroke}{rgb}{0.800000,0.800000,0.800000}%
\pgfsetstrokecolor{currentstroke}%
\pgfsetdash{{1.000000pt}{3.000000pt}}{0.000000pt}%
\pgfpathmoveto{\pgfqpoint{3.389548in}{0.899712in}}%
\pgfpathlineto{\pgfqpoint{3.389548in}{7.099712in}}%
\pgfusepath{stroke}%
\end{pgfscope}%
\begin{pgfscope}%
\definecolor{textcolor}{rgb}{0.150000,0.150000,0.150000}%
\pgfsetstrokecolor{textcolor}%
\pgfsetfillcolor{textcolor}%
\pgftext[x=3.389548in,y=0.821934in,,top]{\color{textcolor}\sffamily\fontsize{24.000000}{28.800000}\selectfont \(\displaystyle 10\)}%
\end{pgfscope}%
\begin{pgfscope}%
\pgfpathrectangle{\pgfqpoint{1.322882in}{0.899712in}}{\pgfqpoint{6.200000in}{6.200000in}} %
\pgfusepath{clip}%
\pgfsetbuttcap%
\pgfsetroundjoin%
\pgfsetlinewidth{0.803000pt}%
\definecolor{currentstroke}{rgb}{0.800000,0.800000,0.800000}%
\pgfsetstrokecolor{currentstroke}%
\pgfsetdash{{1.000000pt}{3.000000pt}}{0.000000pt}%
\pgfpathmoveto{\pgfqpoint{4.422882in}{0.899712in}}%
\pgfpathlineto{\pgfqpoint{4.422882in}{7.099712in}}%
\pgfusepath{stroke}%
\end{pgfscope}%
\begin{pgfscope}%
\definecolor{textcolor}{rgb}{0.150000,0.150000,0.150000}%
\pgfsetstrokecolor{textcolor}%
\pgfsetfillcolor{textcolor}%
\pgftext[x=4.422882in,y=0.821934in,,top]{\color{textcolor}\sffamily\fontsize{24.000000}{28.800000}\selectfont \(\displaystyle 15\)}%
\end{pgfscope}%
\begin{pgfscope}%
\pgfpathrectangle{\pgfqpoint{1.322882in}{0.899712in}}{\pgfqpoint{6.200000in}{6.200000in}} %
\pgfusepath{clip}%
\pgfsetbuttcap%
\pgfsetroundjoin%
\pgfsetlinewidth{0.803000pt}%
\definecolor{currentstroke}{rgb}{0.800000,0.800000,0.800000}%
\pgfsetstrokecolor{currentstroke}%
\pgfsetdash{{1.000000pt}{3.000000pt}}{0.000000pt}%
\pgfpathmoveto{\pgfqpoint{5.456215in}{0.899712in}}%
\pgfpathlineto{\pgfqpoint{5.456215in}{7.099712in}}%
\pgfusepath{stroke}%
\end{pgfscope}%
\begin{pgfscope}%
\definecolor{textcolor}{rgb}{0.150000,0.150000,0.150000}%
\pgfsetstrokecolor{textcolor}%
\pgfsetfillcolor{textcolor}%
\pgftext[x=5.456215in,y=0.821934in,,top]{\color{textcolor}\sffamily\fontsize{24.000000}{28.800000}\selectfont \(\displaystyle 20\)}%
\end{pgfscope}%
\begin{pgfscope}%
\pgfpathrectangle{\pgfqpoint{1.322882in}{0.899712in}}{\pgfqpoint{6.200000in}{6.200000in}} %
\pgfusepath{clip}%
\pgfsetbuttcap%
\pgfsetroundjoin%
\pgfsetlinewidth{0.803000pt}%
\definecolor{currentstroke}{rgb}{0.800000,0.800000,0.800000}%
\pgfsetstrokecolor{currentstroke}%
\pgfsetdash{{1.000000pt}{3.000000pt}}{0.000000pt}%
\pgfpathmoveto{\pgfqpoint{6.489548in}{0.899712in}}%
\pgfpathlineto{\pgfqpoint{6.489548in}{7.099712in}}%
\pgfusepath{stroke}%
\end{pgfscope}%
\begin{pgfscope}%
\definecolor{textcolor}{rgb}{0.150000,0.150000,0.150000}%
\pgfsetstrokecolor{textcolor}%
\pgfsetfillcolor{textcolor}%
\pgftext[x=6.489548in,y=0.821934in,,top]{\color{textcolor}\sffamily\fontsize{24.000000}{28.800000}\selectfont \(\displaystyle 25\)}%
\end{pgfscope}%
\begin{pgfscope}%
\pgfpathrectangle{\pgfqpoint{1.322882in}{0.899712in}}{\pgfqpoint{6.200000in}{6.200000in}} %
\pgfusepath{clip}%
\pgfsetbuttcap%
\pgfsetroundjoin%
\pgfsetlinewidth{0.803000pt}%
\definecolor{currentstroke}{rgb}{0.800000,0.800000,0.800000}%
\pgfsetstrokecolor{currentstroke}%
\pgfsetdash{{1.000000pt}{3.000000pt}}{0.000000pt}%
\pgfpathmoveto{\pgfqpoint{7.522882in}{0.899712in}}%
\pgfpathlineto{\pgfqpoint{7.522882in}{7.099712in}}%
\pgfusepath{stroke}%
\end{pgfscope}%
\begin{pgfscope}%
\definecolor{textcolor}{rgb}{0.150000,0.150000,0.150000}%
\pgfsetstrokecolor{textcolor}%
\pgfsetfillcolor{textcolor}%
\pgftext[x=7.522882in,y=0.821934in,,top]{\color{textcolor}\sffamily\fontsize{24.000000}{28.800000}\selectfont \(\displaystyle 30\)}%
\end{pgfscope}%
\begin{pgfscope}%
\definecolor{textcolor}{rgb}{0.150000,0.150000,0.150000}%
\pgfsetstrokecolor{textcolor}%
\pgfsetfillcolor{textcolor}%
\pgftext[x=4.422882in,y=0.441780in,,top]{\color{textcolor}\sffamily\fontsize{26.400000}{31.680000}\selectfont d}%
\end{pgfscope}%
\begin{pgfscope}%
\pgfpathrectangle{\pgfqpoint{1.322882in}{0.899712in}}{\pgfqpoint{6.200000in}{6.200000in}} %
\pgfusepath{clip}%
\pgfsetbuttcap%
\pgfsetroundjoin%
\pgfsetlinewidth{0.803000pt}%
\definecolor{currentstroke}{rgb}{0.800000,0.800000,0.800000}%
\pgfsetstrokecolor{currentstroke}%
\pgfsetdash{{1.000000pt}{3.000000pt}}{0.000000pt}%
\pgfpathmoveto{\pgfqpoint{1.322882in}{0.899712in}}%
\pgfpathlineto{\pgfqpoint{7.522882in}{0.899712in}}%
\pgfusepath{stroke}%
\end{pgfscope}%
\begin{pgfscope}%
\definecolor{textcolor}{rgb}{0.150000,0.150000,0.150000}%
\pgfsetstrokecolor{textcolor}%
\pgfsetfillcolor{textcolor}%
\pgftext[x=1.245104in,y=0.899712in,right,]{\color{textcolor}\sffamily\fontsize{24.000000}{28.800000}\selectfont \(\displaystyle 10^{-12}\)}%
\end{pgfscope}%
\begin{pgfscope}%
\pgfpathrectangle{\pgfqpoint{1.322882in}{0.899712in}}{\pgfqpoint{6.200000in}{6.200000in}} %
\pgfusepath{clip}%
\pgfsetbuttcap%
\pgfsetroundjoin%
\pgfsetlinewidth{0.803000pt}%
\definecolor{currentstroke}{rgb}{0.800000,0.800000,0.800000}%
\pgfsetstrokecolor{currentstroke}%
\pgfsetdash{{1.000000pt}{3.000000pt}}{0.000000pt}%
\pgfpathmoveto{\pgfqpoint{1.322882in}{1.416379in}}%
\pgfpathlineto{\pgfqpoint{7.522882in}{1.416379in}}%
\pgfusepath{stroke}%
\end{pgfscope}%
\begin{pgfscope}%
\definecolor{textcolor}{rgb}{0.150000,0.150000,0.150000}%
\pgfsetstrokecolor{textcolor}%
\pgfsetfillcolor{textcolor}%
\pgftext[x=1.245104in,y=1.416379in,right,]{\color{textcolor}\sffamily\fontsize{24.000000}{28.800000}\selectfont \(\displaystyle 10^{-11}\)}%
\end{pgfscope}%
\begin{pgfscope}%
\pgfpathrectangle{\pgfqpoint{1.322882in}{0.899712in}}{\pgfqpoint{6.200000in}{6.200000in}} %
\pgfusepath{clip}%
\pgfsetbuttcap%
\pgfsetroundjoin%
\pgfsetlinewidth{0.803000pt}%
\definecolor{currentstroke}{rgb}{0.800000,0.800000,0.800000}%
\pgfsetstrokecolor{currentstroke}%
\pgfsetdash{{1.000000pt}{3.000000pt}}{0.000000pt}%
\pgfpathmoveto{\pgfqpoint{1.322882in}{1.933045in}}%
\pgfpathlineto{\pgfqpoint{7.522882in}{1.933045in}}%
\pgfusepath{stroke}%
\end{pgfscope}%
\begin{pgfscope}%
\definecolor{textcolor}{rgb}{0.150000,0.150000,0.150000}%
\pgfsetstrokecolor{textcolor}%
\pgfsetfillcolor{textcolor}%
\pgftext[x=1.245104in,y=1.933045in,right,]{\color{textcolor}\sffamily\fontsize{24.000000}{28.800000}\selectfont \(\displaystyle 10^{-10}\)}%
\end{pgfscope}%
\begin{pgfscope}%
\pgfpathrectangle{\pgfqpoint{1.322882in}{0.899712in}}{\pgfqpoint{6.200000in}{6.200000in}} %
\pgfusepath{clip}%
\pgfsetbuttcap%
\pgfsetroundjoin%
\pgfsetlinewidth{0.803000pt}%
\definecolor{currentstroke}{rgb}{0.800000,0.800000,0.800000}%
\pgfsetstrokecolor{currentstroke}%
\pgfsetdash{{1.000000pt}{3.000000pt}}{0.000000pt}%
\pgfpathmoveto{\pgfqpoint{1.322882in}{2.449712in}}%
\pgfpathlineto{\pgfqpoint{7.522882in}{2.449712in}}%
\pgfusepath{stroke}%
\end{pgfscope}%
\begin{pgfscope}%
\definecolor{textcolor}{rgb}{0.150000,0.150000,0.150000}%
\pgfsetstrokecolor{textcolor}%
\pgfsetfillcolor{textcolor}%
\pgftext[x=1.245104in,y=2.449712in,right,]{\color{textcolor}\sffamily\fontsize{24.000000}{28.800000}\selectfont \(\displaystyle 10^{-9}\)}%
\end{pgfscope}%
\begin{pgfscope}%
\pgfpathrectangle{\pgfqpoint{1.322882in}{0.899712in}}{\pgfqpoint{6.200000in}{6.200000in}} %
\pgfusepath{clip}%
\pgfsetbuttcap%
\pgfsetroundjoin%
\pgfsetlinewidth{0.803000pt}%
\definecolor{currentstroke}{rgb}{0.800000,0.800000,0.800000}%
\pgfsetstrokecolor{currentstroke}%
\pgfsetdash{{1.000000pt}{3.000000pt}}{0.000000pt}%
\pgfpathmoveto{\pgfqpoint{1.322882in}{2.966379in}}%
\pgfpathlineto{\pgfqpoint{7.522882in}{2.966379in}}%
\pgfusepath{stroke}%
\end{pgfscope}%
\begin{pgfscope}%
\definecolor{textcolor}{rgb}{0.150000,0.150000,0.150000}%
\pgfsetstrokecolor{textcolor}%
\pgfsetfillcolor{textcolor}%
\pgftext[x=1.245104in,y=2.966379in,right,]{\color{textcolor}\sffamily\fontsize{24.000000}{28.800000}\selectfont \(\displaystyle 10^{-8}\)}%
\end{pgfscope}%
\begin{pgfscope}%
\pgfpathrectangle{\pgfqpoint{1.322882in}{0.899712in}}{\pgfqpoint{6.200000in}{6.200000in}} %
\pgfusepath{clip}%
\pgfsetbuttcap%
\pgfsetroundjoin%
\pgfsetlinewidth{0.803000pt}%
\definecolor{currentstroke}{rgb}{0.800000,0.800000,0.800000}%
\pgfsetstrokecolor{currentstroke}%
\pgfsetdash{{1.000000pt}{3.000000pt}}{0.000000pt}%
\pgfpathmoveto{\pgfqpoint{1.322882in}{3.483045in}}%
\pgfpathlineto{\pgfqpoint{7.522882in}{3.483045in}}%
\pgfusepath{stroke}%
\end{pgfscope}%
\begin{pgfscope}%
\definecolor{textcolor}{rgb}{0.150000,0.150000,0.150000}%
\pgfsetstrokecolor{textcolor}%
\pgfsetfillcolor{textcolor}%
\pgftext[x=1.245104in,y=3.483045in,right,]{\color{textcolor}\sffamily\fontsize{24.000000}{28.800000}\selectfont \(\displaystyle 10^{-7}\)}%
\end{pgfscope}%
\begin{pgfscope}%
\pgfpathrectangle{\pgfqpoint{1.322882in}{0.899712in}}{\pgfqpoint{6.200000in}{6.200000in}} %
\pgfusepath{clip}%
\pgfsetbuttcap%
\pgfsetroundjoin%
\pgfsetlinewidth{0.803000pt}%
\definecolor{currentstroke}{rgb}{0.800000,0.800000,0.800000}%
\pgfsetstrokecolor{currentstroke}%
\pgfsetdash{{1.000000pt}{3.000000pt}}{0.000000pt}%
\pgfpathmoveto{\pgfqpoint{1.322882in}{3.999712in}}%
\pgfpathlineto{\pgfqpoint{7.522882in}{3.999712in}}%
\pgfusepath{stroke}%
\end{pgfscope}%
\begin{pgfscope}%
\definecolor{textcolor}{rgb}{0.150000,0.150000,0.150000}%
\pgfsetstrokecolor{textcolor}%
\pgfsetfillcolor{textcolor}%
\pgftext[x=1.245104in,y=3.999712in,right,]{\color{textcolor}\sffamily\fontsize{24.000000}{28.800000}\selectfont \(\displaystyle 10^{-6}\)}%
\end{pgfscope}%
\begin{pgfscope}%
\pgfpathrectangle{\pgfqpoint{1.322882in}{0.899712in}}{\pgfqpoint{6.200000in}{6.200000in}} %
\pgfusepath{clip}%
\pgfsetbuttcap%
\pgfsetroundjoin%
\pgfsetlinewidth{0.803000pt}%
\definecolor{currentstroke}{rgb}{0.800000,0.800000,0.800000}%
\pgfsetstrokecolor{currentstroke}%
\pgfsetdash{{1.000000pt}{3.000000pt}}{0.000000pt}%
\pgfpathmoveto{\pgfqpoint{1.322882in}{4.516379in}}%
\pgfpathlineto{\pgfqpoint{7.522882in}{4.516379in}}%
\pgfusepath{stroke}%
\end{pgfscope}%
\begin{pgfscope}%
\definecolor{textcolor}{rgb}{0.150000,0.150000,0.150000}%
\pgfsetstrokecolor{textcolor}%
\pgfsetfillcolor{textcolor}%
\pgftext[x=1.245104in,y=4.516379in,right,]{\color{textcolor}\sffamily\fontsize{24.000000}{28.800000}\selectfont \(\displaystyle 10^{-5}\)}%
\end{pgfscope}%
\begin{pgfscope}%
\pgfpathrectangle{\pgfqpoint{1.322882in}{0.899712in}}{\pgfqpoint{6.200000in}{6.200000in}} %
\pgfusepath{clip}%
\pgfsetbuttcap%
\pgfsetroundjoin%
\pgfsetlinewidth{0.803000pt}%
\definecolor{currentstroke}{rgb}{0.800000,0.800000,0.800000}%
\pgfsetstrokecolor{currentstroke}%
\pgfsetdash{{1.000000pt}{3.000000pt}}{0.000000pt}%
\pgfpathmoveto{\pgfqpoint{1.322882in}{5.033045in}}%
\pgfpathlineto{\pgfqpoint{7.522882in}{5.033045in}}%
\pgfusepath{stroke}%
\end{pgfscope}%
\begin{pgfscope}%
\definecolor{textcolor}{rgb}{0.150000,0.150000,0.150000}%
\pgfsetstrokecolor{textcolor}%
\pgfsetfillcolor{textcolor}%
\pgftext[x=1.245104in,y=5.033045in,right,]{\color{textcolor}\sffamily\fontsize{24.000000}{28.800000}\selectfont \(\displaystyle 10^{-4}\)}%
\end{pgfscope}%
\begin{pgfscope}%
\pgfpathrectangle{\pgfqpoint{1.322882in}{0.899712in}}{\pgfqpoint{6.200000in}{6.200000in}} %
\pgfusepath{clip}%
\pgfsetbuttcap%
\pgfsetroundjoin%
\pgfsetlinewidth{0.803000pt}%
\definecolor{currentstroke}{rgb}{0.800000,0.800000,0.800000}%
\pgfsetstrokecolor{currentstroke}%
\pgfsetdash{{1.000000pt}{3.000000pt}}{0.000000pt}%
\pgfpathmoveto{\pgfqpoint{1.322882in}{5.549712in}}%
\pgfpathlineto{\pgfqpoint{7.522882in}{5.549712in}}%
\pgfusepath{stroke}%
\end{pgfscope}%
\begin{pgfscope}%
\definecolor{textcolor}{rgb}{0.150000,0.150000,0.150000}%
\pgfsetstrokecolor{textcolor}%
\pgfsetfillcolor{textcolor}%
\pgftext[x=1.245104in,y=5.549712in,right,]{\color{textcolor}\sffamily\fontsize{24.000000}{28.800000}\selectfont \(\displaystyle 10^{-3}\)}%
\end{pgfscope}%
\begin{pgfscope}%
\pgfpathrectangle{\pgfqpoint{1.322882in}{0.899712in}}{\pgfqpoint{6.200000in}{6.200000in}} %
\pgfusepath{clip}%
\pgfsetbuttcap%
\pgfsetroundjoin%
\pgfsetlinewidth{0.803000pt}%
\definecolor{currentstroke}{rgb}{0.800000,0.800000,0.800000}%
\pgfsetstrokecolor{currentstroke}%
\pgfsetdash{{1.000000pt}{3.000000pt}}{0.000000pt}%
\pgfpathmoveto{\pgfqpoint{1.322882in}{6.066379in}}%
\pgfpathlineto{\pgfqpoint{7.522882in}{6.066379in}}%
\pgfusepath{stroke}%
\end{pgfscope}%
\begin{pgfscope}%
\definecolor{textcolor}{rgb}{0.150000,0.150000,0.150000}%
\pgfsetstrokecolor{textcolor}%
\pgfsetfillcolor{textcolor}%
\pgftext[x=1.245104in,y=6.066379in,right,]{\color{textcolor}\sffamily\fontsize{24.000000}{28.800000}\selectfont \(\displaystyle 10^{-2}\)}%
\end{pgfscope}%
\begin{pgfscope}%
\pgfpathrectangle{\pgfqpoint{1.322882in}{0.899712in}}{\pgfqpoint{6.200000in}{6.200000in}} %
\pgfusepath{clip}%
\pgfsetbuttcap%
\pgfsetroundjoin%
\pgfsetlinewidth{0.803000pt}%
\definecolor{currentstroke}{rgb}{0.800000,0.800000,0.800000}%
\pgfsetstrokecolor{currentstroke}%
\pgfsetdash{{1.000000pt}{3.000000pt}}{0.000000pt}%
\pgfpathmoveto{\pgfqpoint{1.322882in}{6.583045in}}%
\pgfpathlineto{\pgfqpoint{7.522882in}{6.583045in}}%
\pgfusepath{stroke}%
\end{pgfscope}%
\begin{pgfscope}%
\definecolor{textcolor}{rgb}{0.150000,0.150000,0.150000}%
\pgfsetstrokecolor{textcolor}%
\pgfsetfillcolor{textcolor}%
\pgftext[x=1.245104in,y=6.583045in,right,]{\color{textcolor}\sffamily\fontsize{24.000000}{28.800000}\selectfont \(\displaystyle 10^{-1}\)}%
\end{pgfscope}%
\begin{pgfscope}%
\pgfpathrectangle{\pgfqpoint{1.322882in}{0.899712in}}{\pgfqpoint{6.200000in}{6.200000in}} %
\pgfusepath{clip}%
\pgfsetbuttcap%
\pgfsetroundjoin%
\pgfsetlinewidth{0.803000pt}%
\definecolor{currentstroke}{rgb}{0.800000,0.800000,0.800000}%
\pgfsetstrokecolor{currentstroke}%
\pgfsetdash{{1.000000pt}{3.000000pt}}{0.000000pt}%
\pgfpathmoveto{\pgfqpoint{1.322882in}{7.099712in}}%
\pgfpathlineto{\pgfqpoint{7.522882in}{7.099712in}}%
\pgfusepath{stroke}%
\end{pgfscope}%
\begin{pgfscope}%
\definecolor{textcolor}{rgb}{0.150000,0.150000,0.150000}%
\pgfsetstrokecolor{textcolor}%
\pgfsetfillcolor{textcolor}%
\pgftext[x=1.245104in,y=7.099712in,right,]{\color{textcolor}\sffamily\fontsize{24.000000}{28.800000}\selectfont \(\displaystyle 10^{0}\)}%
\end{pgfscope}%
\begin{pgfscope}%
\definecolor{textcolor}{rgb}{0.150000,0.150000,0.150000}%
\pgfsetstrokecolor{textcolor}%
\pgfsetfillcolor{textcolor}%
\pgftext[x=0.444955in,y=3.999712in,,bottom,rotate=90.000000]{\color{textcolor}\sffamily\fontsize{26.400000}{31.680000}\selectfont \(\displaystyle (u, f) - (u, f)_{RE}\)}%
\end{pgfscope}%
\begin{pgfscope}%
\pgfpathrectangle{\pgfqpoint{1.322882in}{0.899712in}}{\pgfqpoint{6.200000in}{6.200000in}} %
\pgfusepath{clip}%
\pgfsetroundcap%
\pgfsetroundjoin%
\pgfsetlinewidth{2.007500pt}%
\definecolor{currentstroke}{rgb}{0.298039,0.447059,0.690196}%
\pgfsetstrokecolor{currentstroke}%
\pgfsetdash{}{0pt}%
\pgfpathmoveto{\pgfqpoint{2.149548in}{6.634408in}}%
\pgfpathlineto{\pgfqpoint{2.356215in}{6.324276in}}%
\pgfpathlineto{\pgfqpoint{2.562882in}{6.013428in}}%
\pgfpathlineto{\pgfqpoint{2.769548in}{5.702417in}}%
\pgfpathlineto{\pgfqpoint{2.976215in}{5.391366in}}%
\pgfpathlineto{\pgfqpoint{3.182882in}{5.080305in}}%
\pgfpathlineto{\pgfqpoint{3.389548in}{4.769241in}}%
\pgfpathlineto{\pgfqpoint{3.596215in}{4.458173in}}%
\pgfpathlineto{\pgfqpoint{3.802882in}{4.147095in}}%
\pgfpathlineto{\pgfqpoint{4.009548in}{3.835969in}}%
\pgfpathlineto{\pgfqpoint{4.216215in}{3.524665in}}%
\pgfpathlineto{\pgfqpoint{4.422882in}{3.212732in}}%
\pgfpathlineto{\pgfqpoint{4.629548in}{2.897510in}}%
\pgfpathlineto{\pgfqpoint{4.836215in}{2.565120in}}%
\pgfpathlineto{\pgfqpoint{5.042882in}{2.117068in}}%
\pgfpathlineto{\pgfqpoint{5.249548in}{2.087102in}}%
\pgfpathlineto{\pgfqpoint{5.456215in}{2.156280in}}%
\pgfpathlineto{\pgfqpoint{5.662882in}{2.078424in}}%
\pgfpathlineto{\pgfqpoint{5.869548in}{1.931023in}}%
\pgfpathlineto{\pgfqpoint{6.076215in}{1.919736in}}%
\pgfpathlineto{\pgfqpoint{6.282882in}{1.278931in}}%
\pgfpathlineto{\pgfqpoint{6.489548in}{1.954701in}}%
\pgfpathlineto{\pgfqpoint{6.696215in}{1.921543in}}%
\pgfpathlineto{\pgfqpoint{6.902882in}{2.593943in}}%
\pgfpathlineto{\pgfqpoint{7.109548in}{2.752096in}}%
\pgfpathlineto{\pgfqpoint{7.316215in}{2.784617in}}%
\pgfpathlineto{\pgfqpoint{7.522882in}{3.104313in}}%
\pgfusepath{stroke}%
\end{pgfscope}%
\begin{pgfscope}%
\pgfpathrectangle{\pgfqpoint{1.322882in}{0.899712in}}{\pgfqpoint{6.200000in}{6.200000in}} %
\pgfusepath{clip}%
\pgfsetbuttcap%
\pgfsetroundjoin%
\definecolor{currentfill}{rgb}{0.298039,0.447059,0.690196}%
\pgfsetfillcolor{currentfill}%
\pgfsetlinewidth{0.000000pt}%
\definecolor{currentstroke}{rgb}{0.000000,0.000000,0.000000}%
\pgfsetstrokecolor{currentstroke}%
\pgfsetdash{}{0pt}%
\pgfsys@defobject{currentmarker}{\pgfqpoint{-0.038889in}{-0.038889in}}{\pgfqpoint{0.038889in}{0.038889in}}{%
\pgfpathmoveto{\pgfqpoint{0.000000in}{-0.038889in}}%
\pgfpathcurveto{\pgfqpoint{0.010313in}{-0.038889in}}{\pgfqpoint{0.020206in}{-0.034791in}}{\pgfqpoint{0.027499in}{-0.027499in}}%
\pgfpathcurveto{\pgfqpoint{0.034791in}{-0.020206in}}{\pgfqpoint{0.038889in}{-0.010313in}}{\pgfqpoint{0.038889in}{0.000000in}}%
\pgfpathcurveto{\pgfqpoint{0.038889in}{0.010313in}}{\pgfqpoint{0.034791in}{0.020206in}}{\pgfqpoint{0.027499in}{0.027499in}}%
\pgfpathcurveto{\pgfqpoint{0.020206in}{0.034791in}}{\pgfqpoint{0.010313in}{0.038889in}}{\pgfqpoint{0.000000in}{0.038889in}}%
\pgfpathcurveto{\pgfqpoint{-0.010313in}{0.038889in}}{\pgfqpoint{-0.020206in}{0.034791in}}{\pgfqpoint{-0.027499in}{0.027499in}}%
\pgfpathcurveto{\pgfqpoint{-0.034791in}{0.020206in}}{\pgfqpoint{-0.038889in}{0.010313in}}{\pgfqpoint{-0.038889in}{0.000000in}}%
\pgfpathcurveto{\pgfqpoint{-0.038889in}{-0.010313in}}{\pgfqpoint{-0.034791in}{-0.020206in}}{\pgfqpoint{-0.027499in}{-0.027499in}}%
\pgfpathcurveto{\pgfqpoint{-0.020206in}{-0.034791in}}{\pgfqpoint{-0.010313in}{-0.038889in}}{\pgfqpoint{0.000000in}{-0.038889in}}%
\pgfpathclose%
\pgfusepath{fill}%
}%
\begin{pgfscope}%
\pgfsys@transformshift{2.149548in}{6.634408in}%
\pgfsys@useobject{currentmarker}{}%
\end{pgfscope}%
\begin{pgfscope}%
\pgfsys@transformshift{2.356215in}{6.324276in}%
\pgfsys@useobject{currentmarker}{}%
\end{pgfscope}%
\begin{pgfscope}%
\pgfsys@transformshift{2.562882in}{6.013428in}%
\pgfsys@useobject{currentmarker}{}%
\end{pgfscope}%
\begin{pgfscope}%
\pgfsys@transformshift{2.769548in}{5.702417in}%
\pgfsys@useobject{currentmarker}{}%
\end{pgfscope}%
\begin{pgfscope}%
\pgfsys@transformshift{2.976215in}{5.391366in}%
\pgfsys@useobject{currentmarker}{}%
\end{pgfscope}%
\begin{pgfscope}%
\pgfsys@transformshift{3.182882in}{5.080305in}%
\pgfsys@useobject{currentmarker}{}%
\end{pgfscope}%
\begin{pgfscope}%
\pgfsys@transformshift{3.389548in}{4.769241in}%
\pgfsys@useobject{currentmarker}{}%
\end{pgfscope}%
\begin{pgfscope}%
\pgfsys@transformshift{3.596215in}{4.458173in}%
\pgfsys@useobject{currentmarker}{}%
\end{pgfscope}%
\begin{pgfscope}%
\pgfsys@transformshift{3.802882in}{4.147095in}%
\pgfsys@useobject{currentmarker}{}%
\end{pgfscope}%
\begin{pgfscope}%
\pgfsys@transformshift{4.009548in}{3.835969in}%
\pgfsys@useobject{currentmarker}{}%
\end{pgfscope}%
\begin{pgfscope}%
\pgfsys@transformshift{4.216215in}{3.524665in}%
\pgfsys@useobject{currentmarker}{}%
\end{pgfscope}%
\begin{pgfscope}%
\pgfsys@transformshift{4.422882in}{3.212732in}%
\pgfsys@useobject{currentmarker}{}%
\end{pgfscope}%
\begin{pgfscope}%
\pgfsys@transformshift{4.629548in}{2.897510in}%
\pgfsys@useobject{currentmarker}{}%
\end{pgfscope}%
\begin{pgfscope}%
\pgfsys@transformshift{4.836215in}{2.565120in}%
\pgfsys@useobject{currentmarker}{}%
\end{pgfscope}%
\begin{pgfscope}%
\pgfsys@transformshift{5.042882in}{2.117068in}%
\pgfsys@useobject{currentmarker}{}%
\end{pgfscope}%
\begin{pgfscope}%
\pgfsys@transformshift{5.249548in}{2.087102in}%
\pgfsys@useobject{currentmarker}{}%
\end{pgfscope}%
\begin{pgfscope}%
\pgfsys@transformshift{5.456215in}{2.156280in}%
\pgfsys@useobject{currentmarker}{}%
\end{pgfscope}%
\begin{pgfscope}%
\pgfsys@transformshift{5.662882in}{2.078424in}%
\pgfsys@useobject{currentmarker}{}%
\end{pgfscope}%
\begin{pgfscope}%
\pgfsys@transformshift{5.869548in}{1.931023in}%
\pgfsys@useobject{currentmarker}{}%
\end{pgfscope}%
\begin{pgfscope}%
\pgfsys@transformshift{6.076215in}{1.919736in}%
\pgfsys@useobject{currentmarker}{}%
\end{pgfscope}%
\begin{pgfscope}%
\pgfsys@transformshift{6.282882in}{1.278931in}%
\pgfsys@useobject{currentmarker}{}%
\end{pgfscope}%
\begin{pgfscope}%
\pgfsys@transformshift{6.489548in}{1.954701in}%
\pgfsys@useobject{currentmarker}{}%
\end{pgfscope}%
\begin{pgfscope}%
\pgfsys@transformshift{6.696215in}{1.921543in}%
\pgfsys@useobject{currentmarker}{}%
\end{pgfscope}%
\begin{pgfscope}%
\pgfsys@transformshift{6.902882in}{2.593943in}%
\pgfsys@useobject{currentmarker}{}%
\end{pgfscope}%
\begin{pgfscope}%
\pgfsys@transformshift{7.109548in}{2.752096in}%
\pgfsys@useobject{currentmarker}{}%
\end{pgfscope}%
\begin{pgfscope}%
\pgfsys@transformshift{7.316215in}{2.784617in}%
\pgfsys@useobject{currentmarker}{}%
\end{pgfscope}%
\begin{pgfscope}%
\pgfsys@transformshift{7.522882in}{3.104313in}%
\pgfsys@useobject{currentmarker}{}%
\end{pgfscope}%
\end{pgfscope}%
\begin{pgfscope}%
\pgfpathrectangle{\pgfqpoint{1.322882in}{0.899712in}}{\pgfqpoint{6.200000in}{6.200000in}} %
\pgfusepath{clip}%
\pgfsetroundcap%
\pgfsetroundjoin%
\pgfsetlinewidth{2.007500pt}%
\definecolor{currentstroke}{rgb}{0.333333,0.658824,0.407843}%
\pgfsetstrokecolor{currentstroke}%
\pgfsetdash{}{0pt}%
\pgfpathmoveto{\pgfqpoint{2.149548in}{6.634408in}}%
\pgfpathlineto{\pgfqpoint{2.356215in}{6.324276in}}%
\pgfpathlineto{\pgfqpoint{2.562882in}{6.013427in}}%
\pgfpathlineto{\pgfqpoint{2.769548in}{5.702413in}}%
\pgfpathlineto{\pgfqpoint{2.976215in}{5.391349in}}%
\pgfpathlineto{\pgfqpoint{3.182882in}{5.080220in}}%
\pgfpathlineto{\pgfqpoint{3.389548in}{4.768758in}}%
\pgfpathlineto{\pgfqpoint{3.596215in}{4.452410in}}%
\pgfpathlineto{\pgfqpoint{3.802882in}{4.087784in}}%
\pgfpathlineto{\pgfqpoint{4.009548in}{4.009265in}}%
\pgfpathlineto{\pgfqpoint{4.216215in}{4.393025in}}%
\pgfpathlineto{\pgfqpoint{4.422882in}{4.302041in}}%
\pgfusepath{stroke}%
\end{pgfscope}%
\begin{pgfscope}%
\pgfpathrectangle{\pgfqpoint{1.322882in}{0.899712in}}{\pgfqpoint{6.200000in}{6.200000in}} %
\pgfusepath{clip}%
\pgfsetbuttcap%
\pgfsetmiterjoin%
\definecolor{currentfill}{rgb}{0.333333,0.658824,0.407843}%
\pgfsetfillcolor{currentfill}%
\pgfsetlinewidth{0.000000pt}%
\definecolor{currentstroke}{rgb}{0.000000,0.000000,0.000000}%
\pgfsetstrokecolor{currentstroke}%
\pgfsetdash{}{0pt}%
\pgfsys@defobject{currentmarker}{\pgfqpoint{-0.038889in}{-0.038889in}}{\pgfqpoint{0.038889in}{0.038889in}}{%
\pgfpathmoveto{\pgfqpoint{-0.038889in}{-0.038889in}}%
\pgfpathlineto{\pgfqpoint{0.038889in}{-0.038889in}}%
\pgfpathlineto{\pgfqpoint{0.038889in}{0.038889in}}%
\pgfpathlineto{\pgfqpoint{-0.038889in}{0.038889in}}%
\pgfpathclose%
\pgfusepath{fill}%
}%
\begin{pgfscope}%
\pgfsys@transformshift{2.149548in}{6.634408in}%
\pgfsys@useobject{currentmarker}{}%
\end{pgfscope}%
\begin{pgfscope}%
\pgfsys@transformshift{2.356215in}{6.324276in}%
\pgfsys@useobject{currentmarker}{}%
\end{pgfscope}%
\begin{pgfscope}%
\pgfsys@transformshift{2.562882in}{6.013427in}%
\pgfsys@useobject{currentmarker}{}%
\end{pgfscope}%
\begin{pgfscope}%
\pgfsys@transformshift{2.769548in}{5.702413in}%
\pgfsys@useobject{currentmarker}{}%
\end{pgfscope}%
\begin{pgfscope}%
\pgfsys@transformshift{2.976215in}{5.391349in}%
\pgfsys@useobject{currentmarker}{}%
\end{pgfscope}%
\begin{pgfscope}%
\pgfsys@transformshift{3.182882in}{5.080220in}%
\pgfsys@useobject{currentmarker}{}%
\end{pgfscope}%
\begin{pgfscope}%
\pgfsys@transformshift{3.389548in}{4.768758in}%
\pgfsys@useobject{currentmarker}{}%
\end{pgfscope}%
\begin{pgfscope}%
\pgfsys@transformshift{3.596215in}{4.452410in}%
\pgfsys@useobject{currentmarker}{}%
\end{pgfscope}%
\begin{pgfscope}%
\pgfsys@transformshift{3.802882in}{4.087784in}%
\pgfsys@useobject{currentmarker}{}%
\end{pgfscope}%
\begin{pgfscope}%
\pgfsys@transformshift{4.009548in}{4.009265in}%
\pgfsys@useobject{currentmarker}{}%
\end{pgfscope}%
\begin{pgfscope}%
\pgfsys@transformshift{4.216215in}{4.393025in}%
\pgfsys@useobject{currentmarker}{}%
\end{pgfscope}%
\begin{pgfscope}%
\pgfsys@transformshift{4.422882in}{4.302041in}%
\pgfsys@useobject{currentmarker}{}%
\end{pgfscope}%
\end{pgfscope}%
\begin{pgfscope}%
\pgfpathrectangle{\pgfqpoint{1.322882in}{0.899712in}}{\pgfqpoint{6.200000in}{6.200000in}} %
\pgfusepath{clip}%
\pgfsetroundcap%
\pgfsetroundjoin%
\pgfsetlinewidth{2.007500pt}%
\definecolor{currentstroke}{rgb}{0.768627,0.305882,0.321569}%
\pgfsetstrokecolor{currentstroke}%
\pgfsetdash{}{0pt}%
\pgfpathmoveto{\pgfqpoint{2.149548in}{6.634408in}}%
\pgfpathlineto{\pgfqpoint{2.356215in}{6.324276in}}%
\pgfpathlineto{\pgfqpoint{2.562882in}{6.013427in}}%
\pgfpathlineto{\pgfqpoint{2.769548in}{5.702413in}}%
\pgfpathlineto{\pgfqpoint{2.976215in}{5.391349in}}%
\pgfpathlineto{\pgfqpoint{3.182882in}{5.080235in}}%
\pgfpathlineto{\pgfqpoint{3.389548in}{4.768959in}}%
\pgfusepath{stroke}%
\end{pgfscope}%
\begin{pgfscope}%
\pgfpathrectangle{\pgfqpoint{1.322882in}{0.899712in}}{\pgfqpoint{6.200000in}{6.200000in}} %
\pgfusepath{clip}%
\pgfsetbuttcap%
\pgfsetmiterjoin%
\definecolor{currentfill}{rgb}{0.768627,0.305882,0.321569}%
\pgfsetfillcolor{currentfill}%
\pgfsetlinewidth{0.000000pt}%
\definecolor{currentstroke}{rgb}{0.000000,0.000000,0.000000}%
\pgfsetstrokecolor{currentstroke}%
\pgfsetdash{}{0pt}%
\pgfsys@defobject{currentmarker}{\pgfqpoint{-0.038889in}{-0.038889in}}{\pgfqpoint{0.038889in}{0.038889in}}{%
\pgfpathmoveto{\pgfqpoint{0.000000in}{0.038889in}}%
\pgfpathlineto{\pgfqpoint{-0.038889in}{-0.038889in}}%
\pgfpathlineto{\pgfqpoint{0.038889in}{-0.038889in}}%
\pgfpathclose%
\pgfusepath{fill}%
}%
\begin{pgfscope}%
\pgfsys@transformshift{2.149548in}{6.634408in}%
\pgfsys@useobject{currentmarker}{}%
\end{pgfscope}%
\begin{pgfscope}%
\pgfsys@transformshift{2.356215in}{6.324276in}%
\pgfsys@useobject{currentmarker}{}%
\end{pgfscope}%
\begin{pgfscope}%
\pgfsys@transformshift{2.562882in}{6.013427in}%
\pgfsys@useobject{currentmarker}{}%
\end{pgfscope}%
\begin{pgfscope}%
\pgfsys@transformshift{2.769548in}{5.702413in}%
\pgfsys@useobject{currentmarker}{}%
\end{pgfscope}%
\begin{pgfscope}%
\pgfsys@transformshift{2.976215in}{5.391349in}%
\pgfsys@useobject{currentmarker}{}%
\end{pgfscope}%
\begin{pgfscope}%
\pgfsys@transformshift{3.182882in}{5.080235in}%
\pgfsys@useobject{currentmarker}{}%
\end{pgfscope}%
\begin{pgfscope}%
\pgfsys@transformshift{3.389548in}{4.768959in}%
\pgfsys@useobject{currentmarker}{}%
\end{pgfscope}%
\end{pgfscope}%
\begin{pgfscope}%
\pgfsetrectcap%
\pgfsetmiterjoin%
\pgfsetlinewidth{1.254687pt}%
\definecolor{currentstroke}{rgb}{0.150000,0.150000,0.150000}%
\pgfsetstrokecolor{currentstroke}%
\pgfsetdash{}{0pt}%
\pgfpathmoveto{\pgfqpoint{1.322882in}{0.899712in}}%
\pgfpathlineto{\pgfqpoint{7.522882in}{0.899712in}}%
\pgfusepath{stroke}%
\end{pgfscope}%
\begin{pgfscope}%
\pgfsetrectcap%
\pgfsetmiterjoin%
\pgfsetlinewidth{1.254687pt}%
\definecolor{currentstroke}{rgb}{0.150000,0.150000,0.150000}%
\pgfsetstrokecolor{currentstroke}%
\pgfsetdash{}{0pt}%
\pgfpathmoveto{\pgfqpoint{1.322882in}{0.899712in}}%
\pgfpathlineto{\pgfqpoint{1.322882in}{7.099712in}}%
\pgfusepath{stroke}%
\end{pgfscope}%
\begin{pgfscope}%
\pgfsetbuttcap%
\pgfsetmiterjoin%
\definecolor{currentfill}{rgb}{1.000000,1.000000,1.000000}%
\pgfsetfillcolor{currentfill}%
\pgfsetlinewidth{0.240900pt}%
\definecolor{currentstroke}{rgb}{0.150000,0.150000,0.150000}%
\pgfsetstrokecolor{currentstroke}%
\pgfsetdash{}{0pt}%
\pgfpathmoveto{\pgfqpoint{3.952602in}{5.400917in}}%
\pgfpathlineto{\pgfqpoint{7.356215in}{5.400917in}}%
\pgfpathlineto{\pgfqpoint{7.356215in}{6.933045in}}%
\pgfpathlineto{\pgfqpoint{3.952602in}{6.933045in}}%
\pgfpathclose%
\pgfusepath{stroke,fill}%
\end{pgfscope}%
\begin{pgfscope}%
\pgfsetroundcap%
\pgfsetroundjoin%
\pgfsetlinewidth{2.007500pt}%
\definecolor{currentstroke}{rgb}{0.298039,0.447059,0.690196}%
\pgfsetstrokecolor{currentstroke}%
\pgfsetdash{}{0pt}%
\pgfpathmoveto{\pgfqpoint{4.085935in}{6.674842in}}%
\pgfpathlineto{\pgfqpoint{4.752602in}{6.674842in}}%
\pgfusepath{stroke}%
\end{pgfscope}%
\begin{pgfscope}%
\pgfsetbuttcap%
\pgfsetroundjoin%
\definecolor{currentfill}{rgb}{0.298039,0.447059,0.690196}%
\pgfsetfillcolor{currentfill}%
\pgfsetlinewidth{0.000000pt}%
\definecolor{currentstroke}{rgb}{0.000000,0.000000,0.000000}%
\pgfsetstrokecolor{currentstroke}%
\pgfsetdash{}{0pt}%
\pgfsys@defobject{currentmarker}{\pgfqpoint{-0.038889in}{-0.038889in}}{\pgfqpoint{0.038889in}{0.038889in}}{%
\pgfpathmoveto{\pgfqpoint{0.000000in}{-0.038889in}}%
\pgfpathcurveto{\pgfqpoint{0.010313in}{-0.038889in}}{\pgfqpoint{0.020206in}{-0.034791in}}{\pgfqpoint{0.027499in}{-0.027499in}}%
\pgfpathcurveto{\pgfqpoint{0.034791in}{-0.020206in}}{\pgfqpoint{0.038889in}{-0.010313in}}{\pgfqpoint{0.038889in}{0.000000in}}%
\pgfpathcurveto{\pgfqpoint{0.038889in}{0.010313in}}{\pgfqpoint{0.034791in}{0.020206in}}{\pgfqpoint{0.027499in}{0.027499in}}%
\pgfpathcurveto{\pgfqpoint{0.020206in}{0.034791in}}{\pgfqpoint{0.010313in}{0.038889in}}{\pgfqpoint{0.000000in}{0.038889in}}%
\pgfpathcurveto{\pgfqpoint{-0.010313in}{0.038889in}}{\pgfqpoint{-0.020206in}{0.034791in}}{\pgfqpoint{-0.027499in}{0.027499in}}%
\pgfpathcurveto{\pgfqpoint{-0.034791in}{0.020206in}}{\pgfqpoint{-0.038889in}{0.010313in}}{\pgfqpoint{-0.038889in}{0.000000in}}%
\pgfpathcurveto{\pgfqpoint{-0.038889in}{-0.010313in}}{\pgfqpoint{-0.034791in}{-0.020206in}}{\pgfqpoint{-0.027499in}{-0.027499in}}%
\pgfpathcurveto{\pgfqpoint{-0.020206in}{-0.034791in}}{\pgfqpoint{-0.010313in}{-0.038889in}}{\pgfqpoint{0.000000in}{-0.038889in}}%
\pgfpathclose%
\pgfusepath{fill}%
}%
\begin{pgfscope}%
\pgfsys@transformshift{4.419268in}{6.674842in}%
\pgfsys@useobject{currentmarker}{}%
\end{pgfscope}%
\end{pgfscope}%
\begin{pgfscope}%
\definecolor{textcolor}{rgb}{0.150000,0.150000,0.150000}%
\pgfsetstrokecolor{textcolor}%
\pgfsetfillcolor{textcolor}%
\pgftext[x=5.019268in,y=6.558175in,left,base]{\color{textcolor}\sffamily\fontsize{24.000000}{28.800000}\selectfont FS-QTT-solver}%
\end{pgfscope}%
\begin{pgfscope}%
\pgfsetroundcap%
\pgfsetroundjoin%
\pgfsetlinewidth{2.007500pt}%
\definecolor{currentstroke}{rgb}{0.333333,0.658824,0.407843}%
\pgfsetstrokecolor{currentstroke}%
\pgfsetdash{}{0pt}%
\pgfpathmoveto{\pgfqpoint{4.085935in}{6.197466in}}%
\pgfpathlineto{\pgfqpoint{4.752602in}{6.197466in}}%
\pgfusepath{stroke}%
\end{pgfscope}%
\begin{pgfscope}%
\pgfsetbuttcap%
\pgfsetmiterjoin%
\definecolor{currentfill}{rgb}{0.333333,0.658824,0.407843}%
\pgfsetfillcolor{currentfill}%
\pgfsetlinewidth{0.000000pt}%
\definecolor{currentstroke}{rgb}{0.000000,0.000000,0.000000}%
\pgfsetstrokecolor{currentstroke}%
\pgfsetdash{}{0pt}%
\pgfsys@defobject{currentmarker}{\pgfqpoint{-0.038889in}{-0.038889in}}{\pgfqpoint{0.038889in}{0.038889in}}{%
\pgfpathmoveto{\pgfqpoint{-0.038889in}{-0.038889in}}%
\pgfpathlineto{\pgfqpoint{0.038889in}{-0.038889in}}%
\pgfpathlineto{\pgfqpoint{0.038889in}{0.038889in}}%
\pgfpathlineto{\pgfqpoint{-0.038889in}{0.038889in}}%
\pgfpathclose%
\pgfusepath{fill}%
}%
\begin{pgfscope}%
\pgfsys@transformshift{4.419268in}{6.197466in}%
\pgfsys@useobject{currentmarker}{}%
\end{pgfscope}%
\end{pgfscope}%
\begin{pgfscope}%
\definecolor{textcolor}{rgb}{0.150000,0.150000,0.150000}%
\pgfsetstrokecolor{textcolor}%
\pgfsetfillcolor{textcolor}%
\pgftext[x=5.019268in,y=6.080799in,left,base]{\color{textcolor}\sffamily\fontsize{24.000000}{28.800000}\selectfont FD-QTT-solver}%
\end{pgfscope}%
\begin{pgfscope}%
\pgfsetroundcap%
\pgfsetroundjoin%
\pgfsetlinewidth{2.007500pt}%
\definecolor{currentstroke}{rgb}{0.768627,0.305882,0.321569}%
\pgfsetstrokecolor{currentstroke}%
\pgfsetdash{}{0pt}%
\pgfpathmoveto{\pgfqpoint{4.085935in}{5.720090in}}%
\pgfpathlineto{\pgfqpoint{4.752602in}{5.720090in}}%
\pgfusepath{stroke}%
\end{pgfscope}%
\begin{pgfscope}%
\pgfsetbuttcap%
\pgfsetmiterjoin%
\definecolor{currentfill}{rgb}{0.768627,0.305882,0.321569}%
\pgfsetfillcolor{currentfill}%
\pgfsetlinewidth{0.000000pt}%
\definecolor{currentstroke}{rgb}{0.000000,0.000000,0.000000}%
\pgfsetstrokecolor{currentstroke}%
\pgfsetdash{}{0pt}%
\pgfsys@defobject{currentmarker}{\pgfqpoint{-0.038889in}{-0.038889in}}{\pgfqpoint{0.038889in}{0.038889in}}{%
\pgfpathmoveto{\pgfqpoint{0.000000in}{0.038889in}}%
\pgfpathlineto{\pgfqpoint{-0.038889in}{-0.038889in}}%
\pgfpathlineto{\pgfqpoint{0.038889in}{-0.038889in}}%
\pgfpathclose%
\pgfusepath{fill}%
}%
\begin{pgfscope}%
\pgfsys@transformshift{4.419268in}{5.720090in}%
\pgfsys@useobject{currentmarker}{}%
\end{pgfscope}%
\end{pgfscope}%
\begin{pgfscope}%
\definecolor{textcolor}{rgb}{0.150000,0.150000,0.150000}%
\pgfsetstrokecolor{textcolor}%
\pgfsetfillcolor{textcolor}%
\pgftext[x=5.019268in,y=5.603423in,left,base]{\color{textcolor}\sffamily\fontsize{24.000000}{28.800000}\selectfont FD-solver}%
\end{pgfscope}%
\end{pgfpicture}%
\makeatother%
\endgroup%

%% file: res_analyt_all_u_err.pgf
\begingroup%
\makeatletter%
\begin{pgfpicture}%
\pgfpathrectangle{\pgfpointorigin}{\pgfqpoint{7.778185in}{7.355067in}}%
\pgfusepath{use as bounding box, clip}%
\begin{pgfscope}%
\pgfsetbuttcap%
\pgfsetmiterjoin%
\definecolor{currentfill}{rgb}{1.000000,1.000000,1.000000}%
\pgfsetfillcolor{currentfill}%
\pgfsetlinewidth{0.000000pt}%
\definecolor{currentstroke}{rgb}{1.000000,1.000000,1.000000}%
\pgfsetstrokecolor{currentstroke}%
\pgfsetdash{}{0pt}%
\pgfpathmoveto{\pgfqpoint{0.000000in}{0.000000in}}%
\pgfpathlineto{\pgfqpoint{7.778185in}{0.000000in}}%
\pgfpathlineto{\pgfqpoint{7.778185in}{7.355067in}}%
\pgfpathlineto{\pgfqpoint{0.000000in}{7.355067in}}%
\pgfpathclose%
\pgfusepath{fill}%
\end{pgfscope}%
\begin{pgfscope}%
\pgfsetbuttcap%
\pgfsetmiterjoin%
\definecolor{currentfill}{rgb}{1.000000,1.000000,1.000000}%
\pgfsetfillcolor{currentfill}%
\pgfsetlinewidth{0.000000pt}%
\definecolor{currentstroke}{rgb}{0.000000,0.000000,0.000000}%
\pgfsetstrokecolor{currentstroke}%
\pgfsetstrokeopacity{0.000000}%
\pgfsetdash{}{0pt}%
\pgfpathmoveto{\pgfqpoint{1.319707in}{0.899712in}}%
\pgfpathlineto{\pgfqpoint{7.519707in}{0.899712in}}%
\pgfpathlineto{\pgfqpoint{7.519707in}{7.099712in}}%
\pgfpathlineto{\pgfqpoint{1.319707in}{7.099712in}}%
\pgfpathclose%
\pgfusepath{fill}%
\end{pgfscope}%
\begin{pgfscope}%
\pgfpathrectangle{\pgfqpoint{1.319707in}{0.899712in}}{\pgfqpoint{6.200000in}{6.200000in}} %
\pgfusepath{clip}%
\pgfsetbuttcap%
\pgfsetroundjoin%
\pgfsetlinewidth{0.803000pt}%
\definecolor{currentstroke}{rgb}{0.800000,0.800000,0.800000}%
\pgfsetstrokecolor{currentstroke}%
\pgfsetdash{{1.000000pt}{3.000000pt}}{0.000000pt}%
\pgfpathmoveto{\pgfqpoint{1.319707in}{0.899712in}}%
\pgfpathlineto{\pgfqpoint{1.319707in}{7.099712in}}%
\pgfusepath{stroke}%
\end{pgfscope}%
\begin{pgfscope}%
\definecolor{textcolor}{rgb}{0.150000,0.150000,0.150000}%
\pgfsetstrokecolor{textcolor}%
\pgfsetfillcolor{textcolor}%
\pgftext[x=1.319707in,y=0.821934in,,top]{\color{textcolor}\sffamily\fontsize{24.000000}{28.800000}\selectfont \(\displaystyle 0\)}%
\end{pgfscope}%
\begin{pgfscope}%
\pgfpathrectangle{\pgfqpoint{1.319707in}{0.899712in}}{\pgfqpoint{6.200000in}{6.200000in}} %
\pgfusepath{clip}%
\pgfsetbuttcap%
\pgfsetroundjoin%
\pgfsetlinewidth{0.803000pt}%
\definecolor{currentstroke}{rgb}{0.800000,0.800000,0.800000}%
\pgfsetstrokecolor{currentstroke}%
\pgfsetdash{{1.000000pt}{3.000000pt}}{0.000000pt}%
\pgfpathmoveto{\pgfqpoint{2.353040in}{0.899712in}}%
\pgfpathlineto{\pgfqpoint{2.353040in}{7.099712in}}%
\pgfusepath{stroke}%
\end{pgfscope}%
\begin{pgfscope}%
\definecolor{textcolor}{rgb}{0.150000,0.150000,0.150000}%
\pgfsetstrokecolor{textcolor}%
\pgfsetfillcolor{textcolor}%
\pgftext[x=2.353040in,y=0.821934in,,top]{\color{textcolor}\sffamily\fontsize{24.000000}{28.800000}\selectfont \(\displaystyle 5\)}%
\end{pgfscope}%
\begin{pgfscope}%
\pgfpathrectangle{\pgfqpoint{1.319707in}{0.899712in}}{\pgfqpoint{6.200000in}{6.200000in}} %
\pgfusepath{clip}%
\pgfsetbuttcap%
\pgfsetroundjoin%
\pgfsetlinewidth{0.803000pt}%
\definecolor{currentstroke}{rgb}{0.800000,0.800000,0.800000}%
\pgfsetstrokecolor{currentstroke}%
\pgfsetdash{{1.000000pt}{3.000000pt}}{0.000000pt}%
\pgfpathmoveto{\pgfqpoint{3.386373in}{0.899712in}}%
\pgfpathlineto{\pgfqpoint{3.386373in}{7.099712in}}%
\pgfusepath{stroke}%
\end{pgfscope}%
\begin{pgfscope}%
\definecolor{textcolor}{rgb}{0.150000,0.150000,0.150000}%
\pgfsetstrokecolor{textcolor}%
\pgfsetfillcolor{textcolor}%
\pgftext[x=3.386373in,y=0.821934in,,top]{\color{textcolor}\sffamily\fontsize{24.000000}{28.800000}\selectfont \(\displaystyle 10\)}%
\end{pgfscope}%
\begin{pgfscope}%
\pgfpathrectangle{\pgfqpoint{1.319707in}{0.899712in}}{\pgfqpoint{6.200000in}{6.200000in}} %
\pgfusepath{clip}%
\pgfsetbuttcap%
\pgfsetroundjoin%
\pgfsetlinewidth{0.803000pt}%
\definecolor{currentstroke}{rgb}{0.800000,0.800000,0.800000}%
\pgfsetstrokecolor{currentstroke}%
\pgfsetdash{{1.000000pt}{3.000000pt}}{0.000000pt}%
\pgfpathmoveto{\pgfqpoint{4.419707in}{0.899712in}}%
\pgfpathlineto{\pgfqpoint{4.419707in}{7.099712in}}%
\pgfusepath{stroke}%
\end{pgfscope}%
\begin{pgfscope}%
\definecolor{textcolor}{rgb}{0.150000,0.150000,0.150000}%
\pgfsetstrokecolor{textcolor}%
\pgfsetfillcolor{textcolor}%
\pgftext[x=4.419707in,y=0.821934in,,top]{\color{textcolor}\sffamily\fontsize{24.000000}{28.800000}\selectfont \(\displaystyle 15\)}%
\end{pgfscope}%
\begin{pgfscope}%
\pgfpathrectangle{\pgfqpoint{1.319707in}{0.899712in}}{\pgfqpoint{6.200000in}{6.200000in}} %
\pgfusepath{clip}%
\pgfsetbuttcap%
\pgfsetroundjoin%
\pgfsetlinewidth{0.803000pt}%
\definecolor{currentstroke}{rgb}{0.800000,0.800000,0.800000}%
\pgfsetstrokecolor{currentstroke}%
\pgfsetdash{{1.000000pt}{3.000000pt}}{0.000000pt}%
\pgfpathmoveto{\pgfqpoint{5.453040in}{0.899712in}}%
\pgfpathlineto{\pgfqpoint{5.453040in}{7.099712in}}%
\pgfusepath{stroke}%
\end{pgfscope}%
\begin{pgfscope}%
\definecolor{textcolor}{rgb}{0.150000,0.150000,0.150000}%
\pgfsetstrokecolor{textcolor}%
\pgfsetfillcolor{textcolor}%
\pgftext[x=5.453040in,y=0.821934in,,top]{\color{textcolor}\sffamily\fontsize{24.000000}{28.800000}\selectfont \(\displaystyle 20\)}%
\end{pgfscope}%
\begin{pgfscope}%
\pgfpathrectangle{\pgfqpoint{1.319707in}{0.899712in}}{\pgfqpoint{6.200000in}{6.200000in}} %
\pgfusepath{clip}%
\pgfsetbuttcap%
\pgfsetroundjoin%
\pgfsetlinewidth{0.803000pt}%
\definecolor{currentstroke}{rgb}{0.800000,0.800000,0.800000}%
\pgfsetstrokecolor{currentstroke}%
\pgfsetdash{{1.000000pt}{3.000000pt}}{0.000000pt}%
\pgfpathmoveto{\pgfqpoint{6.486373in}{0.899712in}}%
\pgfpathlineto{\pgfqpoint{6.486373in}{7.099712in}}%
\pgfusepath{stroke}%
\end{pgfscope}%
\begin{pgfscope}%
\definecolor{textcolor}{rgb}{0.150000,0.150000,0.150000}%
\pgfsetstrokecolor{textcolor}%
\pgfsetfillcolor{textcolor}%
\pgftext[x=6.486373in,y=0.821934in,,top]{\color{textcolor}\sffamily\fontsize{24.000000}{28.800000}\selectfont \(\displaystyle 25\)}%
\end{pgfscope}%
\begin{pgfscope}%
\pgfpathrectangle{\pgfqpoint{1.319707in}{0.899712in}}{\pgfqpoint{6.200000in}{6.200000in}} %
\pgfusepath{clip}%
\pgfsetbuttcap%
\pgfsetroundjoin%
\pgfsetlinewidth{0.803000pt}%
\definecolor{currentstroke}{rgb}{0.800000,0.800000,0.800000}%
\pgfsetstrokecolor{currentstroke}%
\pgfsetdash{{1.000000pt}{3.000000pt}}{0.000000pt}%
\pgfpathmoveto{\pgfqpoint{7.519707in}{0.899712in}}%
\pgfpathlineto{\pgfqpoint{7.519707in}{7.099712in}}%
\pgfusepath{stroke}%
\end{pgfscope}%
\begin{pgfscope}%
\definecolor{textcolor}{rgb}{0.150000,0.150000,0.150000}%
\pgfsetstrokecolor{textcolor}%
\pgfsetfillcolor{textcolor}%
\pgftext[x=7.519707in,y=0.821934in,,top]{\color{textcolor}\sffamily\fontsize{24.000000}{28.800000}\selectfont \(\displaystyle 30\)}%
\end{pgfscope}%
\begin{pgfscope}%
\definecolor{textcolor}{rgb}{0.150000,0.150000,0.150000}%
\pgfsetstrokecolor{textcolor}%
\pgfsetfillcolor{textcolor}%
\pgftext[x=4.419707in,y=0.441780in,,top]{\color{textcolor}\sffamily\fontsize{26.400000}{31.680000}\selectfont d}%
\end{pgfscope}%
\begin{pgfscope}%
\pgfpathrectangle{\pgfqpoint{1.319707in}{0.899712in}}{\pgfqpoint{6.200000in}{6.200000in}} %
\pgfusepath{clip}%
\pgfsetbuttcap%
\pgfsetroundjoin%
\pgfsetlinewidth{0.803000pt}%
\definecolor{currentstroke}{rgb}{0.800000,0.800000,0.800000}%
\pgfsetstrokecolor{currentstroke}%
\pgfsetdash{{1.000000pt}{3.000000pt}}{0.000000pt}%
\pgfpathmoveto{\pgfqpoint{1.319707in}{0.899712in}}%
\pgfpathlineto{\pgfqpoint{7.519707in}{0.899712in}}%
\pgfusepath{stroke}%
\end{pgfscope}%
\begin{pgfscope}%
\definecolor{textcolor}{rgb}{0.150000,0.150000,0.150000}%
\pgfsetstrokecolor{textcolor}%
\pgfsetfillcolor{textcolor}%
\pgftext[x=1.241929in,y=0.899712in,right,]{\color{textcolor}\sffamily\fontsize{24.000000}{28.800000}\selectfont \(\displaystyle 10^{-11}\)}%
\end{pgfscope}%
\begin{pgfscope}%
\pgfpathrectangle{\pgfqpoint{1.319707in}{0.899712in}}{\pgfqpoint{6.200000in}{6.200000in}} %
\pgfusepath{clip}%
\pgfsetbuttcap%
\pgfsetroundjoin%
\pgfsetlinewidth{0.803000pt}%
\definecolor{currentstroke}{rgb}{0.800000,0.800000,0.800000}%
\pgfsetstrokecolor{currentstroke}%
\pgfsetdash{{1.000000pt}{3.000000pt}}{0.000000pt}%
\pgfpathmoveto{\pgfqpoint{1.319707in}{1.519712in}}%
\pgfpathlineto{\pgfqpoint{7.519707in}{1.519712in}}%
\pgfusepath{stroke}%
\end{pgfscope}%
\begin{pgfscope}%
\definecolor{textcolor}{rgb}{0.150000,0.150000,0.150000}%
\pgfsetstrokecolor{textcolor}%
\pgfsetfillcolor{textcolor}%
\pgftext[x=1.241929in,y=1.519712in,right,]{\color{textcolor}\sffamily\fontsize{24.000000}{28.800000}\selectfont \(\displaystyle 10^{-10}\)}%
\end{pgfscope}%
\begin{pgfscope}%
\pgfpathrectangle{\pgfqpoint{1.319707in}{0.899712in}}{\pgfqpoint{6.200000in}{6.200000in}} %
\pgfusepath{clip}%
\pgfsetbuttcap%
\pgfsetroundjoin%
\pgfsetlinewidth{0.803000pt}%
\definecolor{currentstroke}{rgb}{0.800000,0.800000,0.800000}%
\pgfsetstrokecolor{currentstroke}%
\pgfsetdash{{1.000000pt}{3.000000pt}}{0.000000pt}%
\pgfpathmoveto{\pgfqpoint{1.319707in}{2.139712in}}%
\pgfpathlineto{\pgfqpoint{7.519707in}{2.139712in}}%
\pgfusepath{stroke}%
\end{pgfscope}%
\begin{pgfscope}%
\definecolor{textcolor}{rgb}{0.150000,0.150000,0.150000}%
\pgfsetstrokecolor{textcolor}%
\pgfsetfillcolor{textcolor}%
\pgftext[x=1.241929in,y=2.139712in,right,]{\color{textcolor}\sffamily\fontsize{24.000000}{28.800000}\selectfont \(\displaystyle 10^{-9}\)}%
\end{pgfscope}%
\begin{pgfscope}%
\pgfpathrectangle{\pgfqpoint{1.319707in}{0.899712in}}{\pgfqpoint{6.200000in}{6.200000in}} %
\pgfusepath{clip}%
\pgfsetbuttcap%
\pgfsetroundjoin%
\pgfsetlinewidth{0.803000pt}%
\definecolor{currentstroke}{rgb}{0.800000,0.800000,0.800000}%
\pgfsetstrokecolor{currentstroke}%
\pgfsetdash{{1.000000pt}{3.000000pt}}{0.000000pt}%
\pgfpathmoveto{\pgfqpoint{1.319707in}{2.759712in}}%
\pgfpathlineto{\pgfqpoint{7.519707in}{2.759712in}}%
\pgfusepath{stroke}%
\end{pgfscope}%
\begin{pgfscope}%
\definecolor{textcolor}{rgb}{0.150000,0.150000,0.150000}%
\pgfsetstrokecolor{textcolor}%
\pgfsetfillcolor{textcolor}%
\pgftext[x=1.241929in,y=2.759712in,right,]{\color{textcolor}\sffamily\fontsize{24.000000}{28.800000}\selectfont \(\displaystyle 10^{-8}\)}%
\end{pgfscope}%
\begin{pgfscope}%
\pgfpathrectangle{\pgfqpoint{1.319707in}{0.899712in}}{\pgfqpoint{6.200000in}{6.200000in}} %
\pgfusepath{clip}%
\pgfsetbuttcap%
\pgfsetroundjoin%
\pgfsetlinewidth{0.803000pt}%
\definecolor{currentstroke}{rgb}{0.800000,0.800000,0.800000}%
\pgfsetstrokecolor{currentstroke}%
\pgfsetdash{{1.000000pt}{3.000000pt}}{0.000000pt}%
\pgfpathmoveto{\pgfqpoint{1.319707in}{3.379712in}}%
\pgfpathlineto{\pgfqpoint{7.519707in}{3.379712in}}%
\pgfusepath{stroke}%
\end{pgfscope}%
\begin{pgfscope}%
\definecolor{textcolor}{rgb}{0.150000,0.150000,0.150000}%
\pgfsetstrokecolor{textcolor}%
\pgfsetfillcolor{textcolor}%
\pgftext[x=1.241929in,y=3.379712in,right,]{\color{textcolor}\sffamily\fontsize{24.000000}{28.800000}\selectfont \(\displaystyle 10^{-7}\)}%
\end{pgfscope}%
\begin{pgfscope}%
\pgfpathrectangle{\pgfqpoint{1.319707in}{0.899712in}}{\pgfqpoint{6.200000in}{6.200000in}} %
\pgfusepath{clip}%
\pgfsetbuttcap%
\pgfsetroundjoin%
\pgfsetlinewidth{0.803000pt}%
\definecolor{currentstroke}{rgb}{0.800000,0.800000,0.800000}%
\pgfsetstrokecolor{currentstroke}%
\pgfsetdash{{1.000000pt}{3.000000pt}}{0.000000pt}%
\pgfpathmoveto{\pgfqpoint{1.319707in}{3.999712in}}%
\pgfpathlineto{\pgfqpoint{7.519707in}{3.999712in}}%
\pgfusepath{stroke}%
\end{pgfscope}%
\begin{pgfscope}%
\definecolor{textcolor}{rgb}{0.150000,0.150000,0.150000}%
\pgfsetstrokecolor{textcolor}%
\pgfsetfillcolor{textcolor}%
\pgftext[x=1.241929in,y=3.999712in,right,]{\color{textcolor}\sffamily\fontsize{24.000000}{28.800000}\selectfont \(\displaystyle 10^{-6}\)}%
\end{pgfscope}%
\begin{pgfscope}%
\pgfpathrectangle{\pgfqpoint{1.319707in}{0.899712in}}{\pgfqpoint{6.200000in}{6.200000in}} %
\pgfusepath{clip}%
\pgfsetbuttcap%
\pgfsetroundjoin%
\pgfsetlinewidth{0.803000pt}%
\definecolor{currentstroke}{rgb}{0.800000,0.800000,0.800000}%
\pgfsetstrokecolor{currentstroke}%
\pgfsetdash{{1.000000pt}{3.000000pt}}{0.000000pt}%
\pgfpathmoveto{\pgfqpoint{1.319707in}{4.619712in}}%
\pgfpathlineto{\pgfqpoint{7.519707in}{4.619712in}}%
\pgfusepath{stroke}%
\end{pgfscope}%
\begin{pgfscope}%
\definecolor{textcolor}{rgb}{0.150000,0.150000,0.150000}%
\pgfsetstrokecolor{textcolor}%
\pgfsetfillcolor{textcolor}%
\pgftext[x=1.241929in,y=4.619712in,right,]{\color{textcolor}\sffamily\fontsize{24.000000}{28.800000}\selectfont \(\displaystyle 10^{-5}\)}%
\end{pgfscope}%
\begin{pgfscope}%
\pgfpathrectangle{\pgfqpoint{1.319707in}{0.899712in}}{\pgfqpoint{6.200000in}{6.200000in}} %
\pgfusepath{clip}%
\pgfsetbuttcap%
\pgfsetroundjoin%
\pgfsetlinewidth{0.803000pt}%
\definecolor{currentstroke}{rgb}{0.800000,0.800000,0.800000}%
\pgfsetstrokecolor{currentstroke}%
\pgfsetdash{{1.000000pt}{3.000000pt}}{0.000000pt}%
\pgfpathmoveto{\pgfqpoint{1.319707in}{5.239712in}}%
\pgfpathlineto{\pgfqpoint{7.519707in}{5.239712in}}%
\pgfusepath{stroke}%
\end{pgfscope}%
\begin{pgfscope}%
\definecolor{textcolor}{rgb}{0.150000,0.150000,0.150000}%
\pgfsetstrokecolor{textcolor}%
\pgfsetfillcolor{textcolor}%
\pgftext[x=1.241929in,y=5.239712in,right,]{\color{textcolor}\sffamily\fontsize{24.000000}{28.800000}\selectfont \(\displaystyle 10^{-4}\)}%
\end{pgfscope}%
\begin{pgfscope}%
\pgfpathrectangle{\pgfqpoint{1.319707in}{0.899712in}}{\pgfqpoint{6.200000in}{6.200000in}} %
\pgfusepath{clip}%
\pgfsetbuttcap%
\pgfsetroundjoin%
\pgfsetlinewidth{0.803000pt}%
\definecolor{currentstroke}{rgb}{0.800000,0.800000,0.800000}%
\pgfsetstrokecolor{currentstroke}%
\pgfsetdash{{1.000000pt}{3.000000pt}}{0.000000pt}%
\pgfpathmoveto{\pgfqpoint{1.319707in}{5.859712in}}%
\pgfpathlineto{\pgfqpoint{7.519707in}{5.859712in}}%
\pgfusepath{stroke}%
\end{pgfscope}%
\begin{pgfscope}%
\definecolor{textcolor}{rgb}{0.150000,0.150000,0.150000}%
\pgfsetstrokecolor{textcolor}%
\pgfsetfillcolor{textcolor}%
\pgftext[x=1.241929in,y=5.859712in,right,]{\color{textcolor}\sffamily\fontsize{24.000000}{28.800000}\selectfont \(\displaystyle 10^{-3}\)}%
\end{pgfscope}%
\begin{pgfscope}%
\pgfpathrectangle{\pgfqpoint{1.319707in}{0.899712in}}{\pgfqpoint{6.200000in}{6.200000in}} %
\pgfusepath{clip}%
\pgfsetbuttcap%
\pgfsetroundjoin%
\pgfsetlinewidth{0.803000pt}%
\definecolor{currentstroke}{rgb}{0.800000,0.800000,0.800000}%
\pgfsetstrokecolor{currentstroke}%
\pgfsetdash{{1.000000pt}{3.000000pt}}{0.000000pt}%
\pgfpathmoveto{\pgfqpoint{1.319707in}{6.479712in}}%
\pgfpathlineto{\pgfqpoint{7.519707in}{6.479712in}}%
\pgfusepath{stroke}%
\end{pgfscope}%
\begin{pgfscope}%
\definecolor{textcolor}{rgb}{0.150000,0.150000,0.150000}%
\pgfsetstrokecolor{textcolor}%
\pgfsetfillcolor{textcolor}%
\pgftext[x=1.241929in,y=6.479712in,right,]{\color{textcolor}\sffamily\fontsize{24.000000}{28.800000}\selectfont \(\displaystyle 10^{-2}\)}%
\end{pgfscope}%
\begin{pgfscope}%
\pgfpathrectangle{\pgfqpoint{1.319707in}{0.899712in}}{\pgfqpoint{6.200000in}{6.200000in}} %
\pgfusepath{clip}%
\pgfsetbuttcap%
\pgfsetroundjoin%
\pgfsetlinewidth{0.803000pt}%
\definecolor{currentstroke}{rgb}{0.800000,0.800000,0.800000}%
\pgfsetstrokecolor{currentstroke}%
\pgfsetdash{{1.000000pt}{3.000000pt}}{0.000000pt}%
\pgfpathmoveto{\pgfqpoint{1.319707in}{7.099712in}}%
\pgfpathlineto{\pgfqpoint{7.519707in}{7.099712in}}%
\pgfusepath{stroke}%
\end{pgfscope}%
\begin{pgfscope}%
\definecolor{textcolor}{rgb}{0.150000,0.150000,0.150000}%
\pgfsetstrokecolor{textcolor}%
\pgfsetfillcolor{textcolor}%
\pgftext[x=1.241929in,y=7.099712in,right,]{\color{textcolor}\sffamily\fontsize{24.000000}{28.800000}\selectfont \(\displaystyle 10^{-1}\)}%
\end{pgfscope}%
\begin{pgfscope}%
\definecolor{textcolor}{rgb}{0.150000,0.150000,0.150000}%
\pgfsetstrokecolor{textcolor}%
\pgfsetfillcolor{textcolor}%
\pgftext[x=0.441780in,y=3.999712in,,bottom,rotate=90.000000]{\color{textcolor}\sffamily\fontsize{26.400000}{31.680000}\selectfont Solution relative error}%
\end{pgfscope}%
\begin{pgfscope}%
\pgfpathrectangle{\pgfqpoint{1.319707in}{0.899712in}}{\pgfqpoint{6.200000in}{6.200000in}} %
\pgfusepath{clip}%
\pgfsetroundcap%
\pgfsetroundjoin%
\pgfsetlinewidth{2.007500pt}%
\definecolor{currentstroke}{rgb}{0.298039,0.447059,0.690196}%
\pgfsetstrokecolor{currentstroke}%
\pgfsetdash{}{0pt}%
\pgfpathmoveto{\pgfqpoint{2.146373in}{6.502303in}}%
\pgfpathlineto{\pgfqpoint{2.353040in}{6.128285in}}%
\pgfpathlineto{\pgfqpoint{2.559707in}{5.754819in}}%
\pgfpathlineto{\pgfqpoint{2.766373in}{5.381495in}}%
\pgfpathlineto{\pgfqpoint{2.973040in}{5.008206in}}%
\pgfpathlineto{\pgfqpoint{3.179707in}{4.634925in}}%
\pgfpathlineto{\pgfqpoint{3.386373in}{4.261648in}}%
\pgfpathlineto{\pgfqpoint{3.593040in}{3.888370in}}%
\pgfpathlineto{\pgfqpoint{3.799707in}{3.515092in}}%
\pgfpathlineto{\pgfqpoint{4.006373in}{3.141812in}}%
\pgfpathlineto{\pgfqpoint{4.213040in}{2.768514in}}%
\pgfpathlineto{\pgfqpoint{4.419707in}{2.395201in}}%
\pgfpathlineto{\pgfqpoint{4.626373in}{2.021409in}}%
\pgfpathlineto{\pgfqpoint{4.833040in}{1.644294in}}%
\pgfpathlineto{\pgfqpoint{5.039707in}{1.287058in}}%
\pgfpathlineto{\pgfqpoint{5.246373in}{1.180357in}}%
\pgfpathlineto{\pgfqpoint{5.453040in}{1.183895in}}%
\pgfpathlineto{\pgfqpoint{5.659707in}{1.266807in}}%
\pgfpathlineto{\pgfqpoint{5.866373in}{1.263572in}}%
\pgfpathlineto{\pgfqpoint{6.073040in}{1.419686in}}%
\pgfpathlineto{\pgfqpoint{6.279707in}{1.435460in}}%
\pgfpathlineto{\pgfqpoint{6.486373in}{1.503690in}}%
\pgfpathlineto{\pgfqpoint{6.693040in}{1.724458in}}%
\pgfpathlineto{\pgfqpoint{6.899707in}{1.822289in}}%
\pgfpathlineto{\pgfqpoint{7.106373in}{1.946664in}}%
\pgfpathlineto{\pgfqpoint{7.313040in}{1.985377in}}%
\pgfpathlineto{\pgfqpoint{7.519707in}{2.307356in}}%
\pgfusepath{stroke}%
\end{pgfscope}%
\begin{pgfscope}%
\pgfpathrectangle{\pgfqpoint{1.319707in}{0.899712in}}{\pgfqpoint{6.200000in}{6.200000in}} %
\pgfusepath{clip}%
\pgfsetbuttcap%
\pgfsetroundjoin%
\definecolor{currentfill}{rgb}{0.298039,0.447059,0.690196}%
\pgfsetfillcolor{currentfill}%
\pgfsetlinewidth{0.000000pt}%
\definecolor{currentstroke}{rgb}{0.000000,0.000000,0.000000}%
\pgfsetstrokecolor{currentstroke}%
\pgfsetdash{}{0pt}%
\pgfsys@defobject{currentmarker}{\pgfqpoint{-0.038889in}{-0.038889in}}{\pgfqpoint{0.038889in}{0.038889in}}{%
\pgfpathmoveto{\pgfqpoint{0.000000in}{-0.038889in}}%
\pgfpathcurveto{\pgfqpoint{0.010313in}{-0.038889in}}{\pgfqpoint{0.020206in}{-0.034791in}}{\pgfqpoint{0.027499in}{-0.027499in}}%
\pgfpathcurveto{\pgfqpoint{0.034791in}{-0.020206in}}{\pgfqpoint{0.038889in}{-0.010313in}}{\pgfqpoint{0.038889in}{0.000000in}}%
\pgfpathcurveto{\pgfqpoint{0.038889in}{0.010313in}}{\pgfqpoint{0.034791in}{0.020206in}}{\pgfqpoint{0.027499in}{0.027499in}}%
\pgfpathcurveto{\pgfqpoint{0.020206in}{0.034791in}}{\pgfqpoint{0.010313in}{0.038889in}}{\pgfqpoint{0.000000in}{0.038889in}}%
\pgfpathcurveto{\pgfqpoint{-0.010313in}{0.038889in}}{\pgfqpoint{-0.020206in}{0.034791in}}{\pgfqpoint{-0.027499in}{0.027499in}}%
\pgfpathcurveto{\pgfqpoint{-0.034791in}{0.020206in}}{\pgfqpoint{-0.038889in}{0.010313in}}{\pgfqpoint{-0.038889in}{0.000000in}}%
\pgfpathcurveto{\pgfqpoint{-0.038889in}{-0.010313in}}{\pgfqpoint{-0.034791in}{-0.020206in}}{\pgfqpoint{-0.027499in}{-0.027499in}}%
\pgfpathcurveto{\pgfqpoint{-0.020206in}{-0.034791in}}{\pgfqpoint{-0.010313in}{-0.038889in}}{\pgfqpoint{0.000000in}{-0.038889in}}%
\pgfpathclose%
\pgfusepath{fill}%
}%
\begin{pgfscope}%
\pgfsys@transformshift{2.146373in}{6.502303in}%
\pgfsys@useobject{currentmarker}{}%
\end{pgfscope}%
\begin{pgfscope}%
\pgfsys@transformshift{2.353040in}{6.128285in}%
\pgfsys@useobject{currentmarker}{}%
\end{pgfscope}%
\begin{pgfscope}%
\pgfsys@transformshift{2.559707in}{5.754819in}%
\pgfsys@useobject{currentmarker}{}%
\end{pgfscope}%
\begin{pgfscope}%
\pgfsys@transformshift{2.766373in}{5.381495in}%
\pgfsys@useobject{currentmarker}{}%
\end{pgfscope}%
\begin{pgfscope}%
\pgfsys@transformshift{2.973040in}{5.008206in}%
\pgfsys@useobject{currentmarker}{}%
\end{pgfscope}%
\begin{pgfscope}%
\pgfsys@transformshift{3.179707in}{4.634925in}%
\pgfsys@useobject{currentmarker}{}%
\end{pgfscope}%
\begin{pgfscope}%
\pgfsys@transformshift{3.386373in}{4.261648in}%
\pgfsys@useobject{currentmarker}{}%
\end{pgfscope}%
\begin{pgfscope}%
\pgfsys@transformshift{3.593040in}{3.888370in}%
\pgfsys@useobject{currentmarker}{}%
\end{pgfscope}%
\begin{pgfscope}%
\pgfsys@transformshift{3.799707in}{3.515092in}%
\pgfsys@useobject{currentmarker}{}%
\end{pgfscope}%
\begin{pgfscope}%
\pgfsys@transformshift{4.006373in}{3.141812in}%
\pgfsys@useobject{currentmarker}{}%
\end{pgfscope}%
\begin{pgfscope}%
\pgfsys@transformshift{4.213040in}{2.768514in}%
\pgfsys@useobject{currentmarker}{}%
\end{pgfscope}%
\begin{pgfscope}%
\pgfsys@transformshift{4.419707in}{2.395201in}%
\pgfsys@useobject{currentmarker}{}%
\end{pgfscope}%
\begin{pgfscope}%
\pgfsys@transformshift{4.626373in}{2.021409in}%
\pgfsys@useobject{currentmarker}{}%
\end{pgfscope}%
\begin{pgfscope}%
\pgfsys@transformshift{4.833040in}{1.644294in}%
\pgfsys@useobject{currentmarker}{}%
\end{pgfscope}%
\begin{pgfscope}%
\pgfsys@transformshift{5.039707in}{1.287058in}%
\pgfsys@useobject{currentmarker}{}%
\end{pgfscope}%
\begin{pgfscope}%
\pgfsys@transformshift{5.246373in}{1.180357in}%
\pgfsys@useobject{currentmarker}{}%
\end{pgfscope}%
\begin{pgfscope}%
\pgfsys@transformshift{5.453040in}{1.183895in}%
\pgfsys@useobject{currentmarker}{}%
\end{pgfscope}%
\begin{pgfscope}%
\pgfsys@transformshift{5.659707in}{1.266807in}%
\pgfsys@useobject{currentmarker}{}%
\end{pgfscope}%
\begin{pgfscope}%
\pgfsys@transformshift{5.866373in}{1.263572in}%
\pgfsys@useobject{currentmarker}{}%
\end{pgfscope}%
\begin{pgfscope}%
\pgfsys@transformshift{6.073040in}{1.419686in}%
\pgfsys@useobject{currentmarker}{}%
\end{pgfscope}%
\begin{pgfscope}%
\pgfsys@transformshift{6.279707in}{1.435460in}%
\pgfsys@useobject{currentmarker}{}%
\end{pgfscope}%
\begin{pgfscope}%
\pgfsys@transformshift{6.486373in}{1.503690in}%
\pgfsys@useobject{currentmarker}{}%
\end{pgfscope}%
\begin{pgfscope}%
\pgfsys@transformshift{6.693040in}{1.724458in}%
\pgfsys@useobject{currentmarker}{}%
\end{pgfscope}%
\begin{pgfscope}%
\pgfsys@transformshift{6.899707in}{1.822289in}%
\pgfsys@useobject{currentmarker}{}%
\end{pgfscope}%
\begin{pgfscope}%
\pgfsys@transformshift{7.106373in}{1.946664in}%
\pgfsys@useobject{currentmarker}{}%
\end{pgfscope}%
\begin{pgfscope}%
\pgfsys@transformshift{7.313040in}{1.985377in}%
\pgfsys@useobject{currentmarker}{}%
\end{pgfscope}%
\begin{pgfscope}%
\pgfsys@transformshift{7.519707in}{2.307356in}%
\pgfsys@useobject{currentmarker}{}%
\end{pgfscope}%
\end{pgfscope}%
\begin{pgfscope}%
\pgfpathrectangle{\pgfqpoint{1.319707in}{0.899712in}}{\pgfqpoint{6.200000in}{6.200000in}} %
\pgfusepath{clip}%
\pgfsetroundcap%
\pgfsetroundjoin%
\pgfsetlinewidth{2.007500pt}%
\definecolor{currentstroke}{rgb}{0.333333,0.658824,0.407843}%
\pgfsetstrokecolor{currentstroke}%
\pgfsetdash{}{0pt}%
\pgfpathmoveto{\pgfqpoint{2.146373in}{6.502303in}}%
\pgfpathlineto{\pgfqpoint{2.353040in}{6.128285in}}%
\pgfpathlineto{\pgfqpoint{2.559707in}{5.754819in}}%
\pgfpathlineto{\pgfqpoint{2.766373in}{5.381495in}}%
\pgfpathlineto{\pgfqpoint{2.973040in}{5.008205in}}%
\pgfpathlineto{\pgfqpoint{3.179707in}{4.634908in}}%
\pgfpathlineto{\pgfqpoint{3.386373in}{4.261423in}}%
\pgfpathlineto{\pgfqpoint{3.593040in}{3.883276in}}%
\pgfpathlineto{\pgfqpoint{3.799707in}{3.457358in}}%
\pgfpathlineto{\pgfqpoint{4.006373in}{3.367184in}}%
\pgfpathlineto{\pgfqpoint{4.213040in}{3.900600in}}%
\pgfpathlineto{\pgfqpoint{4.419707in}{4.397270in}}%
\pgfusepath{stroke}%
\end{pgfscope}%
\begin{pgfscope}%
\pgfpathrectangle{\pgfqpoint{1.319707in}{0.899712in}}{\pgfqpoint{6.200000in}{6.200000in}} %
\pgfusepath{clip}%
\pgfsetbuttcap%
\pgfsetmiterjoin%
\definecolor{currentfill}{rgb}{0.333333,0.658824,0.407843}%
\pgfsetfillcolor{currentfill}%
\pgfsetlinewidth{0.000000pt}%
\definecolor{currentstroke}{rgb}{0.000000,0.000000,0.000000}%
\pgfsetstrokecolor{currentstroke}%
\pgfsetdash{}{0pt}%
\pgfsys@defobject{currentmarker}{\pgfqpoint{-0.038889in}{-0.038889in}}{\pgfqpoint{0.038889in}{0.038889in}}{%
\pgfpathmoveto{\pgfqpoint{-0.038889in}{-0.038889in}}%
\pgfpathlineto{\pgfqpoint{0.038889in}{-0.038889in}}%
\pgfpathlineto{\pgfqpoint{0.038889in}{0.038889in}}%
\pgfpathlineto{\pgfqpoint{-0.038889in}{0.038889in}}%
\pgfpathclose%
\pgfusepath{fill}%
}%
\begin{pgfscope}%
\pgfsys@transformshift{2.146373in}{6.502303in}%
\pgfsys@useobject{currentmarker}{}%
\end{pgfscope}%
\begin{pgfscope}%
\pgfsys@transformshift{2.353040in}{6.128285in}%
\pgfsys@useobject{currentmarker}{}%
\end{pgfscope}%
\begin{pgfscope}%
\pgfsys@transformshift{2.559707in}{5.754819in}%
\pgfsys@useobject{currentmarker}{}%
\end{pgfscope}%
\begin{pgfscope}%
\pgfsys@transformshift{2.766373in}{5.381495in}%
\pgfsys@useobject{currentmarker}{}%
\end{pgfscope}%
\begin{pgfscope}%
\pgfsys@transformshift{2.973040in}{5.008205in}%
\pgfsys@useobject{currentmarker}{}%
\end{pgfscope}%
\begin{pgfscope}%
\pgfsys@transformshift{3.179707in}{4.634908in}%
\pgfsys@useobject{currentmarker}{}%
\end{pgfscope}%
\begin{pgfscope}%
\pgfsys@transformshift{3.386373in}{4.261423in}%
\pgfsys@useobject{currentmarker}{}%
\end{pgfscope}%
\begin{pgfscope}%
\pgfsys@transformshift{3.593040in}{3.883276in}%
\pgfsys@useobject{currentmarker}{}%
\end{pgfscope}%
\begin{pgfscope}%
\pgfsys@transformshift{3.799707in}{3.457358in}%
\pgfsys@useobject{currentmarker}{}%
\end{pgfscope}%
\begin{pgfscope}%
\pgfsys@transformshift{4.006373in}{3.367184in}%
\pgfsys@useobject{currentmarker}{}%
\end{pgfscope}%
\begin{pgfscope}%
\pgfsys@transformshift{4.213040in}{3.900600in}%
\pgfsys@useobject{currentmarker}{}%
\end{pgfscope}%
\begin{pgfscope}%
\pgfsys@transformshift{4.419707in}{4.397270in}%
\pgfsys@useobject{currentmarker}{}%
\end{pgfscope}%
\end{pgfscope}%
\begin{pgfscope}%
\pgfpathrectangle{\pgfqpoint{1.319707in}{0.899712in}}{\pgfqpoint{6.200000in}{6.200000in}} %
\pgfusepath{clip}%
\pgfsetroundcap%
\pgfsetroundjoin%
\pgfsetlinewidth{2.007500pt}%
\definecolor{currentstroke}{rgb}{0.768627,0.305882,0.321569}%
\pgfsetstrokecolor{currentstroke}%
\pgfsetdash{}{0pt}%
\pgfpathmoveto{\pgfqpoint{2.146373in}{6.502303in}}%
\pgfpathlineto{\pgfqpoint{2.353040in}{6.128285in}}%
\pgfpathlineto{\pgfqpoint{2.559707in}{5.754819in}}%
\pgfpathlineto{\pgfqpoint{2.766373in}{5.381495in}}%
\pgfpathlineto{\pgfqpoint{2.973040in}{5.008206in}}%
\pgfpathlineto{\pgfqpoint{3.179707in}{4.634925in}}%
\pgfpathlineto{\pgfqpoint{3.386373in}{4.261648in}}%
\pgfusepath{stroke}%
\end{pgfscope}%
\begin{pgfscope}%
\pgfpathrectangle{\pgfqpoint{1.319707in}{0.899712in}}{\pgfqpoint{6.200000in}{6.200000in}} %
\pgfusepath{clip}%
\pgfsetbuttcap%
\pgfsetmiterjoin%
\definecolor{currentfill}{rgb}{0.768627,0.305882,0.321569}%
\pgfsetfillcolor{currentfill}%
\pgfsetlinewidth{0.000000pt}%
\definecolor{currentstroke}{rgb}{0.000000,0.000000,0.000000}%
\pgfsetstrokecolor{currentstroke}%
\pgfsetdash{}{0pt}%
\pgfsys@defobject{currentmarker}{\pgfqpoint{-0.038889in}{-0.038889in}}{\pgfqpoint{0.038889in}{0.038889in}}{%
\pgfpathmoveto{\pgfqpoint{0.000000in}{0.038889in}}%
\pgfpathlineto{\pgfqpoint{-0.038889in}{-0.038889in}}%
\pgfpathlineto{\pgfqpoint{0.038889in}{-0.038889in}}%
\pgfpathclose%
\pgfusepath{fill}%
}%
\begin{pgfscope}%
\pgfsys@transformshift{2.146373in}{6.502303in}%
\pgfsys@useobject{currentmarker}{}%
\end{pgfscope}%
\begin{pgfscope}%
\pgfsys@transformshift{2.353040in}{6.128285in}%
\pgfsys@useobject{currentmarker}{}%
\end{pgfscope}%
\begin{pgfscope}%
\pgfsys@transformshift{2.559707in}{5.754819in}%
\pgfsys@useobject{currentmarker}{}%
\end{pgfscope}%
\begin{pgfscope}%
\pgfsys@transformshift{2.766373in}{5.381495in}%
\pgfsys@useobject{currentmarker}{}%
\end{pgfscope}%
\begin{pgfscope}%
\pgfsys@transformshift{2.973040in}{5.008206in}%
\pgfsys@useobject{currentmarker}{}%
\end{pgfscope}%
\begin{pgfscope}%
\pgfsys@transformshift{3.179707in}{4.634925in}%
\pgfsys@useobject{currentmarker}{}%
\end{pgfscope}%
\begin{pgfscope}%
\pgfsys@transformshift{3.386373in}{4.261648in}%
\pgfsys@useobject{currentmarker}{}%
\end{pgfscope}%
\end{pgfscope}%
\begin{pgfscope}%
\pgfsetrectcap%
\pgfsetmiterjoin%
\pgfsetlinewidth{1.254687pt}%
\definecolor{currentstroke}{rgb}{0.150000,0.150000,0.150000}%
\pgfsetstrokecolor{currentstroke}%
\pgfsetdash{}{0pt}%
\pgfpathmoveto{\pgfqpoint{1.319707in}{0.899712in}}%
\pgfpathlineto{\pgfqpoint{7.519707in}{0.899712in}}%
\pgfusepath{stroke}%
\end{pgfscope}%
\begin{pgfscope}%
\pgfsetrectcap%
\pgfsetmiterjoin%
\pgfsetlinewidth{1.254687pt}%
\definecolor{currentstroke}{rgb}{0.150000,0.150000,0.150000}%
\pgfsetstrokecolor{currentstroke}%
\pgfsetdash{}{0pt}%
\pgfpathmoveto{\pgfqpoint{1.319707in}{0.899712in}}%
\pgfpathlineto{\pgfqpoint{1.319707in}{7.099712in}}%
\pgfusepath{stroke}%
\end{pgfscope}%
\begin{pgfscope}%
\pgfsetbuttcap%
\pgfsetmiterjoin%
\definecolor{currentfill}{rgb}{1.000000,1.000000,1.000000}%
\pgfsetfillcolor{currentfill}%
\pgfsetlinewidth{0.240900pt}%
\definecolor{currentstroke}{rgb}{0.150000,0.150000,0.150000}%
\pgfsetstrokecolor{currentstroke}%
\pgfsetdash{}{0pt}%
\pgfpathmoveto{\pgfqpoint{3.949427in}{5.400917in}}%
\pgfpathlineto{\pgfqpoint{7.353040in}{5.400917in}}%
\pgfpathlineto{\pgfqpoint{7.353040in}{6.933045in}}%
\pgfpathlineto{\pgfqpoint{3.949427in}{6.933045in}}%
\pgfpathclose%
\pgfusepath{stroke,fill}%
\end{pgfscope}%
\begin{pgfscope}%
\pgfsetroundcap%
\pgfsetroundjoin%
\pgfsetlinewidth{2.007500pt}%
\definecolor{currentstroke}{rgb}{0.298039,0.447059,0.690196}%
\pgfsetstrokecolor{currentstroke}%
\pgfsetdash{}{0pt}%
\pgfpathmoveto{\pgfqpoint{4.082760in}{6.674842in}}%
\pgfpathlineto{\pgfqpoint{4.749427in}{6.674842in}}%
\pgfusepath{stroke}%
\end{pgfscope}%
\begin{pgfscope}%
\pgfsetbuttcap%
\pgfsetroundjoin%
\definecolor{currentfill}{rgb}{0.298039,0.447059,0.690196}%
\pgfsetfillcolor{currentfill}%
\pgfsetlinewidth{0.000000pt}%
\definecolor{currentstroke}{rgb}{0.000000,0.000000,0.000000}%
\pgfsetstrokecolor{currentstroke}%
\pgfsetdash{}{0pt}%
\pgfsys@defobject{currentmarker}{\pgfqpoint{-0.038889in}{-0.038889in}}{\pgfqpoint{0.038889in}{0.038889in}}{%
\pgfpathmoveto{\pgfqpoint{0.000000in}{-0.038889in}}%
\pgfpathcurveto{\pgfqpoint{0.010313in}{-0.038889in}}{\pgfqpoint{0.020206in}{-0.034791in}}{\pgfqpoint{0.027499in}{-0.027499in}}%
\pgfpathcurveto{\pgfqpoint{0.034791in}{-0.020206in}}{\pgfqpoint{0.038889in}{-0.010313in}}{\pgfqpoint{0.038889in}{0.000000in}}%
\pgfpathcurveto{\pgfqpoint{0.038889in}{0.010313in}}{\pgfqpoint{0.034791in}{0.020206in}}{\pgfqpoint{0.027499in}{0.027499in}}%
\pgfpathcurveto{\pgfqpoint{0.020206in}{0.034791in}}{\pgfqpoint{0.010313in}{0.038889in}}{\pgfqpoint{0.000000in}{0.038889in}}%
\pgfpathcurveto{\pgfqpoint{-0.010313in}{0.038889in}}{\pgfqpoint{-0.020206in}{0.034791in}}{\pgfqpoint{-0.027499in}{0.027499in}}%
\pgfpathcurveto{\pgfqpoint{-0.034791in}{0.020206in}}{\pgfqpoint{-0.038889in}{0.010313in}}{\pgfqpoint{-0.038889in}{0.000000in}}%
\pgfpathcurveto{\pgfqpoint{-0.038889in}{-0.010313in}}{\pgfqpoint{-0.034791in}{-0.020206in}}{\pgfqpoint{-0.027499in}{-0.027499in}}%
\pgfpathcurveto{\pgfqpoint{-0.020206in}{-0.034791in}}{\pgfqpoint{-0.010313in}{-0.038889in}}{\pgfqpoint{0.000000in}{-0.038889in}}%
\pgfpathclose%
\pgfusepath{fill}%
}%
\begin{pgfscope}%
\pgfsys@transformshift{4.416093in}{6.674842in}%
\pgfsys@useobject{currentmarker}{}%
\end{pgfscope}%
\end{pgfscope}%
\begin{pgfscope}%
\definecolor{textcolor}{rgb}{0.150000,0.150000,0.150000}%
\pgfsetstrokecolor{textcolor}%
\pgfsetfillcolor{textcolor}%
\pgftext[x=5.016093in,y=6.558175in,left,base]{\color{textcolor}\sffamily\fontsize{24.000000}{28.800000}\selectfont FS-QTT-solver}%
\end{pgfscope}%
\begin{pgfscope}%
\pgfsetroundcap%
\pgfsetroundjoin%
\pgfsetlinewidth{2.007500pt}%
\definecolor{currentstroke}{rgb}{0.333333,0.658824,0.407843}%
\pgfsetstrokecolor{currentstroke}%
\pgfsetdash{}{0pt}%
\pgfpathmoveto{\pgfqpoint{4.082760in}{6.197466in}}%
\pgfpathlineto{\pgfqpoint{4.749427in}{6.197466in}}%
\pgfusepath{stroke}%
\end{pgfscope}%
\begin{pgfscope}%
\pgfsetbuttcap%
\pgfsetmiterjoin%
\definecolor{currentfill}{rgb}{0.333333,0.658824,0.407843}%
\pgfsetfillcolor{currentfill}%
\pgfsetlinewidth{0.000000pt}%
\definecolor{currentstroke}{rgb}{0.000000,0.000000,0.000000}%
\pgfsetstrokecolor{currentstroke}%
\pgfsetdash{}{0pt}%
\pgfsys@defobject{currentmarker}{\pgfqpoint{-0.038889in}{-0.038889in}}{\pgfqpoint{0.038889in}{0.038889in}}{%
\pgfpathmoveto{\pgfqpoint{-0.038889in}{-0.038889in}}%
\pgfpathlineto{\pgfqpoint{0.038889in}{-0.038889in}}%
\pgfpathlineto{\pgfqpoint{0.038889in}{0.038889in}}%
\pgfpathlineto{\pgfqpoint{-0.038889in}{0.038889in}}%
\pgfpathclose%
\pgfusepath{fill}%
}%
\begin{pgfscope}%
\pgfsys@transformshift{4.416093in}{6.197466in}%
\pgfsys@useobject{currentmarker}{}%
\end{pgfscope}%
\end{pgfscope}%
\begin{pgfscope}%
\definecolor{textcolor}{rgb}{0.150000,0.150000,0.150000}%
\pgfsetstrokecolor{textcolor}%
\pgfsetfillcolor{textcolor}%
\pgftext[x=5.016093in,y=6.080799in,left,base]{\color{textcolor}\sffamily\fontsize{24.000000}{28.800000}\selectfont FD-QTT-solver}%
\end{pgfscope}%
\begin{pgfscope}%
\pgfsetroundcap%
\pgfsetroundjoin%
\pgfsetlinewidth{2.007500pt}%
\definecolor{currentstroke}{rgb}{0.768627,0.305882,0.321569}%
\pgfsetstrokecolor{currentstroke}%
\pgfsetdash{}{0pt}%
\pgfpathmoveto{\pgfqpoint{4.082760in}{5.720090in}}%
\pgfpathlineto{\pgfqpoint{4.749427in}{5.720090in}}%
\pgfusepath{stroke}%
\end{pgfscope}%
\begin{pgfscope}%
\pgfsetbuttcap%
\pgfsetmiterjoin%
\definecolor{currentfill}{rgb}{0.768627,0.305882,0.321569}%
\pgfsetfillcolor{currentfill}%
\pgfsetlinewidth{0.000000pt}%
\definecolor{currentstroke}{rgb}{0.000000,0.000000,0.000000}%
\pgfsetstrokecolor{currentstroke}%
\pgfsetdash{}{0pt}%
\pgfsys@defobject{currentmarker}{\pgfqpoint{-0.038889in}{-0.038889in}}{\pgfqpoint{0.038889in}{0.038889in}}{%
\pgfpathmoveto{\pgfqpoint{0.000000in}{0.038889in}}%
\pgfpathlineto{\pgfqpoint{-0.038889in}{-0.038889in}}%
\pgfpathlineto{\pgfqpoint{0.038889in}{-0.038889in}}%
\pgfpathclose%
\pgfusepath{fill}%
}%
\begin{pgfscope}%
\pgfsys@transformshift{4.416093in}{5.720090in}%
\pgfsys@useobject{currentmarker}{}%
\end{pgfscope}%
\end{pgfscope}%
\begin{pgfscope}%
\definecolor{textcolor}{rgb}{0.150000,0.150000,0.150000}%
\pgfsetstrokecolor{textcolor}%
\pgfsetfillcolor{textcolor}%
\pgftext[x=5.016093in,y=5.603423in,left,base]{\color{textcolor}\sffamily\fontsize{24.000000}{28.800000}\selectfont FD-solver}%
\end{pgfscope}%
\end{pgfpicture}%
\makeatother%
\endgroup%

%% file: res_analyt_all_u_der_err.pgf
\begingroup%
\makeatletter%
\begin{pgfpicture}%
\pgfpathrectangle{\pgfpointorigin}{\pgfqpoint{7.778185in}{7.355067in}}%
\pgfusepath{use as bounding box, clip}%
\begin{pgfscope}%
\pgfsetbuttcap%
\pgfsetmiterjoin%
\definecolor{currentfill}{rgb}{1.000000,1.000000,1.000000}%
\pgfsetfillcolor{currentfill}%
\pgfsetlinewidth{0.000000pt}%
\definecolor{currentstroke}{rgb}{1.000000,1.000000,1.000000}%
\pgfsetstrokecolor{currentstroke}%
\pgfsetdash{}{0pt}%
\pgfpathmoveto{\pgfqpoint{0.000000in}{0.000000in}}%
\pgfpathlineto{\pgfqpoint{7.778185in}{0.000000in}}%
\pgfpathlineto{\pgfqpoint{7.778185in}{7.355067in}}%
\pgfpathlineto{\pgfqpoint{0.000000in}{7.355067in}}%
\pgfpathclose%
\pgfusepath{fill}%
\end{pgfscope}%
\begin{pgfscope}%
\pgfsetbuttcap%
\pgfsetmiterjoin%
\definecolor{currentfill}{rgb}{1.000000,1.000000,1.000000}%
\pgfsetfillcolor{currentfill}%
\pgfsetlinewidth{0.000000pt}%
\definecolor{currentstroke}{rgb}{0.000000,0.000000,0.000000}%
\pgfsetstrokecolor{currentstroke}%
\pgfsetstrokeopacity{0.000000}%
\pgfsetdash{}{0pt}%
\pgfpathmoveto{\pgfqpoint{1.319707in}{0.899712in}}%
\pgfpathlineto{\pgfqpoint{7.519707in}{0.899712in}}%
\pgfpathlineto{\pgfqpoint{7.519707in}{7.099712in}}%
\pgfpathlineto{\pgfqpoint{1.319707in}{7.099712in}}%
\pgfpathclose%
\pgfusepath{fill}%
\end{pgfscope}%
\begin{pgfscope}%
\pgfpathrectangle{\pgfqpoint{1.319707in}{0.899712in}}{\pgfqpoint{6.200000in}{6.200000in}} %
\pgfusepath{clip}%
\pgfsetbuttcap%
\pgfsetroundjoin%
\pgfsetlinewidth{0.803000pt}%
\definecolor{currentstroke}{rgb}{0.800000,0.800000,0.800000}%
\pgfsetstrokecolor{currentstroke}%
\pgfsetdash{{1.000000pt}{3.000000pt}}{0.000000pt}%
\pgfpathmoveto{\pgfqpoint{1.319707in}{0.899712in}}%
\pgfpathlineto{\pgfqpoint{1.319707in}{7.099712in}}%
\pgfusepath{stroke}%
\end{pgfscope}%
\begin{pgfscope}%
\definecolor{textcolor}{rgb}{0.150000,0.150000,0.150000}%
\pgfsetstrokecolor{textcolor}%
\pgfsetfillcolor{textcolor}%
\pgftext[x=1.319707in,y=0.821934in,,top]{\color{textcolor}\sffamily\fontsize{24.000000}{28.800000}\selectfont \(\displaystyle 0\)}%
\end{pgfscope}%
\begin{pgfscope}%
\pgfpathrectangle{\pgfqpoint{1.319707in}{0.899712in}}{\pgfqpoint{6.200000in}{6.200000in}} %
\pgfusepath{clip}%
\pgfsetbuttcap%
\pgfsetroundjoin%
\pgfsetlinewidth{0.803000pt}%
\definecolor{currentstroke}{rgb}{0.800000,0.800000,0.800000}%
\pgfsetstrokecolor{currentstroke}%
\pgfsetdash{{1.000000pt}{3.000000pt}}{0.000000pt}%
\pgfpathmoveto{\pgfqpoint{2.353040in}{0.899712in}}%
\pgfpathlineto{\pgfqpoint{2.353040in}{7.099712in}}%
\pgfusepath{stroke}%
\end{pgfscope}%
\begin{pgfscope}%
\definecolor{textcolor}{rgb}{0.150000,0.150000,0.150000}%
\pgfsetstrokecolor{textcolor}%
\pgfsetfillcolor{textcolor}%
\pgftext[x=2.353040in,y=0.821934in,,top]{\color{textcolor}\sffamily\fontsize{24.000000}{28.800000}\selectfont \(\displaystyle 5\)}%
\end{pgfscope}%
\begin{pgfscope}%
\pgfpathrectangle{\pgfqpoint{1.319707in}{0.899712in}}{\pgfqpoint{6.200000in}{6.200000in}} %
\pgfusepath{clip}%
\pgfsetbuttcap%
\pgfsetroundjoin%
\pgfsetlinewidth{0.803000pt}%
\definecolor{currentstroke}{rgb}{0.800000,0.800000,0.800000}%
\pgfsetstrokecolor{currentstroke}%
\pgfsetdash{{1.000000pt}{3.000000pt}}{0.000000pt}%
\pgfpathmoveto{\pgfqpoint{3.386373in}{0.899712in}}%
\pgfpathlineto{\pgfqpoint{3.386373in}{7.099712in}}%
\pgfusepath{stroke}%
\end{pgfscope}%
\begin{pgfscope}%
\definecolor{textcolor}{rgb}{0.150000,0.150000,0.150000}%
\pgfsetstrokecolor{textcolor}%
\pgfsetfillcolor{textcolor}%
\pgftext[x=3.386373in,y=0.821934in,,top]{\color{textcolor}\sffamily\fontsize{24.000000}{28.800000}\selectfont \(\displaystyle 10\)}%
\end{pgfscope}%
\begin{pgfscope}%
\pgfpathrectangle{\pgfqpoint{1.319707in}{0.899712in}}{\pgfqpoint{6.200000in}{6.200000in}} %
\pgfusepath{clip}%
\pgfsetbuttcap%
\pgfsetroundjoin%
\pgfsetlinewidth{0.803000pt}%
\definecolor{currentstroke}{rgb}{0.800000,0.800000,0.800000}%
\pgfsetstrokecolor{currentstroke}%
\pgfsetdash{{1.000000pt}{3.000000pt}}{0.000000pt}%
\pgfpathmoveto{\pgfqpoint{4.419707in}{0.899712in}}%
\pgfpathlineto{\pgfqpoint{4.419707in}{7.099712in}}%
\pgfusepath{stroke}%
\end{pgfscope}%
\begin{pgfscope}%
\definecolor{textcolor}{rgb}{0.150000,0.150000,0.150000}%
\pgfsetstrokecolor{textcolor}%
\pgfsetfillcolor{textcolor}%
\pgftext[x=4.419707in,y=0.821934in,,top]{\color{textcolor}\sffamily\fontsize{24.000000}{28.800000}\selectfont \(\displaystyle 15\)}%
\end{pgfscope}%
\begin{pgfscope}%
\pgfpathrectangle{\pgfqpoint{1.319707in}{0.899712in}}{\pgfqpoint{6.200000in}{6.200000in}} %
\pgfusepath{clip}%
\pgfsetbuttcap%
\pgfsetroundjoin%
\pgfsetlinewidth{0.803000pt}%
\definecolor{currentstroke}{rgb}{0.800000,0.800000,0.800000}%
\pgfsetstrokecolor{currentstroke}%
\pgfsetdash{{1.000000pt}{3.000000pt}}{0.000000pt}%
\pgfpathmoveto{\pgfqpoint{5.453040in}{0.899712in}}%
\pgfpathlineto{\pgfqpoint{5.453040in}{7.099712in}}%
\pgfusepath{stroke}%
\end{pgfscope}%
\begin{pgfscope}%
\definecolor{textcolor}{rgb}{0.150000,0.150000,0.150000}%
\pgfsetstrokecolor{textcolor}%
\pgfsetfillcolor{textcolor}%
\pgftext[x=5.453040in,y=0.821934in,,top]{\color{textcolor}\sffamily\fontsize{24.000000}{28.800000}\selectfont \(\displaystyle 20\)}%
\end{pgfscope}%
\begin{pgfscope}%
\pgfpathrectangle{\pgfqpoint{1.319707in}{0.899712in}}{\pgfqpoint{6.200000in}{6.200000in}} %
\pgfusepath{clip}%
\pgfsetbuttcap%
\pgfsetroundjoin%
\pgfsetlinewidth{0.803000pt}%
\definecolor{currentstroke}{rgb}{0.800000,0.800000,0.800000}%
\pgfsetstrokecolor{currentstroke}%
\pgfsetdash{{1.000000pt}{3.000000pt}}{0.000000pt}%
\pgfpathmoveto{\pgfqpoint{6.486373in}{0.899712in}}%
\pgfpathlineto{\pgfqpoint{6.486373in}{7.099712in}}%
\pgfusepath{stroke}%
\end{pgfscope}%
\begin{pgfscope}%
\definecolor{textcolor}{rgb}{0.150000,0.150000,0.150000}%
\pgfsetstrokecolor{textcolor}%
\pgfsetfillcolor{textcolor}%
\pgftext[x=6.486373in,y=0.821934in,,top]{\color{textcolor}\sffamily\fontsize{24.000000}{28.800000}\selectfont \(\displaystyle 25\)}%
\end{pgfscope}%
\begin{pgfscope}%
\pgfpathrectangle{\pgfqpoint{1.319707in}{0.899712in}}{\pgfqpoint{6.200000in}{6.200000in}} %
\pgfusepath{clip}%
\pgfsetbuttcap%
\pgfsetroundjoin%
\pgfsetlinewidth{0.803000pt}%
\definecolor{currentstroke}{rgb}{0.800000,0.800000,0.800000}%
\pgfsetstrokecolor{currentstroke}%
\pgfsetdash{{1.000000pt}{3.000000pt}}{0.000000pt}%
\pgfpathmoveto{\pgfqpoint{7.519707in}{0.899712in}}%
\pgfpathlineto{\pgfqpoint{7.519707in}{7.099712in}}%
\pgfusepath{stroke}%
\end{pgfscope}%
\begin{pgfscope}%
\definecolor{textcolor}{rgb}{0.150000,0.150000,0.150000}%
\pgfsetstrokecolor{textcolor}%
\pgfsetfillcolor{textcolor}%
\pgftext[x=7.519707in,y=0.821934in,,top]{\color{textcolor}\sffamily\fontsize{24.000000}{28.800000}\selectfont \(\displaystyle 30\)}%
\end{pgfscope}%
\begin{pgfscope}%
\definecolor{textcolor}{rgb}{0.150000,0.150000,0.150000}%
\pgfsetstrokecolor{textcolor}%
\pgfsetfillcolor{textcolor}%
\pgftext[x=4.419707in,y=0.441780in,,top]{\color{textcolor}\sffamily\fontsize{26.400000}{31.680000}\selectfont d}%
\end{pgfscope}%
\begin{pgfscope}%
\pgfpathrectangle{\pgfqpoint{1.319707in}{0.899712in}}{\pgfqpoint{6.200000in}{6.200000in}} %
\pgfusepath{clip}%
\pgfsetbuttcap%
\pgfsetroundjoin%
\pgfsetlinewidth{0.803000pt}%
\definecolor{currentstroke}{rgb}{0.800000,0.800000,0.800000}%
\pgfsetstrokecolor{currentstroke}%
\pgfsetdash{{1.000000pt}{3.000000pt}}{0.000000pt}%
\pgfpathmoveto{\pgfqpoint{1.319707in}{0.899712in}}%
\pgfpathlineto{\pgfqpoint{7.519707in}{0.899712in}}%
\pgfusepath{stroke}%
\end{pgfscope}%
\begin{pgfscope}%
\definecolor{textcolor}{rgb}{0.150000,0.150000,0.150000}%
\pgfsetstrokecolor{textcolor}%
\pgfsetfillcolor{textcolor}%
\pgftext[x=1.241929in,y=0.899712in,right,]{\color{textcolor}\sffamily\fontsize{24.000000}{28.800000}\selectfont \(\displaystyle 10^{-10}\)}%
\end{pgfscope}%
\begin{pgfscope}%
\pgfpathrectangle{\pgfqpoint{1.319707in}{0.899712in}}{\pgfqpoint{6.200000in}{6.200000in}} %
\pgfusepath{clip}%
\pgfsetbuttcap%
\pgfsetroundjoin%
\pgfsetlinewidth{0.803000pt}%
\definecolor{currentstroke}{rgb}{0.800000,0.800000,0.800000}%
\pgfsetstrokecolor{currentstroke}%
\pgfsetdash{{1.000000pt}{3.000000pt}}{0.000000pt}%
\pgfpathmoveto{\pgfqpoint{1.319707in}{1.674712in}}%
\pgfpathlineto{\pgfqpoint{7.519707in}{1.674712in}}%
\pgfusepath{stroke}%
\end{pgfscope}%
\begin{pgfscope}%
\definecolor{textcolor}{rgb}{0.150000,0.150000,0.150000}%
\pgfsetstrokecolor{textcolor}%
\pgfsetfillcolor{textcolor}%
\pgftext[x=1.241929in,y=1.674712in,right,]{\color{textcolor}\sffamily\fontsize{24.000000}{28.800000}\selectfont \(\displaystyle 10^{-9}\)}%
\end{pgfscope}%
\begin{pgfscope}%
\pgfpathrectangle{\pgfqpoint{1.319707in}{0.899712in}}{\pgfqpoint{6.200000in}{6.200000in}} %
\pgfusepath{clip}%
\pgfsetbuttcap%
\pgfsetroundjoin%
\pgfsetlinewidth{0.803000pt}%
\definecolor{currentstroke}{rgb}{0.800000,0.800000,0.800000}%
\pgfsetstrokecolor{currentstroke}%
\pgfsetdash{{1.000000pt}{3.000000pt}}{0.000000pt}%
\pgfpathmoveto{\pgfqpoint{1.319707in}{2.449712in}}%
\pgfpathlineto{\pgfqpoint{7.519707in}{2.449712in}}%
\pgfusepath{stroke}%
\end{pgfscope}%
\begin{pgfscope}%
\definecolor{textcolor}{rgb}{0.150000,0.150000,0.150000}%
\pgfsetstrokecolor{textcolor}%
\pgfsetfillcolor{textcolor}%
\pgftext[x=1.241929in,y=2.449712in,right,]{\color{textcolor}\sffamily\fontsize{24.000000}{28.800000}\selectfont \(\displaystyle 10^{-8}\)}%
\end{pgfscope}%
\begin{pgfscope}%
\pgfpathrectangle{\pgfqpoint{1.319707in}{0.899712in}}{\pgfqpoint{6.200000in}{6.200000in}} %
\pgfusepath{clip}%
\pgfsetbuttcap%
\pgfsetroundjoin%
\pgfsetlinewidth{0.803000pt}%
\definecolor{currentstroke}{rgb}{0.800000,0.800000,0.800000}%
\pgfsetstrokecolor{currentstroke}%
\pgfsetdash{{1.000000pt}{3.000000pt}}{0.000000pt}%
\pgfpathmoveto{\pgfqpoint{1.319707in}{3.224712in}}%
\pgfpathlineto{\pgfqpoint{7.519707in}{3.224712in}}%
\pgfusepath{stroke}%
\end{pgfscope}%
\begin{pgfscope}%
\definecolor{textcolor}{rgb}{0.150000,0.150000,0.150000}%
\pgfsetstrokecolor{textcolor}%
\pgfsetfillcolor{textcolor}%
\pgftext[x=1.241929in,y=3.224712in,right,]{\color{textcolor}\sffamily\fontsize{24.000000}{28.800000}\selectfont \(\displaystyle 10^{-7}\)}%
\end{pgfscope}%
\begin{pgfscope}%
\pgfpathrectangle{\pgfqpoint{1.319707in}{0.899712in}}{\pgfqpoint{6.200000in}{6.200000in}} %
\pgfusepath{clip}%
\pgfsetbuttcap%
\pgfsetroundjoin%
\pgfsetlinewidth{0.803000pt}%
\definecolor{currentstroke}{rgb}{0.800000,0.800000,0.800000}%
\pgfsetstrokecolor{currentstroke}%
\pgfsetdash{{1.000000pt}{3.000000pt}}{0.000000pt}%
\pgfpathmoveto{\pgfqpoint{1.319707in}{3.999712in}}%
\pgfpathlineto{\pgfqpoint{7.519707in}{3.999712in}}%
\pgfusepath{stroke}%
\end{pgfscope}%
\begin{pgfscope}%
\definecolor{textcolor}{rgb}{0.150000,0.150000,0.150000}%
\pgfsetstrokecolor{textcolor}%
\pgfsetfillcolor{textcolor}%
\pgftext[x=1.241929in,y=3.999712in,right,]{\color{textcolor}\sffamily\fontsize{24.000000}{28.800000}\selectfont \(\displaystyle 10^{-6}\)}%
\end{pgfscope}%
\begin{pgfscope}%
\pgfpathrectangle{\pgfqpoint{1.319707in}{0.899712in}}{\pgfqpoint{6.200000in}{6.200000in}} %
\pgfusepath{clip}%
\pgfsetbuttcap%
\pgfsetroundjoin%
\pgfsetlinewidth{0.803000pt}%
\definecolor{currentstroke}{rgb}{0.800000,0.800000,0.800000}%
\pgfsetstrokecolor{currentstroke}%
\pgfsetdash{{1.000000pt}{3.000000pt}}{0.000000pt}%
\pgfpathmoveto{\pgfqpoint{1.319707in}{4.774712in}}%
\pgfpathlineto{\pgfqpoint{7.519707in}{4.774712in}}%
\pgfusepath{stroke}%
\end{pgfscope}%
\begin{pgfscope}%
\definecolor{textcolor}{rgb}{0.150000,0.150000,0.150000}%
\pgfsetstrokecolor{textcolor}%
\pgfsetfillcolor{textcolor}%
\pgftext[x=1.241929in,y=4.774712in,right,]{\color{textcolor}\sffamily\fontsize{24.000000}{28.800000}\selectfont \(\displaystyle 10^{-5}\)}%
\end{pgfscope}%
\begin{pgfscope}%
\pgfpathrectangle{\pgfqpoint{1.319707in}{0.899712in}}{\pgfqpoint{6.200000in}{6.200000in}} %
\pgfusepath{clip}%
\pgfsetbuttcap%
\pgfsetroundjoin%
\pgfsetlinewidth{0.803000pt}%
\definecolor{currentstroke}{rgb}{0.800000,0.800000,0.800000}%
\pgfsetstrokecolor{currentstroke}%
\pgfsetdash{{1.000000pt}{3.000000pt}}{0.000000pt}%
\pgfpathmoveto{\pgfqpoint{1.319707in}{5.549712in}}%
\pgfpathlineto{\pgfqpoint{7.519707in}{5.549712in}}%
\pgfusepath{stroke}%
\end{pgfscope}%
\begin{pgfscope}%
\definecolor{textcolor}{rgb}{0.150000,0.150000,0.150000}%
\pgfsetstrokecolor{textcolor}%
\pgfsetfillcolor{textcolor}%
\pgftext[x=1.241929in,y=5.549712in,right,]{\color{textcolor}\sffamily\fontsize{24.000000}{28.800000}\selectfont \(\displaystyle 10^{-4}\)}%
\end{pgfscope}%
\begin{pgfscope}%
\pgfpathrectangle{\pgfqpoint{1.319707in}{0.899712in}}{\pgfqpoint{6.200000in}{6.200000in}} %
\pgfusepath{clip}%
\pgfsetbuttcap%
\pgfsetroundjoin%
\pgfsetlinewidth{0.803000pt}%
\definecolor{currentstroke}{rgb}{0.800000,0.800000,0.800000}%
\pgfsetstrokecolor{currentstroke}%
\pgfsetdash{{1.000000pt}{3.000000pt}}{0.000000pt}%
\pgfpathmoveto{\pgfqpoint{1.319707in}{6.324712in}}%
\pgfpathlineto{\pgfqpoint{7.519707in}{6.324712in}}%
\pgfusepath{stroke}%
\end{pgfscope}%
\begin{pgfscope}%
\definecolor{textcolor}{rgb}{0.150000,0.150000,0.150000}%
\pgfsetstrokecolor{textcolor}%
\pgfsetfillcolor{textcolor}%
\pgftext[x=1.241929in,y=6.324712in,right,]{\color{textcolor}\sffamily\fontsize{24.000000}{28.800000}\selectfont \(\displaystyle 10^{-3}\)}%
\end{pgfscope}%
\begin{pgfscope}%
\pgfpathrectangle{\pgfqpoint{1.319707in}{0.899712in}}{\pgfqpoint{6.200000in}{6.200000in}} %
\pgfusepath{clip}%
\pgfsetbuttcap%
\pgfsetroundjoin%
\pgfsetlinewidth{0.803000pt}%
\definecolor{currentstroke}{rgb}{0.800000,0.800000,0.800000}%
\pgfsetstrokecolor{currentstroke}%
\pgfsetdash{{1.000000pt}{3.000000pt}}{0.000000pt}%
\pgfpathmoveto{\pgfqpoint{1.319707in}{7.099712in}}%
\pgfpathlineto{\pgfqpoint{7.519707in}{7.099712in}}%
\pgfusepath{stroke}%
\end{pgfscope}%
\begin{pgfscope}%
\definecolor{textcolor}{rgb}{0.150000,0.150000,0.150000}%
\pgfsetstrokecolor{textcolor}%
\pgfsetfillcolor{textcolor}%
\pgftext[x=1.241929in,y=7.099712in,right,]{\color{textcolor}\sffamily\fontsize{24.000000}{28.800000}\selectfont \(\displaystyle 10^{-2}\)}%
\end{pgfscope}%
\begin{pgfscope}%
\definecolor{textcolor}{rgb}{0.150000,0.150000,0.150000}%
\pgfsetstrokecolor{textcolor}%
\pgfsetfillcolor{textcolor}%
\pgftext[x=0.441780in,y=3.999712in,,bottom,rotate=90.000000]{\color{textcolor}\sffamily\fontsize{26.400000}{31.680000}\selectfont Derivative relative error}%
\end{pgfscope}%
\begin{pgfscope}%
\pgfpathrectangle{\pgfqpoint{1.319707in}{0.899712in}}{\pgfqpoint{6.200000in}{6.200000in}} %
\pgfusepath{clip}%
\pgfsetroundcap%
\pgfsetroundjoin%
\pgfsetlinewidth{2.007500pt}%
\definecolor{currentstroke}{rgb}{0.298039,0.447059,0.690196}%
\pgfsetstrokecolor{currentstroke}%
\pgfsetdash{}{0pt}%
\pgfpathmoveto{\pgfqpoint{2.146373in}{7.019030in}}%
\pgfpathlineto{\pgfqpoint{2.353040in}{6.551855in}}%
\pgfpathlineto{\pgfqpoint{2.559707in}{6.085116in}}%
\pgfpathlineto{\pgfqpoint{2.766373in}{5.618484in}}%
\pgfpathlineto{\pgfqpoint{2.973040in}{5.151879in}}%
\pgfpathlineto{\pgfqpoint{3.179707in}{4.685280in}}%
\pgfpathlineto{\pgfqpoint{3.386373in}{4.218683in}}%
\pgfpathlineto{\pgfqpoint{3.593040in}{3.752087in}}%
\pgfpathlineto{\pgfqpoint{3.799707in}{3.285493in}}%
\pgfpathlineto{\pgfqpoint{4.006373in}{2.818902in}}%
\pgfpathlineto{\pgfqpoint{4.213040in}{2.352269in}}%
\pgfpathlineto{\pgfqpoint{4.419707in}{1.885893in}}%
\pgfpathlineto{\pgfqpoint{4.626373in}{1.427153in}}%
\pgfpathlineto{\pgfqpoint{4.833040in}{1.086961in}}%
\pgfpathlineto{\pgfqpoint{5.039707in}{1.103008in}}%
\pgfpathlineto{\pgfqpoint{5.246373in}{1.131013in}}%
\pgfpathlineto{\pgfqpoint{5.453040in}{1.143716in}}%
\pgfpathlineto{\pgfqpoint{5.659707in}{1.388924in}}%
\pgfpathlineto{\pgfqpoint{5.866373in}{1.517225in}}%
\pgfpathlineto{\pgfqpoint{6.073040in}{1.655855in}}%
\pgfpathlineto{\pgfqpoint{6.279707in}{1.658327in}}%
\pgfpathlineto{\pgfqpoint{6.486373in}{1.827359in}}%
\pgfpathlineto{\pgfqpoint{6.693040in}{1.991315in}}%
\pgfpathlineto{\pgfqpoint{6.899707in}{2.185509in}}%
\pgfpathlineto{\pgfqpoint{7.106373in}{2.332614in}}%
\pgfpathlineto{\pgfqpoint{7.313040in}{2.347493in}}%
\pgfpathlineto{\pgfqpoint{7.519707in}{2.570298in}}%
\pgfusepath{stroke}%
\end{pgfscope}%
\begin{pgfscope}%
\pgfpathrectangle{\pgfqpoint{1.319707in}{0.899712in}}{\pgfqpoint{6.200000in}{6.200000in}} %
\pgfusepath{clip}%
\pgfsetbuttcap%
\pgfsetroundjoin%
\definecolor{currentfill}{rgb}{0.298039,0.447059,0.690196}%
\pgfsetfillcolor{currentfill}%
\pgfsetlinewidth{0.000000pt}%
\definecolor{currentstroke}{rgb}{0.000000,0.000000,0.000000}%
\pgfsetstrokecolor{currentstroke}%
\pgfsetdash{}{0pt}%
\pgfsys@defobject{currentmarker}{\pgfqpoint{-0.038889in}{-0.038889in}}{\pgfqpoint{0.038889in}{0.038889in}}{%
\pgfpathmoveto{\pgfqpoint{0.000000in}{-0.038889in}}%
\pgfpathcurveto{\pgfqpoint{0.010313in}{-0.038889in}}{\pgfqpoint{0.020206in}{-0.034791in}}{\pgfqpoint{0.027499in}{-0.027499in}}%
\pgfpathcurveto{\pgfqpoint{0.034791in}{-0.020206in}}{\pgfqpoint{0.038889in}{-0.010313in}}{\pgfqpoint{0.038889in}{0.000000in}}%
\pgfpathcurveto{\pgfqpoint{0.038889in}{0.010313in}}{\pgfqpoint{0.034791in}{0.020206in}}{\pgfqpoint{0.027499in}{0.027499in}}%
\pgfpathcurveto{\pgfqpoint{0.020206in}{0.034791in}}{\pgfqpoint{0.010313in}{0.038889in}}{\pgfqpoint{0.000000in}{0.038889in}}%
\pgfpathcurveto{\pgfqpoint{-0.010313in}{0.038889in}}{\pgfqpoint{-0.020206in}{0.034791in}}{\pgfqpoint{-0.027499in}{0.027499in}}%
\pgfpathcurveto{\pgfqpoint{-0.034791in}{0.020206in}}{\pgfqpoint{-0.038889in}{0.010313in}}{\pgfqpoint{-0.038889in}{0.000000in}}%
\pgfpathcurveto{\pgfqpoint{-0.038889in}{-0.010313in}}{\pgfqpoint{-0.034791in}{-0.020206in}}{\pgfqpoint{-0.027499in}{-0.027499in}}%
\pgfpathcurveto{\pgfqpoint{-0.020206in}{-0.034791in}}{\pgfqpoint{-0.010313in}{-0.038889in}}{\pgfqpoint{0.000000in}{-0.038889in}}%
\pgfpathclose%
\pgfusepath{fill}%
}%
\begin{pgfscope}%
\pgfsys@transformshift{2.146373in}{7.019030in}%
\pgfsys@useobject{currentmarker}{}%
\end{pgfscope}%
\begin{pgfscope}%
\pgfsys@transformshift{2.353040in}{6.551855in}%
\pgfsys@useobject{currentmarker}{}%
\end{pgfscope}%
\begin{pgfscope}%
\pgfsys@transformshift{2.559707in}{6.085116in}%
\pgfsys@useobject{currentmarker}{}%
\end{pgfscope}%
\begin{pgfscope}%
\pgfsys@transformshift{2.766373in}{5.618484in}%
\pgfsys@useobject{currentmarker}{}%
\end{pgfscope}%
\begin{pgfscope}%
\pgfsys@transformshift{2.973040in}{5.151879in}%
\pgfsys@useobject{currentmarker}{}%
\end{pgfscope}%
\begin{pgfscope}%
\pgfsys@transformshift{3.179707in}{4.685280in}%
\pgfsys@useobject{currentmarker}{}%
\end{pgfscope}%
\begin{pgfscope}%
\pgfsys@transformshift{3.386373in}{4.218683in}%
\pgfsys@useobject{currentmarker}{}%
\end{pgfscope}%
\begin{pgfscope}%
\pgfsys@transformshift{3.593040in}{3.752087in}%
\pgfsys@useobject{currentmarker}{}%
\end{pgfscope}%
\begin{pgfscope}%
\pgfsys@transformshift{3.799707in}{3.285493in}%
\pgfsys@useobject{currentmarker}{}%
\end{pgfscope}%
\begin{pgfscope}%
\pgfsys@transformshift{4.006373in}{2.818902in}%
\pgfsys@useobject{currentmarker}{}%
\end{pgfscope}%
\begin{pgfscope}%
\pgfsys@transformshift{4.213040in}{2.352269in}%
\pgfsys@useobject{currentmarker}{}%
\end{pgfscope}%
\begin{pgfscope}%
\pgfsys@transformshift{4.419707in}{1.885893in}%
\pgfsys@useobject{currentmarker}{}%
\end{pgfscope}%
\begin{pgfscope}%
\pgfsys@transformshift{4.626373in}{1.427153in}%
\pgfsys@useobject{currentmarker}{}%
\end{pgfscope}%
\begin{pgfscope}%
\pgfsys@transformshift{4.833040in}{1.086961in}%
\pgfsys@useobject{currentmarker}{}%
\end{pgfscope}%
\begin{pgfscope}%
\pgfsys@transformshift{5.039707in}{1.103008in}%
\pgfsys@useobject{currentmarker}{}%
\end{pgfscope}%
\begin{pgfscope}%
\pgfsys@transformshift{5.246373in}{1.131013in}%
\pgfsys@useobject{currentmarker}{}%
\end{pgfscope}%
\begin{pgfscope}%
\pgfsys@transformshift{5.453040in}{1.143716in}%
\pgfsys@useobject{currentmarker}{}%
\end{pgfscope}%
\begin{pgfscope}%
\pgfsys@transformshift{5.659707in}{1.388924in}%
\pgfsys@useobject{currentmarker}{}%
\end{pgfscope}%
\begin{pgfscope}%
\pgfsys@transformshift{5.866373in}{1.517225in}%
\pgfsys@useobject{currentmarker}{}%
\end{pgfscope}%
\begin{pgfscope}%
\pgfsys@transformshift{6.073040in}{1.655855in}%
\pgfsys@useobject{currentmarker}{}%
\end{pgfscope}%
\begin{pgfscope}%
\pgfsys@transformshift{6.279707in}{1.658327in}%
\pgfsys@useobject{currentmarker}{}%
\end{pgfscope}%
\begin{pgfscope}%
\pgfsys@transformshift{6.486373in}{1.827359in}%
\pgfsys@useobject{currentmarker}{}%
\end{pgfscope}%
\begin{pgfscope}%
\pgfsys@transformshift{6.693040in}{1.991315in}%
\pgfsys@useobject{currentmarker}{}%
\end{pgfscope}%
\begin{pgfscope}%
\pgfsys@transformshift{6.899707in}{2.185509in}%
\pgfsys@useobject{currentmarker}{}%
\end{pgfscope}%
\begin{pgfscope}%
\pgfsys@transformshift{7.106373in}{2.332614in}%
\pgfsys@useobject{currentmarker}{}%
\end{pgfscope}%
\begin{pgfscope}%
\pgfsys@transformshift{7.313040in}{2.347493in}%
\pgfsys@useobject{currentmarker}{}%
\end{pgfscope}%
\begin{pgfscope}%
\pgfsys@transformshift{7.519707in}{2.570298in}%
\pgfsys@useobject{currentmarker}{}%
\end{pgfscope}%
\end{pgfscope}%
\begin{pgfscope}%
\pgfpathrectangle{\pgfqpoint{1.319707in}{0.899712in}}{\pgfqpoint{6.200000in}{6.200000in}} %
\pgfusepath{clip}%
\pgfsetbuttcap%
\pgfsetroundjoin%
\pgfsetlinewidth{2.007500pt}%
\definecolor{currentstroke}{rgb}{0.298039,0.447059,0.690196}%
\pgfsetstrokecolor{currentstroke}%
\pgfsetdash{{6.000000pt}{6.000000pt}}{0.000000pt}%
\pgfpathmoveto{\pgfqpoint{2.146373in}{6.830420in}}%
\pgfpathlineto{\pgfqpoint{2.353040in}{6.364855in}}%
\pgfpathlineto{\pgfqpoint{2.559707in}{5.898500in}}%
\pgfpathlineto{\pgfqpoint{2.766373in}{5.431963in}}%
\pgfpathlineto{\pgfqpoint{2.973040in}{4.965382in}}%
\pgfpathlineto{\pgfqpoint{3.179707in}{4.498789in}}%
\pgfpathlineto{\pgfqpoint{3.386373in}{4.032193in}}%
\pgfpathlineto{\pgfqpoint{3.593040in}{3.565596in}}%
\pgfpathlineto{\pgfqpoint{3.799707in}{3.098998in}}%
\pgfpathlineto{\pgfqpoint{4.006373in}{2.632397in}}%
\pgfpathlineto{\pgfqpoint{4.213040in}{2.165840in}}%
\pgfpathlineto{\pgfqpoint{4.419707in}{1.699504in}}%
\pgfpathlineto{\pgfqpoint{4.626373in}{1.238421in}}%
\pgfpathlineto{\pgfqpoint{4.833040in}{0.927497in}}%
\pgfpathlineto{\pgfqpoint{5.039707in}{0.995561in}}%
\pgfpathlineto{\pgfqpoint{5.246373in}{1.091426in}}%
\pgfpathlineto{\pgfqpoint{5.453040in}{1.199825in}}%
\pgfpathlineto{\pgfqpoint{5.659707in}{1.349419in}}%
\pgfpathlineto{\pgfqpoint{5.866373in}{1.458810in}}%
\pgfpathlineto{\pgfqpoint{6.073040in}{1.600358in}}%
\pgfpathlineto{\pgfqpoint{6.279707in}{1.686544in}}%
\pgfpathlineto{\pgfqpoint{6.486373in}{1.793730in}}%
\pgfpathlineto{\pgfqpoint{6.693040in}{1.897157in}}%
\pgfpathlineto{\pgfqpoint{6.899707in}{2.105121in}}%
\pgfpathlineto{\pgfqpoint{7.106373in}{2.189662in}}%
\pgfpathlineto{\pgfqpoint{7.313040in}{2.256597in}}%
\pgfpathlineto{\pgfqpoint{7.519707in}{2.412199in}}%
\pgfusepath{stroke}%
\end{pgfscope}%
\begin{pgfscope}%
\pgfpathrectangle{\pgfqpoint{1.319707in}{0.899712in}}{\pgfqpoint{6.200000in}{6.200000in}} %
\pgfusepath{clip}%
\pgfsetbuttcap%
\pgfsetroundjoin%
\definecolor{currentfill}{rgb}{0.298039,0.447059,0.690196}%
\pgfsetfillcolor{currentfill}%
\pgfsetlinewidth{0.000000pt}%
\definecolor{currentstroke}{rgb}{0.000000,0.000000,0.000000}%
\pgfsetstrokecolor{currentstroke}%
\pgfsetdash{}{0pt}%
\pgfsys@defobject{currentmarker}{\pgfqpoint{-0.038889in}{-0.038889in}}{\pgfqpoint{0.038889in}{0.038889in}}{%
\pgfpathmoveto{\pgfqpoint{0.000000in}{-0.038889in}}%
\pgfpathcurveto{\pgfqpoint{0.010313in}{-0.038889in}}{\pgfqpoint{0.020206in}{-0.034791in}}{\pgfqpoint{0.027499in}{-0.027499in}}%
\pgfpathcurveto{\pgfqpoint{0.034791in}{-0.020206in}}{\pgfqpoint{0.038889in}{-0.010313in}}{\pgfqpoint{0.038889in}{0.000000in}}%
\pgfpathcurveto{\pgfqpoint{0.038889in}{0.010313in}}{\pgfqpoint{0.034791in}{0.020206in}}{\pgfqpoint{0.027499in}{0.027499in}}%
\pgfpathcurveto{\pgfqpoint{0.020206in}{0.034791in}}{\pgfqpoint{0.010313in}{0.038889in}}{\pgfqpoint{0.000000in}{0.038889in}}%
\pgfpathcurveto{\pgfqpoint{-0.010313in}{0.038889in}}{\pgfqpoint{-0.020206in}{0.034791in}}{\pgfqpoint{-0.027499in}{0.027499in}}%
\pgfpathcurveto{\pgfqpoint{-0.034791in}{0.020206in}}{\pgfqpoint{-0.038889in}{0.010313in}}{\pgfqpoint{-0.038889in}{0.000000in}}%
\pgfpathcurveto{\pgfqpoint{-0.038889in}{-0.010313in}}{\pgfqpoint{-0.034791in}{-0.020206in}}{\pgfqpoint{-0.027499in}{-0.027499in}}%
\pgfpathcurveto{\pgfqpoint{-0.020206in}{-0.034791in}}{\pgfqpoint{-0.010313in}{-0.038889in}}{\pgfqpoint{0.000000in}{-0.038889in}}%
\pgfpathclose%
\pgfusepath{fill}%
}%
\begin{pgfscope}%
\pgfsys@transformshift{2.146373in}{6.830420in}%
\pgfsys@useobject{currentmarker}{}%
\end{pgfscope}%
\begin{pgfscope}%
\pgfsys@transformshift{2.353040in}{6.364855in}%
\pgfsys@useobject{currentmarker}{}%
\end{pgfscope}%
\begin{pgfscope}%
\pgfsys@transformshift{2.559707in}{5.898500in}%
\pgfsys@useobject{currentmarker}{}%
\end{pgfscope}%
\begin{pgfscope}%
\pgfsys@transformshift{2.766373in}{5.431963in}%
\pgfsys@useobject{currentmarker}{}%
\end{pgfscope}%
\begin{pgfscope}%
\pgfsys@transformshift{2.973040in}{4.965382in}%
\pgfsys@useobject{currentmarker}{}%
\end{pgfscope}%
\begin{pgfscope}%
\pgfsys@transformshift{3.179707in}{4.498789in}%
\pgfsys@useobject{currentmarker}{}%
\end{pgfscope}%
\begin{pgfscope}%
\pgfsys@transformshift{3.386373in}{4.032193in}%
\pgfsys@useobject{currentmarker}{}%
\end{pgfscope}%
\begin{pgfscope}%
\pgfsys@transformshift{3.593040in}{3.565596in}%
\pgfsys@useobject{currentmarker}{}%
\end{pgfscope}%
\begin{pgfscope}%
\pgfsys@transformshift{3.799707in}{3.098998in}%
\pgfsys@useobject{currentmarker}{}%
\end{pgfscope}%
\begin{pgfscope}%
\pgfsys@transformshift{4.006373in}{2.632397in}%
\pgfsys@useobject{currentmarker}{}%
\end{pgfscope}%
\begin{pgfscope}%
\pgfsys@transformshift{4.213040in}{2.165840in}%
\pgfsys@useobject{currentmarker}{}%
\end{pgfscope}%
\begin{pgfscope}%
\pgfsys@transformshift{4.419707in}{1.699504in}%
\pgfsys@useobject{currentmarker}{}%
\end{pgfscope}%
\begin{pgfscope}%
\pgfsys@transformshift{4.626373in}{1.238421in}%
\pgfsys@useobject{currentmarker}{}%
\end{pgfscope}%
\begin{pgfscope}%
\pgfsys@transformshift{4.833040in}{0.927497in}%
\pgfsys@useobject{currentmarker}{}%
\end{pgfscope}%
\begin{pgfscope}%
\pgfsys@transformshift{5.039707in}{0.995561in}%
\pgfsys@useobject{currentmarker}{}%
\end{pgfscope}%
\begin{pgfscope}%
\pgfsys@transformshift{5.246373in}{1.091426in}%
\pgfsys@useobject{currentmarker}{}%
\end{pgfscope}%
\begin{pgfscope}%
\pgfsys@transformshift{5.453040in}{1.199825in}%
\pgfsys@useobject{currentmarker}{}%
\end{pgfscope}%
\begin{pgfscope}%
\pgfsys@transformshift{5.659707in}{1.349419in}%
\pgfsys@useobject{currentmarker}{}%
\end{pgfscope}%
\begin{pgfscope}%
\pgfsys@transformshift{5.866373in}{1.458810in}%
\pgfsys@useobject{currentmarker}{}%
\end{pgfscope}%
\begin{pgfscope}%
\pgfsys@transformshift{6.073040in}{1.600358in}%
\pgfsys@useobject{currentmarker}{}%
\end{pgfscope}%
\begin{pgfscope}%
\pgfsys@transformshift{6.279707in}{1.686544in}%
\pgfsys@useobject{currentmarker}{}%
\end{pgfscope}%
\begin{pgfscope}%
\pgfsys@transformshift{6.486373in}{1.793730in}%
\pgfsys@useobject{currentmarker}{}%
\end{pgfscope}%
\begin{pgfscope}%
\pgfsys@transformshift{6.693040in}{1.897157in}%
\pgfsys@useobject{currentmarker}{}%
\end{pgfscope}%
\begin{pgfscope}%
\pgfsys@transformshift{6.899707in}{2.105121in}%
\pgfsys@useobject{currentmarker}{}%
\end{pgfscope}%
\begin{pgfscope}%
\pgfsys@transformshift{7.106373in}{2.189662in}%
\pgfsys@useobject{currentmarker}{}%
\end{pgfscope}%
\begin{pgfscope}%
\pgfsys@transformshift{7.313040in}{2.256597in}%
\pgfsys@useobject{currentmarker}{}%
\end{pgfscope}%
\begin{pgfscope}%
\pgfsys@transformshift{7.519707in}{2.412199in}%
\pgfsys@useobject{currentmarker}{}%
\end{pgfscope}%
\end{pgfscope}%
\begin{pgfscope}%
\pgfpathrectangle{\pgfqpoint{1.319707in}{0.899712in}}{\pgfqpoint{6.200000in}{6.200000in}} %
\pgfusepath{clip}%
\pgfsetroundcap%
\pgfsetroundjoin%
\pgfsetlinewidth{2.007500pt}%
\definecolor{currentstroke}{rgb}{0.333333,0.658824,0.407843}%
\pgfsetstrokecolor{currentstroke}%
\pgfsetdash{}{0pt}%
\pgfpathmoveto{\pgfqpoint{2.146373in}{7.019030in}}%
\pgfpathlineto{\pgfqpoint{2.353040in}{6.551855in}}%
\pgfpathlineto{\pgfqpoint{2.559707in}{6.085116in}}%
\pgfpathlineto{\pgfqpoint{2.766373in}{5.618484in}}%
\pgfpathlineto{\pgfqpoint{2.973040in}{5.151878in}}%
\pgfpathlineto{\pgfqpoint{3.179707in}{4.685251in}}%
\pgfpathlineto{\pgfqpoint{3.386373in}{4.218381in}}%
\pgfpathlineto{\pgfqpoint{3.593040in}{3.792011in}}%
\pgfpathlineto{\pgfqpoint{3.799707in}{3.960709in}}%
\pgfpathlineto{\pgfqpoint{4.006373in}{5.259524in}}%
\pgfpathlineto{\pgfqpoint{4.213040in}{6.132671in}}%
\pgfpathlineto{\pgfqpoint{4.419707in}{6.718053in}}%
\pgfusepath{stroke}%
\end{pgfscope}%
\begin{pgfscope}%
\pgfpathrectangle{\pgfqpoint{1.319707in}{0.899712in}}{\pgfqpoint{6.200000in}{6.200000in}} %
\pgfusepath{clip}%
\pgfsetbuttcap%
\pgfsetmiterjoin%
\definecolor{currentfill}{rgb}{0.333333,0.658824,0.407843}%
\pgfsetfillcolor{currentfill}%
\pgfsetlinewidth{0.000000pt}%
\definecolor{currentstroke}{rgb}{0.000000,0.000000,0.000000}%
\pgfsetstrokecolor{currentstroke}%
\pgfsetdash{}{0pt}%
\pgfsys@defobject{currentmarker}{\pgfqpoint{-0.038889in}{-0.038889in}}{\pgfqpoint{0.038889in}{0.038889in}}{%
\pgfpathmoveto{\pgfqpoint{-0.038889in}{-0.038889in}}%
\pgfpathlineto{\pgfqpoint{0.038889in}{-0.038889in}}%
\pgfpathlineto{\pgfqpoint{0.038889in}{0.038889in}}%
\pgfpathlineto{\pgfqpoint{-0.038889in}{0.038889in}}%
\pgfpathclose%
\pgfusepath{fill}%
}%
\begin{pgfscope}%
\pgfsys@transformshift{2.146373in}{7.019030in}%
\pgfsys@useobject{currentmarker}{}%
\end{pgfscope}%
\begin{pgfscope}%
\pgfsys@transformshift{2.353040in}{6.551855in}%
\pgfsys@useobject{currentmarker}{}%
\end{pgfscope}%
\begin{pgfscope}%
\pgfsys@transformshift{2.559707in}{6.085116in}%
\pgfsys@useobject{currentmarker}{}%
\end{pgfscope}%
\begin{pgfscope}%
\pgfsys@transformshift{2.766373in}{5.618484in}%
\pgfsys@useobject{currentmarker}{}%
\end{pgfscope}%
\begin{pgfscope}%
\pgfsys@transformshift{2.973040in}{5.151878in}%
\pgfsys@useobject{currentmarker}{}%
\end{pgfscope}%
\begin{pgfscope}%
\pgfsys@transformshift{3.179707in}{4.685251in}%
\pgfsys@useobject{currentmarker}{}%
\end{pgfscope}%
\begin{pgfscope}%
\pgfsys@transformshift{3.386373in}{4.218381in}%
\pgfsys@useobject{currentmarker}{}%
\end{pgfscope}%
\begin{pgfscope}%
\pgfsys@transformshift{3.593040in}{3.792011in}%
\pgfsys@useobject{currentmarker}{}%
\end{pgfscope}%
\begin{pgfscope}%
\pgfsys@transformshift{3.799707in}{3.960709in}%
\pgfsys@useobject{currentmarker}{}%
\end{pgfscope}%
\begin{pgfscope}%
\pgfsys@transformshift{4.006373in}{5.259524in}%
\pgfsys@useobject{currentmarker}{}%
\end{pgfscope}%
\begin{pgfscope}%
\pgfsys@transformshift{4.213040in}{6.132671in}%
\pgfsys@useobject{currentmarker}{}%
\end{pgfscope}%
\begin{pgfscope}%
\pgfsys@transformshift{4.419707in}{6.718053in}%
\pgfsys@useobject{currentmarker}{}%
\end{pgfscope}%
\end{pgfscope}%
\begin{pgfscope}%
\pgfpathrectangle{\pgfqpoint{1.319707in}{0.899712in}}{\pgfqpoint{6.200000in}{6.200000in}} %
\pgfusepath{clip}%
\pgfsetbuttcap%
\pgfsetroundjoin%
\pgfsetlinewidth{2.007500pt}%
\definecolor{currentstroke}{rgb}{0.333333,0.658824,0.407843}%
\pgfsetstrokecolor{currentstroke}%
\pgfsetdash{{6.000000pt}{6.000000pt}}{0.000000pt}%
\pgfpathmoveto{\pgfqpoint{2.146373in}{6.830420in}}%
\pgfpathlineto{\pgfqpoint{2.353040in}{6.364855in}}%
\pgfpathlineto{\pgfqpoint{2.559707in}{5.898500in}}%
\pgfpathlineto{\pgfqpoint{2.766373in}{5.431963in}}%
\pgfpathlineto{\pgfqpoint{2.973040in}{4.965381in}}%
\pgfpathlineto{\pgfqpoint{3.179707in}{4.498737in}}%
\pgfpathlineto{\pgfqpoint{3.386373in}{4.031591in}}%
\pgfpathlineto{\pgfqpoint{3.593040in}{3.552593in}}%
\pgfpathlineto{\pgfqpoint{3.799707in}{3.028034in}}%
\pgfpathlineto{\pgfqpoint{4.006373in}{3.594744in}}%
\pgfpathlineto{\pgfqpoint{4.213040in}{4.159540in}}%
\pgfpathlineto{\pgfqpoint{4.419707in}{4.681567in}}%
\pgfusepath{stroke}%
\end{pgfscope}%
\begin{pgfscope}%
\pgfpathrectangle{\pgfqpoint{1.319707in}{0.899712in}}{\pgfqpoint{6.200000in}{6.200000in}} %
\pgfusepath{clip}%
\pgfsetbuttcap%
\pgfsetmiterjoin%
\definecolor{currentfill}{rgb}{0.333333,0.658824,0.407843}%
\pgfsetfillcolor{currentfill}%
\pgfsetlinewidth{0.000000pt}%
\definecolor{currentstroke}{rgb}{0.000000,0.000000,0.000000}%
\pgfsetstrokecolor{currentstroke}%
\pgfsetdash{}{0pt}%
\pgfsys@defobject{currentmarker}{\pgfqpoint{-0.038889in}{-0.038889in}}{\pgfqpoint{0.038889in}{0.038889in}}{%
\pgfpathmoveto{\pgfqpoint{-0.038889in}{-0.038889in}}%
\pgfpathlineto{\pgfqpoint{0.038889in}{-0.038889in}}%
\pgfpathlineto{\pgfqpoint{0.038889in}{0.038889in}}%
\pgfpathlineto{\pgfqpoint{-0.038889in}{0.038889in}}%
\pgfpathclose%
\pgfusepath{fill}%
}%
\begin{pgfscope}%
\pgfsys@transformshift{2.146373in}{6.830420in}%
\pgfsys@useobject{currentmarker}{}%
\end{pgfscope}%
\begin{pgfscope}%
\pgfsys@transformshift{2.353040in}{6.364855in}%
\pgfsys@useobject{currentmarker}{}%
\end{pgfscope}%
\begin{pgfscope}%
\pgfsys@transformshift{2.559707in}{5.898500in}%
\pgfsys@useobject{currentmarker}{}%
\end{pgfscope}%
\begin{pgfscope}%
\pgfsys@transformshift{2.766373in}{5.431963in}%
\pgfsys@useobject{currentmarker}{}%
\end{pgfscope}%
\begin{pgfscope}%
\pgfsys@transformshift{2.973040in}{4.965381in}%
\pgfsys@useobject{currentmarker}{}%
\end{pgfscope}%
\begin{pgfscope}%
\pgfsys@transformshift{3.179707in}{4.498737in}%
\pgfsys@useobject{currentmarker}{}%
\end{pgfscope}%
\begin{pgfscope}%
\pgfsys@transformshift{3.386373in}{4.031591in}%
\pgfsys@useobject{currentmarker}{}%
\end{pgfscope}%
\begin{pgfscope}%
\pgfsys@transformshift{3.593040in}{3.552593in}%
\pgfsys@useobject{currentmarker}{}%
\end{pgfscope}%
\begin{pgfscope}%
\pgfsys@transformshift{3.799707in}{3.028034in}%
\pgfsys@useobject{currentmarker}{}%
\end{pgfscope}%
\begin{pgfscope}%
\pgfsys@transformshift{4.006373in}{3.594744in}%
\pgfsys@useobject{currentmarker}{}%
\end{pgfscope}%
\begin{pgfscope}%
\pgfsys@transformshift{4.213040in}{4.159540in}%
\pgfsys@useobject{currentmarker}{}%
\end{pgfscope}%
\begin{pgfscope}%
\pgfsys@transformshift{4.419707in}{4.681567in}%
\pgfsys@useobject{currentmarker}{}%
\end{pgfscope}%
\end{pgfscope}%
\begin{pgfscope}%
\pgfpathrectangle{\pgfqpoint{1.319707in}{0.899712in}}{\pgfqpoint{6.200000in}{6.200000in}} %
\pgfusepath{clip}%
\pgfsetroundcap%
\pgfsetroundjoin%
\pgfsetlinewidth{2.007500pt}%
\definecolor{currentstroke}{rgb}{0.768627,0.305882,0.321569}%
\pgfsetstrokecolor{currentstroke}%
\pgfsetdash{}{0pt}%
\pgfpathmoveto{\pgfqpoint{2.146373in}{7.019030in}}%
\pgfpathlineto{\pgfqpoint{2.353040in}{6.551855in}}%
\pgfpathlineto{\pgfqpoint{2.559707in}{6.085116in}}%
\pgfpathlineto{\pgfqpoint{2.766373in}{5.618484in}}%
\pgfpathlineto{\pgfqpoint{2.973040in}{5.151879in}}%
\pgfpathlineto{\pgfqpoint{3.179707in}{4.685280in}}%
\pgfpathlineto{\pgfqpoint{3.386373in}{4.218683in}}%
\pgfusepath{stroke}%
\end{pgfscope}%
\begin{pgfscope}%
\pgfpathrectangle{\pgfqpoint{1.319707in}{0.899712in}}{\pgfqpoint{6.200000in}{6.200000in}} %
\pgfusepath{clip}%
\pgfsetbuttcap%
\pgfsetmiterjoin%
\definecolor{currentfill}{rgb}{0.768627,0.305882,0.321569}%
\pgfsetfillcolor{currentfill}%
\pgfsetlinewidth{0.000000pt}%
\definecolor{currentstroke}{rgb}{0.000000,0.000000,0.000000}%
\pgfsetstrokecolor{currentstroke}%
\pgfsetdash{}{0pt}%
\pgfsys@defobject{currentmarker}{\pgfqpoint{-0.038889in}{-0.038889in}}{\pgfqpoint{0.038889in}{0.038889in}}{%
\pgfpathmoveto{\pgfqpoint{0.000000in}{0.038889in}}%
\pgfpathlineto{\pgfqpoint{-0.038889in}{-0.038889in}}%
\pgfpathlineto{\pgfqpoint{0.038889in}{-0.038889in}}%
\pgfpathclose%
\pgfusepath{fill}%
}%
\begin{pgfscope}%
\pgfsys@transformshift{2.146373in}{7.019030in}%
\pgfsys@useobject{currentmarker}{}%
\end{pgfscope}%
\begin{pgfscope}%
\pgfsys@transformshift{2.353040in}{6.551855in}%
\pgfsys@useobject{currentmarker}{}%
\end{pgfscope}%
\begin{pgfscope}%
\pgfsys@transformshift{2.559707in}{6.085116in}%
\pgfsys@useobject{currentmarker}{}%
\end{pgfscope}%
\begin{pgfscope}%
\pgfsys@transformshift{2.766373in}{5.618484in}%
\pgfsys@useobject{currentmarker}{}%
\end{pgfscope}%
\begin{pgfscope}%
\pgfsys@transformshift{2.973040in}{5.151879in}%
\pgfsys@useobject{currentmarker}{}%
\end{pgfscope}%
\begin{pgfscope}%
\pgfsys@transformshift{3.179707in}{4.685280in}%
\pgfsys@useobject{currentmarker}{}%
\end{pgfscope}%
\begin{pgfscope}%
\pgfsys@transformshift{3.386373in}{4.218683in}%
\pgfsys@useobject{currentmarker}{}%
\end{pgfscope}%
\end{pgfscope}%
\begin{pgfscope}%
\pgfpathrectangle{\pgfqpoint{1.319707in}{0.899712in}}{\pgfqpoint{6.200000in}{6.200000in}} %
\pgfusepath{clip}%
\pgfsetbuttcap%
\pgfsetroundjoin%
\pgfsetlinewidth{2.007500pt}%
\definecolor{currentstroke}{rgb}{0.768627,0.305882,0.321569}%
\pgfsetstrokecolor{currentstroke}%
\pgfsetdash{{6.000000pt}{6.000000pt}}{0.000000pt}%
\pgfpathmoveto{\pgfqpoint{2.146373in}{6.830420in}}%
\pgfpathlineto{\pgfqpoint{2.353040in}{6.364855in}}%
\pgfpathlineto{\pgfqpoint{2.559707in}{5.898500in}}%
\pgfpathlineto{\pgfqpoint{2.766373in}{5.431963in}}%
\pgfpathlineto{\pgfqpoint{2.973040in}{4.965382in}}%
\pgfpathlineto{\pgfqpoint{3.179707in}{4.498789in}}%
\pgfpathlineto{\pgfqpoint{3.386373in}{4.032193in}}%
\pgfusepath{stroke}%
\end{pgfscope}%
\begin{pgfscope}%
\pgfpathrectangle{\pgfqpoint{1.319707in}{0.899712in}}{\pgfqpoint{6.200000in}{6.200000in}} %
\pgfusepath{clip}%
\pgfsetbuttcap%
\pgfsetmiterjoin%
\definecolor{currentfill}{rgb}{0.768627,0.305882,0.321569}%
\pgfsetfillcolor{currentfill}%
\pgfsetlinewidth{0.000000pt}%
\definecolor{currentstroke}{rgb}{0.000000,0.000000,0.000000}%
\pgfsetstrokecolor{currentstroke}%
\pgfsetdash{}{0pt}%
\pgfsys@defobject{currentmarker}{\pgfqpoint{-0.038889in}{-0.038889in}}{\pgfqpoint{0.038889in}{0.038889in}}{%
\pgfpathmoveto{\pgfqpoint{0.000000in}{0.038889in}}%
\pgfpathlineto{\pgfqpoint{-0.038889in}{-0.038889in}}%
\pgfpathlineto{\pgfqpoint{0.038889in}{-0.038889in}}%
\pgfpathclose%
\pgfusepath{fill}%
}%
\begin{pgfscope}%
\pgfsys@transformshift{2.146373in}{6.830420in}%
\pgfsys@useobject{currentmarker}{}%
\end{pgfscope}%
\begin{pgfscope}%
\pgfsys@transformshift{2.353040in}{6.364855in}%
\pgfsys@useobject{currentmarker}{}%
\end{pgfscope}%
\begin{pgfscope}%
\pgfsys@transformshift{2.559707in}{5.898500in}%
\pgfsys@useobject{currentmarker}{}%
\end{pgfscope}%
\begin{pgfscope}%
\pgfsys@transformshift{2.766373in}{5.431963in}%
\pgfsys@useobject{currentmarker}{}%
\end{pgfscope}%
\begin{pgfscope}%
\pgfsys@transformshift{2.973040in}{4.965382in}%
\pgfsys@useobject{currentmarker}{}%
\end{pgfscope}%
\begin{pgfscope}%
\pgfsys@transformshift{3.179707in}{4.498789in}%
\pgfsys@useobject{currentmarker}{}%
\end{pgfscope}%
\begin{pgfscope}%
\pgfsys@transformshift{3.386373in}{4.032193in}%
\pgfsys@useobject{currentmarker}{}%
\end{pgfscope}%
\end{pgfscope}%
\begin{pgfscope}%
\pgfsetrectcap%
\pgfsetmiterjoin%
\pgfsetlinewidth{1.254687pt}%
\definecolor{currentstroke}{rgb}{0.150000,0.150000,0.150000}%
\pgfsetstrokecolor{currentstroke}%
\pgfsetdash{}{0pt}%
\pgfpathmoveto{\pgfqpoint{1.319707in}{0.899712in}}%
\pgfpathlineto{\pgfqpoint{7.519707in}{0.899712in}}%
\pgfusepath{stroke}%
\end{pgfscope}%
\begin{pgfscope}%
\pgfsetrectcap%
\pgfsetmiterjoin%
\pgfsetlinewidth{1.254687pt}%
\definecolor{currentstroke}{rgb}{0.150000,0.150000,0.150000}%
\pgfsetstrokecolor{currentstroke}%
\pgfsetdash{}{0pt}%
\pgfpathmoveto{\pgfqpoint{1.319707in}{0.899712in}}%
\pgfpathlineto{\pgfqpoint{1.319707in}{7.099712in}}%
\pgfusepath{stroke}%
\end{pgfscope}%
\begin{pgfscope}%
\pgfsetbuttcap%
\pgfsetmiterjoin%
\definecolor{currentfill}{rgb}{1.000000,1.000000,1.000000}%
\pgfsetfillcolor{currentfill}%
\pgfsetlinewidth{0.240900pt}%
\definecolor{currentstroke}{rgb}{0.150000,0.150000,0.150000}%
\pgfsetstrokecolor{currentstroke}%
\pgfsetdash{}{0pt}%
\pgfpathmoveto{\pgfqpoint{2.844448in}{5.298179in}}%
\pgfpathlineto{\pgfqpoint{7.353040in}{5.298179in}}%
\pgfpathlineto{\pgfqpoint{7.353040in}{6.933045in}}%
\pgfpathlineto{\pgfqpoint{2.844448in}{6.933045in}}%
\pgfpathclose%
\pgfusepath{stroke,fill}%
\end{pgfscope}%
\begin{pgfscope}%
\pgfsetroundcap%
\pgfsetroundjoin%
\pgfsetlinewidth{2.007500pt}%
\definecolor{currentstroke}{rgb}{0.298039,0.447059,0.690196}%
\pgfsetstrokecolor{currentstroke}%
\pgfsetdash{}{0pt}%
\pgfpathmoveto{\pgfqpoint{2.977781in}{6.657523in}}%
\pgfpathlineto{\pgfqpoint{3.644448in}{6.657523in}}%
\pgfusepath{stroke}%
\end{pgfscope}%
\begin{pgfscope}%
\pgfsetbuttcap%
\pgfsetroundjoin%
\definecolor{currentfill}{rgb}{0.298039,0.447059,0.690196}%
\pgfsetfillcolor{currentfill}%
\pgfsetlinewidth{0.000000pt}%
\definecolor{currentstroke}{rgb}{0.000000,0.000000,0.000000}%
\pgfsetstrokecolor{currentstroke}%
\pgfsetdash{}{0pt}%
\pgfsys@defobject{currentmarker}{\pgfqpoint{-0.038889in}{-0.038889in}}{\pgfqpoint{0.038889in}{0.038889in}}{%
\pgfpathmoveto{\pgfqpoint{0.000000in}{-0.038889in}}%
\pgfpathcurveto{\pgfqpoint{0.010313in}{-0.038889in}}{\pgfqpoint{0.020206in}{-0.034791in}}{\pgfqpoint{0.027499in}{-0.027499in}}%
\pgfpathcurveto{\pgfqpoint{0.034791in}{-0.020206in}}{\pgfqpoint{0.038889in}{-0.010313in}}{\pgfqpoint{0.038889in}{0.000000in}}%
\pgfpathcurveto{\pgfqpoint{0.038889in}{0.010313in}}{\pgfqpoint{0.034791in}{0.020206in}}{\pgfqpoint{0.027499in}{0.027499in}}%
\pgfpathcurveto{\pgfqpoint{0.020206in}{0.034791in}}{\pgfqpoint{0.010313in}{0.038889in}}{\pgfqpoint{0.000000in}{0.038889in}}%
\pgfpathcurveto{\pgfqpoint{-0.010313in}{0.038889in}}{\pgfqpoint{-0.020206in}{0.034791in}}{\pgfqpoint{-0.027499in}{0.027499in}}%
\pgfpathcurveto{\pgfqpoint{-0.034791in}{0.020206in}}{\pgfqpoint{-0.038889in}{0.010313in}}{\pgfqpoint{-0.038889in}{0.000000in}}%
\pgfpathcurveto{\pgfqpoint{-0.038889in}{-0.010313in}}{\pgfqpoint{-0.034791in}{-0.020206in}}{\pgfqpoint{-0.027499in}{-0.027499in}}%
\pgfpathcurveto{\pgfqpoint{-0.020206in}{-0.034791in}}{\pgfqpoint{-0.010313in}{-0.038889in}}{\pgfqpoint{0.000000in}{-0.038889in}}%
\pgfpathclose%
\pgfusepath{fill}%
}%
\begin{pgfscope}%
\pgfsys@transformshift{3.311115in}{6.657523in}%
\pgfsys@useobject{currentmarker}{}%
\end{pgfscope}%
\end{pgfscope}%
\begin{pgfscope}%
\definecolor{textcolor}{rgb}{0.150000,0.150000,0.150000}%
\pgfsetstrokecolor{textcolor}%
\pgfsetfillcolor{textcolor}%
\pgftext[x=3.911115in,y=6.540856in,left,base]{\color{textcolor}\sffamily\fontsize{24.000000}{28.800000}\selectfont FS-QTT-solver, \(\displaystyle \partial u / \partial x\)}%
\end{pgfscope}%
\begin{pgfscope}%
\pgfsetroundcap%
\pgfsetroundjoin%
\pgfsetlinewidth{2.007500pt}%
\definecolor{currentstroke}{rgb}{0.333333,0.658824,0.407843}%
\pgfsetstrokecolor{currentstroke}%
\pgfsetdash{}{0pt}%
\pgfpathmoveto{\pgfqpoint{2.977781in}{6.145901in}}%
\pgfpathlineto{\pgfqpoint{3.644448in}{6.145901in}}%
\pgfusepath{stroke}%
\end{pgfscope}%
\begin{pgfscope}%
\pgfsetbuttcap%
\pgfsetmiterjoin%
\definecolor{currentfill}{rgb}{0.333333,0.658824,0.407843}%
\pgfsetfillcolor{currentfill}%
\pgfsetlinewidth{0.000000pt}%
\definecolor{currentstroke}{rgb}{0.000000,0.000000,0.000000}%
\pgfsetstrokecolor{currentstroke}%
\pgfsetdash{}{0pt}%
\pgfsys@defobject{currentmarker}{\pgfqpoint{-0.038889in}{-0.038889in}}{\pgfqpoint{0.038889in}{0.038889in}}{%
\pgfpathmoveto{\pgfqpoint{-0.038889in}{-0.038889in}}%
\pgfpathlineto{\pgfqpoint{0.038889in}{-0.038889in}}%
\pgfpathlineto{\pgfqpoint{0.038889in}{0.038889in}}%
\pgfpathlineto{\pgfqpoint{-0.038889in}{0.038889in}}%
\pgfpathclose%
\pgfusepath{fill}%
}%
\begin{pgfscope}%
\pgfsys@transformshift{3.311115in}{6.145901in}%
\pgfsys@useobject{currentmarker}{}%
\end{pgfscope}%
\end{pgfscope}%
\begin{pgfscope}%
\definecolor{textcolor}{rgb}{0.150000,0.150000,0.150000}%
\pgfsetstrokecolor{textcolor}%
\pgfsetfillcolor{textcolor}%
\pgftext[x=3.911115in,y=6.029234in,left,base]{\color{textcolor}\sffamily\fontsize{24.000000}{28.800000}\selectfont FD-QTT-solver, \(\displaystyle \partial u / \partial x\)}%
\end{pgfscope}%
\begin{pgfscope}%
\pgfsetroundcap%
\pgfsetroundjoin%
\pgfsetlinewidth{2.007500pt}%
\definecolor{currentstroke}{rgb}{0.768627,0.305882,0.321569}%
\pgfsetstrokecolor{currentstroke}%
\pgfsetdash{}{0pt}%
\pgfpathmoveto{\pgfqpoint{2.977781in}{5.634279in}}%
\pgfpathlineto{\pgfqpoint{3.644448in}{5.634279in}}%
\pgfusepath{stroke}%
\end{pgfscope}%
\begin{pgfscope}%
\pgfsetbuttcap%
\pgfsetmiterjoin%
\definecolor{currentfill}{rgb}{0.768627,0.305882,0.321569}%
\pgfsetfillcolor{currentfill}%
\pgfsetlinewidth{0.000000pt}%
\definecolor{currentstroke}{rgb}{0.000000,0.000000,0.000000}%
\pgfsetstrokecolor{currentstroke}%
\pgfsetdash{}{0pt}%
\pgfsys@defobject{currentmarker}{\pgfqpoint{-0.038889in}{-0.038889in}}{\pgfqpoint{0.038889in}{0.038889in}}{%
\pgfpathmoveto{\pgfqpoint{0.000000in}{0.038889in}}%
\pgfpathlineto{\pgfqpoint{-0.038889in}{-0.038889in}}%
\pgfpathlineto{\pgfqpoint{0.038889in}{-0.038889in}}%
\pgfpathclose%
\pgfusepath{fill}%
}%
\begin{pgfscope}%
\pgfsys@transformshift{3.311115in}{5.634279in}%
\pgfsys@useobject{currentmarker}{}%
\end{pgfscope}%
\end{pgfscope}%
\begin{pgfscope}%
\definecolor{textcolor}{rgb}{0.150000,0.150000,0.150000}%
\pgfsetstrokecolor{textcolor}%
\pgfsetfillcolor{textcolor}%
\pgftext[x=3.911115in,y=5.517612in,left,base]{\color{textcolor}\sffamily\fontsize{24.000000}{28.800000}\selectfont FD-solver, \(\displaystyle \partial u / \partial x\)}%
\end{pgfscope}%
\end{pgfpicture}%
\makeatother%
\endgroup%

%% file: res_rhs1_all_u_calc_erank.pgf
\begingroup%
\makeatletter%
\begin{pgfpicture}%
\pgfpathrectangle{\pgfpointorigin}{\pgfqpoint{7.364436in}{7.355067in}}%
\pgfusepath{use as bounding box, clip}%
\begin{pgfscope}%
\pgfsetbuttcap%
\pgfsetmiterjoin%
\definecolor{currentfill}{rgb}{1.000000,1.000000,1.000000}%
\pgfsetfillcolor{currentfill}%
\pgfsetlinewidth{0.000000pt}%
\definecolor{currentstroke}{rgb}{1.000000,1.000000,1.000000}%
\pgfsetstrokecolor{currentstroke}%
\pgfsetdash{}{0pt}%
\pgfpathmoveto{\pgfqpoint{0.000000in}{0.000000in}}%
\pgfpathlineto{\pgfqpoint{7.364436in}{0.000000in}}%
\pgfpathlineto{\pgfqpoint{7.364436in}{7.355067in}}%
\pgfpathlineto{\pgfqpoint{0.000000in}{7.355067in}}%
\pgfpathclose%
\pgfusepath{fill}%
\end{pgfscope}%
\begin{pgfscope}%
\pgfsetbuttcap%
\pgfsetmiterjoin%
\definecolor{currentfill}{rgb}{1.000000,1.000000,1.000000}%
\pgfsetfillcolor{currentfill}%
\pgfsetlinewidth{0.000000pt}%
\definecolor{currentstroke}{rgb}{0.000000,0.000000,0.000000}%
\pgfsetstrokecolor{currentstroke}%
\pgfsetstrokeopacity{0.000000}%
\pgfsetdash{}{0pt}%
\pgfpathmoveto{\pgfqpoint{0.905958in}{0.899712in}}%
\pgfpathlineto{\pgfqpoint{7.105958in}{0.899712in}}%
\pgfpathlineto{\pgfqpoint{7.105958in}{7.099712in}}%
\pgfpathlineto{\pgfqpoint{0.905958in}{7.099712in}}%
\pgfpathclose%
\pgfusepath{fill}%
\end{pgfscope}%
\begin{pgfscope}%
\pgfpathrectangle{\pgfqpoint{0.905958in}{0.899712in}}{\pgfqpoint{6.200000in}{6.200000in}} %
\pgfusepath{clip}%
\pgfsetbuttcap%
\pgfsetroundjoin%
\pgfsetlinewidth{0.803000pt}%
\definecolor{currentstroke}{rgb}{0.800000,0.800000,0.800000}%
\pgfsetstrokecolor{currentstroke}%
\pgfsetdash{{1.000000pt}{3.000000pt}}{0.000000pt}%
\pgfpathmoveto{\pgfqpoint{0.905958in}{0.899712in}}%
\pgfpathlineto{\pgfqpoint{0.905958in}{7.099712in}}%
\pgfusepath{stroke}%
\end{pgfscope}%
\begin{pgfscope}%
\pgfsetbuttcap%
\pgfsetroundjoin%
\definecolor{currentfill}{rgb}{0.150000,0.150000,0.150000}%
\pgfsetfillcolor{currentfill}%
\pgfsetlinewidth{0.803000pt}%
\definecolor{currentstroke}{rgb}{0.150000,0.150000,0.150000}%
\pgfsetstrokecolor{currentstroke}%
\pgfsetdash{}{0pt}%
\pgfsys@defobject{currentmarker}{\pgfqpoint{0.000000in}{0.000000in}}{\pgfqpoint{0.000000in}{0.000000in}}{%
\pgfpathmoveto{\pgfqpoint{0.000000in}{0.000000in}}%
\pgfpathlineto{\pgfqpoint{0.000000in}{0.000000in}}%
\pgfusepath{stroke,fill}%
}%
\begin{pgfscope}%
\pgfsys@transformshift{0.905958in}{0.899712in}%
\pgfsys@useobject{currentmarker}{}%
\end{pgfscope}%
\end{pgfscope}%
\begin{pgfscope}%
\definecolor{textcolor}{rgb}{0.150000,0.150000,0.150000}%
\pgfsetstrokecolor{textcolor}%
\pgfsetfillcolor{textcolor}%
\pgftext[x=0.905958in,y=0.821934in,,top]{\color{textcolor}\sffamily\fontsize{24.000000}{28.800000}\selectfont \(\displaystyle 0\)}%
\end{pgfscope}%
\begin{pgfscope}%
\pgfpathrectangle{\pgfqpoint{0.905958in}{0.899712in}}{\pgfqpoint{6.200000in}{6.200000in}} %
\pgfusepath{clip}%
\pgfsetbuttcap%
\pgfsetroundjoin%
\pgfsetlinewidth{0.803000pt}%
\definecolor{currentstroke}{rgb}{0.800000,0.800000,0.800000}%
\pgfsetstrokecolor{currentstroke}%
\pgfsetdash{{1.000000pt}{3.000000pt}}{0.000000pt}%
\pgfpathmoveto{\pgfqpoint{1.939292in}{0.899712in}}%
\pgfpathlineto{\pgfqpoint{1.939292in}{7.099712in}}%
\pgfusepath{stroke}%
\end{pgfscope}%
\begin{pgfscope}%
\pgfsetbuttcap%
\pgfsetroundjoin%
\definecolor{currentfill}{rgb}{0.150000,0.150000,0.150000}%
\pgfsetfillcolor{currentfill}%
\pgfsetlinewidth{0.803000pt}%
\definecolor{currentstroke}{rgb}{0.150000,0.150000,0.150000}%
\pgfsetstrokecolor{currentstroke}%
\pgfsetdash{}{0pt}%
\pgfsys@defobject{currentmarker}{\pgfqpoint{0.000000in}{0.000000in}}{\pgfqpoint{0.000000in}{0.000000in}}{%
\pgfpathmoveto{\pgfqpoint{0.000000in}{0.000000in}}%
\pgfpathlineto{\pgfqpoint{0.000000in}{0.000000in}}%
\pgfusepath{stroke,fill}%
}%
\begin{pgfscope}%
\pgfsys@transformshift{1.939292in}{0.899712in}%
\pgfsys@useobject{currentmarker}{}%
\end{pgfscope}%
\end{pgfscope}%
\begin{pgfscope}%
\definecolor{textcolor}{rgb}{0.150000,0.150000,0.150000}%
\pgfsetstrokecolor{textcolor}%
\pgfsetfillcolor{textcolor}%
\pgftext[x=1.939292in,y=0.821934in,,top]{\color{textcolor}\sffamily\fontsize{24.000000}{28.800000}\selectfont \(\displaystyle 5\)}%
\end{pgfscope}%
\begin{pgfscope}%
\pgfpathrectangle{\pgfqpoint{0.905958in}{0.899712in}}{\pgfqpoint{6.200000in}{6.200000in}} %
\pgfusepath{clip}%
\pgfsetbuttcap%
\pgfsetroundjoin%
\pgfsetlinewidth{0.803000pt}%
\definecolor{currentstroke}{rgb}{0.800000,0.800000,0.800000}%
\pgfsetstrokecolor{currentstroke}%
\pgfsetdash{{1.000000pt}{3.000000pt}}{0.000000pt}%
\pgfpathmoveto{\pgfqpoint{2.972625in}{0.899712in}}%
\pgfpathlineto{\pgfqpoint{2.972625in}{7.099712in}}%
\pgfusepath{stroke}%
\end{pgfscope}%
\begin{pgfscope}%
\pgfsetbuttcap%
\pgfsetroundjoin%
\definecolor{currentfill}{rgb}{0.150000,0.150000,0.150000}%
\pgfsetfillcolor{currentfill}%
\pgfsetlinewidth{0.803000pt}%
\definecolor{currentstroke}{rgb}{0.150000,0.150000,0.150000}%
\pgfsetstrokecolor{currentstroke}%
\pgfsetdash{}{0pt}%
\pgfsys@defobject{currentmarker}{\pgfqpoint{0.000000in}{0.000000in}}{\pgfqpoint{0.000000in}{0.000000in}}{%
\pgfpathmoveto{\pgfqpoint{0.000000in}{0.000000in}}%
\pgfpathlineto{\pgfqpoint{0.000000in}{0.000000in}}%
\pgfusepath{stroke,fill}%
}%
\begin{pgfscope}%
\pgfsys@transformshift{2.972625in}{0.899712in}%
\pgfsys@useobject{currentmarker}{}%
\end{pgfscope}%
\end{pgfscope}%
\begin{pgfscope}%
\definecolor{textcolor}{rgb}{0.150000,0.150000,0.150000}%
\pgfsetstrokecolor{textcolor}%
\pgfsetfillcolor{textcolor}%
\pgftext[x=2.972625in,y=0.821934in,,top]{\color{textcolor}\sffamily\fontsize{24.000000}{28.800000}\selectfont \(\displaystyle 10\)}%
\end{pgfscope}%
\begin{pgfscope}%
\pgfpathrectangle{\pgfqpoint{0.905958in}{0.899712in}}{\pgfqpoint{6.200000in}{6.200000in}} %
\pgfusepath{clip}%
\pgfsetbuttcap%
\pgfsetroundjoin%
\pgfsetlinewidth{0.803000pt}%
\definecolor{currentstroke}{rgb}{0.800000,0.800000,0.800000}%
\pgfsetstrokecolor{currentstroke}%
\pgfsetdash{{1.000000pt}{3.000000pt}}{0.000000pt}%
\pgfpathmoveto{\pgfqpoint{4.005958in}{0.899712in}}%
\pgfpathlineto{\pgfqpoint{4.005958in}{7.099712in}}%
\pgfusepath{stroke}%
\end{pgfscope}%
\begin{pgfscope}%
\pgfsetbuttcap%
\pgfsetroundjoin%
\definecolor{currentfill}{rgb}{0.150000,0.150000,0.150000}%
\pgfsetfillcolor{currentfill}%
\pgfsetlinewidth{0.803000pt}%
\definecolor{currentstroke}{rgb}{0.150000,0.150000,0.150000}%
\pgfsetstrokecolor{currentstroke}%
\pgfsetdash{}{0pt}%
\pgfsys@defobject{currentmarker}{\pgfqpoint{0.000000in}{0.000000in}}{\pgfqpoint{0.000000in}{0.000000in}}{%
\pgfpathmoveto{\pgfqpoint{0.000000in}{0.000000in}}%
\pgfpathlineto{\pgfqpoint{0.000000in}{0.000000in}}%
\pgfusepath{stroke,fill}%
}%
\begin{pgfscope}%
\pgfsys@transformshift{4.005958in}{0.899712in}%
\pgfsys@useobject{currentmarker}{}%
\end{pgfscope}%
\end{pgfscope}%
\begin{pgfscope}%
\definecolor{textcolor}{rgb}{0.150000,0.150000,0.150000}%
\pgfsetstrokecolor{textcolor}%
\pgfsetfillcolor{textcolor}%
\pgftext[x=4.005958in,y=0.821934in,,top]{\color{textcolor}\sffamily\fontsize{24.000000}{28.800000}\selectfont \(\displaystyle 15\)}%
\end{pgfscope}%
\begin{pgfscope}%
\pgfpathrectangle{\pgfqpoint{0.905958in}{0.899712in}}{\pgfqpoint{6.200000in}{6.200000in}} %
\pgfusepath{clip}%
\pgfsetbuttcap%
\pgfsetroundjoin%
\pgfsetlinewidth{0.803000pt}%
\definecolor{currentstroke}{rgb}{0.800000,0.800000,0.800000}%
\pgfsetstrokecolor{currentstroke}%
\pgfsetdash{{1.000000pt}{3.000000pt}}{0.000000pt}%
\pgfpathmoveto{\pgfqpoint{5.039292in}{0.899712in}}%
\pgfpathlineto{\pgfqpoint{5.039292in}{7.099712in}}%
\pgfusepath{stroke}%
\end{pgfscope}%
\begin{pgfscope}%
\pgfsetbuttcap%
\pgfsetroundjoin%
\definecolor{currentfill}{rgb}{0.150000,0.150000,0.150000}%
\pgfsetfillcolor{currentfill}%
\pgfsetlinewidth{0.803000pt}%
\definecolor{currentstroke}{rgb}{0.150000,0.150000,0.150000}%
\pgfsetstrokecolor{currentstroke}%
\pgfsetdash{}{0pt}%
\pgfsys@defobject{currentmarker}{\pgfqpoint{0.000000in}{0.000000in}}{\pgfqpoint{0.000000in}{0.000000in}}{%
\pgfpathmoveto{\pgfqpoint{0.000000in}{0.000000in}}%
\pgfpathlineto{\pgfqpoint{0.000000in}{0.000000in}}%
\pgfusepath{stroke,fill}%
}%
\begin{pgfscope}%
\pgfsys@transformshift{5.039292in}{0.899712in}%
\pgfsys@useobject{currentmarker}{}%
\end{pgfscope}%
\end{pgfscope}%
\begin{pgfscope}%
\definecolor{textcolor}{rgb}{0.150000,0.150000,0.150000}%
\pgfsetstrokecolor{textcolor}%
\pgfsetfillcolor{textcolor}%
\pgftext[x=5.039292in,y=0.821934in,,top]{\color{textcolor}\sffamily\fontsize{24.000000}{28.800000}\selectfont \(\displaystyle 20\)}%
\end{pgfscope}%
\begin{pgfscope}%
\pgfpathrectangle{\pgfqpoint{0.905958in}{0.899712in}}{\pgfqpoint{6.200000in}{6.200000in}} %
\pgfusepath{clip}%
\pgfsetbuttcap%
\pgfsetroundjoin%
\pgfsetlinewidth{0.803000pt}%
\definecolor{currentstroke}{rgb}{0.800000,0.800000,0.800000}%
\pgfsetstrokecolor{currentstroke}%
\pgfsetdash{{1.000000pt}{3.000000pt}}{0.000000pt}%
\pgfpathmoveto{\pgfqpoint{6.072625in}{0.899712in}}%
\pgfpathlineto{\pgfqpoint{6.072625in}{7.099712in}}%
\pgfusepath{stroke}%
\end{pgfscope}%
\begin{pgfscope}%
\pgfsetbuttcap%
\pgfsetroundjoin%
\definecolor{currentfill}{rgb}{0.150000,0.150000,0.150000}%
\pgfsetfillcolor{currentfill}%
\pgfsetlinewidth{0.803000pt}%
\definecolor{currentstroke}{rgb}{0.150000,0.150000,0.150000}%
\pgfsetstrokecolor{currentstroke}%
\pgfsetdash{}{0pt}%
\pgfsys@defobject{currentmarker}{\pgfqpoint{0.000000in}{0.000000in}}{\pgfqpoint{0.000000in}{0.000000in}}{%
\pgfpathmoveto{\pgfqpoint{0.000000in}{0.000000in}}%
\pgfpathlineto{\pgfqpoint{0.000000in}{0.000000in}}%
\pgfusepath{stroke,fill}%
}%
\begin{pgfscope}%
\pgfsys@transformshift{6.072625in}{0.899712in}%
\pgfsys@useobject{currentmarker}{}%
\end{pgfscope}%
\end{pgfscope}%
\begin{pgfscope}%
\definecolor{textcolor}{rgb}{0.150000,0.150000,0.150000}%
\pgfsetstrokecolor{textcolor}%
\pgfsetfillcolor{textcolor}%
\pgftext[x=6.072625in,y=0.821934in,,top]{\color{textcolor}\sffamily\fontsize{24.000000}{28.800000}\selectfont \(\displaystyle 25\)}%
\end{pgfscope}%
\begin{pgfscope}%
\pgfpathrectangle{\pgfqpoint{0.905958in}{0.899712in}}{\pgfqpoint{6.200000in}{6.200000in}} %
\pgfusepath{clip}%
\pgfsetbuttcap%
\pgfsetroundjoin%
\pgfsetlinewidth{0.803000pt}%
\definecolor{currentstroke}{rgb}{0.800000,0.800000,0.800000}%
\pgfsetstrokecolor{currentstroke}%
\pgfsetdash{{1.000000pt}{3.000000pt}}{0.000000pt}%
\pgfpathmoveto{\pgfqpoint{7.105958in}{0.899712in}}%
\pgfpathlineto{\pgfqpoint{7.105958in}{7.099712in}}%
\pgfusepath{stroke}%
\end{pgfscope}%
\begin{pgfscope}%
\pgfsetbuttcap%
\pgfsetroundjoin%
\definecolor{currentfill}{rgb}{0.150000,0.150000,0.150000}%
\pgfsetfillcolor{currentfill}%
\pgfsetlinewidth{0.803000pt}%
\definecolor{currentstroke}{rgb}{0.150000,0.150000,0.150000}%
\pgfsetstrokecolor{currentstroke}%
\pgfsetdash{}{0pt}%
\pgfsys@defobject{currentmarker}{\pgfqpoint{0.000000in}{0.000000in}}{\pgfqpoint{0.000000in}{0.000000in}}{%
\pgfpathmoveto{\pgfqpoint{0.000000in}{0.000000in}}%
\pgfpathlineto{\pgfqpoint{0.000000in}{0.000000in}}%
\pgfusepath{stroke,fill}%
}%
\begin{pgfscope}%
\pgfsys@transformshift{7.105958in}{0.899712in}%
\pgfsys@useobject{currentmarker}{}%
\end{pgfscope}%
\end{pgfscope}%
\begin{pgfscope}%
\definecolor{textcolor}{rgb}{0.150000,0.150000,0.150000}%
\pgfsetstrokecolor{textcolor}%
\pgfsetfillcolor{textcolor}%
\pgftext[x=7.105958in,y=0.821934in,,top]{\color{textcolor}\sffamily\fontsize{24.000000}{28.800000}\selectfont \(\displaystyle 30\)}%
\end{pgfscope}%
\begin{pgfscope}%
\definecolor{textcolor}{rgb}{0.150000,0.150000,0.150000}%
\pgfsetstrokecolor{textcolor}%
\pgfsetfillcolor{textcolor}%
\pgftext[x=4.005958in,y=0.441780in,,top]{\color{textcolor}\sffamily\fontsize{26.400000}{31.680000}\selectfont d}%
\end{pgfscope}%
\begin{pgfscope}%
\pgfpathrectangle{\pgfqpoint{0.905958in}{0.899712in}}{\pgfqpoint{6.200000in}{6.200000in}} %
\pgfusepath{clip}%
\pgfsetbuttcap%
\pgfsetroundjoin%
\pgfsetlinewidth{0.803000pt}%
\definecolor{currentstroke}{rgb}{0.800000,0.800000,0.800000}%
\pgfsetstrokecolor{currentstroke}%
\pgfsetdash{{1.000000pt}{3.000000pt}}{0.000000pt}%
\pgfpathmoveto{\pgfqpoint{0.905958in}{0.899712in}}%
\pgfpathlineto{\pgfqpoint{7.105958in}{0.899712in}}%
\pgfusepath{stroke}%
\end{pgfscope}%
\begin{pgfscope}%
\pgfsetbuttcap%
\pgfsetroundjoin%
\definecolor{currentfill}{rgb}{0.150000,0.150000,0.150000}%
\pgfsetfillcolor{currentfill}%
\pgfsetlinewidth{0.803000pt}%
\definecolor{currentstroke}{rgb}{0.150000,0.150000,0.150000}%
\pgfsetstrokecolor{currentstroke}%
\pgfsetdash{}{0pt}%
\pgfsys@defobject{currentmarker}{\pgfqpoint{0.000000in}{0.000000in}}{\pgfqpoint{0.000000in}{0.000000in}}{%
\pgfpathmoveto{\pgfqpoint{0.000000in}{0.000000in}}%
\pgfpathlineto{\pgfqpoint{0.000000in}{0.000000in}}%
\pgfusepath{stroke,fill}%
}%
\begin{pgfscope}%
\pgfsys@transformshift{0.905958in}{0.899712in}%
\pgfsys@useobject{currentmarker}{}%
\end{pgfscope}%
\end{pgfscope}%
\begin{pgfscope}%
\definecolor{textcolor}{rgb}{0.150000,0.150000,0.150000}%
\pgfsetstrokecolor{textcolor}%
\pgfsetfillcolor{textcolor}%
\pgftext[x=0.828181in,y=0.899712in,right,]{\color{textcolor}\sffamily\fontsize{24.000000}{28.800000}\selectfont \(\displaystyle 6\)}%
\end{pgfscope}%
\begin{pgfscope}%
\pgfpathrectangle{\pgfqpoint{0.905958in}{0.899712in}}{\pgfqpoint{6.200000in}{6.200000in}} %
\pgfusepath{clip}%
\pgfsetbuttcap%
\pgfsetroundjoin%
\pgfsetlinewidth{0.803000pt}%
\definecolor{currentstroke}{rgb}{0.800000,0.800000,0.800000}%
\pgfsetstrokecolor{currentstroke}%
\pgfsetdash{{1.000000pt}{3.000000pt}}{0.000000pt}%
\pgfpathmoveto{\pgfqpoint{0.905958in}{1.674712in}}%
\pgfpathlineto{\pgfqpoint{7.105958in}{1.674712in}}%
\pgfusepath{stroke}%
\end{pgfscope}%
\begin{pgfscope}%
\pgfsetbuttcap%
\pgfsetroundjoin%
\definecolor{currentfill}{rgb}{0.150000,0.150000,0.150000}%
\pgfsetfillcolor{currentfill}%
\pgfsetlinewidth{0.803000pt}%
\definecolor{currentstroke}{rgb}{0.150000,0.150000,0.150000}%
\pgfsetstrokecolor{currentstroke}%
\pgfsetdash{}{0pt}%
\pgfsys@defobject{currentmarker}{\pgfqpoint{0.000000in}{0.000000in}}{\pgfqpoint{0.000000in}{0.000000in}}{%
\pgfpathmoveto{\pgfqpoint{0.000000in}{0.000000in}}%
\pgfpathlineto{\pgfqpoint{0.000000in}{0.000000in}}%
\pgfusepath{stroke,fill}%
}%
\begin{pgfscope}%
\pgfsys@transformshift{0.905958in}{1.674712in}%
\pgfsys@useobject{currentmarker}{}%
\end{pgfscope}%
\end{pgfscope}%
\begin{pgfscope}%
\definecolor{textcolor}{rgb}{0.150000,0.150000,0.150000}%
\pgfsetstrokecolor{textcolor}%
\pgfsetfillcolor{textcolor}%
\pgftext[x=0.828181in,y=1.674712in,right,]{\color{textcolor}\sffamily\fontsize{24.000000}{28.800000}\selectfont \(\displaystyle 7\)}%
\end{pgfscope}%
\begin{pgfscope}%
\pgfpathrectangle{\pgfqpoint{0.905958in}{0.899712in}}{\pgfqpoint{6.200000in}{6.200000in}} %
\pgfusepath{clip}%
\pgfsetbuttcap%
\pgfsetroundjoin%
\pgfsetlinewidth{0.803000pt}%
\definecolor{currentstroke}{rgb}{0.800000,0.800000,0.800000}%
\pgfsetstrokecolor{currentstroke}%
\pgfsetdash{{1.000000pt}{3.000000pt}}{0.000000pt}%
\pgfpathmoveto{\pgfqpoint{0.905958in}{2.449712in}}%
\pgfpathlineto{\pgfqpoint{7.105958in}{2.449712in}}%
\pgfusepath{stroke}%
\end{pgfscope}%
\begin{pgfscope}%
\pgfsetbuttcap%
\pgfsetroundjoin%
\definecolor{currentfill}{rgb}{0.150000,0.150000,0.150000}%
\pgfsetfillcolor{currentfill}%
\pgfsetlinewidth{0.803000pt}%
\definecolor{currentstroke}{rgb}{0.150000,0.150000,0.150000}%
\pgfsetstrokecolor{currentstroke}%
\pgfsetdash{}{0pt}%
\pgfsys@defobject{currentmarker}{\pgfqpoint{0.000000in}{0.000000in}}{\pgfqpoint{0.000000in}{0.000000in}}{%
\pgfpathmoveto{\pgfqpoint{0.000000in}{0.000000in}}%
\pgfpathlineto{\pgfqpoint{0.000000in}{0.000000in}}%
\pgfusepath{stroke,fill}%
}%
\begin{pgfscope}%
\pgfsys@transformshift{0.905958in}{2.449712in}%
\pgfsys@useobject{currentmarker}{}%
\end{pgfscope}%
\end{pgfscope}%
\begin{pgfscope}%
\definecolor{textcolor}{rgb}{0.150000,0.150000,0.150000}%
\pgfsetstrokecolor{textcolor}%
\pgfsetfillcolor{textcolor}%
\pgftext[x=0.828181in,y=2.449712in,right,]{\color{textcolor}\sffamily\fontsize{24.000000}{28.800000}\selectfont \(\displaystyle 8\)}%
\end{pgfscope}%
\begin{pgfscope}%
\pgfpathrectangle{\pgfqpoint{0.905958in}{0.899712in}}{\pgfqpoint{6.200000in}{6.200000in}} %
\pgfusepath{clip}%
\pgfsetbuttcap%
\pgfsetroundjoin%
\pgfsetlinewidth{0.803000pt}%
\definecolor{currentstroke}{rgb}{0.800000,0.800000,0.800000}%
\pgfsetstrokecolor{currentstroke}%
\pgfsetdash{{1.000000pt}{3.000000pt}}{0.000000pt}%
\pgfpathmoveto{\pgfqpoint{0.905958in}{3.224712in}}%
\pgfpathlineto{\pgfqpoint{7.105958in}{3.224712in}}%
\pgfusepath{stroke}%
\end{pgfscope}%
\begin{pgfscope}%
\pgfsetbuttcap%
\pgfsetroundjoin%
\definecolor{currentfill}{rgb}{0.150000,0.150000,0.150000}%
\pgfsetfillcolor{currentfill}%
\pgfsetlinewidth{0.803000pt}%
\definecolor{currentstroke}{rgb}{0.150000,0.150000,0.150000}%
\pgfsetstrokecolor{currentstroke}%
\pgfsetdash{}{0pt}%
\pgfsys@defobject{currentmarker}{\pgfqpoint{0.000000in}{0.000000in}}{\pgfqpoint{0.000000in}{0.000000in}}{%
\pgfpathmoveto{\pgfqpoint{0.000000in}{0.000000in}}%
\pgfpathlineto{\pgfqpoint{0.000000in}{0.000000in}}%
\pgfusepath{stroke,fill}%
}%
\begin{pgfscope}%
\pgfsys@transformshift{0.905958in}{3.224712in}%
\pgfsys@useobject{currentmarker}{}%
\end{pgfscope}%
\end{pgfscope}%
\begin{pgfscope}%
\definecolor{textcolor}{rgb}{0.150000,0.150000,0.150000}%
\pgfsetstrokecolor{textcolor}%
\pgfsetfillcolor{textcolor}%
\pgftext[x=0.828181in,y=3.224712in,right,]{\color{textcolor}\sffamily\fontsize{24.000000}{28.800000}\selectfont \(\displaystyle 9\)}%
\end{pgfscope}%
\begin{pgfscope}%
\pgfpathrectangle{\pgfqpoint{0.905958in}{0.899712in}}{\pgfqpoint{6.200000in}{6.200000in}} %
\pgfusepath{clip}%
\pgfsetbuttcap%
\pgfsetroundjoin%
\pgfsetlinewidth{0.803000pt}%
\definecolor{currentstroke}{rgb}{0.800000,0.800000,0.800000}%
\pgfsetstrokecolor{currentstroke}%
\pgfsetdash{{1.000000pt}{3.000000pt}}{0.000000pt}%
\pgfpathmoveto{\pgfqpoint{0.905958in}{3.999712in}}%
\pgfpathlineto{\pgfqpoint{7.105958in}{3.999712in}}%
\pgfusepath{stroke}%
\end{pgfscope}%
\begin{pgfscope}%
\pgfsetbuttcap%
\pgfsetroundjoin%
\definecolor{currentfill}{rgb}{0.150000,0.150000,0.150000}%
\pgfsetfillcolor{currentfill}%
\pgfsetlinewidth{0.803000pt}%
\definecolor{currentstroke}{rgb}{0.150000,0.150000,0.150000}%
\pgfsetstrokecolor{currentstroke}%
\pgfsetdash{}{0pt}%
\pgfsys@defobject{currentmarker}{\pgfqpoint{0.000000in}{0.000000in}}{\pgfqpoint{0.000000in}{0.000000in}}{%
\pgfpathmoveto{\pgfqpoint{0.000000in}{0.000000in}}%
\pgfpathlineto{\pgfqpoint{0.000000in}{0.000000in}}%
\pgfusepath{stroke,fill}%
}%
\begin{pgfscope}%
\pgfsys@transformshift{0.905958in}{3.999712in}%
\pgfsys@useobject{currentmarker}{}%
\end{pgfscope}%
\end{pgfscope}%
\begin{pgfscope}%
\definecolor{textcolor}{rgb}{0.150000,0.150000,0.150000}%
\pgfsetstrokecolor{textcolor}%
\pgfsetfillcolor{textcolor}%
\pgftext[x=0.828181in,y=3.999712in,right,]{\color{textcolor}\sffamily\fontsize{24.000000}{28.800000}\selectfont \(\displaystyle 10\)}%
\end{pgfscope}%
\begin{pgfscope}%
\pgfpathrectangle{\pgfqpoint{0.905958in}{0.899712in}}{\pgfqpoint{6.200000in}{6.200000in}} %
\pgfusepath{clip}%
\pgfsetbuttcap%
\pgfsetroundjoin%
\pgfsetlinewidth{0.803000pt}%
\definecolor{currentstroke}{rgb}{0.800000,0.800000,0.800000}%
\pgfsetstrokecolor{currentstroke}%
\pgfsetdash{{1.000000pt}{3.000000pt}}{0.000000pt}%
\pgfpathmoveto{\pgfqpoint{0.905958in}{4.774712in}}%
\pgfpathlineto{\pgfqpoint{7.105958in}{4.774712in}}%
\pgfusepath{stroke}%
\end{pgfscope}%
\begin{pgfscope}%
\pgfsetbuttcap%
\pgfsetroundjoin%
\definecolor{currentfill}{rgb}{0.150000,0.150000,0.150000}%
\pgfsetfillcolor{currentfill}%
\pgfsetlinewidth{0.803000pt}%
\definecolor{currentstroke}{rgb}{0.150000,0.150000,0.150000}%
\pgfsetstrokecolor{currentstroke}%
\pgfsetdash{}{0pt}%
\pgfsys@defobject{currentmarker}{\pgfqpoint{0.000000in}{0.000000in}}{\pgfqpoint{0.000000in}{0.000000in}}{%
\pgfpathmoveto{\pgfqpoint{0.000000in}{0.000000in}}%
\pgfpathlineto{\pgfqpoint{0.000000in}{0.000000in}}%
\pgfusepath{stroke,fill}%
}%
\begin{pgfscope}%
\pgfsys@transformshift{0.905958in}{4.774712in}%
\pgfsys@useobject{currentmarker}{}%
\end{pgfscope}%
\end{pgfscope}%
\begin{pgfscope}%
\definecolor{textcolor}{rgb}{0.150000,0.150000,0.150000}%
\pgfsetstrokecolor{textcolor}%
\pgfsetfillcolor{textcolor}%
\pgftext[x=0.828181in,y=4.774712in,right,]{\color{textcolor}\sffamily\fontsize{24.000000}{28.800000}\selectfont \(\displaystyle 11\)}%
\end{pgfscope}%
\begin{pgfscope}%
\pgfpathrectangle{\pgfqpoint{0.905958in}{0.899712in}}{\pgfqpoint{6.200000in}{6.200000in}} %
\pgfusepath{clip}%
\pgfsetbuttcap%
\pgfsetroundjoin%
\pgfsetlinewidth{0.803000pt}%
\definecolor{currentstroke}{rgb}{0.800000,0.800000,0.800000}%
\pgfsetstrokecolor{currentstroke}%
\pgfsetdash{{1.000000pt}{3.000000pt}}{0.000000pt}%
\pgfpathmoveto{\pgfqpoint{0.905958in}{5.549712in}}%
\pgfpathlineto{\pgfqpoint{7.105958in}{5.549712in}}%
\pgfusepath{stroke}%
\end{pgfscope}%
\begin{pgfscope}%
\pgfsetbuttcap%
\pgfsetroundjoin%
\definecolor{currentfill}{rgb}{0.150000,0.150000,0.150000}%
\pgfsetfillcolor{currentfill}%
\pgfsetlinewidth{0.803000pt}%
\definecolor{currentstroke}{rgb}{0.150000,0.150000,0.150000}%
\pgfsetstrokecolor{currentstroke}%
\pgfsetdash{}{0pt}%
\pgfsys@defobject{currentmarker}{\pgfqpoint{0.000000in}{0.000000in}}{\pgfqpoint{0.000000in}{0.000000in}}{%
\pgfpathmoveto{\pgfqpoint{0.000000in}{0.000000in}}%
\pgfpathlineto{\pgfqpoint{0.000000in}{0.000000in}}%
\pgfusepath{stroke,fill}%
}%
\begin{pgfscope}%
\pgfsys@transformshift{0.905958in}{5.549712in}%
\pgfsys@useobject{currentmarker}{}%
\end{pgfscope}%
\end{pgfscope}%
\begin{pgfscope}%
\definecolor{textcolor}{rgb}{0.150000,0.150000,0.150000}%
\pgfsetstrokecolor{textcolor}%
\pgfsetfillcolor{textcolor}%
\pgftext[x=0.828181in,y=5.549712in,right,]{\color{textcolor}\sffamily\fontsize{24.000000}{28.800000}\selectfont \(\displaystyle 12\)}%
\end{pgfscope}%
\begin{pgfscope}%
\pgfpathrectangle{\pgfqpoint{0.905958in}{0.899712in}}{\pgfqpoint{6.200000in}{6.200000in}} %
\pgfusepath{clip}%
\pgfsetbuttcap%
\pgfsetroundjoin%
\pgfsetlinewidth{0.803000pt}%
\definecolor{currentstroke}{rgb}{0.800000,0.800000,0.800000}%
\pgfsetstrokecolor{currentstroke}%
\pgfsetdash{{1.000000pt}{3.000000pt}}{0.000000pt}%
\pgfpathmoveto{\pgfqpoint{0.905958in}{6.324712in}}%
\pgfpathlineto{\pgfqpoint{7.105958in}{6.324712in}}%
\pgfusepath{stroke}%
\end{pgfscope}%
\begin{pgfscope}%
\pgfsetbuttcap%
\pgfsetroundjoin%
\definecolor{currentfill}{rgb}{0.150000,0.150000,0.150000}%
\pgfsetfillcolor{currentfill}%
\pgfsetlinewidth{0.803000pt}%
\definecolor{currentstroke}{rgb}{0.150000,0.150000,0.150000}%
\pgfsetstrokecolor{currentstroke}%
\pgfsetdash{}{0pt}%
\pgfsys@defobject{currentmarker}{\pgfqpoint{0.000000in}{0.000000in}}{\pgfqpoint{0.000000in}{0.000000in}}{%
\pgfpathmoveto{\pgfqpoint{0.000000in}{0.000000in}}%
\pgfpathlineto{\pgfqpoint{0.000000in}{0.000000in}}%
\pgfusepath{stroke,fill}%
}%
\begin{pgfscope}%
\pgfsys@transformshift{0.905958in}{6.324712in}%
\pgfsys@useobject{currentmarker}{}%
\end{pgfscope}%
\end{pgfscope}%
\begin{pgfscope}%
\definecolor{textcolor}{rgb}{0.150000,0.150000,0.150000}%
\pgfsetstrokecolor{textcolor}%
\pgfsetfillcolor{textcolor}%
\pgftext[x=0.828181in,y=6.324712in,right,]{\color{textcolor}\sffamily\fontsize{24.000000}{28.800000}\selectfont \(\displaystyle 13\)}%
\end{pgfscope}%
\begin{pgfscope}%
\pgfpathrectangle{\pgfqpoint{0.905958in}{0.899712in}}{\pgfqpoint{6.200000in}{6.200000in}} %
\pgfusepath{clip}%
\pgfsetbuttcap%
\pgfsetroundjoin%
\pgfsetlinewidth{0.803000pt}%
\definecolor{currentstroke}{rgb}{0.800000,0.800000,0.800000}%
\pgfsetstrokecolor{currentstroke}%
\pgfsetdash{{1.000000pt}{3.000000pt}}{0.000000pt}%
\pgfpathmoveto{\pgfqpoint{0.905958in}{7.099712in}}%
\pgfpathlineto{\pgfqpoint{7.105958in}{7.099712in}}%
\pgfusepath{stroke}%
\end{pgfscope}%
\begin{pgfscope}%
\pgfsetbuttcap%
\pgfsetroundjoin%
\definecolor{currentfill}{rgb}{0.150000,0.150000,0.150000}%
\pgfsetfillcolor{currentfill}%
\pgfsetlinewidth{0.803000pt}%
\definecolor{currentstroke}{rgb}{0.150000,0.150000,0.150000}%
\pgfsetstrokecolor{currentstroke}%
\pgfsetdash{}{0pt}%
\pgfsys@defobject{currentmarker}{\pgfqpoint{0.000000in}{0.000000in}}{\pgfqpoint{0.000000in}{0.000000in}}{%
\pgfpathmoveto{\pgfqpoint{0.000000in}{0.000000in}}%
\pgfpathlineto{\pgfqpoint{0.000000in}{0.000000in}}%
\pgfusepath{stroke,fill}%
}%
\begin{pgfscope}%
\pgfsys@transformshift{0.905958in}{7.099712in}%
\pgfsys@useobject{currentmarker}{}%
\end{pgfscope}%
\end{pgfscope}%
\begin{pgfscope}%
\definecolor{textcolor}{rgb}{0.150000,0.150000,0.150000}%
\pgfsetstrokecolor{textcolor}%
\pgfsetfillcolor{textcolor}%
\pgftext[x=0.828181in,y=7.099712in,right,]{\color{textcolor}\sffamily\fontsize{24.000000}{28.800000}\selectfont \(\displaystyle 14\)}%
\end{pgfscope}%
\begin{pgfscope}%
\definecolor{textcolor}{rgb}{0.150000,0.150000,0.150000}%
\pgfsetstrokecolor{textcolor}%
\pgfsetfillcolor{textcolor}%
\pgftext[x=0.441780in,y=3.999712in,,bottom,rotate=90.000000]{\color{textcolor}\sffamily\fontsize{26.400000}{31.680000}\selectfont Solution effective TT-rank}%
\end{pgfscope}%
\begin{pgfscope}%
\pgfpathrectangle{\pgfqpoint{0.905958in}{0.899712in}}{\pgfqpoint{6.200000in}{6.200000in}} %
\pgfusepath{clip}%
\pgfsetroundcap%
\pgfsetroundjoin%
\pgfsetlinewidth{2.007500pt}%
\definecolor{currentstroke}{rgb}{0.298039,0.447059,0.690196}%
\pgfsetstrokecolor{currentstroke}%
\pgfsetdash{}{0pt}%
\pgfpathmoveto{\pgfqpoint{1.732625in}{0.899712in}}%
\pgfpathlineto{\pgfqpoint{1.939292in}{2.407842in}}%
\pgfpathlineto{\pgfqpoint{2.145958in}{3.330376in}}%
\pgfpathlineto{\pgfqpoint{2.352625in}{3.755717in}}%
\pgfpathlineto{\pgfqpoint{2.559292in}{4.930792in}}%
\pgfpathlineto{\pgfqpoint{2.765958in}{5.389729in}}%
\pgfpathlineto{\pgfqpoint{2.972625in}{4.796105in}}%
\pgfpathlineto{\pgfqpoint{3.179292in}{5.331711in}}%
\pgfpathlineto{\pgfqpoint{3.385958in}{5.473362in}}%
\pgfpathlineto{\pgfqpoint{3.592625in}{4.283941in}}%
\pgfpathlineto{\pgfqpoint{3.799292in}{4.361007in}}%
\pgfpathlineto{\pgfqpoint{4.005958in}{6.260811in}}%
\pgfpathlineto{\pgfqpoint{4.212625in}{6.448822in}}%
\pgfpathlineto{\pgfqpoint{4.419292in}{4.806487in}}%
\pgfpathlineto{\pgfqpoint{4.625958in}{6.431735in}}%
\pgfpathlineto{\pgfqpoint{4.832625in}{5.217133in}}%
\pgfpathlineto{\pgfqpoint{5.039292in}{6.156569in}}%
\pgfpathlineto{\pgfqpoint{5.245958in}{6.243314in}}%
\pgfpathlineto{\pgfqpoint{5.452625in}{4.947056in}}%
\pgfpathlineto{\pgfqpoint{5.659292in}{5.337054in}}%
\pgfpathlineto{\pgfqpoint{5.865958in}{4.853833in}}%
\pgfpathlineto{\pgfqpoint{6.072625in}{4.855638in}}%
\pgfpathlineto{\pgfqpoint{6.279292in}{4.867044in}}%
\pgfpathlineto{\pgfqpoint{6.485958in}{6.377234in}}%
\pgfpathlineto{\pgfqpoint{6.692625in}{5.126318in}}%
\pgfpathlineto{\pgfqpoint{6.899292in}{4.539732in}}%
\pgfpathlineto{\pgfqpoint{7.105958in}{4.734594in}}%
\pgfusepath{stroke}%
\end{pgfscope}%
\begin{pgfscope}%
\pgfpathrectangle{\pgfqpoint{0.905958in}{0.899712in}}{\pgfqpoint{6.200000in}{6.200000in}} %
\pgfusepath{clip}%
\pgfsetbuttcap%
\pgfsetroundjoin%
\definecolor{currentfill}{rgb}{0.298039,0.447059,0.690196}%
\pgfsetfillcolor{currentfill}%
\pgfsetlinewidth{0.000000pt}%
\definecolor{currentstroke}{rgb}{0.000000,0.000000,0.000000}%
\pgfsetstrokecolor{currentstroke}%
\pgfsetdash{}{0pt}%
\pgfsys@defobject{currentmarker}{\pgfqpoint{-0.038889in}{-0.038889in}}{\pgfqpoint{0.038889in}{0.038889in}}{%
\pgfpathmoveto{\pgfqpoint{0.000000in}{-0.038889in}}%
\pgfpathcurveto{\pgfqpoint{0.010313in}{-0.038889in}}{\pgfqpoint{0.020206in}{-0.034791in}}{\pgfqpoint{0.027499in}{-0.027499in}}%
\pgfpathcurveto{\pgfqpoint{0.034791in}{-0.020206in}}{\pgfqpoint{0.038889in}{-0.010313in}}{\pgfqpoint{0.038889in}{0.000000in}}%
\pgfpathcurveto{\pgfqpoint{0.038889in}{0.010313in}}{\pgfqpoint{0.034791in}{0.020206in}}{\pgfqpoint{0.027499in}{0.027499in}}%
\pgfpathcurveto{\pgfqpoint{0.020206in}{0.034791in}}{\pgfqpoint{0.010313in}{0.038889in}}{\pgfqpoint{0.000000in}{0.038889in}}%
\pgfpathcurveto{\pgfqpoint{-0.010313in}{0.038889in}}{\pgfqpoint{-0.020206in}{0.034791in}}{\pgfqpoint{-0.027499in}{0.027499in}}%
\pgfpathcurveto{\pgfqpoint{-0.034791in}{0.020206in}}{\pgfqpoint{-0.038889in}{0.010313in}}{\pgfqpoint{-0.038889in}{0.000000in}}%
\pgfpathcurveto{\pgfqpoint{-0.038889in}{-0.010313in}}{\pgfqpoint{-0.034791in}{-0.020206in}}{\pgfqpoint{-0.027499in}{-0.027499in}}%
\pgfpathcurveto{\pgfqpoint{-0.020206in}{-0.034791in}}{\pgfqpoint{-0.010313in}{-0.038889in}}{\pgfqpoint{0.000000in}{-0.038889in}}%
\pgfpathclose%
\pgfusepath{fill}%
}%
\begin{pgfscope}%
\pgfsys@transformshift{1.732625in}{0.899712in}%
\pgfsys@useobject{currentmarker}{}%
\end{pgfscope}%
\begin{pgfscope}%
\pgfsys@transformshift{1.939292in}{2.407842in}%
\pgfsys@useobject{currentmarker}{}%
\end{pgfscope}%
\begin{pgfscope}%
\pgfsys@transformshift{2.145958in}{3.330376in}%
\pgfsys@useobject{currentmarker}{}%
\end{pgfscope}%
\begin{pgfscope}%
\pgfsys@transformshift{2.352625in}{3.755717in}%
\pgfsys@useobject{currentmarker}{}%
\end{pgfscope}%
\begin{pgfscope}%
\pgfsys@transformshift{2.559292in}{4.930792in}%
\pgfsys@useobject{currentmarker}{}%
\end{pgfscope}%
\begin{pgfscope}%
\pgfsys@transformshift{2.765958in}{5.389729in}%
\pgfsys@useobject{currentmarker}{}%
\end{pgfscope}%
\begin{pgfscope}%
\pgfsys@transformshift{2.972625in}{4.796105in}%
\pgfsys@useobject{currentmarker}{}%
\end{pgfscope}%
\begin{pgfscope}%
\pgfsys@transformshift{3.179292in}{5.331711in}%
\pgfsys@useobject{currentmarker}{}%
\end{pgfscope}%
\begin{pgfscope}%
\pgfsys@transformshift{3.385958in}{5.473362in}%
\pgfsys@useobject{currentmarker}{}%
\end{pgfscope}%
\begin{pgfscope}%
\pgfsys@transformshift{3.592625in}{4.283941in}%
\pgfsys@useobject{currentmarker}{}%
\end{pgfscope}%
\begin{pgfscope}%
\pgfsys@transformshift{3.799292in}{4.361007in}%
\pgfsys@useobject{currentmarker}{}%
\end{pgfscope}%
\begin{pgfscope}%
\pgfsys@transformshift{4.005958in}{6.260811in}%
\pgfsys@useobject{currentmarker}{}%
\end{pgfscope}%
\begin{pgfscope}%
\pgfsys@transformshift{4.212625in}{6.448822in}%
\pgfsys@useobject{currentmarker}{}%
\end{pgfscope}%
\begin{pgfscope}%
\pgfsys@transformshift{4.419292in}{4.806487in}%
\pgfsys@useobject{currentmarker}{}%
\end{pgfscope}%
\begin{pgfscope}%
\pgfsys@transformshift{4.625958in}{6.431735in}%
\pgfsys@useobject{currentmarker}{}%
\end{pgfscope}%
\begin{pgfscope}%
\pgfsys@transformshift{4.832625in}{5.217133in}%
\pgfsys@useobject{currentmarker}{}%
\end{pgfscope}%
\begin{pgfscope}%
\pgfsys@transformshift{5.039292in}{6.156569in}%
\pgfsys@useobject{currentmarker}{}%
\end{pgfscope}%
\begin{pgfscope}%
\pgfsys@transformshift{5.245958in}{6.243314in}%
\pgfsys@useobject{currentmarker}{}%
\end{pgfscope}%
\begin{pgfscope}%
\pgfsys@transformshift{5.452625in}{4.947056in}%
\pgfsys@useobject{currentmarker}{}%
\end{pgfscope}%
\begin{pgfscope}%
\pgfsys@transformshift{5.659292in}{5.337054in}%
\pgfsys@useobject{currentmarker}{}%
\end{pgfscope}%
\begin{pgfscope}%
\pgfsys@transformshift{5.865958in}{4.853833in}%
\pgfsys@useobject{currentmarker}{}%
\end{pgfscope}%
\begin{pgfscope}%
\pgfsys@transformshift{6.072625in}{4.855638in}%
\pgfsys@useobject{currentmarker}{}%
\end{pgfscope}%
\begin{pgfscope}%
\pgfsys@transformshift{6.279292in}{4.867044in}%
\pgfsys@useobject{currentmarker}{}%
\end{pgfscope}%
\begin{pgfscope}%
\pgfsys@transformshift{6.485958in}{6.377234in}%
\pgfsys@useobject{currentmarker}{}%
\end{pgfscope}%
\begin{pgfscope}%
\pgfsys@transformshift{6.692625in}{5.126318in}%
\pgfsys@useobject{currentmarker}{}%
\end{pgfscope}%
\begin{pgfscope}%
\pgfsys@transformshift{6.899292in}{4.539732in}%
\pgfsys@useobject{currentmarker}{}%
\end{pgfscope}%
\begin{pgfscope}%
\pgfsys@transformshift{7.105958in}{4.734594in}%
\pgfsys@useobject{currentmarker}{}%
\end{pgfscope}%
\end{pgfscope}%
\begin{pgfscope}%
\pgfpathrectangle{\pgfqpoint{0.905958in}{0.899712in}}{\pgfqpoint{6.200000in}{6.200000in}} %
\pgfusepath{clip}%
\pgfsetroundcap%
\pgfsetroundjoin%
\pgfsetlinewidth{2.007500pt}%
\definecolor{currentstroke}{rgb}{0.333333,0.658824,0.407843}%
\pgfsetstrokecolor{currentstroke}%
\pgfsetdash{}{0pt}%
\pgfpathmoveto{\pgfqpoint{1.732625in}{0.899712in}}%
\pgfpathlineto{\pgfqpoint{1.939292in}{1.953215in}}%
\pgfpathlineto{\pgfqpoint{2.145958in}{2.767002in}}%
\pgfpathlineto{\pgfqpoint{2.352625in}{3.221156in}}%
\pgfpathlineto{\pgfqpoint{2.559292in}{3.642514in}}%
\pgfpathlineto{\pgfqpoint{2.765958in}{3.772569in}}%
\pgfpathlineto{\pgfqpoint{2.972625in}{3.747281in}}%
\pgfpathlineto{\pgfqpoint{3.179292in}{4.003567in}}%
\pgfpathlineto{\pgfqpoint{3.385958in}{3.931027in}}%
\pgfpathlineto{\pgfqpoint{3.592625in}{3.741379in}}%
\pgfpathlineto{\pgfqpoint{3.799292in}{3.943088in}}%
\pgfpathlineto{\pgfqpoint{4.005958in}{3.808436in}}%
\pgfusepath{stroke}%
\end{pgfscope}%
\begin{pgfscope}%
\pgfpathrectangle{\pgfqpoint{0.905958in}{0.899712in}}{\pgfqpoint{6.200000in}{6.200000in}} %
\pgfusepath{clip}%
\pgfsetbuttcap%
\pgfsetmiterjoin%
\definecolor{currentfill}{rgb}{0.333333,0.658824,0.407843}%
\pgfsetfillcolor{currentfill}%
\pgfsetlinewidth{0.000000pt}%
\definecolor{currentstroke}{rgb}{0.000000,0.000000,0.000000}%
\pgfsetstrokecolor{currentstroke}%
\pgfsetdash{}{0pt}%
\pgfsys@defobject{currentmarker}{\pgfqpoint{-0.038889in}{-0.038889in}}{\pgfqpoint{0.038889in}{0.038889in}}{%
\pgfpathmoveto{\pgfqpoint{-0.038889in}{-0.038889in}}%
\pgfpathlineto{\pgfqpoint{0.038889in}{-0.038889in}}%
\pgfpathlineto{\pgfqpoint{0.038889in}{0.038889in}}%
\pgfpathlineto{\pgfqpoint{-0.038889in}{0.038889in}}%
\pgfpathclose%
\pgfusepath{fill}%
}%
\begin{pgfscope}%
\pgfsys@transformshift{1.732625in}{0.899712in}%
\pgfsys@useobject{currentmarker}{}%
\end{pgfscope}%
\begin{pgfscope}%
\pgfsys@transformshift{1.939292in}{1.953215in}%
\pgfsys@useobject{currentmarker}{}%
\end{pgfscope}%
\begin{pgfscope}%
\pgfsys@transformshift{2.145958in}{2.767002in}%
\pgfsys@useobject{currentmarker}{}%
\end{pgfscope}%
\begin{pgfscope}%
\pgfsys@transformshift{2.352625in}{3.221156in}%
\pgfsys@useobject{currentmarker}{}%
\end{pgfscope}%
\begin{pgfscope}%
\pgfsys@transformshift{2.559292in}{3.642514in}%
\pgfsys@useobject{currentmarker}{}%
\end{pgfscope}%
\begin{pgfscope}%
\pgfsys@transformshift{2.765958in}{3.772569in}%
\pgfsys@useobject{currentmarker}{}%
\end{pgfscope}%
\begin{pgfscope}%
\pgfsys@transformshift{2.972625in}{3.747281in}%
\pgfsys@useobject{currentmarker}{}%
\end{pgfscope}%
\begin{pgfscope}%
\pgfsys@transformshift{3.179292in}{4.003567in}%
\pgfsys@useobject{currentmarker}{}%
\end{pgfscope}%
\begin{pgfscope}%
\pgfsys@transformshift{3.385958in}{3.931027in}%
\pgfsys@useobject{currentmarker}{}%
\end{pgfscope}%
\begin{pgfscope}%
\pgfsys@transformshift{3.592625in}{3.741379in}%
\pgfsys@useobject{currentmarker}{}%
\end{pgfscope}%
\begin{pgfscope}%
\pgfsys@transformshift{3.799292in}{3.943088in}%
\pgfsys@useobject{currentmarker}{}%
\end{pgfscope}%
\begin{pgfscope}%
\pgfsys@transformshift{4.005958in}{3.808436in}%
\pgfsys@useobject{currentmarker}{}%
\end{pgfscope}%
\end{pgfscope}%
\begin{pgfscope}%
\pgfsetrectcap%
\pgfsetmiterjoin%
\pgfsetlinewidth{1.254687pt}%
\definecolor{currentstroke}{rgb}{0.150000,0.150000,0.150000}%
\pgfsetstrokecolor{currentstroke}%
\pgfsetdash{}{0pt}%
\pgfpathmoveto{\pgfqpoint{0.905958in}{0.899712in}}%
\pgfpathlineto{\pgfqpoint{7.105958in}{0.899712in}}%
\pgfusepath{stroke}%
\end{pgfscope}%
\begin{pgfscope}%
\pgfsetrectcap%
\pgfsetmiterjoin%
\pgfsetlinewidth{1.254687pt}%
\definecolor{currentstroke}{rgb}{0.150000,0.150000,0.150000}%
\pgfsetstrokecolor{currentstroke}%
\pgfsetdash{}{0pt}%
\pgfpathmoveto{\pgfqpoint{0.905958in}{0.899712in}}%
\pgfpathlineto{\pgfqpoint{0.905958in}{7.099712in}}%
\pgfusepath{stroke}%
\end{pgfscope}%
\begin{pgfscope}%
\pgfsetbuttcap%
\pgfsetmiterjoin%
\definecolor{currentfill}{rgb}{1.000000,1.000000,1.000000}%
\pgfsetfillcolor{currentfill}%
\pgfsetlinewidth{0.240900pt}%
\definecolor{currentstroke}{rgb}{0.150000,0.150000,0.150000}%
\pgfsetstrokecolor{currentstroke}%
\pgfsetdash{}{0pt}%
\pgfpathmoveto{\pgfqpoint{3.535678in}{1.066379in}}%
\pgfpathlineto{\pgfqpoint{6.939292in}{1.066379in}}%
\pgfpathlineto{\pgfqpoint{6.939292in}{2.121131in}}%
\pgfpathlineto{\pgfqpoint{3.535678in}{2.121131in}}%
\pgfpathclose%
\pgfusepath{stroke,fill}%
\end{pgfscope}%
\begin{pgfscope}%
\pgfsetroundcap%
\pgfsetroundjoin%
\pgfsetlinewidth{2.007500pt}%
\definecolor{currentstroke}{rgb}{0.298039,0.447059,0.690196}%
\pgfsetstrokecolor{currentstroke}%
\pgfsetdash{}{0pt}%
\pgfpathmoveto{\pgfqpoint{3.669012in}{1.862928in}}%
\pgfpathlineto{\pgfqpoint{4.335678in}{1.862928in}}%
\pgfusepath{stroke}%
\end{pgfscope}%
\begin{pgfscope}%
\pgfsetbuttcap%
\pgfsetroundjoin%
\definecolor{currentfill}{rgb}{0.298039,0.447059,0.690196}%
\pgfsetfillcolor{currentfill}%
\pgfsetlinewidth{0.000000pt}%
\definecolor{currentstroke}{rgb}{0.000000,0.000000,0.000000}%
\pgfsetstrokecolor{currentstroke}%
\pgfsetdash{}{0pt}%
\pgfsys@defobject{currentmarker}{\pgfqpoint{-0.038889in}{-0.038889in}}{\pgfqpoint{0.038889in}{0.038889in}}{%
\pgfpathmoveto{\pgfqpoint{0.000000in}{-0.038889in}}%
\pgfpathcurveto{\pgfqpoint{0.010313in}{-0.038889in}}{\pgfqpoint{0.020206in}{-0.034791in}}{\pgfqpoint{0.027499in}{-0.027499in}}%
\pgfpathcurveto{\pgfqpoint{0.034791in}{-0.020206in}}{\pgfqpoint{0.038889in}{-0.010313in}}{\pgfqpoint{0.038889in}{0.000000in}}%
\pgfpathcurveto{\pgfqpoint{0.038889in}{0.010313in}}{\pgfqpoint{0.034791in}{0.020206in}}{\pgfqpoint{0.027499in}{0.027499in}}%
\pgfpathcurveto{\pgfqpoint{0.020206in}{0.034791in}}{\pgfqpoint{0.010313in}{0.038889in}}{\pgfqpoint{0.000000in}{0.038889in}}%
\pgfpathcurveto{\pgfqpoint{-0.010313in}{0.038889in}}{\pgfqpoint{-0.020206in}{0.034791in}}{\pgfqpoint{-0.027499in}{0.027499in}}%
\pgfpathcurveto{\pgfqpoint{-0.034791in}{0.020206in}}{\pgfqpoint{-0.038889in}{0.010313in}}{\pgfqpoint{-0.038889in}{0.000000in}}%
\pgfpathcurveto{\pgfqpoint{-0.038889in}{-0.010313in}}{\pgfqpoint{-0.034791in}{-0.020206in}}{\pgfqpoint{-0.027499in}{-0.027499in}}%
\pgfpathcurveto{\pgfqpoint{-0.020206in}{-0.034791in}}{\pgfqpoint{-0.010313in}{-0.038889in}}{\pgfqpoint{0.000000in}{-0.038889in}}%
\pgfpathclose%
\pgfusepath{fill}%
}%
\begin{pgfscope}%
\pgfsys@transformshift{4.002345in}{1.862928in}%
\pgfsys@useobject{currentmarker}{}%
\end{pgfscope}%
\end{pgfscope}%
\begin{pgfscope}%
\definecolor{textcolor}{rgb}{0.150000,0.150000,0.150000}%
\pgfsetstrokecolor{textcolor}%
\pgfsetfillcolor{textcolor}%
\pgftext[x=4.602345in,y=1.746261in,left,base]{\color{textcolor}\sffamily\fontsize{24.000000}{28.800000}\selectfont FS-QTT-solver}%
\end{pgfscope}%
\begin{pgfscope}%
\pgfsetroundcap%
\pgfsetroundjoin%
\pgfsetlinewidth{2.007500pt}%
\definecolor{currentstroke}{rgb}{0.333333,0.658824,0.407843}%
\pgfsetstrokecolor{currentstroke}%
\pgfsetdash{}{0pt}%
\pgfpathmoveto{\pgfqpoint{3.669012in}{1.385552in}}%
\pgfpathlineto{\pgfqpoint{4.335678in}{1.385552in}}%
\pgfusepath{stroke}%
\end{pgfscope}%
\begin{pgfscope}%
\pgfsetbuttcap%
\pgfsetmiterjoin%
\definecolor{currentfill}{rgb}{0.333333,0.658824,0.407843}%
\pgfsetfillcolor{currentfill}%
\pgfsetlinewidth{0.000000pt}%
\definecolor{currentstroke}{rgb}{0.000000,0.000000,0.000000}%
\pgfsetstrokecolor{currentstroke}%
\pgfsetdash{}{0pt}%
\pgfsys@defobject{currentmarker}{\pgfqpoint{-0.038889in}{-0.038889in}}{\pgfqpoint{0.038889in}{0.038889in}}{%
\pgfpathmoveto{\pgfqpoint{-0.038889in}{-0.038889in}}%
\pgfpathlineto{\pgfqpoint{0.038889in}{-0.038889in}}%
\pgfpathlineto{\pgfqpoint{0.038889in}{0.038889in}}%
\pgfpathlineto{\pgfqpoint{-0.038889in}{0.038889in}}%
\pgfpathclose%
\pgfusepath{fill}%
}%
\begin{pgfscope}%
\pgfsys@transformshift{4.002345in}{1.385552in}%
\pgfsys@useobject{currentmarker}{}%
\end{pgfscope}%
\end{pgfscope}%
\begin{pgfscope}%
\definecolor{textcolor}{rgb}{0.150000,0.150000,0.150000}%
\pgfsetstrokecolor{textcolor}%
\pgfsetfillcolor{textcolor}%
\pgftext[x=4.602345in,y=1.268885in,left,base]{\color{textcolor}\sffamily\fontsize{24.000000}{28.800000}\selectfont FD-QTT-solver}%
\end{pgfscope}%
\end{pgfpicture}%
\makeatother%
\endgroup%

%% file: res_ps4_all_u_calc_erank.pgf
\begingroup%
\makeatletter%
\begin{pgfpicture}%
\pgfpathrectangle{\pgfpointorigin}{\pgfqpoint{7.364436in}{7.355067in}}%
\pgfusepath{use as bounding box, clip}%
\begin{pgfscope}%
\pgfsetbuttcap%
\pgfsetmiterjoin%
\definecolor{currentfill}{rgb}{1.000000,1.000000,1.000000}%
\pgfsetfillcolor{currentfill}%
\pgfsetlinewidth{0.000000pt}%
\definecolor{currentstroke}{rgb}{1.000000,1.000000,1.000000}%
\pgfsetstrokecolor{currentstroke}%
\pgfsetdash{}{0pt}%
\pgfpathmoveto{\pgfqpoint{0.000000in}{0.000000in}}%
\pgfpathlineto{\pgfqpoint{7.364436in}{0.000000in}}%
\pgfpathlineto{\pgfqpoint{7.364436in}{7.355067in}}%
\pgfpathlineto{\pgfqpoint{0.000000in}{7.355067in}}%
\pgfpathclose%
\pgfusepath{fill}%
\end{pgfscope}%
\begin{pgfscope}%
\pgfsetbuttcap%
\pgfsetmiterjoin%
\definecolor{currentfill}{rgb}{1.000000,1.000000,1.000000}%
\pgfsetfillcolor{currentfill}%
\pgfsetlinewidth{0.000000pt}%
\definecolor{currentstroke}{rgb}{0.000000,0.000000,0.000000}%
\pgfsetstrokecolor{currentstroke}%
\pgfsetstrokeopacity{0.000000}%
\pgfsetdash{}{0pt}%
\pgfpathmoveto{\pgfqpoint{0.905958in}{0.899712in}}%
\pgfpathlineto{\pgfqpoint{7.105958in}{0.899712in}}%
\pgfpathlineto{\pgfqpoint{7.105958in}{7.099712in}}%
\pgfpathlineto{\pgfqpoint{0.905958in}{7.099712in}}%
\pgfpathclose%
\pgfusepath{fill}%
\end{pgfscope}%
\begin{pgfscope}%
\pgfpathrectangle{\pgfqpoint{0.905958in}{0.899712in}}{\pgfqpoint{6.200000in}{6.200000in}} %
\pgfusepath{clip}%
\pgfsetbuttcap%
\pgfsetroundjoin%
\pgfsetlinewidth{0.803000pt}%
\definecolor{currentstroke}{rgb}{0.800000,0.800000,0.800000}%
\pgfsetstrokecolor{currentstroke}%
\pgfsetdash{{1.000000pt}{3.000000pt}}{0.000000pt}%
\pgfpathmoveto{\pgfqpoint{0.905958in}{0.899712in}}%
\pgfpathlineto{\pgfqpoint{0.905958in}{7.099712in}}%
\pgfusepath{stroke}%
\end{pgfscope}%
\begin{pgfscope}%
\pgfsetbuttcap%
\pgfsetroundjoin%
\definecolor{currentfill}{rgb}{0.150000,0.150000,0.150000}%
\pgfsetfillcolor{currentfill}%
\pgfsetlinewidth{0.803000pt}%
\definecolor{currentstroke}{rgb}{0.150000,0.150000,0.150000}%
\pgfsetstrokecolor{currentstroke}%
\pgfsetdash{}{0pt}%
\pgfsys@defobject{currentmarker}{\pgfqpoint{0.000000in}{0.000000in}}{\pgfqpoint{0.000000in}{0.000000in}}{%
\pgfpathmoveto{\pgfqpoint{0.000000in}{0.000000in}}%
\pgfpathlineto{\pgfqpoint{0.000000in}{0.000000in}}%
\pgfusepath{stroke,fill}%
}%
\begin{pgfscope}%
\pgfsys@transformshift{0.905958in}{0.899712in}%
\pgfsys@useobject{currentmarker}{}%
\end{pgfscope}%
\end{pgfscope}%
\begin{pgfscope}%
\definecolor{textcolor}{rgb}{0.150000,0.150000,0.150000}%
\pgfsetstrokecolor{textcolor}%
\pgfsetfillcolor{textcolor}%
\pgftext[x=0.905958in,y=0.821934in,,top]{\color{textcolor}\sffamily\fontsize{24.000000}{28.800000}\selectfont \(\displaystyle 4\)}%
\end{pgfscope}%
\begin{pgfscope}%
\pgfpathrectangle{\pgfqpoint{0.905958in}{0.899712in}}{\pgfqpoint{6.200000in}{6.200000in}} %
\pgfusepath{clip}%
\pgfsetbuttcap%
\pgfsetroundjoin%
\pgfsetlinewidth{0.803000pt}%
\definecolor{currentstroke}{rgb}{0.800000,0.800000,0.800000}%
\pgfsetstrokecolor{currentstroke}%
\pgfsetdash{{1.000000pt}{3.000000pt}}{0.000000pt}%
\pgfpathmoveto{\pgfqpoint{1.680958in}{0.899712in}}%
\pgfpathlineto{\pgfqpoint{1.680958in}{7.099712in}}%
\pgfusepath{stroke}%
\end{pgfscope}%
\begin{pgfscope}%
\pgfsetbuttcap%
\pgfsetroundjoin%
\definecolor{currentfill}{rgb}{0.150000,0.150000,0.150000}%
\pgfsetfillcolor{currentfill}%
\pgfsetlinewidth{0.803000pt}%
\definecolor{currentstroke}{rgb}{0.150000,0.150000,0.150000}%
\pgfsetstrokecolor{currentstroke}%
\pgfsetdash{}{0pt}%
\pgfsys@defobject{currentmarker}{\pgfqpoint{0.000000in}{0.000000in}}{\pgfqpoint{0.000000in}{0.000000in}}{%
\pgfpathmoveto{\pgfqpoint{0.000000in}{0.000000in}}%
\pgfpathlineto{\pgfqpoint{0.000000in}{0.000000in}}%
\pgfusepath{stroke,fill}%
}%
\begin{pgfscope}%
\pgfsys@transformshift{1.680958in}{0.899712in}%
\pgfsys@useobject{currentmarker}{}%
\end{pgfscope}%
\end{pgfscope}%
\begin{pgfscope}%
\definecolor{textcolor}{rgb}{0.150000,0.150000,0.150000}%
\pgfsetstrokecolor{textcolor}%
\pgfsetfillcolor{textcolor}%
\pgftext[x=1.680958in,y=0.821934in,,top]{\color{textcolor}\sffamily\fontsize{24.000000}{28.800000}\selectfont \(\displaystyle 6\)}%
\end{pgfscope}%
\begin{pgfscope}%
\pgfpathrectangle{\pgfqpoint{0.905958in}{0.899712in}}{\pgfqpoint{6.200000in}{6.200000in}} %
\pgfusepath{clip}%
\pgfsetbuttcap%
\pgfsetroundjoin%
\pgfsetlinewidth{0.803000pt}%
\definecolor{currentstroke}{rgb}{0.800000,0.800000,0.800000}%
\pgfsetstrokecolor{currentstroke}%
\pgfsetdash{{1.000000pt}{3.000000pt}}{0.000000pt}%
\pgfpathmoveto{\pgfqpoint{2.455958in}{0.899712in}}%
\pgfpathlineto{\pgfqpoint{2.455958in}{7.099712in}}%
\pgfusepath{stroke}%
\end{pgfscope}%
\begin{pgfscope}%
\pgfsetbuttcap%
\pgfsetroundjoin%
\definecolor{currentfill}{rgb}{0.150000,0.150000,0.150000}%
\pgfsetfillcolor{currentfill}%
\pgfsetlinewidth{0.803000pt}%
\definecolor{currentstroke}{rgb}{0.150000,0.150000,0.150000}%
\pgfsetstrokecolor{currentstroke}%
\pgfsetdash{}{0pt}%
\pgfsys@defobject{currentmarker}{\pgfqpoint{0.000000in}{0.000000in}}{\pgfqpoint{0.000000in}{0.000000in}}{%
\pgfpathmoveto{\pgfqpoint{0.000000in}{0.000000in}}%
\pgfpathlineto{\pgfqpoint{0.000000in}{0.000000in}}%
\pgfusepath{stroke,fill}%
}%
\begin{pgfscope}%
\pgfsys@transformshift{2.455958in}{0.899712in}%
\pgfsys@useobject{currentmarker}{}%
\end{pgfscope}%
\end{pgfscope}%
\begin{pgfscope}%
\definecolor{textcolor}{rgb}{0.150000,0.150000,0.150000}%
\pgfsetstrokecolor{textcolor}%
\pgfsetfillcolor{textcolor}%
\pgftext[x=2.455958in,y=0.821934in,,top]{\color{textcolor}\sffamily\fontsize{24.000000}{28.800000}\selectfont \(\displaystyle 8\)}%
\end{pgfscope}%
\begin{pgfscope}%
\pgfpathrectangle{\pgfqpoint{0.905958in}{0.899712in}}{\pgfqpoint{6.200000in}{6.200000in}} %
\pgfusepath{clip}%
\pgfsetbuttcap%
\pgfsetroundjoin%
\pgfsetlinewidth{0.803000pt}%
\definecolor{currentstroke}{rgb}{0.800000,0.800000,0.800000}%
\pgfsetstrokecolor{currentstroke}%
\pgfsetdash{{1.000000pt}{3.000000pt}}{0.000000pt}%
\pgfpathmoveto{\pgfqpoint{3.230958in}{0.899712in}}%
\pgfpathlineto{\pgfqpoint{3.230958in}{7.099712in}}%
\pgfusepath{stroke}%
\end{pgfscope}%
\begin{pgfscope}%
\pgfsetbuttcap%
\pgfsetroundjoin%
\definecolor{currentfill}{rgb}{0.150000,0.150000,0.150000}%
\pgfsetfillcolor{currentfill}%
\pgfsetlinewidth{0.803000pt}%
\definecolor{currentstroke}{rgb}{0.150000,0.150000,0.150000}%
\pgfsetstrokecolor{currentstroke}%
\pgfsetdash{}{0pt}%
\pgfsys@defobject{currentmarker}{\pgfqpoint{0.000000in}{0.000000in}}{\pgfqpoint{0.000000in}{0.000000in}}{%
\pgfpathmoveto{\pgfqpoint{0.000000in}{0.000000in}}%
\pgfpathlineto{\pgfqpoint{0.000000in}{0.000000in}}%
\pgfusepath{stroke,fill}%
}%
\begin{pgfscope}%
\pgfsys@transformshift{3.230958in}{0.899712in}%
\pgfsys@useobject{currentmarker}{}%
\end{pgfscope}%
\end{pgfscope}%
\begin{pgfscope}%
\definecolor{textcolor}{rgb}{0.150000,0.150000,0.150000}%
\pgfsetstrokecolor{textcolor}%
\pgfsetfillcolor{textcolor}%
\pgftext[x=3.230958in,y=0.821934in,,top]{\color{textcolor}\sffamily\fontsize{24.000000}{28.800000}\selectfont \(\displaystyle 10\)}%
\end{pgfscope}%
\begin{pgfscope}%
\pgfpathrectangle{\pgfqpoint{0.905958in}{0.899712in}}{\pgfqpoint{6.200000in}{6.200000in}} %
\pgfusepath{clip}%
\pgfsetbuttcap%
\pgfsetroundjoin%
\pgfsetlinewidth{0.803000pt}%
\definecolor{currentstroke}{rgb}{0.800000,0.800000,0.800000}%
\pgfsetstrokecolor{currentstroke}%
\pgfsetdash{{1.000000pt}{3.000000pt}}{0.000000pt}%
\pgfpathmoveto{\pgfqpoint{4.005958in}{0.899712in}}%
\pgfpathlineto{\pgfqpoint{4.005958in}{7.099712in}}%
\pgfusepath{stroke}%
\end{pgfscope}%
\begin{pgfscope}%
\pgfsetbuttcap%
\pgfsetroundjoin%
\definecolor{currentfill}{rgb}{0.150000,0.150000,0.150000}%
\pgfsetfillcolor{currentfill}%
\pgfsetlinewidth{0.803000pt}%
\definecolor{currentstroke}{rgb}{0.150000,0.150000,0.150000}%
\pgfsetstrokecolor{currentstroke}%
\pgfsetdash{}{0pt}%
\pgfsys@defobject{currentmarker}{\pgfqpoint{0.000000in}{0.000000in}}{\pgfqpoint{0.000000in}{0.000000in}}{%
\pgfpathmoveto{\pgfqpoint{0.000000in}{0.000000in}}%
\pgfpathlineto{\pgfqpoint{0.000000in}{0.000000in}}%
\pgfusepath{stroke,fill}%
}%
\begin{pgfscope}%
\pgfsys@transformshift{4.005958in}{0.899712in}%
\pgfsys@useobject{currentmarker}{}%
\end{pgfscope}%
\end{pgfscope}%
\begin{pgfscope}%
\definecolor{textcolor}{rgb}{0.150000,0.150000,0.150000}%
\pgfsetstrokecolor{textcolor}%
\pgfsetfillcolor{textcolor}%
\pgftext[x=4.005958in,y=0.821934in,,top]{\color{textcolor}\sffamily\fontsize{24.000000}{28.800000}\selectfont \(\displaystyle 12\)}%
\end{pgfscope}%
\begin{pgfscope}%
\pgfpathrectangle{\pgfqpoint{0.905958in}{0.899712in}}{\pgfqpoint{6.200000in}{6.200000in}} %
\pgfusepath{clip}%
\pgfsetbuttcap%
\pgfsetroundjoin%
\pgfsetlinewidth{0.803000pt}%
\definecolor{currentstroke}{rgb}{0.800000,0.800000,0.800000}%
\pgfsetstrokecolor{currentstroke}%
\pgfsetdash{{1.000000pt}{3.000000pt}}{0.000000pt}%
\pgfpathmoveto{\pgfqpoint{4.780958in}{0.899712in}}%
\pgfpathlineto{\pgfqpoint{4.780958in}{7.099712in}}%
\pgfusepath{stroke}%
\end{pgfscope}%
\begin{pgfscope}%
\pgfsetbuttcap%
\pgfsetroundjoin%
\definecolor{currentfill}{rgb}{0.150000,0.150000,0.150000}%
\pgfsetfillcolor{currentfill}%
\pgfsetlinewidth{0.803000pt}%
\definecolor{currentstroke}{rgb}{0.150000,0.150000,0.150000}%
\pgfsetstrokecolor{currentstroke}%
\pgfsetdash{}{0pt}%
\pgfsys@defobject{currentmarker}{\pgfqpoint{0.000000in}{0.000000in}}{\pgfqpoint{0.000000in}{0.000000in}}{%
\pgfpathmoveto{\pgfqpoint{0.000000in}{0.000000in}}%
\pgfpathlineto{\pgfqpoint{0.000000in}{0.000000in}}%
\pgfusepath{stroke,fill}%
}%
\begin{pgfscope}%
\pgfsys@transformshift{4.780958in}{0.899712in}%
\pgfsys@useobject{currentmarker}{}%
\end{pgfscope}%
\end{pgfscope}%
\begin{pgfscope}%
\definecolor{textcolor}{rgb}{0.150000,0.150000,0.150000}%
\pgfsetstrokecolor{textcolor}%
\pgfsetfillcolor{textcolor}%
\pgftext[x=4.780958in,y=0.821934in,,top]{\color{textcolor}\sffamily\fontsize{24.000000}{28.800000}\selectfont \(\displaystyle 14\)}%
\end{pgfscope}%
\begin{pgfscope}%
\pgfpathrectangle{\pgfqpoint{0.905958in}{0.899712in}}{\pgfqpoint{6.200000in}{6.200000in}} %
\pgfusepath{clip}%
\pgfsetbuttcap%
\pgfsetroundjoin%
\pgfsetlinewidth{0.803000pt}%
\definecolor{currentstroke}{rgb}{0.800000,0.800000,0.800000}%
\pgfsetstrokecolor{currentstroke}%
\pgfsetdash{{1.000000pt}{3.000000pt}}{0.000000pt}%
\pgfpathmoveto{\pgfqpoint{5.555958in}{0.899712in}}%
\pgfpathlineto{\pgfqpoint{5.555958in}{7.099712in}}%
\pgfusepath{stroke}%
\end{pgfscope}%
\begin{pgfscope}%
\pgfsetbuttcap%
\pgfsetroundjoin%
\definecolor{currentfill}{rgb}{0.150000,0.150000,0.150000}%
\pgfsetfillcolor{currentfill}%
\pgfsetlinewidth{0.803000pt}%
\definecolor{currentstroke}{rgb}{0.150000,0.150000,0.150000}%
\pgfsetstrokecolor{currentstroke}%
\pgfsetdash{}{0pt}%
\pgfsys@defobject{currentmarker}{\pgfqpoint{0.000000in}{0.000000in}}{\pgfqpoint{0.000000in}{0.000000in}}{%
\pgfpathmoveto{\pgfqpoint{0.000000in}{0.000000in}}%
\pgfpathlineto{\pgfqpoint{0.000000in}{0.000000in}}%
\pgfusepath{stroke,fill}%
}%
\begin{pgfscope}%
\pgfsys@transformshift{5.555958in}{0.899712in}%
\pgfsys@useobject{currentmarker}{}%
\end{pgfscope}%
\end{pgfscope}%
\begin{pgfscope}%
\definecolor{textcolor}{rgb}{0.150000,0.150000,0.150000}%
\pgfsetstrokecolor{textcolor}%
\pgfsetfillcolor{textcolor}%
\pgftext[x=5.555958in,y=0.821934in,,top]{\color{textcolor}\sffamily\fontsize{24.000000}{28.800000}\selectfont \(\displaystyle 16\)}%
\end{pgfscope}%
\begin{pgfscope}%
\pgfpathrectangle{\pgfqpoint{0.905958in}{0.899712in}}{\pgfqpoint{6.200000in}{6.200000in}} %
\pgfusepath{clip}%
\pgfsetbuttcap%
\pgfsetroundjoin%
\pgfsetlinewidth{0.803000pt}%
\definecolor{currentstroke}{rgb}{0.800000,0.800000,0.800000}%
\pgfsetstrokecolor{currentstroke}%
\pgfsetdash{{1.000000pt}{3.000000pt}}{0.000000pt}%
\pgfpathmoveto{\pgfqpoint{6.330958in}{0.899712in}}%
\pgfpathlineto{\pgfqpoint{6.330958in}{7.099712in}}%
\pgfusepath{stroke}%
\end{pgfscope}%
\begin{pgfscope}%
\pgfsetbuttcap%
\pgfsetroundjoin%
\definecolor{currentfill}{rgb}{0.150000,0.150000,0.150000}%
\pgfsetfillcolor{currentfill}%
\pgfsetlinewidth{0.803000pt}%
\definecolor{currentstroke}{rgb}{0.150000,0.150000,0.150000}%
\pgfsetstrokecolor{currentstroke}%
\pgfsetdash{}{0pt}%
\pgfsys@defobject{currentmarker}{\pgfqpoint{0.000000in}{0.000000in}}{\pgfqpoint{0.000000in}{0.000000in}}{%
\pgfpathmoveto{\pgfqpoint{0.000000in}{0.000000in}}%
\pgfpathlineto{\pgfqpoint{0.000000in}{0.000000in}}%
\pgfusepath{stroke,fill}%
}%
\begin{pgfscope}%
\pgfsys@transformshift{6.330958in}{0.899712in}%
\pgfsys@useobject{currentmarker}{}%
\end{pgfscope}%
\end{pgfscope}%
\begin{pgfscope}%
\definecolor{textcolor}{rgb}{0.150000,0.150000,0.150000}%
\pgfsetstrokecolor{textcolor}%
\pgfsetfillcolor{textcolor}%
\pgftext[x=6.330958in,y=0.821934in,,top]{\color{textcolor}\sffamily\fontsize{24.000000}{28.800000}\selectfont \(\displaystyle 18\)}%
\end{pgfscope}%
\begin{pgfscope}%
\pgfpathrectangle{\pgfqpoint{0.905958in}{0.899712in}}{\pgfqpoint{6.200000in}{6.200000in}} %
\pgfusepath{clip}%
\pgfsetbuttcap%
\pgfsetroundjoin%
\pgfsetlinewidth{0.803000pt}%
\definecolor{currentstroke}{rgb}{0.800000,0.800000,0.800000}%
\pgfsetstrokecolor{currentstroke}%
\pgfsetdash{{1.000000pt}{3.000000pt}}{0.000000pt}%
\pgfpathmoveto{\pgfqpoint{7.105958in}{0.899712in}}%
\pgfpathlineto{\pgfqpoint{7.105958in}{7.099712in}}%
\pgfusepath{stroke}%
\end{pgfscope}%
\begin{pgfscope}%
\pgfsetbuttcap%
\pgfsetroundjoin%
\definecolor{currentfill}{rgb}{0.150000,0.150000,0.150000}%
\pgfsetfillcolor{currentfill}%
\pgfsetlinewidth{0.803000pt}%
\definecolor{currentstroke}{rgb}{0.150000,0.150000,0.150000}%
\pgfsetstrokecolor{currentstroke}%
\pgfsetdash{}{0pt}%
\pgfsys@defobject{currentmarker}{\pgfqpoint{0.000000in}{0.000000in}}{\pgfqpoint{0.000000in}{0.000000in}}{%
\pgfpathmoveto{\pgfqpoint{0.000000in}{0.000000in}}%
\pgfpathlineto{\pgfqpoint{0.000000in}{0.000000in}}%
\pgfusepath{stroke,fill}%
}%
\begin{pgfscope}%
\pgfsys@transformshift{7.105958in}{0.899712in}%
\pgfsys@useobject{currentmarker}{}%
\end{pgfscope}%
\end{pgfscope}%
\begin{pgfscope}%
\definecolor{textcolor}{rgb}{0.150000,0.150000,0.150000}%
\pgfsetstrokecolor{textcolor}%
\pgfsetfillcolor{textcolor}%
\pgftext[x=7.105958in,y=0.821934in,,top]{\color{textcolor}\sffamily\fontsize{24.000000}{28.800000}\selectfont \(\displaystyle 20\)}%
\end{pgfscope}%
\begin{pgfscope}%
\definecolor{textcolor}{rgb}{0.150000,0.150000,0.150000}%
\pgfsetstrokecolor{textcolor}%
\pgfsetfillcolor{textcolor}%
\pgftext[x=4.005958in,y=0.441780in,,top]{\color{textcolor}\sffamily\fontsize{26.400000}{31.680000}\selectfont d}%
\end{pgfscope}%
\begin{pgfscope}%
\pgfpathrectangle{\pgfqpoint{0.905958in}{0.899712in}}{\pgfqpoint{6.200000in}{6.200000in}} %
\pgfusepath{clip}%
\pgfsetbuttcap%
\pgfsetroundjoin%
\pgfsetlinewidth{0.803000pt}%
\definecolor{currentstroke}{rgb}{0.800000,0.800000,0.800000}%
\pgfsetstrokecolor{currentstroke}%
\pgfsetdash{{1.000000pt}{3.000000pt}}{0.000000pt}%
\pgfpathmoveto{\pgfqpoint{0.905958in}{0.899712in}}%
\pgfpathlineto{\pgfqpoint{7.105958in}{0.899712in}}%
\pgfusepath{stroke}%
\end{pgfscope}%
\begin{pgfscope}%
\pgfsetbuttcap%
\pgfsetroundjoin%
\definecolor{currentfill}{rgb}{0.150000,0.150000,0.150000}%
\pgfsetfillcolor{currentfill}%
\pgfsetlinewidth{0.803000pt}%
\definecolor{currentstroke}{rgb}{0.150000,0.150000,0.150000}%
\pgfsetstrokecolor{currentstroke}%
\pgfsetdash{}{0pt}%
\pgfsys@defobject{currentmarker}{\pgfqpoint{0.000000in}{0.000000in}}{\pgfqpoint{0.000000in}{0.000000in}}{%
\pgfpathmoveto{\pgfqpoint{0.000000in}{0.000000in}}%
\pgfpathlineto{\pgfqpoint{0.000000in}{0.000000in}}%
\pgfusepath{stroke,fill}%
}%
\begin{pgfscope}%
\pgfsys@transformshift{0.905958in}{0.899712in}%
\pgfsys@useobject{currentmarker}{}%
\end{pgfscope}%
\end{pgfscope}%
\begin{pgfscope}%
\definecolor{textcolor}{rgb}{0.150000,0.150000,0.150000}%
\pgfsetstrokecolor{textcolor}%
\pgfsetfillcolor{textcolor}%
\pgftext[x=0.828181in,y=0.899712in,right,]{\color{textcolor}\sffamily\fontsize{24.000000}{28.800000}\selectfont \(\displaystyle 0\)}%
\end{pgfscope}%
\begin{pgfscope}%
\pgfpathrectangle{\pgfqpoint{0.905958in}{0.899712in}}{\pgfqpoint{6.200000in}{6.200000in}} %
\pgfusepath{clip}%
\pgfsetbuttcap%
\pgfsetroundjoin%
\pgfsetlinewidth{0.803000pt}%
\definecolor{currentstroke}{rgb}{0.800000,0.800000,0.800000}%
\pgfsetstrokecolor{currentstroke}%
\pgfsetdash{{1.000000pt}{3.000000pt}}{0.000000pt}%
\pgfpathmoveto{\pgfqpoint{0.905958in}{1.933045in}}%
\pgfpathlineto{\pgfqpoint{7.105958in}{1.933045in}}%
\pgfusepath{stroke}%
\end{pgfscope}%
\begin{pgfscope}%
\pgfsetbuttcap%
\pgfsetroundjoin%
\definecolor{currentfill}{rgb}{0.150000,0.150000,0.150000}%
\pgfsetfillcolor{currentfill}%
\pgfsetlinewidth{0.803000pt}%
\definecolor{currentstroke}{rgb}{0.150000,0.150000,0.150000}%
\pgfsetstrokecolor{currentstroke}%
\pgfsetdash{}{0pt}%
\pgfsys@defobject{currentmarker}{\pgfqpoint{0.000000in}{0.000000in}}{\pgfqpoint{0.000000in}{0.000000in}}{%
\pgfpathmoveto{\pgfqpoint{0.000000in}{0.000000in}}%
\pgfpathlineto{\pgfqpoint{0.000000in}{0.000000in}}%
\pgfusepath{stroke,fill}%
}%
\begin{pgfscope}%
\pgfsys@transformshift{0.905958in}{1.933045in}%
\pgfsys@useobject{currentmarker}{}%
\end{pgfscope}%
\end{pgfscope}%
\begin{pgfscope}%
\definecolor{textcolor}{rgb}{0.150000,0.150000,0.150000}%
\pgfsetstrokecolor{textcolor}%
\pgfsetfillcolor{textcolor}%
\pgftext[x=0.828181in,y=1.933045in,right,]{\color{textcolor}\sffamily\fontsize{24.000000}{28.800000}\selectfont \(\displaystyle 10\)}%
\end{pgfscope}%
\begin{pgfscope}%
\pgfpathrectangle{\pgfqpoint{0.905958in}{0.899712in}}{\pgfqpoint{6.200000in}{6.200000in}} %
\pgfusepath{clip}%
\pgfsetbuttcap%
\pgfsetroundjoin%
\pgfsetlinewidth{0.803000pt}%
\definecolor{currentstroke}{rgb}{0.800000,0.800000,0.800000}%
\pgfsetstrokecolor{currentstroke}%
\pgfsetdash{{1.000000pt}{3.000000pt}}{0.000000pt}%
\pgfpathmoveto{\pgfqpoint{0.905958in}{2.966379in}}%
\pgfpathlineto{\pgfqpoint{7.105958in}{2.966379in}}%
\pgfusepath{stroke}%
\end{pgfscope}%
\begin{pgfscope}%
\pgfsetbuttcap%
\pgfsetroundjoin%
\definecolor{currentfill}{rgb}{0.150000,0.150000,0.150000}%
\pgfsetfillcolor{currentfill}%
\pgfsetlinewidth{0.803000pt}%
\definecolor{currentstroke}{rgb}{0.150000,0.150000,0.150000}%
\pgfsetstrokecolor{currentstroke}%
\pgfsetdash{}{0pt}%
\pgfsys@defobject{currentmarker}{\pgfqpoint{0.000000in}{0.000000in}}{\pgfqpoint{0.000000in}{0.000000in}}{%
\pgfpathmoveto{\pgfqpoint{0.000000in}{0.000000in}}%
\pgfpathlineto{\pgfqpoint{0.000000in}{0.000000in}}%
\pgfusepath{stroke,fill}%
}%
\begin{pgfscope}%
\pgfsys@transformshift{0.905958in}{2.966379in}%
\pgfsys@useobject{currentmarker}{}%
\end{pgfscope}%
\end{pgfscope}%
\begin{pgfscope}%
\definecolor{textcolor}{rgb}{0.150000,0.150000,0.150000}%
\pgfsetstrokecolor{textcolor}%
\pgfsetfillcolor{textcolor}%
\pgftext[x=0.828181in,y=2.966379in,right,]{\color{textcolor}\sffamily\fontsize{24.000000}{28.800000}\selectfont \(\displaystyle 20\)}%
\end{pgfscope}%
\begin{pgfscope}%
\pgfpathrectangle{\pgfqpoint{0.905958in}{0.899712in}}{\pgfqpoint{6.200000in}{6.200000in}} %
\pgfusepath{clip}%
\pgfsetbuttcap%
\pgfsetroundjoin%
\pgfsetlinewidth{0.803000pt}%
\definecolor{currentstroke}{rgb}{0.800000,0.800000,0.800000}%
\pgfsetstrokecolor{currentstroke}%
\pgfsetdash{{1.000000pt}{3.000000pt}}{0.000000pt}%
\pgfpathmoveto{\pgfqpoint{0.905958in}{3.999712in}}%
\pgfpathlineto{\pgfqpoint{7.105958in}{3.999712in}}%
\pgfusepath{stroke}%
\end{pgfscope}%
\begin{pgfscope}%
\pgfsetbuttcap%
\pgfsetroundjoin%
\definecolor{currentfill}{rgb}{0.150000,0.150000,0.150000}%
\pgfsetfillcolor{currentfill}%
\pgfsetlinewidth{0.803000pt}%
\definecolor{currentstroke}{rgb}{0.150000,0.150000,0.150000}%
\pgfsetstrokecolor{currentstroke}%
\pgfsetdash{}{0pt}%
\pgfsys@defobject{currentmarker}{\pgfqpoint{0.000000in}{0.000000in}}{\pgfqpoint{0.000000in}{0.000000in}}{%
\pgfpathmoveto{\pgfqpoint{0.000000in}{0.000000in}}%
\pgfpathlineto{\pgfqpoint{0.000000in}{0.000000in}}%
\pgfusepath{stroke,fill}%
}%
\begin{pgfscope}%
\pgfsys@transformshift{0.905958in}{3.999712in}%
\pgfsys@useobject{currentmarker}{}%
\end{pgfscope}%
\end{pgfscope}%
\begin{pgfscope}%
\definecolor{textcolor}{rgb}{0.150000,0.150000,0.150000}%
\pgfsetstrokecolor{textcolor}%
\pgfsetfillcolor{textcolor}%
\pgftext[x=0.828181in,y=3.999712in,right,]{\color{textcolor}\sffamily\fontsize{24.000000}{28.800000}\selectfont \(\displaystyle 30\)}%
\end{pgfscope}%
\begin{pgfscope}%
\pgfpathrectangle{\pgfqpoint{0.905958in}{0.899712in}}{\pgfqpoint{6.200000in}{6.200000in}} %
\pgfusepath{clip}%
\pgfsetbuttcap%
\pgfsetroundjoin%
\pgfsetlinewidth{0.803000pt}%
\definecolor{currentstroke}{rgb}{0.800000,0.800000,0.800000}%
\pgfsetstrokecolor{currentstroke}%
\pgfsetdash{{1.000000pt}{3.000000pt}}{0.000000pt}%
\pgfpathmoveto{\pgfqpoint{0.905958in}{5.033045in}}%
\pgfpathlineto{\pgfqpoint{7.105958in}{5.033045in}}%
\pgfusepath{stroke}%
\end{pgfscope}%
\begin{pgfscope}%
\pgfsetbuttcap%
\pgfsetroundjoin%
\definecolor{currentfill}{rgb}{0.150000,0.150000,0.150000}%
\pgfsetfillcolor{currentfill}%
\pgfsetlinewidth{0.803000pt}%
\definecolor{currentstroke}{rgb}{0.150000,0.150000,0.150000}%
\pgfsetstrokecolor{currentstroke}%
\pgfsetdash{}{0pt}%
\pgfsys@defobject{currentmarker}{\pgfqpoint{0.000000in}{0.000000in}}{\pgfqpoint{0.000000in}{0.000000in}}{%
\pgfpathmoveto{\pgfqpoint{0.000000in}{0.000000in}}%
\pgfpathlineto{\pgfqpoint{0.000000in}{0.000000in}}%
\pgfusepath{stroke,fill}%
}%
\begin{pgfscope}%
\pgfsys@transformshift{0.905958in}{5.033045in}%
\pgfsys@useobject{currentmarker}{}%
\end{pgfscope}%
\end{pgfscope}%
\begin{pgfscope}%
\definecolor{textcolor}{rgb}{0.150000,0.150000,0.150000}%
\pgfsetstrokecolor{textcolor}%
\pgfsetfillcolor{textcolor}%
\pgftext[x=0.828181in,y=5.033045in,right,]{\color{textcolor}\sffamily\fontsize{24.000000}{28.800000}\selectfont \(\displaystyle 40\)}%
\end{pgfscope}%
\begin{pgfscope}%
\pgfpathrectangle{\pgfqpoint{0.905958in}{0.899712in}}{\pgfqpoint{6.200000in}{6.200000in}} %
\pgfusepath{clip}%
\pgfsetbuttcap%
\pgfsetroundjoin%
\pgfsetlinewidth{0.803000pt}%
\definecolor{currentstroke}{rgb}{0.800000,0.800000,0.800000}%
\pgfsetstrokecolor{currentstroke}%
\pgfsetdash{{1.000000pt}{3.000000pt}}{0.000000pt}%
\pgfpathmoveto{\pgfqpoint{0.905958in}{6.066379in}}%
\pgfpathlineto{\pgfqpoint{7.105958in}{6.066379in}}%
\pgfusepath{stroke}%
\end{pgfscope}%
\begin{pgfscope}%
\pgfsetbuttcap%
\pgfsetroundjoin%
\definecolor{currentfill}{rgb}{0.150000,0.150000,0.150000}%
\pgfsetfillcolor{currentfill}%
\pgfsetlinewidth{0.803000pt}%
\definecolor{currentstroke}{rgb}{0.150000,0.150000,0.150000}%
\pgfsetstrokecolor{currentstroke}%
\pgfsetdash{}{0pt}%
\pgfsys@defobject{currentmarker}{\pgfqpoint{0.000000in}{0.000000in}}{\pgfqpoint{0.000000in}{0.000000in}}{%
\pgfpathmoveto{\pgfqpoint{0.000000in}{0.000000in}}%
\pgfpathlineto{\pgfqpoint{0.000000in}{0.000000in}}%
\pgfusepath{stroke,fill}%
}%
\begin{pgfscope}%
\pgfsys@transformshift{0.905958in}{6.066379in}%
\pgfsys@useobject{currentmarker}{}%
\end{pgfscope}%
\end{pgfscope}%
\begin{pgfscope}%
\definecolor{textcolor}{rgb}{0.150000,0.150000,0.150000}%
\pgfsetstrokecolor{textcolor}%
\pgfsetfillcolor{textcolor}%
\pgftext[x=0.828181in,y=6.066379in,right,]{\color{textcolor}\sffamily\fontsize{24.000000}{28.800000}\selectfont \(\displaystyle 50\)}%
\end{pgfscope}%
\begin{pgfscope}%
\pgfpathrectangle{\pgfqpoint{0.905958in}{0.899712in}}{\pgfqpoint{6.200000in}{6.200000in}} %
\pgfusepath{clip}%
\pgfsetbuttcap%
\pgfsetroundjoin%
\pgfsetlinewidth{0.803000pt}%
\definecolor{currentstroke}{rgb}{0.800000,0.800000,0.800000}%
\pgfsetstrokecolor{currentstroke}%
\pgfsetdash{{1.000000pt}{3.000000pt}}{0.000000pt}%
\pgfpathmoveto{\pgfqpoint{0.905958in}{7.099712in}}%
\pgfpathlineto{\pgfqpoint{7.105958in}{7.099712in}}%
\pgfusepath{stroke}%
\end{pgfscope}%
\begin{pgfscope}%
\pgfsetbuttcap%
\pgfsetroundjoin%
\definecolor{currentfill}{rgb}{0.150000,0.150000,0.150000}%
\pgfsetfillcolor{currentfill}%
\pgfsetlinewidth{0.803000pt}%
\definecolor{currentstroke}{rgb}{0.150000,0.150000,0.150000}%
\pgfsetstrokecolor{currentstroke}%
\pgfsetdash{}{0pt}%
\pgfsys@defobject{currentmarker}{\pgfqpoint{0.000000in}{0.000000in}}{\pgfqpoint{0.000000in}{0.000000in}}{%
\pgfpathmoveto{\pgfqpoint{0.000000in}{0.000000in}}%
\pgfpathlineto{\pgfqpoint{0.000000in}{0.000000in}}%
\pgfusepath{stroke,fill}%
}%
\begin{pgfscope}%
\pgfsys@transformshift{0.905958in}{7.099712in}%
\pgfsys@useobject{currentmarker}{}%
\end{pgfscope}%
\end{pgfscope}%
\begin{pgfscope}%
\definecolor{textcolor}{rgb}{0.150000,0.150000,0.150000}%
\pgfsetstrokecolor{textcolor}%
\pgfsetfillcolor{textcolor}%
\pgftext[x=0.828181in,y=7.099712in,right,]{\color{textcolor}\sffamily\fontsize{24.000000}{28.800000}\selectfont \(\displaystyle 60\)}%
\end{pgfscope}%
\begin{pgfscope}%
\definecolor{textcolor}{rgb}{0.150000,0.150000,0.150000}%
\pgfsetstrokecolor{textcolor}%
\pgfsetfillcolor{textcolor}%
\pgftext[x=0.441780in,y=3.999712in,,bottom,rotate=90.000000]{\color{textcolor}\sffamily\fontsize{26.400000}{31.680000}\selectfont Solution effective TT-rank}%
\end{pgfscope}%
\begin{pgfscope}%
\pgfpathrectangle{\pgfqpoint{0.905958in}{0.899712in}}{\pgfqpoint{6.200000in}{6.200000in}} %
\pgfusepath{clip}%
\pgfsetroundcap%
\pgfsetroundjoin%
\pgfsetlinewidth{2.007500pt}%
\definecolor{currentstroke}{rgb}{0.298039,0.447059,0.690196}%
\pgfsetstrokecolor{currentstroke}%
\pgfsetdash{}{0pt}%
\pgfpathmoveto{\pgfqpoint{0.905958in}{1.642027in}}%
\pgfpathlineto{\pgfqpoint{1.293458in}{2.049080in}}%
\pgfpathlineto{\pgfqpoint{1.680958in}{2.453130in}}%
\pgfpathlineto{\pgfqpoint{2.068458in}{2.883694in}}%
\pgfpathlineto{\pgfqpoint{2.455958in}{3.212499in}}%
\pgfpathlineto{\pgfqpoint{2.843458in}{3.090464in}}%
\pgfpathlineto{\pgfqpoint{3.230958in}{3.409719in}}%
\pgfpathlineto{\pgfqpoint{3.618458in}{3.306670in}}%
\pgfpathlineto{\pgfqpoint{4.005958in}{3.754181in}}%
\pgfpathlineto{\pgfqpoint{4.393458in}{4.188729in}}%
\pgfpathlineto{\pgfqpoint{4.780958in}{4.603908in}}%
\pgfpathlineto{\pgfqpoint{5.168458in}{4.852152in}}%
\pgfpathlineto{\pgfqpoint{5.555958in}{4.918091in}}%
\pgfpathlineto{\pgfqpoint{5.943458in}{5.007679in}}%
\pgfpathlineto{\pgfqpoint{6.330958in}{5.327267in}}%
\pgfpathlineto{\pgfqpoint{6.718458in}{5.626718in}}%
\pgfpathlineto{\pgfqpoint{7.105958in}{6.436941in}}%
\pgfusepath{stroke}%
\end{pgfscope}%
\begin{pgfscope}%
\pgfpathrectangle{\pgfqpoint{0.905958in}{0.899712in}}{\pgfqpoint{6.200000in}{6.200000in}} %
\pgfusepath{clip}%
\pgfsetbuttcap%
\pgfsetroundjoin%
\definecolor{currentfill}{rgb}{0.298039,0.447059,0.690196}%
\pgfsetfillcolor{currentfill}%
\pgfsetlinewidth{0.000000pt}%
\definecolor{currentstroke}{rgb}{0.000000,0.000000,0.000000}%
\pgfsetstrokecolor{currentstroke}%
\pgfsetdash{}{0pt}%
\pgfsys@defobject{currentmarker}{\pgfqpoint{-0.038889in}{-0.038889in}}{\pgfqpoint{0.038889in}{0.038889in}}{%
\pgfpathmoveto{\pgfqpoint{0.000000in}{-0.038889in}}%
\pgfpathcurveto{\pgfqpoint{0.010313in}{-0.038889in}}{\pgfqpoint{0.020206in}{-0.034791in}}{\pgfqpoint{0.027499in}{-0.027499in}}%
\pgfpathcurveto{\pgfqpoint{0.034791in}{-0.020206in}}{\pgfqpoint{0.038889in}{-0.010313in}}{\pgfqpoint{0.038889in}{0.000000in}}%
\pgfpathcurveto{\pgfqpoint{0.038889in}{0.010313in}}{\pgfqpoint{0.034791in}{0.020206in}}{\pgfqpoint{0.027499in}{0.027499in}}%
\pgfpathcurveto{\pgfqpoint{0.020206in}{0.034791in}}{\pgfqpoint{0.010313in}{0.038889in}}{\pgfqpoint{0.000000in}{0.038889in}}%
\pgfpathcurveto{\pgfqpoint{-0.010313in}{0.038889in}}{\pgfqpoint{-0.020206in}{0.034791in}}{\pgfqpoint{-0.027499in}{0.027499in}}%
\pgfpathcurveto{\pgfqpoint{-0.034791in}{0.020206in}}{\pgfqpoint{-0.038889in}{0.010313in}}{\pgfqpoint{-0.038889in}{0.000000in}}%
\pgfpathcurveto{\pgfqpoint{-0.038889in}{-0.010313in}}{\pgfqpoint{-0.034791in}{-0.020206in}}{\pgfqpoint{-0.027499in}{-0.027499in}}%
\pgfpathcurveto{\pgfqpoint{-0.020206in}{-0.034791in}}{\pgfqpoint{-0.010313in}{-0.038889in}}{\pgfqpoint{0.000000in}{-0.038889in}}%
\pgfpathclose%
\pgfusepath{fill}%
}%
\begin{pgfscope}%
\pgfsys@transformshift{0.905958in}{1.642027in}%
\pgfsys@useobject{currentmarker}{}%
\end{pgfscope}%
\begin{pgfscope}%
\pgfsys@transformshift{1.293458in}{2.049080in}%
\pgfsys@useobject{currentmarker}{}%
\end{pgfscope}%
\begin{pgfscope}%
\pgfsys@transformshift{1.680958in}{2.453130in}%
\pgfsys@useobject{currentmarker}{}%
\end{pgfscope}%
\begin{pgfscope}%
\pgfsys@transformshift{2.068458in}{2.883694in}%
\pgfsys@useobject{currentmarker}{}%
\end{pgfscope}%
\begin{pgfscope}%
\pgfsys@transformshift{2.455958in}{3.212499in}%
\pgfsys@useobject{currentmarker}{}%
\end{pgfscope}%
\begin{pgfscope}%
\pgfsys@transformshift{2.843458in}{3.090464in}%
\pgfsys@useobject{currentmarker}{}%
\end{pgfscope}%
\begin{pgfscope}%
\pgfsys@transformshift{3.230958in}{3.409719in}%
\pgfsys@useobject{currentmarker}{}%
\end{pgfscope}%
\begin{pgfscope}%
\pgfsys@transformshift{3.618458in}{3.306670in}%
\pgfsys@useobject{currentmarker}{}%
\end{pgfscope}%
\begin{pgfscope}%
\pgfsys@transformshift{4.005958in}{3.754181in}%
\pgfsys@useobject{currentmarker}{}%
\end{pgfscope}%
\begin{pgfscope}%
\pgfsys@transformshift{4.393458in}{4.188729in}%
\pgfsys@useobject{currentmarker}{}%
\end{pgfscope}%
\begin{pgfscope}%
\pgfsys@transformshift{4.780958in}{4.603908in}%
\pgfsys@useobject{currentmarker}{}%
\end{pgfscope}%
\begin{pgfscope}%
\pgfsys@transformshift{5.168458in}{4.852152in}%
\pgfsys@useobject{currentmarker}{}%
\end{pgfscope}%
\begin{pgfscope}%
\pgfsys@transformshift{5.555958in}{4.918091in}%
\pgfsys@useobject{currentmarker}{}%
\end{pgfscope}%
\begin{pgfscope}%
\pgfsys@transformshift{5.943458in}{5.007679in}%
\pgfsys@useobject{currentmarker}{}%
\end{pgfscope}%
\begin{pgfscope}%
\pgfsys@transformshift{6.330958in}{5.327267in}%
\pgfsys@useobject{currentmarker}{}%
\end{pgfscope}%
\begin{pgfscope}%
\pgfsys@transformshift{6.718458in}{5.626718in}%
\pgfsys@useobject{currentmarker}{}%
\end{pgfscope}%
\begin{pgfscope}%
\pgfsys@transformshift{7.105958in}{6.436941in}%
\pgfsys@useobject{currentmarker}{}%
\end{pgfscope}%
\end{pgfscope}%
\begin{pgfscope}%
\pgfpathrectangle{\pgfqpoint{0.905958in}{0.899712in}}{\pgfqpoint{6.200000in}{6.200000in}} %
\pgfusepath{clip}%
\pgfsetroundcap%
\pgfsetroundjoin%
\pgfsetlinewidth{2.007500pt}%
\definecolor{currentstroke}{rgb}{0.333333,0.658824,0.407843}%
\pgfsetstrokecolor{currentstroke}%
\pgfsetdash{}{0pt}%
\pgfpathmoveto{\pgfqpoint{0.905958in}{1.603564in}}%
\pgfpathlineto{\pgfqpoint{1.293458in}{1.992586in}}%
\pgfpathlineto{\pgfqpoint{1.680958in}{2.422791in}}%
\pgfpathlineto{\pgfqpoint{2.068458in}{2.823159in}}%
\pgfpathlineto{\pgfqpoint{2.455958in}{3.214963in}}%
\pgfpathlineto{\pgfqpoint{2.843458in}{3.481885in}}%
\pgfpathlineto{\pgfqpoint{3.230958in}{4.183556in}}%
\pgfpathlineto{\pgfqpoint{3.618458in}{4.768515in}}%
\pgfpathlineto{\pgfqpoint{4.005958in}{5.571662in}}%
\pgfpathlineto{\pgfqpoint{4.393458in}{5.554919in}}%
\pgfpathlineto{\pgfqpoint{4.780958in}{5.911399in}}%
\pgfpathlineto{\pgfqpoint{5.168458in}{5.853523in}}%
\pgfusepath{stroke}%
\end{pgfscope}%
\begin{pgfscope}%
\pgfpathrectangle{\pgfqpoint{0.905958in}{0.899712in}}{\pgfqpoint{6.200000in}{6.200000in}} %
\pgfusepath{clip}%
\pgfsetbuttcap%
\pgfsetmiterjoin%
\definecolor{currentfill}{rgb}{0.333333,0.658824,0.407843}%
\pgfsetfillcolor{currentfill}%
\pgfsetlinewidth{0.000000pt}%
\definecolor{currentstroke}{rgb}{0.000000,0.000000,0.000000}%
\pgfsetstrokecolor{currentstroke}%
\pgfsetdash{}{0pt}%
\pgfsys@defobject{currentmarker}{\pgfqpoint{-0.038889in}{-0.038889in}}{\pgfqpoint{0.038889in}{0.038889in}}{%
\pgfpathmoveto{\pgfqpoint{-0.038889in}{-0.038889in}}%
\pgfpathlineto{\pgfqpoint{0.038889in}{-0.038889in}}%
\pgfpathlineto{\pgfqpoint{0.038889in}{0.038889in}}%
\pgfpathlineto{\pgfqpoint{-0.038889in}{0.038889in}}%
\pgfpathclose%
\pgfusepath{fill}%
}%
\begin{pgfscope}%
\pgfsys@transformshift{0.905958in}{1.603564in}%
\pgfsys@useobject{currentmarker}{}%
\end{pgfscope}%
\begin{pgfscope}%
\pgfsys@transformshift{1.293458in}{1.992586in}%
\pgfsys@useobject{currentmarker}{}%
\end{pgfscope}%
\begin{pgfscope}%
\pgfsys@transformshift{1.680958in}{2.422791in}%
\pgfsys@useobject{currentmarker}{}%
\end{pgfscope}%
\begin{pgfscope}%
\pgfsys@transformshift{2.068458in}{2.823159in}%
\pgfsys@useobject{currentmarker}{}%
\end{pgfscope}%
\begin{pgfscope}%
\pgfsys@transformshift{2.455958in}{3.214963in}%
\pgfsys@useobject{currentmarker}{}%
\end{pgfscope}%
\begin{pgfscope}%
\pgfsys@transformshift{2.843458in}{3.481885in}%
\pgfsys@useobject{currentmarker}{}%
\end{pgfscope}%
\begin{pgfscope}%
\pgfsys@transformshift{3.230958in}{4.183556in}%
\pgfsys@useobject{currentmarker}{}%
\end{pgfscope}%
\begin{pgfscope}%
\pgfsys@transformshift{3.618458in}{4.768515in}%
\pgfsys@useobject{currentmarker}{}%
\end{pgfscope}%
\begin{pgfscope}%
\pgfsys@transformshift{4.005958in}{5.571662in}%
\pgfsys@useobject{currentmarker}{}%
\end{pgfscope}%
\begin{pgfscope}%
\pgfsys@transformshift{4.393458in}{5.554919in}%
\pgfsys@useobject{currentmarker}{}%
\end{pgfscope}%
\begin{pgfscope}%
\pgfsys@transformshift{4.780958in}{5.911399in}%
\pgfsys@useobject{currentmarker}{}%
\end{pgfscope}%
\begin{pgfscope}%
\pgfsys@transformshift{5.168458in}{5.853523in}%
\pgfsys@useobject{currentmarker}{}%
\end{pgfscope}%
\end{pgfscope}%
\begin{pgfscope}%
\pgfsetrectcap%
\pgfsetmiterjoin%
\pgfsetlinewidth{1.254687pt}%
\definecolor{currentstroke}{rgb}{0.150000,0.150000,0.150000}%
\pgfsetstrokecolor{currentstroke}%
\pgfsetdash{}{0pt}%
\pgfpathmoveto{\pgfqpoint{0.905958in}{0.899712in}}%
\pgfpathlineto{\pgfqpoint{7.105958in}{0.899712in}}%
\pgfusepath{stroke}%
\end{pgfscope}%
\begin{pgfscope}%
\pgfsetrectcap%
\pgfsetmiterjoin%
\pgfsetlinewidth{1.254687pt}%
\definecolor{currentstroke}{rgb}{0.150000,0.150000,0.150000}%
\pgfsetstrokecolor{currentstroke}%
\pgfsetdash{}{0pt}%
\pgfpathmoveto{\pgfqpoint{0.905958in}{0.899712in}}%
\pgfpathlineto{\pgfqpoint{0.905958in}{7.099712in}}%
\pgfusepath{stroke}%
\end{pgfscope}%
\begin{pgfscope}%
\pgfsetbuttcap%
\pgfsetmiterjoin%
\definecolor{currentfill}{rgb}{1.000000,1.000000,1.000000}%
\pgfsetfillcolor{currentfill}%
\pgfsetlinewidth{0.240900pt}%
\definecolor{currentstroke}{rgb}{0.150000,0.150000,0.150000}%
\pgfsetstrokecolor{currentstroke}%
\pgfsetdash{}{0pt}%
\pgfpathmoveto{\pgfqpoint{1.072625in}{5.878293in}}%
\pgfpathlineto{\pgfqpoint{4.476238in}{5.878293in}}%
\pgfpathlineto{\pgfqpoint{4.476238in}{6.933045in}}%
\pgfpathlineto{\pgfqpoint{1.072625in}{6.933045in}}%
\pgfpathclose%
\pgfusepath{stroke,fill}%
\end{pgfscope}%
\begin{pgfscope}%
\pgfsetroundcap%
\pgfsetroundjoin%
\pgfsetlinewidth{2.007500pt}%
\definecolor{currentstroke}{rgb}{0.298039,0.447059,0.690196}%
\pgfsetstrokecolor{currentstroke}%
\pgfsetdash{}{0pt}%
\pgfpathmoveto{\pgfqpoint{1.205958in}{6.674842in}}%
\pgfpathlineto{\pgfqpoint{1.872625in}{6.674842in}}%
\pgfusepath{stroke}%
\end{pgfscope}%
\begin{pgfscope}%
\pgfsetbuttcap%
\pgfsetroundjoin%
\definecolor{currentfill}{rgb}{0.298039,0.447059,0.690196}%
\pgfsetfillcolor{currentfill}%
\pgfsetlinewidth{0.000000pt}%
\definecolor{currentstroke}{rgb}{0.000000,0.000000,0.000000}%
\pgfsetstrokecolor{currentstroke}%
\pgfsetdash{}{0pt}%
\pgfsys@defobject{currentmarker}{\pgfqpoint{-0.038889in}{-0.038889in}}{\pgfqpoint{0.038889in}{0.038889in}}{%
\pgfpathmoveto{\pgfqpoint{0.000000in}{-0.038889in}}%
\pgfpathcurveto{\pgfqpoint{0.010313in}{-0.038889in}}{\pgfqpoint{0.020206in}{-0.034791in}}{\pgfqpoint{0.027499in}{-0.027499in}}%
\pgfpathcurveto{\pgfqpoint{0.034791in}{-0.020206in}}{\pgfqpoint{0.038889in}{-0.010313in}}{\pgfqpoint{0.038889in}{0.000000in}}%
\pgfpathcurveto{\pgfqpoint{0.038889in}{0.010313in}}{\pgfqpoint{0.034791in}{0.020206in}}{\pgfqpoint{0.027499in}{0.027499in}}%
\pgfpathcurveto{\pgfqpoint{0.020206in}{0.034791in}}{\pgfqpoint{0.010313in}{0.038889in}}{\pgfqpoint{0.000000in}{0.038889in}}%
\pgfpathcurveto{\pgfqpoint{-0.010313in}{0.038889in}}{\pgfqpoint{-0.020206in}{0.034791in}}{\pgfqpoint{-0.027499in}{0.027499in}}%
\pgfpathcurveto{\pgfqpoint{-0.034791in}{0.020206in}}{\pgfqpoint{-0.038889in}{0.010313in}}{\pgfqpoint{-0.038889in}{0.000000in}}%
\pgfpathcurveto{\pgfqpoint{-0.038889in}{-0.010313in}}{\pgfqpoint{-0.034791in}{-0.020206in}}{\pgfqpoint{-0.027499in}{-0.027499in}}%
\pgfpathcurveto{\pgfqpoint{-0.020206in}{-0.034791in}}{\pgfqpoint{-0.010313in}{-0.038889in}}{\pgfqpoint{0.000000in}{-0.038889in}}%
\pgfpathclose%
\pgfusepath{fill}%
}%
\begin{pgfscope}%
\pgfsys@transformshift{1.539292in}{6.674842in}%
\pgfsys@useobject{currentmarker}{}%
\end{pgfscope}%
\end{pgfscope}%
\begin{pgfscope}%
\definecolor{textcolor}{rgb}{0.150000,0.150000,0.150000}%
\pgfsetstrokecolor{textcolor}%
\pgfsetfillcolor{textcolor}%
\pgftext[x=2.139292in,y=6.558175in,left,base]{\color{textcolor}\sffamily\fontsize{24.000000}{28.800000}\selectfont FS-QTT-solver}%
\end{pgfscope}%
\begin{pgfscope}%
\pgfsetroundcap%
\pgfsetroundjoin%
\pgfsetlinewidth{2.007500pt}%
\definecolor{currentstroke}{rgb}{0.333333,0.658824,0.407843}%
\pgfsetstrokecolor{currentstroke}%
\pgfsetdash{}{0pt}%
\pgfpathmoveto{\pgfqpoint{1.205958in}{6.197466in}}%
\pgfpathlineto{\pgfqpoint{1.872625in}{6.197466in}}%
\pgfusepath{stroke}%
\end{pgfscope}%
\begin{pgfscope}%
\pgfsetbuttcap%
\pgfsetmiterjoin%
\definecolor{currentfill}{rgb}{0.333333,0.658824,0.407843}%
\pgfsetfillcolor{currentfill}%
\pgfsetlinewidth{0.000000pt}%
\definecolor{currentstroke}{rgb}{0.000000,0.000000,0.000000}%
\pgfsetstrokecolor{currentstroke}%
\pgfsetdash{}{0pt}%
\pgfsys@defobject{currentmarker}{\pgfqpoint{-0.038889in}{-0.038889in}}{\pgfqpoint{0.038889in}{0.038889in}}{%
\pgfpathmoveto{\pgfqpoint{-0.038889in}{-0.038889in}}%
\pgfpathlineto{\pgfqpoint{0.038889in}{-0.038889in}}%
\pgfpathlineto{\pgfqpoint{0.038889in}{0.038889in}}%
\pgfpathlineto{\pgfqpoint{-0.038889in}{0.038889in}}%
\pgfpathclose%
\pgfusepath{fill}%
}%
\begin{pgfscope}%
\pgfsys@transformshift{1.539292in}{6.197466in}%
\pgfsys@useobject{currentmarker}{}%
\end{pgfscope}%
\end{pgfscope}%
\begin{pgfscope}%
\definecolor{textcolor}{rgb}{0.150000,0.150000,0.150000}%
\pgfsetstrokecolor{textcolor}%
\pgfsetfillcolor{textcolor}%
\pgftext[x=2.139292in,y=6.080799in,left,base]{\color{textcolor}\sffamily\fontsize{24.000000}{28.800000}\selectfont FD-QTT-solver}%
\end{pgfscope}%
\end{pgfpicture}%
\makeatother%
\endgroup%